\documentclass[11pt]{article}

\usepackage[numbers]{natbib}
\usepackage{amsbsy}
\usepackage{amsgen}
\usepackage{diagrams}
\usepackage{amstext}
\usepackage{amsxtra}
\usepackage{amsmath}
\usepackage{amsopn}
\usepackage{amsthm}
\newtheorem{thm}{Theorem}[section]
\theoremstyle{definition}
\newtheorem{dfn}{Definition}[section]
\theoremstyle{plain}
\newtheorem*{mr}{Main Result}
\theoremstyle{remark}
\newtheorem{note}{Remark}[section]
\theoremstyle{plain}

\theoremstyle{plain}

\theoremstyle{plain}
\newtheorem{pro}[thm]{Proposition}
\theoremstyle{conjecture}
\newtheorem{con}{Conjecture}[section]

\begin{document}\date{}

\title{Instanton Correction, Wall Crossing And Mirror Symmetry Of Hitchin's Moduli Spaces }
\author{Wenxuan Lu\footnote{wenxuanl@mit.edu}} \maketitle

\begin{abstract}
We study two instanton correction problems of Hitchin's moduli spaces along with their wall crossing formulas. The hyperkahler metric of a Hitchin's moduli space can be put into an instanton-corrected form according to physicists Gaiotto, Moore and Neitzke. The problem boils down to the construction of a set of special coordinates which can be constructed as Fock-Goncharov coordinates associated with foliations of quadratic differentials on a Riemann surface. A wall crossing formula of Kontsevich and Soibelman arises both as a crucial consistency condition and an effective computational tool. On the other hand Gross and Siebert have succeeded in determining instanton corrections of complex structures of Calabi-Yau varieties in the context of mirror symmetry from a singular affine structure with additional data.   We will show that the two instanton correction problems are equivalent in an appropriate sense via the identification of the wall crossing formulas in the metric problem with consistency conditions in the complex structure problem. This is a nontrivial statement of mirror symmetry of Hitchin's moduli spaces which till now has been mostly studied in the framework of geometric Langlands duality.  This result provides examples of Calabi-Yau varieties where the instanton correction (in the sense of mirror symmetry) of  metrics and complex structures can be determined. This equivalence also relates certain enumerative problems in foliations to some gluing constructions of affine varieties. \\\\\\\\\\
\end{abstract}

 \newpage

 \tableofcontents{\section{Introduction}}

In this paper we describe some sophisticated constructions on Hitchin's moduli spaces which eventually match together to form a quite harmonious picture. It provides us some new insights of  differential geometry and algebraic geometry of Hitchin's moduli spaces in the context of mirror symmetry.

$Instanton$ $correction$ in the title refers to the problem of determining the effect of instantons. This kind of problems is crucial in many geometric problems with physical origins. Sometimes certain geometric object (such as a complex structure or a metric) has a decomposition into possibly infinitely many pieces such that the first one represents the "classical" contribution ignoring the quantum effects. Quantum effects usually have a perturbative part and a nonperturbative part. The nonperturbative part can often be considered as contributions associated to some solitonic solutions of the underlying physical theory known as instantons. If this is the case, we say that the object receives instanton corrections.

 According to the Glossary of Polchinski's book STRING THEORY\citep{P} an instanton is

 "$\ldots$ in a Euclidean path integral, a nonconstant configuration that is a local but not a global minimum of the action. Such configurations are usually localized in a spacetime, are usually topologically nontrivial, and are of interest when they give rise to effects such as tunneling that are not obtained from small fluctuations around a constant configuration. Spacetime instantons are instantons in the effective field theory in the spacetime. Worldsheet instantons are instantons in the worldsheet quantum field theory and correspond to worldsheets wrapping around nontrivial two-cycles of spacetime."

  A famous example of instantons is given by gauge theoretical instantons studied in the Donaldson theory. Geometrically these are (anti-)self-dual connections. Another one is worldsheet instantons of  the supersymmetric sigma model in mirror symmetry. These are Gromov-Witten invariants counting holomorphic curves.  Different kinds of instantons and spatially extended solitonic objects such as branes can be related by taking physical dualities, low energy approximations, dimensional reductions, etc.

  Sometimes instanton corrections can be calculated from first principles directly. This is usually hard and involves handling the moduli space of instantons. It is also very common that they can be determined indirectly by imposing dualities or consistency conditions. A famous example is the calculation of worldsheet instanton corrections of quintic threefolds via periods of the mirror\citep{Ca}.

  Another typical phenomenon is $wall$ $crossing$. It is often the case that our instantons contributions  can depend on some moduli parameters so that they are generically locally constant but can behave discontinuously when the moduli parameters cross some exceptional  locus called stability walls. Then there should be a wall crossing formula which basically is a continuity condition. Roughly speaking the wall crossing formula we are going to study in this paper  tells us that when one crosses a stability wall some naturally defined products (ordered compositions) of  instanton contribution terms do not change although the order and the set of instantons
must change. The wall crossing formula is crucial for the consistency of the instanton correction problem and sometimes can actually determine instanton corrections given some input data.\\

The problem discussed in this paper concerns instanton corrections and wall crossing on the moduli space of solutions of  Hitchin's equations on a Riemann surface. Let $\mathcal{M}(R)$ be the moduli space of pairs $(A,\varphi)$ satisfying Hitchin's equations
$$
\begin{array}{rr}F_{A}+R^{2}[\varphi,\bar{\varphi}]&=0\\
\bar{\partial}_{A}\varphi &=0\\
\end{array}
$$on a fixed Riemann surface $C$ of genus $g$ modulo the gauge equivalence. Here $A$ is a connection on a   bundle $E$ of rank 2 and degree 0 over $C$ whose gauge group is $G = SU(2)$, and $\varphi$ is a holomorphic $ad(E)$-valued one form. $R$ is a positive  parameter. We consider this problem in the context of mirror symmetry as a case which is simpler than general compact Calabi-Yau manifolds but  more nontrivial than Fano varieties or local Calabi-Yau varieties.

Mirror symmetry of moduli spaces of Hitchin's systems has attracted some attentions due to the recent physical interpretation of geometric Langlands program \citep{HT, H3,KW,GW,W}. On the other hand, progress in the research field of mirror symmetry has been blocked by the lack of  understandings of instanton corrections  beyond   Fano and local Calabi-Yau cases. Therefore the author feels that the instanton correction problem (including the associated wall crossing problem) of Hitchin's moduli spaces has not only some potential significance to the  Langlands program but also some chances to deepen our understandings of hard parts of mirror symmetry of Calabi-Yau  manifolds.\\

The strongest hint  that this problem is doable is from some new insights offered by physicists. Gaiotto, Moore and Neitzke in \citep{G1} have argued that due to the  hyperkahler structure of the moduli space, the instanton correction problem can be approached via the twistor method and should be reduced to solutions of certain  infinite-dimensional Riemann-Hilbert problems. In particular, the famous wall crossing problem of instantons can be described directly in this way and Kontsevich-Soibelman's general wall crossing formula is recovered.   They go on and provide a geometric construction of these corrections. The instantons corrections in their work are gauge theoretical but instantons are described indirectly as suggested by some brane constructions. On the other hand, Kontsevich-Soibelman's formula also governs the purely algebraic construction of instanton corrections to complex structures in mirror symmetry given by the ground breaking work of Gross and Siebert \citep{GS1}\footnote{This work uses some crucial ideas from Kontsevich and Soibelman's paper\citep{KS1} which deals with the two dimensional case by a different language.}. More precisely, it governs the consistency of certain data of tropical nature (so-called "scattering diagrams") in Gross-Siebert's construction from which the mirror Calabi-Yau family can be built.

Therefore it seems  that there could be a way of showing that for certain noncompact hyperkahler varieties (our Hitchin's moduli spaces) the two instanton correction problems (the metric problem and the complex structure problem) are essentially equivalent via the identification of wall crossing formulas.  In this paper we  show such an equivalence in a rather strong  sense.\\

The main result in this paper can be summarized as follows. Precise statements are contained in section 9.
\begin{mr} (Imprecise version) \ \ The metric instanton correction problem for a Hitchin moduli space is   equivalent to the complex structure
 instanton correction problem in the following sense.

\begin{thm}Consider an $SU(2)$ Hitchin's moduli space with prescribed singularities.  \begin{itemize}\item One can construct on the base of the Hitchin's fibration  a singular  integral affine structure together with a polyhedral decomposition, a nontrivial polarization\footnote{This is a multi-valued piecewise integral linear function adapted to the polyhedral decomposition.} and a log smooth structure\footnote{The term is used in the sense of logarithmic geometry. Basically  it is some data on the central fiber of a family which encodes some information of the family.}. \item A compatible system of consistent $structures$ in the sense of  Gross and Siebert can be constructed from these data. Instanton corrections associated with log morphisms\footnote{These are morphisms between some canonically defined thickenings (deformations) of affine pieces of the central fiber. We need them to achieve consistency of the gluing of these local deformations.} in this construction correspond to critical trajectories of  quadratic differentials which are instantons in the metric instanton correction problem. \item Systems of consistency conditions   of the complex structure problem can be identified with wall crossing formulas of the metric problem. \end{itemize}\end{thm}

\begin{thm} The Hitchin's moduli space (interpreted as the moduli space of $SL(2,\mathbf{C})$ flat connections)   is isomorphic to a generic fiber of the toric degeneration constructed by Gross and Siebert's algorithm  solving the instanton correction problem of complex structures.\end{thm} \end{mr}

The equivalence of the two instanton problems has the following consequences.
\begin{enumerate}
\item This is a nontrivial statement of mirror symmetry. The problem of instanton correction of Calabi-Yau metrics is very important but often completely ignored in the literature of mirror symmetry. According to Strominger, Yau and Zaslow's speculation\citep{SYZ}, the instanton correction of the metric and the complex structure should be contributed by holomorphic disks wrapping some special Lagrangian fibers of the mirror which is supposed to be the total space of a singular special Lagrangian fibration. It is strongly believed that such a description should be understood in an appropriate limit sense associated to a family instead of a single Calabi-Yau space. It is also speculated that the problem of enumerating disks could be reformulated as a problem of enumerating some tropical objects.  From this perspective, Gross and Siebert have solved the problem of instanton corrections of complex structures on the tropical level\footnote{Actually it still requires some work to relate the instanton data in Gross-Siebert approach to tropical objects. See \citep{GPS} for some results in this direction. See section 3 for further discussions.}. By identifying Gross-Siebert's instanton corrections with Gaiotto-Moore-Neitzke's instanton corrections we solve the metric instanton correction problem (in the sense of mirror symmetry) of Hitchin's moduli spaces to the same degree. The equivalence  also answers the question about the relation between metrics and complex structures in mirror symmetry which is largely unknown in general.

\item Usually a statement of mirror symmetry is  formulated for a family of Calabi-Yau's which approach a $large$ $complex$ $limit$. However it is not a priori clear  how this could be implemented for a Hitchin's moduli space because  there seems to be no natural way to obtain algebraic degenerations   with moduli interpretations. On the other hand, in Gaiotto, Moore and Neitzke's work there is a very natural way of introducing an
additional parameter (the parameter $R$ in Hitchin's equations) to the underlying moduli problem which gives us a substitution of a large complex degeneration. The equivalence that we will show then says that from such a family we can construct an algebraic degeneration with explicitly given equations which is compatible with the requirement of mirror symmetry and our original Hitchin's moduli space is a fiber.
\item The equivalence is also an effective tool of actually calculating instanton corrections in mirror symmetry of Hitchin's moduli spaces --- a problem which seems to be very hard to do directly. Calculating instanton corrections in Gaiotto, Moore and Neitzke's work is in principle a gauge theoretical problem because in their work  Hitchin's moduli spaces arise as target spaces of low energy theories of  gauge theories. Form there Gaiotto, Moore and Neitzke have found a construction  of solutions of the Riemann-Hilbert problems for metric instanton corrections based on the theory of quadratic differentials and Fock-Goncharov's construction of certain coordinates on moduli spaces of flat connections.  In this approach, instantons are realized geometrically as some critical trajectories of quadratic differential foliations.  Physically we  use some brane constructions of gauge theory to transform the instanton problem in gauge theory to a problem of counting some critical strings which are boundaries of M2-branes in the M-theory. Of course it is hard to justify mathematically  these physical ideas in general, but like many previous mathematical works on physics-related geometric problems (e.g.mirror symmetry) we can formulate and prove some  nontrivial consequences ( e.g.  calculating Gromov-Witten invariants (worldsheet instantons) via periods) and the solution of the equivalence problem would give us another such example. Our equivalence  means that the enumerative problem of critical trajectories of quadratic differential foliations has an unexpected relation with gluing deformations of some affine varieties.

 \item The geometric meaning of wall crossing is not  clear\footnote{Even the meaning of the word "(stability) wall" is not  clear.} in Gross and Siebert's solution of the complex structure problem. It is expected to account for the jumping of holomorphic disks when one changes some moduli parameters but this picture has not been understood or even studied for Hitchin's moduli spaces. The equivalence offers a very natural  explanation of geometric meaning from another perspective.  Some examples even suggest that the wall crossing in the complex structure problem can be understood at an elementary level of changes of defining equations of the mirror degenerations.

\item It can also be considered as a check of the compatibility of SYZ's
differential geometric version of mirror conjecture (in some limiting
form) and Gross-Siebert-Kontsevich-Soibelman's algebraic version.\end{enumerate}

There are basically two steps  relating the metric problem to the complex structure problem. In the first step we start from   the metric problem and try to construct the input data of the complex structure problem. The input of Gross-Siebert construction consists of "an integral affine structure with singularities with a polyhedral decomposition, a polarization and a  positive log smooth structure" which must be produced from the metric problem.  Consider a  Hitchin's moduli space and deform it  by varying a parameter in Hitchin's equations. The hyperkahler metrics on the moduli spaces can be reduced via an twistor type ansatz of Gaiotto-Moore-Neitzke to Fock-Goncharov coordinates labeled by one cycles on the spectral curves of Hitchin's moduli spaces. Fock-Goncharov coordinates exhibit discontinuous jumps and wall crossing phenomenon when we vary some underlying parameters and the jumps are Kontsevich-Soibelman  transformations. A singular affine structure is induced by the singular special Kahler structure on the base of the so-called Hitchin's fibration. We can use periods of one cycles to produce a polyhedral decomposition and use discontinuous jumps to produce a log smooth structure. The construction involves several layers of structures and  different languages. To match them we  are guided by two principles: labeling by charges and equivalence of wall crossing formulas and consistency conditions.

In the second step we run Gross-Siebert's algorithm and verify the equivalence. This  step is complicated because so is the construction of Gross and Siebert. The equivalence consists of two parts. The first part is a geometric identification of wall crossing formulas in the metric problem with consistency conditions in the complex structure problem.  This is a consequence of our assignments of log smooth data and is not surprising at all. But the ways that the instanton data are encoded are different.  Then we gluing some canonically constructed affine pieces following Gross and Siebert. Eventually we will see that the same Hitchin's moduli space we started with (considered as the moduli space of flat connections) is embedded as a generic fiber into the degeneration obtained in this way. This means that the algebraic degeneration obtained by incorporating instanton corrections  produces an algebraic  degeneration of the moduli space. Moreover the instanton data for the construction of this  degeneration are induced by and equivalent to the metric  instanton data.

Although some constructions in this paper are strongly motivated by some physical ideas that are not rigorously formulated and proved yet the main result is mathematically rigorous because the two sets of instanton data are well defined. In the metric problem they are critical trajectories of quadratic differential foliations while in the complex structure problem they are some log morphisms between some rings.\\

We review background materials from section 2 to section 8 along with some observations and discussions. Some gaps in the literature are filled. The reader can start from section 9.1 where a summary of section 2-8 is given in the beginning. 9.1 also contains a description of main difficulties and an outline of the proof.  Main theorems are proved in section 9.2 and 9.3. Section 9.4 is important for understanding the meanings of the main theorems. Some examples are carefully computed and interpreted in section 9.5. We only study Hitchin's equations whose gauge group is $SU(2)$ for reasons to be explained at the end of section 7, but we allow possibly irregular singularities.\\

The construction in this work is essentially a synthesis of ideas of several important works in different areas. They are: Seiberg and Witten's work on exact description of low energy effective actions and spectra of $N=2$ supersymmetric gauge theories\footnote{In this paper by Seiberg-Witten theory we always mean this work instead of the study of Seiberg-Witten equations which is more familiar to mathematicians}\citep{SW1},  Strominger,Yau and Zaslow's conjecture of mirror symmetry as a T-duality \citep{SYZ} and its extension to families, Gaiotto, Moore and Neitzke's description of instanton corrections of Hitchin's moduli spaces via a clever ansatz of the associated twistor data \citep{G1}, Gaiotto, Moore and Neitzke's geometric realization of their ansatz \citep{G2}, Fock and Goncharov's work on higher Teichmuller theory \citep{FG1,FG2}, Fomin and Zelevinsky's cluster algebras \citep{FZ1,FZ2,FZ3}, Kontsevich and Soibelman's  general wall crossing formula with respect to the change of stability conditions \citep{KS}, and finally Gross and Siebert's work (partially based on Kontsevich and Soibelman's ideas) which provides an purely algebraic solution of the instanton correction problem  of complex structures in mirror symmetry \citep{GS1, GS2}.

 The equivalence considered in this paper must have been anticipated by some experts. In fact the appearance of \citep{GS1}, \citep{KS} and \citep{G1}  makes a strong equivalence of some sort very plausible. The main idea of the proof is essentially trivial once the natural strategy of relating the two instanton correction problems is clear. It is also clear that once the equivalence is established some very nice consequences   deepening our understanding of the geometry and mirror symmetry of Hitchin's moduli spaces would immediately follow. However to actually formulate and prove an appropriate form of the equivalence turns out to be trickier than the author had expected. Due to the complexity of objects and structures involved and quite different languages used in foundational theories many  not very hard but still nontrivial details become unavoidable. Therefore the author believes it could be helpful to write down all the details. The effort can be further justified by the fact that some interesting  consequences and phenomena reveal themselves  after the details have been handled. Finally the author has tried to  put a large amount of background materials that are scattered in the literature into one place and hopes it would make the paper essentially self-contained.\\\\

\noindent {\bf Acknowledgements}\ \ The author would like to express his gratitude to his advisor Shing-Tung Yau for his valuable guidance and support. The author wants to thank his friends in Yau's seminar Po-Ning Chen, Chen-Yu Chi, Ming-Tao Chuan, Si Li, Yi Li, Chien-Hao Liu, Yu-Shen Lin, Ruifang Song, Valentino Tosatti and Jie Zhou for various discussions. The author is also grateful to  Frederik Denef, Ron Donagi, Maxim Kontsevich, Naichung Conan Leung,  Andrew Neitzke, Yong-Geun Oh, Tony Pantev, Paul Seidel, Bernd Siebert and Yan Soibelman for   helpful comments and discussions via conversations or email correspondences.

\section{Hitchin's Moduli Spaces And Special Kahler Geometry}
In this paper we will consider solutions of Hitchin's equations
\begin{equation}
\begin{array}{rr}F_{A}-\phi\wedge\phi &=0\\
d_{A}\phi=d_{A}\star\phi &=0
\end{array}
\end{equation} on a fixed Riemann surface $C$ of genus $g$. Here $A$ is a connection on a   bundle $E$ of rank 2 and degree 0 over $C$ whose gauge group is $G = SU(2)$, $\star$ is the Hodge star, $d_{A}:=d+A$, $\phi$ is an $ad(E)$-valued one form and $\wedge$ really means that we take the wedge of the one form part and the Lie bracket [ , ] of the bundle valued part. We will also use complex notations. In other words, we study a pair ($E$,$\varphi$) called a Higgs bundle or a Higgs pair. Here $E$ is a $holomorphic$-$G$ bundle and $\varphi$ is a $holomorphic$ one form valued in $ad(E)$. A Higgs bundle is obtained from a solution $(A, \phi)$ of equation (1) in the following way. The (0,1) part of $d_{A}$ denoted as $\bar{\partial}_{A}$ defines the holomorphic structure on $E$ and $\varphi$ is the (1,0) part of $i\phi$, $i\phi =\varphi+\bar{\varphi}$. $\varphi$ is also known as the Higgs field. In terms of Higgs bundles, the equivalent form of the equations are\begin{equation}\begin{array}{rr}F_{A}+[\varphi,\bar{\varphi}]&=0\\
\bar{\partial}_{A}\varphi &=0\\
\end{array}
\end{equation}
Following Seiberg-Witten and Gaiotto-Moore-Neitzke, we introduce an additional parameter $R$ and modify the  equations into \begin{equation}
\begin{array}{rr}F_{A}+R^{2}[\varphi,\bar{\varphi}]&=0\\
\bar{\partial}_{A}\varphi &=0\\
\end{array}
\end{equation}For the meaning and the significance of this parameter, see section 5.

 It is well known after Hitchin \citep{H1} that the moduli space $\mathcal{M}$ of solutions of Hitchin's equations modulo the gauge equivalence is a  noncompact hyperkahler space. It is obtained by an infinite dimensional hyperkahler quotient construction of the moduli space. In fact the tangent space of the pair $(A,\phi)$ is an infinite dimensional affine space endowed with a natural flat hyperkahler metric\footnote{Here the notations are the same as those of \citep{KW}.} \begin{equation}ds^{2}=-{1\over 4\pi}\int_{C}|d^{2}z|\mathrm{Tr}(\delta A_{z}\otimes\delta A_{\bar{z}}+\delta A_{\bar{z}}\otimes\delta A_{z}+\delta \phi_{z}\otimes\delta \phi_{\bar{z}}+\delta \phi_{\bar{z}}\otimes\delta \phi_{z})\end{equation}where $A=A_{z}dz+A_{\bar{z}}d\bar{z}$ and $\phi=\phi_{z}dz+\phi_{\bar{z}}d\bar{z}$ for the holomorphic coordinate $z$ over $C$. $\delta$ denotes the tangent vectors. The group of gauge transformations acts on this flat hyperkahler space and the set of solutions of Hitchin's equations turns out to be the zero level set of associated three moment maps. Therefore $\mathcal{M}$ is obtained as a hyperkahler quotient.

 As a hyperkahler space, it has a set of compatible complex structures parameterized by $\xi\in CP^{1}$ and generated by three independent complex structures $J_{1},J_{2},J_{3}$  (in \citep{KW} these are denoted by $-J, -K, I$. For their explicit descriptions see \citep{KW} or \citep{H1}). The three independent complex structures satisfy the quaternion relations $$J^{2}_{1}=J^{2}_{2}=J^{2}_{3}=J_{1}J_{2}J_{3}=-1$$The set of all compatible complex structures are given by $$J_{\xi}:={i(-\xi+\bar{\xi})J_{1}-(\xi+\bar{\xi})J_{2}+(1-|\xi|^{2})J_{3}\over 1+|\xi|^{2}}$$ where $\xi\in CP^{1}$  is called the twistor parameter. Let $\omega_{i}$ be the Kahler form in $J_{i}$. Define $$\Omega_{1}=\omega_{2}+i\omega_{3}$$ $$\Omega_{2}=\omega_{3}+i\omega_{1}$$ $$\Omega_{3}=\omega_{1}+i\omega_{2}$$ then $\Omega_{i}$ is a holomorphic symplectic form in $J_{i}$. It turns out that $\Omega_{1}$ does not depend on the complex structure of the Riemann surface while the other two do.

 \begin{itemize}\item There are two opposite special complex structures which  are identified with $\pm J_{3}$ (i.e. $\xi$ is 0 or $\infty$). The moduli space $\mathcal{M}$ in $J_{3}$ is identified as the moduli space of semistable Higgs pairs. This moduli space is quasi-projective.\item $\mathcal{M}$ is identified as the moduli space of  $SL(2,\mathbf{C})$ flat connections\footnote{The gauge group is $SL(2,\mathbf{C})$ because it is the complexification of $SU(2)$.} when $\xi\neq 0,\infty$. \citep{C,D,H1}. In fact, Hitchin's equations tell us that
the new connection ${R\over\xi}\varphi+A+R\xi\bar{\varphi}$ is flat.   The moduli space of flat connections is Stein.\item The moduli space of $SL(2,\mathbf{C})$ flat connections is analytically isomorphic the moduli space of representations $\pi_{1}(C)\rightarrow SL(2,\mathbf{C})$ (i.e. the categorical quotient of the conjugate action on the space of such representations.). The latter is also called the character variety. The isomorphism is not algebraic. This is called the Riemann-Hilbert correspondence. The moduli space of fundamental group representations is an affine variety. \item There is a $\mathbf{C}^{\times}$ action on $\mathcal{M}$ given by $\varphi\rightarrow \lambda\varphi$ where $\lambda\in\mathbf{C}^{\times}$. If we restrict it to $U(1)$ (i.e. $|\lambda|=1$) then the action is isometric with respect to the hyperkahler metric. In general it preserves $\pm J_{3}$ and all the other $J_{\xi}$ are in one orbit. More precisely speaking we have a lifted  $\mathbf{C}^{\times}$ action on the so-called twistor space $\mathcal{M}\times CP^{1}$ (see section 5) which covers the natural $\mathbf{C}^{\times}$ action on $CP^{1}$ fixing the north pole and the south pole labeling $\pm J_{3}$. Because of the meaning of the twistor sphere $CP^{1}$ (i.e. it is the set of twistor parameters) this means that all the complex structures except $\pm J_{3}$ are holomorphically equivalent \citep{H1}. Each one of them has the moduli interpretation as the moduli space of flat connections. But note that  only $J_{1}$ and $-J_{1}$ are independent of the complex structure of the Riemann surface. So we view $J_{1}$ as the canonical complex structure of the moduli space of flat connections. \end{itemize}

There is a natural fibration from the moduli space to the space of quadratic differentials (we denote it by $B$) by taking the determinant of $\varphi$ (this is called Hitchin's map). It is called Hitchin's fibration \citep{H1,H2} and has the following properties\begin{itemize} \item The Hitchin's map denoted as $det$ is holomorphic, surjective and proper with respect to the complex structures in which the moduli space $\mathcal{M}$ is the moduli space of semistable Higgs pairs. \item Fibers of this map have nice geometric meanings. Define a curve in the total space of the canonical bundle of $C$ by the characteristic polynomial of $\varphi$
\begin{equation}det(x-\varphi)=0\end{equation}Note that the trace of $\varphi$ is zero. So the equation is $x^{2}+det\ \varphi=0$.  This curve is called the spectral curve or the Seiberg-Witten curve $S$ and there is one such curve associated to each element of $B$. As an abelian variety, the fiber above $u\in B$ is  the Prym variety $\tilde{J}(S_{u})$\footnote{Denote  the Jacobian  of $S_{u}$ by $J(S_{u})$. The Prym variety is the kernel of the natural map  $J(S_{u})\rightarrow J(C)$ induced by $S_{u}\rightarrow C$.} of the projection $S_{u}\rightarrow C$. Therefore we have  realized the moduli space as a family of complex abelian varieties.\item Let $\mathcal{N}$ be the moduli space of semistable bundles of fixed determinant with rank 2 and degree 0. Then  $T^{*}\mathcal{N}$ is a subset of $\mathcal{M}$. In fact $\mathcal{M}$ is a fiberwise compactification of $T^{*}\mathcal{N}$ with respect to the function $det$. Each fiber is compactified by adding a codimensional $g$ subvariety. Moreover the holomorphic symplectic structure of $T^{*}\mathcal{N}$ is the restriction of the holomorphic symplectic structure induced by the hyperkahler structure of $\mathcal{M}$.\item The complex dimension of $\mathcal{M}$ is $6g-6$. In particular we have to assume that $g> 1$. However later we will allow singularities of solutions of Hitchin's equations and this restriction will be removed then. \end{itemize}

These abelian varieties (fibers of the Hitchin's fibration) are not special Lagrangian submanifolds with respect to $J_{3}$. In fact they are complex Lagrangian in this case. But once we rotate the complex structure to $J_{1}$ (this is known as a hyperkahler rotation, see \citep{H5}) then these torus fibers are special Lagrangian. So after a hyperkahler rotation we are in a situation where Strominger, Yau and Zaslow's recipe applies. In other words, the mirror of this noncompact hyperkahler space should be the total space of a special Lagrangian torus fibration whose base is topologically the same as $B$ and whose fibers are dual tori. Since our space is a moduli space, we may wonder if its mirror is also a moduli space. This is the question investigated by Hausel and Thaddeus and they found that the mirror is in fact the Hitchin's moduli space on the same Riemann surface with the gauge group $G$ (in our case it is $SL(2,\mathbf{C})$) replaced by the Langlands dual group $G^{L}$ (in our case it is $PGL(2,\mathbf{C})$). In fact due to the explicitness of Hitchin's fibrations, one can check that the fibers of the  Hitchin fibration of the second moduli space is dual to the fiber of the first moduli space as abelian varieties. This discovery suggests that this mirror symmetry has something to do with the Langlands duality and later it partly inspired  Kapustin, Witten and Gukov's ambitious program of reformulating the geometric Langlands duality using mirror symmetry of Hitchin's moduli spaces \citep{GW,KW,W}. In particular, the mirror relation is generalized to the cases with (possibly irregular) singularities.\\

It is well known that the base of a special Lagrangian fibration has some special geometric structures. When the total space is hyperkahler and fibers are complex abelian varieties in certain complex structures, the geometric structure is even more special and is known as a special Kahler structure. Good references include \citep{F,H4,Sa}.

\begin{dfn}A special Kahler structure on a Kahler manifold with Kahler form $\omega$ is a real flat torsion-free symplectic connection (symplectic means $\nabla\omega=0$) $\nabla$ satisfying$$d_{\nabla}I=0$$where $I$ is the complex structure and $d_{\nabla}$ is the extension of the connection to the de Rham complex valued in $TM$.\end{dfn} There are two kinds of special coordinates that one can introduce on a special Kahler manifold. The first kind is a system of flat Darboux coordinates $(x_{i},y_{i})$. They are flat in the sense that $$\nabla dx_{i}=\nabla dy_{i}=0$$ and they are Darboux in the sense that $$\omega=dx^{i}\wedge dy_{i}$$ Later we will also refer this set of coordinates as affine coordinates because their transition functions are affine transformations. The other set of special coordinates are dual pairs of holomorphic coordinates $(a_{i}, a^{D}_{i})$ which can be chosen to be adapted to a given set of affine coordinates. They are redundant, i.e. one needs only half of them (for example $a_{i}$) to provide local holomorphic coordinates. They are adapted to affine coordinates $(x_{i},y_{i})$ if\begin{equation}\mathrm{Re}(da_{i})=dx_{i}, \mathrm{Re}(da^{D}_{i})=-dy_{i}\end{equation} Given a set of affine coordinates one can always find a set of adapted special holomorphic coordinates and vice versa.  The condition of being dual is $${\partial\over\partial a_{i}}={i\over 2}({\partial\over\partial x_{i}}-\tau_{ij}{\partial\over\partial y_{j}})$$ where\begin{equation}\tau_{ij}={\partial a^{D}_{j}\over\partial a_{i}}\end{equation}To make $\omega$ a type $(1,1)$ form, we must have $\tau_{ij}=\tau_{ji}$. So locally we have a holomorphic function $F$ called the prepotential such that $$a_{i}^{D}={\partial F\over\partial a_{i}}$$ The Kahler form is \begin{equation}\omega={\sqrt{-1}\over 2}\mathrm{Im}(\tau_{ij})da^{i}\wedge d\bar{a}^{j}\end{equation}and the metric is \begin{equation}ds^{2}=\mathrm{Im}(\tau_{ij})da^{i}d\bar{a}^{j}=-{\sqrt{-1}\over 2}(da^{D}_{j}d\bar{a}^{j}-da^{j}d\bar{a}^{D}_{j})\end{equation}

Away from singular fibers of the Hitchin's fibration, our moduli space is an algebraically integrable system which always induces a special Kahler structure on the base. In out context, it can be explicitly described as follows \citep{Do,DoM,DW,HHP}. There is a holomorphic one form known as the Seiberg-Witten differential $\lambda$. It is the restriction to $S$ of the tautological one form. In fact the spectral curve $S$ is a ramified double cover lying in the total space of the cotangent bundle of $C$. Locally, one can choose Darboux coordinates and write the canonical holomorphic symplectic form as $dx\wedge dz$ where $z$ is holomorphic coordinate on $C$ and $x$ is the vertical coordinate. Then we define $\lambda$ to be the restriction of $xdz$ to $S$. Note that because of (5) \begin{equation}\lambda^{2}=-det\ \varphi\end{equation}
\begin{note}The definition of the Seiberg-Witten differential given here is the same as the one in \citep{G2}. However in most literatures a Seiberg-Witten differential is meant to be a one form over $\mathcal{M}$. The second formulation is in fact induced by the first one. This follows from two facts. First, the differential of a Seiberg-Witten differential in both cases is the canonical symplectic form in the context (for the first formulation it is the canonical one on the cotangent bundle, for the second it is the holomorphic symplectic form of $\mathcal{M}$ which extends the canonical holomorphic symplectic form on $T^{*}\mathcal{N}$). Second, the symplectic form in the first formulation induces via the Abel-Jacobi map the symplectic form in the second formulation. See \citep{DW} for details  along with a comparison with Mukai's famous result that the moduli space of simple sheaves of a symplectic surface has a natural symplectic structure.Since we will discuss the foliation of  quadratic differentials  over the Riemann surface $C$  we will be using the first formulation. Therefore unlike most literatures the  special Kahler structure on the base of the Hitchin's fibration will be constructed from periods of the Seiberg-Witten differential defined in this way. A thorough treatment from this point of view for general gauge groups is given in \citep{HHP}.  \end{note}

We introduce some notions which will be used throughout the paper. Recall that a spectral curve $S$ is a hyperelliptic curve defined by (5). Strictly speaking that equation defines a curve with punctures before we complete  it into a projective curve. These punctures are lifts of singularities of $\lambda^{2}$ on $C$ (we shall allow $\lambda^{2}$ to have poles). By its normalization (when we define the Prym variety we actually use this normalization) denoted as $\bar{S}$ we mean this projective completion by filling punctures.
\begin{dfn}\citep{G2} Let $S_{u}$ be a spectral curve of a Hitchin's moduli space $\mathcal{M}$. We consider the odd part of $H_{1}(S_{u},\mathbf{Z})$. Here odd means that the cycle is invariant under the combined operation of exchanging the two sheets  and reversing the orientation. They fit into a local system over the nonsingular part of $B$. The local system is called the charge lattice and  degenerates at the singular locus of $B$ where some cycles become vanishing cycles. A charge is an element of the charge lattice local system. So locally by choosing an trivialization of the local system a charge is just an element of the associated integral lattice and  it has monodromies. The charge lattice is denoted as $\hat{\Gamma}$ and is endowed with the skew-symmetric intersection pairing of integral one cycles. The gauge charge lattice denoted by $\Gamma_{gau}$ is defined to be the local system of odd parts of $H_{1}(\bar{S_{u}},\mathbf{Z})$ together with the intersection paring. The flavor charge lattice $\Gamma_{flavor}$ is the radical of the intersection paring in $\hat{\Gamma}$. It consists of integral combinations of loops around the punctures. Note these lattices  fit into an exact sequence $$0\rightarrow\Gamma_{flavor}\rightarrow\hat{\Gamma}\rightarrow\Gamma_{gau}\rightarrow 0$$\end{dfn}

Let $(A_{i},B_{i}), 1\leq i\leq3g-3$ be a symplectic basis of the gauge charge lattice. The genus of $C$ is $g$ while the genus of  $\bar{S_{u}}$ is $4g-3$. So the rank of $\Gamma_{gau}$ is $6g-6$ which matches the dimension of $B$.  Special holomorphic coordinates of the special Kahler structure are given by period maps\begin{equation}a_{i}(u):={1\over\pi}\int_{A_{i}}\lambda, \ \ a^{D}_{i}(u):={1\over\pi}\int_{B_{i}}\lambda\end{equation} \begin{equation}\tau_{ij}(u):={da^{D}_{j}(u)/du\over da_{i}(u)/du}\end{equation} Note that $a_{i}(u)$ and $a^{D}_{i}(u)$ depend on
the holomorphic coordinate  on $B$    denoted by $u$ canonically defined by the value of the Hitchin's map $$u=det\ \varphi$$

\begin{dfn}The central charge $Z_{\gamma}$ for a charge  $\gamma$ is defined by \begin{equation}Z_{\gamma}={1\over\pi}\int_{\gamma}\lambda\end{equation}The central charge depends on $u$ since $\lambda$ does.\end{dfn}

Of course due to the existence of singular fibers what we really have on the base $B$ is a singular special Kahler structure. We will focus on the induced affine structure.\begin{dfn}An affine structure with singularities on a topological manifold is an (nonsingular) affine structure outside a locally finite union of locally closed submanifolds of codimension greater or equal to 2. It is integral if transition functions of the affine structure over the complement of the singular locus are integral affine transformations. The singular locus is denoted by $\Delta$.\end{dfn}
Going around a component of singular locus gives us a holonomy representation of the corresponding element of the fundamental group of the regular part to the group of affine transformations $Aff(M_{\mathbf{R}})$ where $M:=\mathbf{Z}^{n}$ is an integral lattice and $M_{\mathbf{R}}:=M\otimes \mathbf{R}$. For more information on singular affine structures, see \citep{GS2}. Note that the lattices defining torus fibers fit together to form a (degenerate) local system over $B$ and as such it also has monodromies around the singular locus and theses monodromies are precisely the linear parts of the corresponding holonomy representations. For Hitchin's moduli spaces, singular locus arise when some cycles on spectral curves over generic points of the base degenerate and therefore the monodromies can be read from the Picard-Lefschetz transformations of vanishing cycles.\\

The above theory holds for  more general types of Hitchin moduli spaces. We allow possibly singular solutions of Hitchin's equations. See \citep{BB, Ko,S, GW, W} for a fraction of the huge literature. The exposition below mostly follows  \citep{GW,W}. It is a little bit technical. The reader can safely skip it as long as he or she believes there is a reliable foundational theory for Hitchin's moduli spaces with prescribed singular behaviors at finitely many singularities.

Recall that when the twistor parameter $\xi\in \mathbf{C}^{\times}$, the moduli space of Hitchin's equation over a Riemann surface $C$ is  the moduli space of flat connections. Suppose $(A,\varphi)$ is a solution of Hitchin's equations (2),  then \begin{equation}\mathcal{A}:= {R\over\xi}\varphi+A+R\xi\bar{\varphi}\end{equation} is a flat $SL(2,\mathbf{C})$ connection.

We will assume there are possibly irregular singularities of the Hitchin's equations. Let $p$ be a singularity and take a trivialization of the holomorphic bundle $E$ in a small neighborhood of $p$ with local holomorphic coordinate $z$ such that the (0,1) part of $d_{\mathcal{A}}$ in this coordinate is given as $\bar{\partial}_{\mathcal{A}}=\bar{\partial}_{\bar{z}}$. A covariantly constant section of the $flat$ bundle $E$ is denoted as $\Psi$.  Let $\partial_{\mathcal{A}}$ be the (1,0) part of $d_{\mathcal{A}}$ and define $\mathcal{A}_{z}$ by $\partial_{\mathcal{A}}= \partial_{z}+\mathcal{A}_{z}$. The flatness of $\mathcal{A}$ implies that $\mathcal{A}_{z}$ is  holomorphic away $p$. The singularity of $\mathcal{A}_{z}$ is of the form $$\mathcal{A}_{z}= {T_{n}\over z^{n}}+{T_{n-1}\over z^{n-1}}+\cdots +{T_{1}\over z}+\cdots, n\geq 1$$where the second ellipses represent regular terms. If $n=1$ then it is a regular singularity. Otherwise $p$ is an irregular singularity. The bundle $E$ is both holomorphic and flat away from $p$. It can be extended over $p$ as a holomorphic bundle, but the extension is not unique (depending on the choice of a parabolic structure). The extensions will not be needed for this paper. We will assume that $T_{n}$ is regular and semi-simple which means that it can be diagonalized and has distinct eigenvalues. This assumption is unnecessary and we only use it to simplify the exposition. All the facts concerning the moduli space continue to be true  even if we relax this assumption a little bit. See section 6 in \citep{W} for discussions concerning this point.

Under this assumption we can find a meromorphic gauge transformation to diagonalize all $T_{i}$, i.e. in this gauge\begin{equation}\mathcal{A}_{z}= {T_{n}\over z^{n}}+{T_{n-1}\over z^{n-1}}+\cdots +{T_{1}\over z}+\cdots, n\geq 1\end{equation} with $T_{i}\in \mathrm{t}_{\mathbf{C}}$ where $\mathrm{t}_{\mathbf{C}}$ is the Lie algebra of of a maximal torus $\mathrm{T}_{\mathbf{C}}$. Note that $T_{1}$ is the residue of $\mathcal{A}_{z}dz$.  The diagonalization is not unique. There are additional meromorphic gauge transformations that can change the the eigenvalues of $T_{1}$. Fixing $T_{1}$ would fix this freedom. To formulate the moduli problem, $T_{1}$, which is of topological nature, is fixed. It turns out that picking $T_{1}$ is equivalent to picking a holomorphic extension over the singular point of the flat bundle $E$ on the punctured disk. However there are still additional freedom of permuting eigenvalues of $T_{n}$ which is a Weyl group action and of holomorphic gauge transformations which are diagonal up to order $|z|^{n}$. The gauge group action will be taken care of when we formulate the moduli problem.

Usually it is useful to separate $A$ and $\phi$. If we write $z=re^{i\theta}$, then  convenient form of the singular parts is given as$$A=\alpha d\theta$$\begin{equation}\phi={dz\over 2}({u_{n}\over z^{n}}+\cdots +{u_{1}\over z}+\cdots)+{d\bar{z}\over 2}({\bar{u}_{n}\over \bar{z}^{n}}+\cdots +{\bar{u}_{1}\over \bar{z}}+\cdots)\end{equation}where $\alpha\in \mathrm{t}$, the Lie algebra of a maximal torus of $G$ and $u_{i}\in \mathrm{t}_{\mathbf{C}}$. After a gauge transformation, $\mathcal{A}_{z}$ can be put into the standard form with $T_{1}=-i(\alpha-i\mathrm{Im}\ u_{1}), T_{k}=u_{k}, k>1$. This is called the local model of abelian singularities in \citep{W}. The purpose of the regular semi-simplicity assumption of $T_{n}$ is to get this local model.

The Hitchin's moduli space with possibly irregular singularities (denoted as $\mathcal{M}$) referred in this paper is defined to be the space of pairs $(A,\phi)$ satisfying Hitchin's equations with prescribed local models of abelian singularities modulo the action of the group of $SU(2)$-valued gauge transformations that are $\mathrm{T}$-valued modulo terms of order $|z|^{n}$. Note that the gauge group described here preserves the form of the local model. It has been shown that Hitchin's hyperkahler quotient construction extends to this case and the moduli space $\mathcal{M}$ is hyperkahler.\\

$\mathcal{M}$ is identified with the moduli space of flat connections with prescribed singular parts (15)  at singularities   when the twistor parameter $\xi\neq 0,\infty$.  This means in (15) we fix $T_{1}, T_{2},\cdots, T_{n}$. The underlying complex structures of the moduli space of flat connections for different $\xi$ are holomorphically equivalent.

The complex dimension of the moduli space is $6g-6+2n$ if there is one singularity which is an order $n$ pole. Note that $6g-6$ is the dimension of the moduli space of nonsingular solutions. So the extra freedom introduced by allowing one order $n$ singularity is $2n$. It is obvious how to generalize to the case of several singularities.

The holomorphic symplectic form of $\mathcal{M}$ (which  depends on $\xi$) is given by $$\Omega=-{i\over 4\pi}\int_{C}\mathrm{Tr}\delta\mathcal{A}\wedge\delta\mathcal{A}$$Since all $T_{i}$ are fixed  $\delta\mathcal{A}$ would be nonsingular even though $\mathcal{A}$ is singular. This makes the space  $\mathcal{M}$  a natural place to define this holomorphic symplectic form.

Although it seems that the construction depends on the underlying complex structure of $C$ as well as the asymptotic data $T_{i}, i\geq 2$, it turns out the canonical complex structure  is independent of the choice of complex structure of $C$, positions of singularities and $T_{i}, i\geq 2$. In fact, varying the complex structure of $C$, positions of singularities and $T_{i}, i\geq 2$ gives rise to the so-called isomonodromic deformations, see the discussion below as well as in \citep{W}. It is also worthy pointing out that the symplectic structure of the moduli space does not depend on $T_{i}, i\geq 2$ either\citep{Bo}. On the other hand the complex structure  of the moduli space depends holomorphically on $T_{1}$ and the symplectic structure also depends on it.

The story of the Seiberg-Witten differential is modified accordingly \citep{DW,G2}. The residue of $\varphi$ at a singularity can be  diagonalized as \begin{equation}\varphi={dz\over z}(\left(\begin{array}{cc}m & 0\\0 & -m \end{array}\right)+\cdots)\end{equation} Then $m$ is the residue of $\lambda$ because $$\lambda^{2}=-det\ \varphi\sim {m^{2}\over z^{2}}dz^{2}$$These residues are also called masses or mass parameters by physicists, see section 4.\\

There seems to exist analogous relations between moduli spaces of flat connections and analytically isomorphic moduli spaces of fundamental group representations when we allow singularities (the author does not know if there is  a precise reference for general cases. The results in this paper are fine without it as we do our constructions on the moduli space of flat connections). For example in \citep{W} the author identifies $\mathcal{M}$ (the moduli space of flat connections with prescribed singularities) with the following moduli space denoted by $\mathcal{Y}^{\ast}(T_{1})$ and defined below.

By the flatness of $\mathcal{A}$, $\Psi$ is annihilated by both $\bar{\partial}_{\mathcal{A}}$ and $\partial_{\mathcal{A}}$. Locally near $p$ but away from $p$, the first condition means $\Psi$ is holomorphic (since $\bar{\partial}_{\mathcal{A}}=\bar{\partial}_{\bar{z}}$). The second condition means $\Psi$ is a holomorphic solution of a meromorphic equation $$(\partial_{z}+\mathcal{A}_{z})\Psi=0$$ and as such it exhibits the Stokes phenomenon near $p$ if $p$ is an irregular singularity \citep{Wa,W}.

This means a small open disk containing $p$ is decomposed into $(2n-2)$ sectors. Each sector contains precisely one "Stokes ray". In  sector $\alpha$, up to gauge equivalence there is a unique fundamental solution matrix $Y_{\alpha}$ with appropriately prescribed asymptotic behavior. In the intersection $S_{\alpha}\cap S_{\alpha+1}$
there is a matrix (Stokes matrix) $M_{\alpha}$ such that $Y_{\alpha+1}=Y_{\alpha}M_{\alpha}$.

 The set of $(M_{i}, T_{1})$ for a singularity is called the generalized monodromy data of that singularity. It determines the actual monodromy $\hat{M}$ around the singularity by $$\hat{M}=\exp(-2\pi T_{1})M_{2n-2}M_{2n-1}\cdots M_{1}$$where $\exp(-2\pi T_{1})$ is called the formal monodromy.
We always consider generalize monodromy data up to gauge equivalence. Here by gauge equivalence we mean the data must be considered modulo the action of the maximal torus $\mathrm{T}_{\mathbf{C}}$.

 Without loss of generality suppose there is one singularity. Let $U_{i}$ and $V_{i}$ be respectively the  images of generators $A_{i}$ and $B_{i}$ of the fundamental group of the genus $g$ Riemann surface into $G_{\mathrm{C}}$ under the monodromy representation.  These monodromies together with the generalized monodromy data at the singularity and the connecting matrix defined below are called the generalized monodromy data for the punctured surface. Then we consider the space of generalized monodromy data modulo gauge equivalence. It is denoted as $\mathcal{Y}^{\ast}$ in \citep{W}.

We have the following monodromy relation  $$U_{1}V_{1}U^{-1}_{1}V_{1}^{-1}\cdots U_{g}V_{g}U^{-1}_{g}V_{g}^{-1}W\exp(-2\pi T_{1})M_{2n-2}M_{2n-1}\cdots M_{1}W^{-1}=1$$ where $W$ accounts for a connecting matrix representing the parallel transport from a base point in the first Stokes sector of the singularity to the base point chosen to write the relation among $U$ and $V$.  $\mathcal{Y}^{\ast}$ is the space of solutions of the above identity modulo gauge equivalence. The gauge group here is $G_{\mathbf{C}}\times\mathrm{T}_{\mathbf{C}}$. $g\in G_{\mathbf{C}}$ acts as $U_{i}\rightarrow g U_{i} g^{-1}, V_{i}\rightarrow g V_{i} g^{-1}, W\rightarrow gW$ while $h\in\mathrm{T}_{\mathbf{C}}$ acts as $w\rightarrow W h^{-1}, M_{\alpha}\rightarrow hM_{\alpha}h^{-1}$. $\mathcal{Y}^{\ast}$ is an affine variety.

There is another space denoted as $\mathcal{Y}^{\ast}(T_{1})$. It is the  space of  generalized monodromy data  with fixed $T_{1}$. It is natural to fix $T_{1}$ because of its topological nature (a residue) and its physical interpretation as the mass parameter. If we use  local coordinates to write the local model as before $$\mathcal{A}_{z}= {T_{n}\over z^{n}}+{T_{n-1}\over z^{n-1}}+\cdots +{T_{1}\over z}+\cdots, n\geq 1$$ with $T_{i}\in \mathrm{t}_{\mathbf{C}}$, then when defining $\mathcal{Y}^{\ast}$ we fix $T_{2},\cdots, T_{n}$ while when defining $\mathcal{Y}^{\ast}(T_{1})$ we fix $T_{1}, T_{2},\cdots, T_{n}$. If we change $T_{2},\cdots, T_{n}$ generalized monodromy data do not change. That is why it was called an isomonodromic deformation before.

\section{From Affine Structures To Instanton Corrections of Complex Structures}

\noindent {\bf Motivations}\\
A good reference is \citep{A}. Let us start by recalling the famous proposal of Strominger, Yau and Zaslow \citep{SYZ}.\\

\noindent {\bf Conjecture of Strominger-Yau-Zaslow} (very rough version) \  If $X$ and $\tilde{X}$  is a mirror pair of Calabi-Yau manifolds  then there should be singular special Lagrangian torus fibrations $X\rightarrow B$ and $\tilde{X}\rightarrow\tilde{B}$ which are dual torus fibrations.\\

There are some problems with this version of SYZ conjecture.
\begin{itemize}\item We have to consider families of Calabi-Yau's. In mirror symmetry we usually consider the so-called large complex degeneration. It means that the family approaches a limit point in the complex moduli space with maximally unipotent monodromy. However in this paper we are not going to use this algebraic geometric definition.  Instead  we will work with a differential geometric characterization proposed in section 5.
\item Even for families we may have to understand it only in a limit sense.
\item The problem of instanton corrections by holomorphic disks. It is expected by some heuristic arguments that the complex structures and perhaps the Calabi-Yau metrics should receive some contributions labeled by worldsheet instantons: holomorphic disks wrapping some special Lagrangian fibers of the mirror. \end{itemize}

A refined version which has some chances to be true was proposed by Gross-Wilson \citep{GW1} and Kontsevich-Soibelman \citep{KS2}. \\

\noindent {\bf Conjecture of Limit Form}(Gross-Wilson, Kontsevich-Soibelman)\footnote{This form of SYZ conjecture was proposed by Gross-Wilson and Kontsevich-Soibelman around 2000. But the affine structure and Monge-Ampere structure were already known to Hitchin \citep{H4}\citep{H5}. The idea that one should work with a large complex family to study the SYZ conjecture also appeared in Strominger-Yau-Zaslow's original paper \citep{SYZ}.} :  Consider a large complex degeneration of Calabi-Yau manifolds. If we rescale the Ricci-flat metrics to fix the diameter then there exists a subsequence converging in the Gromov-Hausdorff sense to a limit metric space. This limit space is expected to have a singular affine structure and a singular Monge-Ampere metric.\\

A singular affine structure is an affine structure outside a codimensional 2 subset. A (real) Monge-Ampere metric means that the metric can be locally written as the second order derivative of a convex function $K$ and $K$ satisfies the real Monge-Ampere equation. We expect  that such a structure exists because there is a nonsingular (real) Monge-Ampere metric on the base of a nonsingular special Lagrangian torus fibration and it is expected that in the large complex limit the discriminant locus of the singular fibrations (the Calabi-Yau varieties) on the base will collapse to be codimension 2. If we have a dual special Lagrangian torus fibration then the potential of the corresponding (dual) Monge-Ampere metric should be obtained by taking the Legendre transform of the potential $K$ of the Monge-Ampere metric associated to the original fibration.

So mirror symmetry should be geometrically realized in three steps.\begin{enumerate}
\item Construct the limit structure (singular affine structure plus the potential of the Monge-Ampere metric) from a given large complex degeneration of Calabi-Yau varieties.
\item Do the Legendre transform with respect to the potential  of the limit structure.
\item Solve the instanton correction problem (also known as the reconstruction problem): reconstruct a family of Calabi-Yau manifolds (which is meant to be the mirror family) from the dual limit structure and the construction is supposed to incorporate the contributions of instantons.\end{enumerate}

The idea of Gross and Siebert in \citep{GS1}\citep{GS2} is that we forget about special Lagrangians and metrics and instead focus on some limit structure which does not refer to metrics and try to build everything from algebraic geometry.

The limit structure in their program is the following data \begin{itemize}\item A singular integral affine manifold $B$,

\item A polyhedral decomposition $\mathcal{F}$ ,

 \item A polarization $\varphi$,

\item A log smooth structure. \end{itemize}

The new three steps in mirror symmetry are\begin{enumerate}\item Construct a triple (the first three items given above) from a toric degeneration of Calabi-Yau varieties. Toric degeneration is a notion adapted to Gross-Siebert's program (see below).
\item Perform the Legendre transform to get a dual triple.
\item Reconstruction Problem: reconstruct a toric degeneration (the mirror family) of Calabi-Yau varieties from the dual triple and a log smooth structure which is supposed to encode the information of instantons.\end{enumerate}
The hard part is of course the third one which is the algebraic geometric analogue of the instanton correction problem. The  complex structure instanton correction problems dealt with in this paper is precisely this problem.
\begin{note} There is a related program proposed by Kontsevich and Soibelman \citep{KS1} which is the first one that uses log morphisms (although they did not use the language of log geometry). They only considered the case of the two dimension (K3 surfaces). It seems this work at least partially inspired the algebraic formalism of wall crossing formulas in section 6 although the precise
relation between the two works is not clear to the author.\end{note}

The rest of this section will be a review of of \citep{GS1}. Most definitions are just copied from there. It is important to know that the word $affine$ has two different meanings here. Sometimes it is used for affine structures or charts in the sense of differential geometry while sometimes it is used for affine varieties in the sense of algebraic geometry. Anyway its meaning will be clear from the contexts.\\

\noindent {\bf Review of the complex structure instanton correction problem}\\

 We consider a decomposition of the topological space $B$ into closed cells of various dimensions which induces the affine structure with singularities on $B$. Each cell is identified as a (possibly unbounded) integral  convex polyhedron in $\mathrm{M}_{\mathbf{R}}$ and as such  maximal dimensional cells carry  induced affine charts in the interiors and we demand that cells are glued in integral affine manners. Such a manifold is called an integral polyhedral complex. The reader can find the more formal definition in \citep{GS1,GS2}. To get a singular affine structure one still need to specify the affine transformations in the normal direction of the   gluing of maximal dimensional cells. These are called fan structures. It is easy to see that we only need fan structures for each vertex $v$ (zero dimensional cell).

\begin{dfn}Let $\tau$ be a cell in the polyhedral decomposition denoted by $\mathcal{F}$ and $U_{\tau}$ be the union of interiors of all cells containing $\tau$. A fan structure along $\tau$ is a continuous map $S_{\tau}: U_{\tau} \rightarrow \mathbf{R}^{k}$ with
\begin{enumerate}
\item $S_{\tau}^{-1}(0)=\mathrm{Int} \ \tau$.
\item If $e:\tau \rightarrow\sigma$ is a inclusion morphism of cells then $S_{\tau}\mid_{Int \ \sigma}$ is an integral affine submersion onto its image.
\item The collection of cones $K_{e}:=\mathbf{R}_{\geq0}\cdot S_{\tau}(\sigma\cap U_{\tau})$ defines a fan $\Sigma_{\tau}$ in $\mathbf{R}^{k}$.
\end{enumerate}Two fan structures are equivalent if they differ by an integral linear transformation. If $\tau\subseteq \sigma\in\mathcal{F}$ the fan structure along $\sigma$ induced by $\tau$ is the composition $$U_{\sigma}\rightarrow U_{\tau}\rightarrow\mathbf{R}^{k}\rightarrow\mathbf{R}^{k}/L_{\sigma}\simeq \mathbf{R}^{l}$$where $L_{\sigma}$ is the linear span of $S_{\tau}(Int\ \sigma)$.\end{dfn}
\begin{dfn}An integral polyhedral complex $B$ of dimension $n$ is called an singular integral affine manifold with a polyhedral decomposition (also called an integral tropical manifold) if  there is a fan structure $S_{v}$ for each vertex $v$ such that:\begin{itemize}\item The support $|\Sigma_{v}|:=\bigcup_{C\in\Sigma_{v}}C$ is convex with nonempty interiors.\item If $v,w$ are vertex of a cell $\tau$, then the fan structures along $\tau$ induced by $S_{v}$ and $S_{w}$ are equivalent.\end{itemize}\end{dfn}

With fan structures around all vertices we are able to move from one maximal dimensional cell to the adjacent ones by passing through $U_{v}$. The affine charts in the interiors of maximal dimensional cells and the charts given by the fan structures at vertices provide us an integral affine atlas. However there could be nontrivial monodromies around codimensional two singular locus and this is why we have to allow singularities of the affine structure. In this way we get a singular integral affine structure together with the polyhedral decomposition. In our problem the singular integral affine structure is obtained before we make any polyhedral decomposition. So the decomposition must be chosen to be compatible with the affine structure in the sense that we recover the original affine structure.

For a vertex $v$, let $\Lambda_{v}$ be the free abelian group of integral tangent vector fields in $T_{B,v}$ and for a cell $\tau\in \mathcal{F}$ containing $v$ define
\begin{equation}\Lambda_{\tau}:=\Lambda_{v}\cap T_{\tau,v}\end{equation}The affine structure outside the singular locus induces a natural flat connection on the tangent bundle $T_{B}$. Let $\Lambda_{\mathbf{R}}$ be the local system of flat sections and $\Lambda$ be the lattice in $\Lambda_{\mathbf{R}}$ induced by the embedding $M\rightarrow M_{\mathbf{R}}$. Note that $\Lambda_{\tau}$ is the restriction of $\Lambda$. There is a monodromy representation which is the linear part of the holonomy representation. Consider two vertices $v,v^{'}$ contained in an $(n-1)$- dimensional cell $\rho$. Then the monodromy transformation associated to a loop starting from the fan structure of $v$, going through the fan structure of $v^{'}$ and back to the fan structure of $v$  is shown to have the following form \begin{equation} m\rightarrow m+\langle m, \check{d}_{\rho}\rangle m_{vv^{'}}^{\rho}\end{equation}where $m\in \Lambda_{v}$, $m_{vv^{'}}^{\rho}\in\Lambda_{\rho}$ and $\check{d}_{\rho}\in\Lambda_{\rho}^{\perp}\subseteq\Lambda_{v}^{*}$ is the primitive integral vector evaluating positively in one of the two maximal dimensional cells containing $\rho$ (we choose the one we meet first as we go around a loop).

As for the polarization $\varphi$, it is a convex multi-valued (integral) piecewise linear function on $B$.
\begin{dfn}An integral affine function on an open subset $U$ of $B$ is a continuous function $U\rightarrow\mathbf{R}$ which is integral affine on $U-\Delta$. Integral affine functions define a sheaf $Aff(B,\mathbf{Z})$. An integral piecewise linear function on $U$ is  a continuous function $\varphi: U\rightarrow\mathbf{R}$  such that if $S_{\tau}$ is the fan structure along $\tau\in\mathcal{F}$ then $\varphi\mid_{U\cap U_{\tau}}= \lambda+S_{\tau}^{*}(\bar{\varphi})$ where $\lambda$ is an (integral) affine function on $U_{\tau}$ and $\bar{\varphi}$ is  piecewise integral linear with respect to the fan $\Sigma_{\tau}$. Integral piecewise linear functions define a sheaf $PL_{\mathcal{F}}(B,\mathbf{Z})$. A polarization on $B$ is a section of the sheaf $PL_{\mathcal{F}}(B,\mathbf{Z})/ Aff(B,\mathbf{Z})$\end{dfn}
A local representative of  a polarization induces a strictly convex integral affine function on each fan $\Sigma_{\tau}$.

We can define the Legendre transform (dual) of the a triple $(B,\mathcal{F},\varphi)$ which is an involution. We refer the reader to \citep{GS2} for details. We just want to point out that the transform maps a cell  $\tau\in\mathcal{F}$ to the Newton polyhedron of $\bar{\varphi}$ where $\bar{\varphi}$ is the local representative associated to $\tau$ i.e. $\varphi\mid_{U\cap U_{\tau}}= \lambda+S_{\tau}^{*}(\bar{\varphi})$. For a  strictly convex piecewise linear function on a fan $\Sigma\in N_{\mathbf{R}}$ denoted by $\bar{\varphi}$  its Newton polyhedron is \begin{equation}\Xi:=\{x\in M_{\mathbf{R}}\mid \bar{\varphi} +x\geq 0\}\end{equation} The sum of the dimension of a cell and the dimension of its dual cell is equal to the dimension of $B$.

Gross and Siebert showed how to obtain such a triple $(B,\mathcal{F},\varphi)$ from a  toric degeneration of Calabi-Yau varieties. Toric
degeneration is a notion adapted to Gross-Siebert's
program. Roughly speaking it means the central
fiber is obtained by gluing some toric varieties
along toric strata and outside a codim
2 locus near every point there is a local model
of the degeneration given by a monomial from
an affine toric variety.   It is conjectured that a
large complex degeneration is birationally equivalent
to a toric degeneration.

 For readers' convenience  we copy the precise definition of toric degenerations from \citep{GS1} in the below. However in this paper we do not use the definition. Nor  shall we need the construction of the triple from a toric degeneration. Interested readers can find them in section 1 of \citep{GS1}. What we need instead is an existence theorem (proposition 3.3) and some explicit constructions given in section 9.\\

 \noindent {\bf Definition of Toric Degenerations}\\\\
 A $totally$ $degenerate$ $CY$-$pair$  is a pair\footnote{In this work $D$ is empty because our polyhedral complexes do not have boundaries.} $(X,D)$ where $X$ is a reduced variety\footnote{It may look strange that we did not say directly that $X$ is Calabi-Yau. But that can be deduced from the conditions given here, see theorem 2.39 and definition 4.1 in \citep{GS2}. } and $D$ is a reduced divisor satisfying the following conditions. Let $\nu: \tilde{X}\rightarrow X$ be the normalization and $C\subseteq\tilde{X}$ its conductor locus. Then $\tilde{X}$ is a disjoint union of algebraically convex (this means that there is a proper map to an affine variety) toric varieties, and $C$ is a divisor such that $[C]+\nu^{*}[D]$ is the sum of all  toric prime divisors, $\nu|_{C}: C\rightarrow\nu(C)$ is unramified and generically two-to-one, and the diagram
 \begin{diagram}[labelstyle=\textstyle]
C & & \rTo & &\tilde{X}  \\
\dTo & &   & & \dTo_{\nu} \\
\nu(C) & & \rTo & & X \\
\end{diagram}
is cartesian and cocartesian.

Let $T$ be the spectrum of a discrete valuation $\mathrm{k}$-algebra ($\mathrm{k}$ is the underlying field) with closed point $O$ and uniformizing parameter $t$. Let $\Upsilon$ be a $\mathrm{k}$-scheme and $\mathcal{D}, X$ reduced divisors  of $\Upsilon$. A $log$ $smooth$ $morphism$ $\pi: (\Upsilon,X;\mathcal{D})\rightarrow (T,O)$ is a morphism $\pi: (\Upsilon,X)\rightarrow (T,O)$ satisfying the following properties. For any $x\in\Upsilon$ there is an etale neighborhood $U$ of $x$ such that the following commutative diagram holds
 \begin{diagram}[labelstyle=\textstyle]
U & & \rTo^{\Phi} & &Spec\ \mathrm{k}[P]   \\
\dTo^{\pi|_{U}} & &   & & \dTo_{G} \\
T & & \rTo^{\Psi} & & Spec\ \mathrm{k}[\mathbf{N}] \\
\end{diagram}where $P$ is a toric monoid, $\Psi$ and $G$ are defined by sending the generator $z^{1}\in\mathrm{k}[\mathbf{N}]$ to $t$ and to a nonconstant monomial $z^{\rho}\in\mathrm{k}[P]$ respectively, $\Phi$ is etale with preimage  of the toric boundary divisor equal to the pullback to $U$ of $X\cup \mathcal{D}$.

A $toric$ $degeneration$ $of$ $CY$-$pairs$ $over$ $T$ is a flat morphism $\pi: \Upsilon\rightarrow T$ together with a reduced divisor $\mathcal{D}\subseteq \Upsilon$ with the following properties.
\begin{itemize}\item $\Upsilon$ is normal.\item The central fiber $X:=\pi^{-1}(O)$ together with $D:=\mathcal{D}\cap X$ is a totally degenerate CY-pair.\item Away from a closed subset of $\Upsilon$ of relative codimension two not containing any toric stratum of $X$, the map $\pi: (\Upsilon,X;\mathcal{D})\rightarrow (T,O)$ is log smooth.\end{itemize} There is also a notion of $formal$ $toric$ $degeneration$ $of$ $CY$-$pairs$ whose definition is omitted here.\\\\

Now the problem of instanton corrections of complex structures concerns the inverse problem: suppose we are given a triple, can we construct a toric degeneration which induces this triple? The name of the problems is perhaps a little bit confusing as it is not very clear what the instantons are and what the word "correction" means precisely. We will comment on this issue later in this section.

From the triple, one can first construct the central fiber in the following way. The upper convex hull of the graph of a local representative of the polarization $\varphi_{v}$ near a vertex $v$ defines a convex rational polyhedral cone $C_{v}\in T_{B,v}\bigoplus \mathbf{R}$. Let $P_{v}:=C_{v}\cap (\Lambda_{v}\oplus \mathbf{Z})$ be the corresponding toric monoid\footnote{A monoid is a semigroup with the identity.} (in the sense of toric geometry). In other words $$P_{v}=\{m:=(\bar{m},h) \in \Lambda_{v} \times\mathbf{Z}\mid h\geq \varphi_{v}(\bar{m})\}$$We introduce formal variables $z$ and define the ring $\mathbf{C}[P_{v}]$ to be the ring generated by $z^{m}$ with the multiplicative relations induced by additive relations in $P_{v}$. According to the relation of the triple $(B,\mathcal{F},\varphi)$ and the toric degeneration, the vertex $v$ is supposed to be associated to a zero dimensional toric stratum of the total space and an etale local model for the degeneration near this stratum is given by the map$$\mathbf{C}[t]\rightarrow \mathbf{C}[P_{v}], t\rightarrow z^{(0,1)}$$ The central fiber for this local model is the union of affine toric varieties $$\bigcup_{K}\mathrm{Spec}\ \mathbf{ C }[K\cap(\Lambda_{v}\oplus \mathbf{Z})]$$ where the union is over faces $K$ of $C_{v}$  not containing $(0,1)$ and therefore is indexed by maximal cells $\sigma$ containing $v$. We need to know how they are glued. Recall the definition of a fan structure and $\Sigma_{v}$, we can see that the coordinate ring of  this affine piece of the central fiber is given by   $\mathrm{Spec} \ \mathbf{C}[\Sigma_{v}]$ where $\mathbf{C}[\Sigma_{v}]$ is the monoid ring associated to the monoid defined by $$\bar{m}+\bar{m}^{'}:=
 \begin{array}{ll}\bar{m}+\bar{m}^{'},  & \exists K\in\Sigma : \bar{m},\bar{m}^{'}\in K \\ \infty   & \mathrm{otherwise}\\
 \end{array}
$$and we formally set $z^{\infty}=0$. For $v\subseteq\tau\in\mathcal{F}$, define $$\begin{array}{ll} \tau^{-1}\Sigma_{v}:=\{K_{e}+\Lambda_{\tau,R}\mid K_{e}\in \Sigma_{v}, e:v\rightarrow\sigma & \mathrm{factors\ through} \ \tau\}\end{array}$$
\begin{equation}V(\tau):= \mathrm{Spec} \ \mathbf{C}[\tau^{-1}\Sigma_{v}]\end{equation} Clearly $V(\tau)$ is a union of toric strata labeled by cones in $\tau^{-1}\Sigma_{v}$. We write a toric strata as $V_{e}$ for a cone $e$ in $\tau^{-1}\Sigma_{v}$. Note that $v^{-1}\Sigma_{v}=\Sigma_{v}$. It is not hard to see that if $\omega\subseteq\tau$, then there is a natural embedding $V(\tau)\rightarrow V(\omega)$. In this way we glue affine pieces according to the relation in the polyhedral decomposition and obtain a scheme\footnote{The construction  allows the possibility of composing toric automorphisms of each affine toric piece and the additional data is called an open gluing data. In this work we always choose them to be the trivial ones. Do not confuse these toric automorphisms of affine pieces of the central fiber with log (non-toric) automorphisms of thickenings of affine pieces to be defined later.} denoted by $X_{0}$ which is meant to be the central fiber of the toric degeneration to be constructed. Moreover Gross and Siebert proved that any central fiber of a toric degeneration can be obtained in this way. Details can be found in \citep{GS2}.

To (re)construct the degeneration, we need to "remember" how the central fiber is embedded in the family and this means that we need additional data on $X_{0}$. In algebraic geometry, a convenient way to encode some information of a family in the central fiber is to use the notions of logarithmic geometry and this was also the point of view taken by Gross and Siebert.  The following is the operational definition (which is shown to be equivalent to the ordinary one).
\begin{dfn}Define $LS^{+}_{pre,V(v)}:=\bigoplus_{e} O_{V_{e}}$ where the sum is over all $e: v\rightarrow \rho$ with $\dim \rho=n-1$. Let $Z$ be a closed subset of codimension 2 not containing any toric stratum of $X_{0}$.  A log smooth structure over an open subset $U\in V(v)\setminus Z$ is a rational section
 $(f_{e})$ of $LS^{+}_{pre,V(v)}$ whose zeroes and poles do not contain any toric stratum and which satisfies the following condition:\begin{equation}\Pi_{i=1}^{l}\check{d}_{\rho_{i}} \otimes f_{e_{i}} \mid _{V_{h}}=1\end{equation} where $(\rho_{i}), 1\leq i\leq l$ is a cyclic ordering of all  $(n-1)$-cells containing an $(n-2)$-cell $\tau$ of the polyhedral decomposition which contains $v$. $V_{h}$ is the strata corresponding to the cone $h$ given by $h: v\rightarrow \tau$. $\check{d}_{\rho_{i}}\in \Lambda_{\rho_{i}}^{\perp}\subseteq \Lambda_{v}^{*}$ are generators compatible with the cyclic ordering. In (22) we treat the first factor additively and the second factor multiplicatively. The log smooth structure is positive if it is a section of $\bigoplus_{e}O_{V_{e}}^{\times}$ which extends across $Z$ as a section of $\bigoplus_{e}O_{V_{e}}$. The canonical minimal choice of $Z$ is the union over codimensional one cells $\rho$ of the vanishing locus of $f_{v\rightarrow\rho}$. This defines a log smooth structure locally over a chart. To define it globally we also need a compatible condition (change of vertex formula) associated to the change of charts because of nontrivial monodromies. Let $v,v^{'}$ be two vertices of an $(n-1)$-dimensional cell $\rho$ and $e,e^{'}$ are $e:v\rightarrow\rho$ and $e^{'}:v^{'}\rightarrow\rho$ respectively. The condition is \begin{equation}f_{e^{'}}=z^{m^{\rho}_{v^{'}v}}f_{e}\end{equation}To emphasize the dependence on the vertex, later we may write the section $f_{e}$ as $f_{\rho,v}$.\end{dfn}

 The geometric meaning of a log smooth structure is that at a generic point of $V_{e}$ a local model of the toric degeneration with central fiber $X_{0}$ is given  by the equation $$zw-f_{e}t^{p}=0$$ where $p$ is the integral length of the dual cell (under the Legendre transform) of $\rho$ and $z,w$ are two variables in the monoid associated to the two generators of $\Lambda_{v}/\Lambda_{\rho}$. Hence a log smooth structure does tell us something about the embedding of the central fiber.

 Now we can state the main theorem of Gross and Siebert.
 \begin{thm} Any locally rigid, positive, pre-polarized toric log Calabi-Yau variety with proper irreducible components arises from a formal toric degeneration of Calabi-Yau varieties\footnote{ The theorem is actually true more generally  for pairs, see \citep{GS1}. We only state it for  varieties instead of pairs because the central fiber is only viewed as a variety in this paper.}.\end{thm}

 A pre-polarized toric log Calabi-Yau variety is a scheme $X_{0}$ constructed as above from a triple $(B,\mathcal{F},\varphi)$ together with the polarization $\varphi$ and a log smooth structure. Local rigidity is a technical condition which guarantees uniqueness in certain constructions. In all our examples, explicit constructions will be obtained and  it is not necessary to check this technical condition. In fact we do not use this theorem at all because our spaces, namely Hitchin's moduli spaces, are  noncompact. Moreover we need more than just an existence statement. What is important is the proof of this theorem. We will build the explicit degenerations of a Hitchin's moduli space  by partly following the proof of this theorem. We do not assume the cells of the polyhedral decomposition are bounded.

 The proof of this theorem  is very complicated. The basic idea is that we  deform affine pieces\footnote{Here the word 'affine' is used in the sense of algebraic geometry.} of the central fiber before gluing. This will introduce inconsistencies in general. So we have to modify the manners of the gluing by composing certain (auto)morphisms. These (auto)morphisms will be associated to some codimension one subsets and they must satisfy some consistency conditions as well.

 To be more precise, let us suppose that we want to find a formal deformation of the central fiber over $ \mathrm{Spec} \mathbf{C}[[t]]$ where $t$ is the deformation parameter. We first define canonical $k$-th order thickenings of affine pieces of $X_{0}$ in the following way. For a maximal dimensional cell $\sigma$ in the polyhedral decomposition, a vertex $v\in\sigma$ and an element $m=(\bar{m},h)\in P_{v}$ we define $\mathrm{ord}_{\sigma}(m)$ to be the vanishing order of $z^{m}$ over the strata $V_{v\rightarrow\sigma}$. More precisely, since the Legendre dual $\check{v}$ of a vertex $v$ is a maximal dimensional cell whose vertices are $\check{\sigma}: -\lambda_{\sigma}\in \Lambda_{v}^{*}$ \footnote{See the definition of Newton polyhedrons.} where $ \lambda_{\sigma}$ is the linear function defined by the local representative $\varphi_{v}$ of the polarization on the maximal dimensional cell $\sigma\supseteq v$, we may naturally view  $m$ as an element of the stalk of the sheave of integral affine function $\mathrm{Aff}(\check{B},\mathbf{Z})$\footnote{$\check{B}$ is the Lengdre dual of $B$ and is isomorphic to $B$ topologically. $\Gamma(Int \ \check{v},\mathrm{Aff}(\check{B},\mathbf{Z}))=(\Lambda^{*}_{v})^{*}\oplus\mathbf{Z}=\Lambda_{v}\oplus\mathbf{Z}\ni m$.} over a point $x\in\sigma\subseteq B$ and define its order as \begin{equation}\mathrm{ord}_{\sigma}(m)=\langle\bar{m},-\lambda_{\sigma}\rangle+h\end{equation}where $\bar{m}$ is the image of $m$ under the projection $\mathrm{Aff}(\check{B},\mathbf{Z})_{\check{\sigma}}\rightarrow \Lambda_{\sigma}$ which is induced by the natural exact sequence\begin{equation}0\rightarrow\mathbf{Z}\rightarrow\mathrm{Aff}(\check{B},\mathbf{Z})\rightarrow\Lambda\rightarrow0\end{equation}$m$ is called an exponent at $x$. For a maximal dimensional cell $\sigma$ we define its exponent to be the exponent over any of its interior point away from singularities and this is a well defined element in $\mathrm{Aff}(\check{B},\mathbf{Z})_{\check{\sigma}}$ (note that $\check{\sigma}$ is just a point). The exponent at $\sigma$ is denoted by $m_{\sigma}$ (or just $m$). The order is invariant under the monodromy action on $m$. For a subset $A$ contained in a cell and containing $x$ the order is \begin{equation}\mathrm{ord}_{A}(m):=\max\{\mathrm{ord}_{\sigma}(m)\mid A\subseteq\sigma\in\mathcal{F}_{\max}\}\end{equation}Suppose $\sigma,\sigma^{'}$ are two maximal dimensional cells containing $v$ with nonempty intersection and $m$ is an exponent on $\sigma$ (i.e. on an interior point $x\in\sigma$). Let $m^{'}$ be the parallel transport of $m$ to $\sigma^{'}$ induced by a parallel transport of $\Lambda$, we define the order of $m$ on $\sigma^{'}$ which does not contain $x$ as \begin{equation}\mathrm{ord}_{\sigma^{'}}(m):=\mathrm{ord}_{\sigma^{'}}(m^{'})\end{equation}Note that due to the monodromy to define the parallel transport we must choose a vertex  contained in $\sigma\cap\sigma^{'}$ and use its fan structure.  The definition does not depend on the choice of $v$ because the definition of order is invariant under local monodromy.

 We will need the following proposition in \citep{GS1} later.
 \begin{pro}Let $m$ be an exponent at $x$ and $\tau$ be the minimal cell containing $x$. If $\sigma^{+}$ and $\sigma^{-}$ are maximal dimensional cells containing $\tau$ such that the corresponding maximal cones in $\Sigma_{\tau}$ contain $\bar{m}$ and $-\bar{m}$ respectively, then
 $$\mathrm{ord}_{\sigma^{-}}(m)=\mathrm{max}\{\mathrm{ord}_{\sigma}(m)\mid\sigma\in\mathcal{F}_{max},\tau\subseteq\sigma\}$$
$$\mathrm{ord}_{\sigma^{+}}(m)=\mathrm{min}\{\mathrm{ord}_{\sigma}(m)\mid\sigma\in\mathcal{F}_{max},\tau\subseteq\sigma\}$$\end{pro}
Basically this proposition tells us that if we propagate an exponent $m$ in the direction $-\bar{m}$ its order increases.

 Let us continue and define a generalization of $P_{v}$ for $\omega\subseteq \sigma\in\mathcal{F_{\max}}$
 \begin{equation}P_{\omega,\sigma}:=\{m\in \mathrm{Aff}(\check{B},\mathbf{Z})_{\check{\sigma}}\mid \forall \sigma^{'}\in\mathcal{F_{\max}}, \omega\subseteq \sigma^{'}: \mathrm{ord}_{\sigma^{'}}(m)\geq 0\}\end{equation}For any choice of local representative $\varphi_{v}$ at $v$ of the polarization $\varphi$ and any maximal dimensional cell $\sigma$ containing $v$, we have $P_{v,\sigma}=P_{v}$ canonically. We also define \begin{equation}P_{x}:=\{m\in\mathrm{Aff}(\check{B},\mathbf{Z})_{x}\mid \forall\sigma\in\mathcal{F}_{max},x\in\sigma: \mathrm{ord}_{\sigma}(m)\geq0\}\end{equation} If $\omega\ni x$, $P_{\omega,\sigma}\simeq P_{x}$. As a result if $x,x^{'}\in Int(\omega)-\Delta$ then any maximal dimensional cell $\sigma$ containing $x$ induces an isomorphism $P_{x}\simeq P_{x^{'}}$.

 Let $\sigma^{'}$ be another maximal dimensional cell containing $\omega$ then parallel transport via the fan structure of a vertex $v\in\sigma\cap\sigma^{'}$ induces an isomorphism $P_{\omega,\sigma}\simeq P_{\omega,\sigma^{'}}$.

 The canonical $k$-th order thickenings of affine pieces of $X_{0}$ are
 \begin{equation}R_{g,\sigma}^{k}:=(k[P_{\omega,\sigma}]/I_{g,\sigma}^{>k})_{f_{g,\sigma}}\end{equation}where $g$ is an inclusion morphism $g: \omega\rightarrow \tau, \omega,\tau\in\mathcal{F}$ and $\sigma$ is a maximal dimension cell containing $\tau$. The lower subscript outside the bracket means we take the localization of the quotient ring inside the bracket with respect to elements $f_{g,\sigma}$ and as the symbol has suggested $f_{g,\sigma}$ are constructed from the log smooth structure.  For a codimensional one cell $\rho$ containing $\tau$ denote the composition $v\rightarrow\omega\rightarrow\tau\rightarrow\rho$ by $e_{\rho}$ and denote the associated codimensional one cone in the fan $\Sigma_{v}$ by $K_{e_{\rho}}\in\Sigma_{v}$. By its definition $$f_{e_{\rho}}=\sum_{\bar{m}\in K_{e_{\rho}}\cap\Lambda_{v}}f_{e_{\rho},\bar{m}}z^{\bar{m}}$$ It has a canonical lift to $\mathbf{C}[P_{v,\sigma}]$ which is \begin{equation}\sum_{m\in P_{v,\sigma}}f_{e_{\rho},m}z^{m}\end{equation} where $f_{e_{\rho},m}=f_{e_{\rho},\bar{m}}$ if $\bar{m}\in K_{e_{\rho}}$ and $\mathrm{ord}_{\rho}(m)=0$ and $f_{e_{\rho},m}=0$ otherwise. $f_{g,\sigma}$ is defined to be \begin{equation}f_{g,\sigma}:=\prod_{\rho\supseteq\tau}\sum_{m\in P_{v,\sigma}}f_{e_{\rho},m}z^{m}\end{equation}If $\tau$ is maximal dimensional we define $f_{g,\sigma}=1$. $I_{g,\sigma}^{>k}$ denotes the ideal generated by $P_{g,\sigma}^{>k}$ which is the set of those
 elements in the  monoid $P_{\omega,\sigma}$  such that $\mathrm{ord}_{\tau}(m)> k$. $\mathrm{Spec}R_{g,\sigma}^{k}$ is a thickening of the complement of $(\bigcup_{\rho\supseteq\tau}V(f_{\rho,v}), v\subseteq\omega)$ in  $\mathrm{Spec}(\mathbf{C}[P_{\omega,\sigma}]/I^{>0}_{g,\sigma})\subseteq V(\omega)$ where $V(f_{\rho,v})$ is the locus defined by the vanishing of the   log smooth structure $f_{\rho,v}$, $v\subseteq\omega$. If $g:\omega\rightarrow\tau, g^{'}:\omega^{'}\rightarrow\tau^{'}$ and $\omega\subseteq\omega^{'}, \tau\supseteq\tau^{'}$ then there is a canonical homomorphism $R_{g,\sigma}^{k}\rightarrow R_{g^{'},\sigma}^{k}$.

 As said before, when one tries to glue these deformations of affine pieces one encounters inconsistencies (which are not necessarily associated to the existence of nontrivial monodromy around singularities of the affine manifold) and has to compose some (auto)morphisms called log (auto)morphisms.
 \begin{dfn}A log ring is a ring $R$ together with a monoid homomorphism $\alpha: P\rightarrow (R,\cdot)$ from a monoid $P$. A log morphism between two log rings $\alpha: P\rightarrow R$ and $\alpha^{'}: P^{'}\rightarrow R^{'}$ is a triple $(\psi: R\rightarrow R^{'}, \beta: P\rightarrow P^{'}, \theta:P\rightarrow(R^{'})^{\times})$ satisfying $$\psi\circ\alpha= \theta\cdot(\alpha^{'}\circ\beta)$$In our case, the log ring is $\alpha: P_{\omega,\sigma}\rightarrow R_{g,\sigma}^{k}$ where $\alpha$ send $m$ to $z^{m}$. $\beta$ will be fixed and $\theta$ will factor through $P_{\omega,\sigma}\rightarrow\Lambda_{\sigma}$. So we will use $\theta$ to denote the homomorphism $\Lambda_{\sigma}\rightarrow R_{g,\sigma}^{k}$. We also use $\bar{\theta}$ instead of $\psi$ to denote the ring homomorphism. So\begin{equation}\bar{\theta}(z^{m})=\theta(\bar{m})\cdot z^{\beta(m)}\end{equation}The composition of two log morphisms $\theta_{1},\theta_{2}$ is $$(\theta_{1}\circ\theta_{2})(m)=\theta_{1}(m)\cdot\bar{\theta}_{1}(\theta_{2}(m))$$\end{dfn}Log (auto)morphisms in our problem are associated to codimensional one subsets. There are two kinds of them. The first kind is called slabs which are codimensional one polyhedral subsets  of codimensional one cells of the polyhedra decomposition together with higher order corrections ("attached" to slabs) of the gluing.
  \begin{dfn}A slab is a convex rational $(n-1)$-dimensional polyhedral subset $\mathbf{b}$ of a codimensional one (i.e. $(n-1)$-dimensional) cell $\rho_{\mathbf{b}}\in\mathcal{F}$ together with functions \begin{equation}f_{\mathbf{b},x}:=\sum_{m\in P_{x}.\bar{m}\in\Lambda_{\rho_{\mathbf{b}}}}c_{m}z^{m}\in\mathbf{C}[P_{x}] \end{equation}for each $x\in \mathbf{b}-\Delta$ satisfying\begin{itemize}\item Let $x,x^{'}\in \mathbf{b}-\Delta$ and let $v=v[x]$ be the unique vertex in the same  connected component of $\mathbf{b}-\Delta$ of $x$. Give the analogous meaning to $v^{'}=v^{'}[x^{'}]$. Let $\Pi$ be the parallel transport  $\mathbf{C}[P_{x}]\rightarrow\mathbf{C}[P_{x^{'}}]$ along a path inside the closure of $U_{\rho_{\mathbf{b}}}$ avoiding $\Delta$ from $x$ to $x^{'}$, then \begin{equation}f_{\mathbf{b},x^{'}}=z^{m^{\rho_{\mathbf{b}}}_{v^{'}v}}\Pi(f_{\mathbf{b},x})\end{equation}Consider the two adjacent maximal dimensional cells containing $\rho_{\mathbf{b}}$. Let $\beta$ and $\beta^{'}$ be parallel transports of exponents from one maximal dimensional cell to another via the fan structure of $v$ and $v^{'}$ respectively. Then there is a relation $$\beta^{'}(m)=\beta(m)+\pi(\bar{m})\cdot m^{\rho_{\mathbf{b}}}_{v^{'}v}$$where $\pi(\bar{m})$ is defined in (38) and $m^{\rho_{\mathbf{b}}}_{v^{'}v}\in \Lambda_{\rho_{\mathbf{b}}}$. This relations defines $m^{\rho_{\mathbf{b}}}_{v^{'}v}$. Note that $m^{\rho_{\mathbf{b}}}_{v^{'}v}$ is $m^{\rho_{\mathbf{b}}}_{v^{'}v}$ defined in (19).\item Let $f_{e}$ be the data associated to $e: v\rightarrow\rho_{\mathbf{b}}$ with $v=v[x]$ defining the log smooth structure. Let $\Pi$ be the parallel transport $\mathbf{C}[P_{x}]\rightarrow\mathbf{C}[P_{v}]$. Then \begin{equation}f_{e}=\Pi(\sum_{m\in P_{x}, \mathrm{ord}_{\mathbf{b}}(m)=0}c_{m}z^{m})\end{equation}\end{itemize} \end{dfn}$f_{e}$ in (36) is understood as the lift of the log smooth data to $\mathbf{C}[P_{v,\sigma}]$ defined by (31). Comparing the definition of a slab and the definition of a log smooth structure, we see that a log smooth structure induces a slab structure on each codimensional one cell in the polyhedral decomposition such that the order of each term of the slab function  is zero.

  The other kind is called walls and unlike slabs they are codimensional one polyhedral subsets of maximal dimensional cells.

  \begin{dfn}A wall is a convex rational $(n-1)$-dimensional polyhedral subset $\mathbf{p}$ of a maximal dimensional cell $\sigma_{\mathbf{p}}$ (We require that $\mathbf{p}$  has nonempty intersection with the interior of $\sigma_{\mathbf{p}}$) together with\begin{itemize}\item An $(n-2)$-dimensional face $\mathbf{q}\subseteq\mathbf{p}$ called the base of $\mathbf{p}$.\item An exponent $m_{\mathbf{p}}$ on $\sigma_{\mathbf{p}}$ such that $\mathrm{ord}_{\sigma_{\mathbf{p}}}(m_{\mathbf{p}})>0$, $$\mathbf{p}=(\mathbf{q}-\mathbf{R}_{\geq0}\bar{m}_{\mathbf{p}})\cap\sigma_{\mathbf{p}}$$ and $m_{\mathbf{p},x}\in P_{x}$ for all $x\in\mathbf{p}-\Delta$.\item A number $c_{\mathbf{p}}$ from which we define a function \begin{equation}f_{\mathbf{p},x}:=1+c_{\mathbf{p}}z^{m_{\mathbf{p}},x}\end{equation}\end{itemize} The slides and the top of the wall are defined as $$Slides(\mathbf{p}):=(\partial\mathbf{q}-\mathbf{R}_{\geq0}\bar{m}_{\mathbf{p}})\cap\sigma_{\mathbf{p}}$$ $$Top(\mathbf{p}):=closure(\partial\mathbf{p}-(\mathbf{q}\cup Slides(\mathbf{p})))$$\end{dfn}

  Later in this section we will produce log morphisms from $f_{\mathbf{b},x}$ and $f_{\mathbf{p},x}$.

  \begin{dfn}A locally finite collection of slabs and walls $\aleph$ together with a polyhedral decomposition $\mathcal{F}_{\aleph}$ of the union of  supports of elements of $\aleph$ gives us a $structure$ (also denoted by $\aleph$) if the following conditions are satisfied. \begin{itemize}\item Every codimensional one cell in the original polyhedral decomposition $\mathcal{F}$ (not $\mathcal{F}_{\aleph}$) is considered as the support of a slab of the $structure$. Each slab in the $structure$ defines a codimensional one cell in $\mathcal{F}_{\aleph}$. \item Define  a Gross-Siebert chamber to be  the closure of a connected component of the complement of the union of supports of elements of $\aleph$. Then every Gross-Siebert chamber is convex and its interior is disjoint from any wall. Sometimes we also use the name "chamber" instead of Gross-Siebert chamber. \item Any wall is a union of elements of  $\mathcal{F}_{\aleph}$.\item Each maximal cell in $\mathcal{F}$ contains only finitely many slabs or walls.\end{itemize} \end{dfn}To make sure that after adding slabs and walls the gluing is indeed consistent one must show that going along a loop around a codimensional two cell of the $structure$ the ordered composition of these log morphisms is trivial. If this is true, then we say the $structure$ is consistent to order $k$.

 Let us describe the consistency condition more carefully.  For a chamber $\mathrm{u}$ there is a unique $\sigma_{\mathrm{u}}\in\mathcal{F_{\max}}$ with $\mathrm{u}\subseteq\sigma_{\mathrm{u}}$. So for each pair $(g, \mathrm{u})$ such that  $g$ is a morphism $g:\omega\rightarrow \tau$ in the polyhedral decomposition and $\tau\in\sigma_{\mathrm{u}}$ we have the rings $R_{g,\sigma_{\mathrm{u}}}^{k}$. A joint $j$ of the $structure$ is a codimensional two cell of the $structure$. In the normal (two dimensional) space of $j$ codimensional one cells of the $structure$ containing $j$ are rays cyclicly ordered and numbered. Hence so are the chambers around $j$. Using $\theta_{i}$ to denote the log automorphism from $R_{g,\sigma_{\mathrm{u_{i}}}}^{k}$ to $R_{g,\sigma_{\mathrm{u_{i+1}}}}^{k}$ to be defined below, the consistency condition is
 \begin{dfn}A $structure$ is consistent to order $k$ if it is consistent to order $k$ for all joints. It is consistent  to order $k$ for the joint $j$ if for any $g:\omega\rightarrow \tau$ with $j\cap \omega\neq \emptyset$ and $\tau\in \sigma_{j}$ the composition $$\theta_{k}^{j}:= \theta_{l}\circ\cdots\circ\theta_{1}$$ is the identity. Here $l$ is the number of codimensional one cells in the $structure$ containing $j$ and $\sigma_{j}$ is the minimal cell in $\mathcal{F}$ containing $j$.\end{dfn}\begin{note}The minimal cell $\sigma_{j}$ in the polyhedral decomposition (instead of the $structure$) containing $j$ has codimensional at most two. We call the joint $j$ a codimensional zero, one or two joint according to the codimension of that minimal cell\footnote{But the joint itself is always codimensional two as a subset.}.  \end{note}

 Let us specify the definition of the log morphisms. We glue $R_{g,\sigma_{\mathrm{u}}}^{k}$ together according to the inclusion relations of the $structure$. The gluing process can be decomposed into the following two types of basic gluing "morphisms" $$(g: \omega\rightarrow\tau,\mathrm{u})\rightarrow (g^{'}:\omega^{'}\rightarrow\tau^{'},\mathrm{u}^{'})$$ i.e. we glue $R_{g,\sigma_{\mathrm{u}}}^{k}$ and $R_{g^{'},\sigma_{\mathrm{u}^{'}}}^{k}$ via:\begin{itemize}\item Change of strata: $\omega\subseteq\omega^{'}, \tau\supseteq\tau^{'}, \mathbf{u}=\mathbf{u}^{'}$.\item Change of chambers: $\omega=\omega^{'}, \tau=\tau^{'},\mathrm{dim}\ \mathbf{u}\cap\mathbf{u}^{'}=n-1, \omega\cap\mathbf{u}\cap\mathbf{u}^{'}\neq\emptyset$.\end{itemize}

 For changes of strata the log morphisms are defined to be  the trivial ones induced by the canonical map $\beta:P_{\omega,\sigma_{\mathrm{u}}}\rightarrow P_{\omega^{'},\sigma_{\mathrm{u}^{'}}}$. For changes of chambers, we distinguish two cases.\begin{itemize}\item$\sigma_{\mathrm{u}}=\sigma_{\mathrm{u}^{'}}=\sigma$. In this case there exists  an $(n-1)$-dimensional cell $b$ of the $structure$ contained in the intersection of $\mathbf{u}$ and $\mathbf{u}^{'}$ and having nonempty intersection with $\omega$. Let $\mathbf{p}_{i}, 1\leq i\leq r$ be walls in the $structure$ containing $b$ and let $f_{i}$ be the image of $f_{\mathbf{p}_{i},x}$ in $R_{g,\sigma}^{k}$. The tangent space of the intersection of two chambers $T_{\mathbf{u}\cap\mathbf{u^{'}}}$ is an $(n-1)$-dimensional space and let $\pi: \Lambda_{\sigma}\rightarrow\mathbf{Z}$ be the epimorphism contracting $T_{\mathbf{u}\cap\mathbf{u^{'}}}\cap\Lambda_{\sigma}$ and evaluating positively on vectors pointing from $\mathbf{u}$ to $\mathbf{u}^{'}$. Then we define the log morphism by $\beta=\mathrm{Id}: P_{\omega,\sigma}\rightarrow P_{\omega,\sigma}$,\begin{equation}\theta: \bar{m}\rightarrow(\prod_{i=1}^{r}f_{i})^{-\pi(\bar{m})}\end{equation}It is an log automorphism of $R_{g,\sigma_{\mathrm{u}}}^{k}$.\item$\sigma_{\mathrm{u}}\neq\sigma_{\mathrm{u}^{'}}$. If this is the case then the intersection of the two chambers must be a codimensional one cell of the polyhedral decomposition. Denote by $b$ the codimensional one cell of the $structure$ contained in the intersection and having nonempty intersection with $\omega$. Let $\mathbf{b}$ be the unique slab whose support is $b$. Let $x\in(\omega\cap b)-\Delta$ and $e: v\rightarrow\omega$ where $v=v[x]$. We define a log morphism by defining $\beta$ as the parallel transport through $v$ and \begin{equation}\theta: \bar{m}\rightarrow(f_{\mathbf{b},x})^{-\pi(\bar{m})}\end{equation}where $\pi$ is the epimorphism $\Lambda_{\sigma_{\mathbf{u}}}\rightarrow \mathbf{Z}$ with kernel $\Lambda_{\rho}$ and is positive on vectors pointing from $\mathbf{u}$ to $\mathbf{u}^{'}$. $f_{\mathbf{b},x}$ is viewed as an element of $R_{g^{'},\sigma_{\mathrm{u}^{'}}}^{k}$ via the fan structure at $v$. The condition (35) guarantees that the log morphism does not depend on the choice of $x$ (or equivalently the choice of $v$).\end{itemize}

 Changes of strata commute with changes of chambers. The consistency condition in codimension two formulated above guarantees that if a gluing morphism has two decompositions into basic gluing morphisms then the ordered compositions of log morphisms associated to the two decompositions are the same.

The idea of using log morphisms to
correct gluing construction was due to Kontsevich
and Soibelman in dimension two in the
somewhat different framework of non-archimedean
analytic spaces. Gross and Siebert's construction is adapted to
the theory of logarithmic geometry and works
for any dimensions. Furthermore explicit calculations
of the degeneration can be done.

 We are not done yet because what we really want is a formal degeneration (which in nice cases can be algebraized into a genuine deformation over rings of finite type), i.e. we want to let $k$ go to infinity. So we must  show the consistent $structure$ in order $k$ is compatible to the consistent $structure$ in order $k+1$.
 \begin{dfn}Two $structures$ $\aleph$ and $\aleph^{'}$ are compatible to order $k$ if \begin{itemize}\item If $(\mathbf{p},m,c)$ is a wall in $\aleph$ with $c\neq0, \mathrm{ord}_{\sigma_{\mathbf{p}}}(m)\leq k$, then it is a wall of $\aleph^{'}$ and vice versa.\item If $\mathbf{b}$ and $\mathbf{b^{'}}$ are slab in $\aleph$ and $\aleph^{'}$ respectively and $x\in (Int(b)\cap Int(b^{'}))-\Delta$, then $$f_{\mathbf{b},x}=f_{\mathbf{b^{'}},x}\ \mathrm{mod}\ t^{k+1}$$\end{itemize}If we have a sequence of  $structure$ $\aleph_{k}, k\geq0$ such that  for any $k$ $\aleph_{k}$ is consistent to order $k$ and $\aleph_{k},\aleph_{k+1}$ are compatible to order $k$, then we say we have a compatible system of consistent $structures$.\end{dfn}By the work described above the solution of the complex structure instanton correction problem, i.e. the solution of the (re)construction problem of the degeneration  (theorem 3.1) is a  corollary of the existence of a compatible system of consistent $structures$. In fact
  the following proposition is proved in \citep{GS1}.
  \begin{pro}For a compatible system of consistent $structures$ which is inductively constructed starting from a positive log smooth structure (see the next proposition for its meaning) one can construct a formal toric degeneration of Calabi-Yau varieties with central fiber $X_{0}$  which induces the triple $(B,\mathcal{F},\varphi)$ and the log smooth structure.\end{pro}

 The existence of a compatible system of consistent $structures$ is proved by induction on $k$. When $k=0$  the following proposition is easy to verify by checking definitions.
 \begin{pro}A $structure$ consistent to order zero containing only slabs defines a  log smooth structure and vice versa.\end{pro}

 In fact we can use (36) to define $f_{e}$. Then the consistency in order zero becomes (22) while the condition (35) becomes (23).

  In this sense, a log smooth structure is the initial data which determines all higher order corrections by consistency such that finally a formal deformation can be constructed. In general each order could introduce walls and slabs. The existence of infinitely many walls/slabs is a generic phenomenon but at each order there are only finitely many of them.\\

The construction looks  formidable.  Let us describe the scenario in the real two dimension with only one joint in the polyhedral decomposition \citep{GPS}. It is also important for understanding our examples.

 Let $M:=\mathbf{Z}^{2}$, $N:= \mathrm{Hom}(M,\mathbf{Z})$ and define the group ring $\mathbf{C}[M]\ni z^{m},m\in M$ $$x:=z^{(1,0)}, y:=z^{(0,1)}$$ A log derivation $\xi$ is an element of $(\lim_{\leftarrow}\mathbf{C}[M]\otimes \mathbf{C}[[t]]/t^{k})\bigotimes N$. It is of the form $a\otimes n$ where $n\in N$, we write it as $a\partial_{n}$ and it induces an ordinary derivation $$(a\partial_{n})(z^{m}):=a\langle m,n\rangle z^{m}$$
 \begin{dfn} A ray or a line in $\mathbf{R}^{2}\simeq M_{\mathbf{R}}$ is a pair $(l, f_{l})$ where $l$ is  either a ray with integral slope $m_{0}\in M$  or a line with integral slope $m_{0}\in M$ such that $$f_{l}\in \lim_{\leftarrow}\mathbf{C}[z^{m_{0}}]\otimes \mathbf{C}[[t]]/t^{k}$$ and $f_{l}=1$ mod $z^{m_{0}}t$.
 A universal scattering diagram $D$ is a set of lines and rays such that for each $k$ there are only finitely many pairs $(l, f_{l})$ with $f_{l}\neq 1$ mod $t^{k}$.\end{dfn}
\begin{note} $D$ is called a scattering diagram in \citep{GPS}. We choose  to call it a universal scattering diagram here because it is a "union" of scattering diagrams defined in \citep{GS1}.\end{note}

 For a generic closed loop on the plane we define an ordered product (composition) as follows. Choose an orientation of the loop and we can order the intersections of the loop with rays and lines (viewed as two rays) as $(l_{i})$ with the loop meeting $l_{i}$ before $l_{j}$ if $i<j$. Now define an automorphism \begin{equation}\theta_{i}:=\exp(\log(f_{l_{i}})\partial_{n_{0}})\end{equation}where $n_{0}\in N$ is the primitive normal vector of $l_{i}$ positively oriented along the loop. We can make the composition of these automorphisms in the above defined order. We have the following simple but basic lemma of Kontsevich and Soibelman.
 \begin{thm} Let $\theta_{i}$ be the log automorphism associated to the ray $(l_{i})$ of a universal scattering diagram $D$. There exists a universal scattering diagram $S(D)$ containing $D$ such that the new one is obtained from the old one by adding only rays and such that the ordered product of automorphisms $$\theta_{s}\circ\theta_{s-1}\circ\cdots\circ\theta_{1}$$ is the identity in $S(D)$ for any loop. Moreover there is only one such universal scattering diagram which is minimal in the sense that it does not contain trivial rays or lines and does not have two rays or lines with the same support.\end{thm}

  The theorem is proved by induction on $k$. Let $D_{0}:=D$. Suppose we have already built a universal  scattering diagram $D_{k-1}$ with $$\theta_{\gamma, D_{k-1}}=\mathrm{Id} \ \mathrm{mod} \ t^{k}$$ for any closed loop $\gamma$. We want to build $D_{k}$ such that $\theta_{\gamma, D_{k}}=\mathrm{Id} \ \mathrm{mod} \ t^{k+1}$.

  To this end, we consider $D_{k-1}^{'}$ which consists of all rays and lines $\sigma$ in $D_{k-1}$ with $f_{\sigma}\neq 1\ \mathrm{mod} \ t^{k+1}$. Suppose $p$ is a singularity of $D_{k-1}^{'}$. By this we means that it is either an initial point of a ray or an intersection point of rays/lines. Take a  closed simple loop around $p$ small enough to contain no other singularities. By the definition of $D_{k-1}^{'}$,$$\theta_{\gamma, D_{k-1}}=\theta_{\gamma, D_{k-1}^{'}}\ \ \mathrm{mod}\ t^{k+1}$$The problem is local. By inductive assumption we can assume that $\theta_{\gamma, D_{k-1}^{'}}$ is expanded as $$\theta_{\gamma, D_{k-1}^{'}}=\exp(\sum c_{i}z^{m_{i}}\partial_{n_{i}})$$ This is a finite sum. Set $D[p]:=\{(p+\mathrm{R}_{\geq0}m_{i}, 1\pm c_{i}z^{m_{i}})\}$ and choose the sign such that its associated automorphism is $\exp(- c_{i}z^{m_{i}}\partial_{n_{i}})$ modulo $t^{k+1}$. Then clearly $$\theta_{\gamma, D_{k-1}\cup D[p]}= \mathrm{Id}\ \ \mathrm{mod}\ t^{k+1}$$ Now $D_{k}$ is defined to be $$D_{k}:= D_{k-1}\cup \bigcup_{p}D[p]$$Finally define $S(D)$ to be the union of all $D_{k}$'s. \begin{note}For each $k$, there are only finitely many singularities in $D_{k-1}^{'}$  because $D_{k-1}^{'}$ itself is a finite set. So at each step one only adds finitely many rays.\end{note}
  \begin{note}The proof makes it clear that the possibly infinite ordered compositions should be understood in the sense of taking successive truncations and the projective limit of them. \end{note}
  \begin{dfn}If the ordered compositions are all identities modulo $k$ for all loops we say the universal scattering diagram is consistent to order $k$. If it is consistent for all $k$  we say it is a consistent universal scattering diagram.\end{dfn}

 Here the affine base $B$ is the plane. Slabs and walls are rays and lines in a universal scattering diagram. In fact, initially the polyhedral decomposition provides some lines and rays (in the ordinary sense) which are codimensional one cells. We then decode the definition of the log morphisms and find that they give functions $f_{i}$ as above attached to these lines and rays and this produces a universal scattering diagram. Now the consistency conditions of the $structure$ consisting of slabs and walls become the requirement that the ordered product around any loop is the identity. As explained in the remark 3.5 it is understood  in the projective limit sense and therefore is actually a compatible system of consistency conditions. So a consistent universal sacttering diagram  is really the union of elements of a compatible system  of  consistent $structures$.  \\

 The notion of scattering diagrams of rays and lines actually works in any dimensions simply because the consistency condition is a codimensional two condition.
\begin{dfn}Let $j$ be a joint and $\sigma_{j}$ be the minimal cell in the polyhedral decomposition containing $j$. For a vertex $v\in\sigma_{j}$ we consider the normal space $$Q_{j,\mathbf{R}}^{v}:=\Lambda_{v,\mathbf{R}}/\Lambda_{j,\mathbf{R}}$$Let $\bar{\bar{m}}$ be the image of $m$ in $Q_{j,\mathbf{R}}^{v}$ via projection. If $\tau$ is a cell containing $j$ then let $\bar{\bar{\tau}}\in Q_{j,\mathbf{R}}^{v}$ be the image of the tangent wedge of $\tau$ along $j$. By a cut in $Q_{j,\mathbf{R}}^{v}$ we mean a half line starting from the origin which is contained in $\bar{\bar{\rho}}$ for some $\rho$ which is a codimensional one cell containing $j$. $Q_{j,\mathbf{R}}^{v}$ is divided into chambers by cuts. A ray is a triple $(\mathbf{t},m_{\mathbf{t}},c_{\mathbf{t}})$ (sometimes also denoted simply as $\mathbf{t}$) where $\mathbf{t}$ is a one-dimensional rational cone generated by $\bar{\bar{q}}$ for $q\in(\Lambda_{v}-\Lambda_{j})$. $m_{\mathbf{t}}$ is an exponent on a maximal dimensional cell $\sigma$ such that $\pm\bar{\bar{m_{\mathbf{t}}}}\in \mathbf{t}\cap\bar{\bar{\sigma}}$ and $m_{\mathbf{t}}\in P_{x}$ for all $x\in (j-\Delta)$. $c_{\mathbf{t}}$ is a constant.

A scattering diagram for $j$ at the vertex $v$ is the following data \begin{itemize}\item A cell $\omega\in\mathcal{F}$ whose interior has nonempty intersection with $j$ and $v\in\omega\subseteq\sigma_{j}$.\item A finite set of rays $(\mathbf{t}_{i})$.\item For each cut $\mathbf{c}$ and any $(j\cap Int\ \omega)-\Delta$ a function $f_{\mathbf{c},x}\in \mathbf{C}[P_{x}]$ having the same properties of slab functions $f_{\mathbf{b},x}$ in definition .\item An orientation of $Q_{j,\mathbf{R}}^{v}$.\end{itemize}A scattering diagram is denoted by $\mathcal{D}=\{\mathbf{t},f_{\mathbf{c}}\}$. A cyclic ordering of maximal dimensional cells $\sigma_{1},\cdots \sigma_{r}=\sigma_{0}$ containing $j$ induces a cyclic ordering of of cuts $\mathbf{c}_{i}\subseteq\bar{\bar{\sigma_{i-1}}}\cap\bar{\bar{\sigma_{i}}}$. Note that $\bar{\bar{\sigma_{i}}}$ are "chambers" in $Q_{j,\mathbf{R}}^{v}$ divided by $\bar{\bar{\sigma_{i-1}}}\cap\bar{\bar{\sigma_{i}}}$.

For a ray $\mathbf{t}_{i}\subseteq\bar{\bar{\sigma_{s}}}$ we can define a log automorphism $\theta_{\mathbf{t}_{i}}:\Lambda_{\sigma_{s}}\rightarrow (R_{g,\sigma_{s}}^{k})^{\times}$ \begin{equation}\theta_{\mathbf{t}_{i}}:=\exp(-\log(1+c_{\mathbf{t}_{i}}z^{m_{\mathbf{t}_{i}}})\partial_{n_{\mathbf{t}_{i}}})\end{equation}Explicitly $\theta$ is \begin{equation}\theta_{\mathbf{t}_{i}}:m\rightarrow(1+c_{\mathbf{t}_{i}}z^{m_{\mathbf{t}_{i}}})^{-\langle\bar{m},n_{\mathbf{t}_{i}}\rangle}\end{equation}where $n_{\mathbf{t}_{i}}$ is the generator of normal vectors of $\mathbf{t}_{i}$ in $Q_{j,\mathbf{R}}^{v}$ oriented positively with respect to the orientation of $Q_{j,\mathbf{R}}^{v}$.

For a cut $\mathbf{c}$ one can similarly define a log morphism $\theta_{\mathbf{c}_{i}}:P_{\omega,\sigma_{i-1}}\rightarrow(R_{g,\sigma_{i}}^{k})^{\times}$ \begin{equation}\theta_{\mathbf{c}_{i}}:m\rightarrow(f_{\mathbf{c}_{i},x})^{-\langle\bar{m},n_{\mathbf{c}_{i}}\rangle}\end{equation}Finally we let $\theta_{\mathcal{D},g}^{k} $ be the ordered composition along a loop (around the origin) of log morphisms associated to all rays and cuts. It is an log automorphism. A scattering diagram is consistent to order $k$ if modulo $t^{k+1}$ $$\theta_{\mathcal{D},\mathrm{Id}_{\sigma_{j}}}^{k}=1$$ \end{dfn}
\begin{note}It is easy to show the consistency does not depend on the choice of the vertex $v$ and from now on we can drop $v$ when discussing scattering diagrams.\end{note}
A $structure$ induces a collection of scattering diagrams labeled by joints in the following way. For every joint $j$ and cell $\omega\in\mathcal{F}$ with $v\in\omega\in\mathcal{F}$ and $\omega\cap j\neq\emptyset$ the projection of a slab  together with the function attached to it gives us a cut $\mathbf{c}$ with $f_{\mathbf{c},x}$. For a wall $\mathbf{p}$ containing $j$ there are two cases.\begin{itemize}\item $j\in\partial\mathbf{p}$. We obtain a ray $(\bar{\bar{\mathbf{p}}},m_{\mathbf{p}},c_{\mathbf{p}})$.\item $j\cap Int\ \mathbf{p}\neq\emptyset$. In that case $\bar{\bar{\mathbf{p}}}$ is a line containing the origin and we add a pair of rays with opposite directions.\end{itemize}  The consistency condition of a $structure$ now becomes the consistency conditions of the associated scattering diagrams.

\begin{note}Our previous discussion of universal scattering diagrams in the two dimension is special in the sense that the normal spaces of all joints are the same two dimensional space. A sufficient condition for this to be true is that the underlying singular affine space is topologically $\mathbf{C}\simeq\mathbf{R}^{2}$. Since in our case the affine base is the space of quadratic differentials this condition is satisfied if the dimension is correct.\end{note}

Although a general scattering diagram looks like a two dimensional thing, the problem in higher dimensions is actually much more complicated  than its counterpart in the two dimension due to two reasons. First when we try to make the order $k$ consistent scattering diagram from the order $k-1$ by induction we encounter nontrivial contributions from the nonzero dimensionality of $j$. Second, what we really need to construct the degeneration are cells of the $structure$. So although we use scattering diagrams to discuss the consistency condition we still have to build codimensional one slabs and walls.

In general a collection of consistent scattering diagrams and the associated compatible system of consistent $structures$ are constructed together by induction. Starting from a scattering diagram which is consistent to order $k-1$ at $j$  and is obtained from a $structure$ $\aleph_{k-1}$ it is shown by very difficult arguments in \citep{GS1} that after adding only rays  $\theta_{\mathcal{D^{'}}}^{k}$ for the new scattering diagram $D^{'}$ can be put into a canonical form. If the codimension of the joint is two then we may have to modify functions $f_{\mathbf{c}}$ attached to cuts (slabs). From this new scattering diagram one can build a new structure $\aleph_{k}$ by adding wall and changing slab function. The modifications of slabs functions of different joints do not interact. If a ray $\mathbf{t}$ is added to the scattering diagram then a wall $(\mathbf{p}_{\mathbf{t}},m_{\mathbf{t}},c_{\mathbf{t}})$ is added to the structure with $$\mathbf{p}_{\mathbf{t}}:=(j-\mathbf{R}_{\geq 0}\bar{m}_{\mathbf{t}})\cap\sigma$$where $\sigma$ is the unique cell with $\mathbf{t}\subseteq\bar{\bar{\sigma}}$. We do this for all joints. This structure is already consistent for any codimensional zero joint $j$ if $\sigma_{j}$ is bounded. For codimensional one joints to achieve consistency one has to modify slab functions again (this is called homological modifications) because the modifications of slabs functions given before may not be consistent for different joints. Finally we need a normalization procedure for slab functions to obtain consistency around codimensional two joints. The $structure$ obtained is consistent to order $k$.

Note that when there are are several joints the effects of adding rays (walls) for different joints will interact. The interaction has two consequences. One is the necessity of homological modifications of slab functions at a given order. The other is that the intersections of the added  rays (walls)  produce new joints.  In the discussion of universal scattering diagrams in two dimensional problems we have  taken care of them because  the consistency conditions are formulated for all singularities of the (universal) scattering diagram instead of just the origin (one joint).\\

 In what sense can we call the above construction  a procedure of computing instanton corrections? According to the philosophy of Strominger, Yau, Zaslow and many others, one should count holomorphic disks in the mirror family wrapping some special Lagrangian fibers. This kind of instanton corrections has been understood (although perhaps not completely) in some cases of mirror symmetry such as some so-called Landau-Ginzburg models and some toric (noncompact) Calabi-Yau varieties. But it has never been understood for any compact Calabi-Yau manifolds or for noncompact hyperkahler manifolds such as Hitchin's moduli spaces. On the other hand, it is a trend in recent years to replace enumerative problems in holomorphic geometry by enumerative problems in tropical geometry which are more or less combinatorial\footnote{See \citep{Mi,NS} for some examples.}. This is not a good place to explain the ideas of tropical geometry. It is enough to know that Gross and Siebert conjecture that corrections given in their method are essentially tropical data and should be eventually equivalent to the corrections by holomorphic disks. There are some works in this direction \citep{GP,GPS,Gr,Ni}. Also note that instead of studying the affine base of the mirror, their construction stays in the same affine manifold. So the data here should be the dual of the tropical data. Since we will use their method, we will also stay on one side of the mirror symmetry which is why we do not really have to consider the Hitchin's moduli space with gauge group $PGL(2, \mathbf{C})$. That is also why we do not need to consider the dualization (Legendre transform) of the triple at all in this paper. In this sense, the instanton problem for mirror symmetry has two steps\begin{itemize}\item Find the corrections of the complex structures. This step has been completely solved by Gross and Siebert and is a major breakthrough in this area. \item Identify the corrections as given by instantons of the mirror.\end{itemize}In this work, we will take another route\begin{itemize}\item Find the corrections of the complex structures for the Hitchin's moduli space.\item Show that these corrections are indeed given by "instantons". But these instantons are neither objects on the mirror nor objects on the moduli space itself. They are objects on the underlying Riemann surface (they are in fact some critical trajectories on the surface).\item Moreover, the instanton correction problem here is equivalent to another instanton problem : the instanton correction problem of hyperkahler metrics.\end{itemize}

 \section{Seiberg Witten Theory And Wall Crossing}

 We have considered the problem of instanton corrections of complex structures, now let us turn to the problem of instanton corrections of Calabi-Yau metrics.  The story starts in the seemly unrelated context of determining the exact form of the low energy effective theory of $N=2$ four dimensional supersymmetric gauge theory. This is the famous Seiberg Witten theory \citep{SW1}. For introductions see \citep{Bi,Kl,Ler}. It is relevant here because the special Kahler base $B$ of the Hitchin's moduli space under Hitchin's fibration appears in this theory as the "quantum moduli space". This section also describes the physical background of the phenomenon of wall crossings. Mathematicians who are not interested in physics can ignore this section except for the definition 4.1.

 Let us illustrate  basic ideas in this theory which are relevant to this work by considering the so-called $N=2$ pure $SU(2)$ theory which is the first model solved by Seiberg and Witten. It turns out that in low energy, the theory is effectively (partly) described  by a sigma model  whose target is the so-called "quantum moduli space" and therefore it is important to know the metric on this moduli space. This metric has to be Kahler but is allowed to have singularities. As a complex manifold, this space is just the complex plane $C$ with holomorphic coordinate $u$. Seiberg and Witten showed that the metric is a singular special Kahler metric with two singularities. They determined the monodromies of these two singularities. The special holomorphic coordinates $a, a_{D}$ are functions of $u$ $$a=a(u), a_{D}=a_{D}(u)$$The so-called electric-magnetic duality then requires that $(a(u), a_{D}(u))$ is a section of an $SL(2,\mathbf{Z})\simeq\mathrm{Sp}(2,\mathbf{Z})$ local system over the complement of singularities. So the problem has been reduced to a Riemann-Hilbert problem of determining a section of a local system with prescribed monodromies. This problem was also solved geometrically by Seiberg and Witten by uniformization theory. They explicitly wrote down an elliptic curve whose moduli is parameterized by $u$. Define \begin{equation}a:=\int_{A}\lambda, \ \ a^{D}:=\int_{B}\lambda\end{equation} where $\lambda$ is a canonically defined one form (the Seiberg Witten differential) and $(A, B)$ is a basis of one cycles of the elliptic curve. The $SL(2,\mathbf{Z})$ action is obtained by transforming one cycles. The two singularities are the locus where the elliptic curves are singular and the monodromies are Picard-Lefschetz transformations. This will be a part of our second example and details will be provided in section 9.

Soon after Seiberg and Witten's breakthrough Donagi, Witten and other people realized that the quantum moduli space is actually the base of a Hitchin's moduli space and the fibers are Prym abelian varieties (which happen to be elliptic curves in this example) of spectral curves \citep{Do,DW}. In our example this is clear if one thinks of the curve defined by $x^{2}=\lambda^{2}$ as the spectral curve covering $C=CP^{1}$ and identify $\lambda^{2}$ with $-det\ \varphi$. Under the action of a monodromy matrix $\mathbf{M}\in SL(2,\mathbf{Z})$, \begin{equation}\left(\begin{array}{c}a(u)\\a_{D}(u)\end{array}\right)\rightarrow \mathbf{M}\left(\begin{array}{c}a(u)\\a_{D}(u) \end{array}\right)\end{equation}

So far we have been discussing metrics on $B$ instead of on the Hitchin's moduli space. One needs only one more step to reach it. In another paper \citep{SW2}, Seiberg and Witten considered the compactification of a four dimensional gauge theory to three dimensions. In other words, one replaces the (Euclidean) spacetime $\mathbf{R}^{4}$ by $\mathbf{R}^{3}\times S^{1}$ where $S^{1}$ is a circle of radius $R$\footnote{It turns out that this $R$ is the same $R$ in  equation (3).}. It turns out the low energy effective theory can be formulated as a three dimensional sigma model with spacetime $R^{3}$ but with the target space (quantum moduli space) replaced by the total space of the Hitchin's fibration, i.e. the Hitchin' moduli space itself (\citep{SW2,G2,CK})! Supersymmetry requires that the metric on the target space must be hyperkahler.\\

What do all these facts have to do with instanton correction problem? Well, although it may not be obvious the solution of the four dimensional theory actually has the form of an exact solution obtained by incorporating all instanton effects. The prepotential of the special Kahler metric  can be calculated order by order from the Picard-Fuchs equations satisfied by periods and it has a form of a summation over infinitely many instantons. Later, in a tour de force, Nekrasov and Okounkov \citep{NO} calculated the instanton contributions directly according to rules of instanton calculus in quantum field theory (and hence it could be considered as a first principle verification of Seiberg Witten theory) and derived the Seiberg Witten solution. One may wonder if there is a similar story for the three dimensional theory such that the hyperkahler metric on the Hitchin's moduli space can be exactly determined by calculating all instanton corrections. A direct attack in the spirit of Nekrasov and Okounkov is absent at present, but as we will show in the next section there is an indirect way to do that\footnote{There are some partial first principle calculations which are consistent with it. The references are \citep{CDP,CP}}.\\

One of the basic ingredients of that approach is the determination of BPS spectra of the gauge theory. This is also one of the most important consequences of Seiberg Witten theory and so is described in this section. We consider "particles" with electric and magnetic charges. Note that the gauge charge lattice in this example form a rank two integral lattice $\mathbf{Z}^{2}$ with a symplectic pairing which is nothing but the intersection pairing. Locally we choose a split of a basis (for example: $\gamma_{1}:= e:=(1,0)=A,\ \  \gamma_{2}:= m:=(0,1)=B)$ and call half of the one cycles the electric charges $\gamma_{e}$ and the other half magnetic charges $\gamma_{m}$. So a general cycle (charge) has a decomposition $$\gamma= \gamma_{e}+\gamma_{m}=n_{e}e +n_{m}m$$ where $n_{e},n_{m}$ are integers. Then the central charge defined in section 2 for a charge $\gamma$ is given by \begin{equation}Z_{\gamma}(u)= {1\over \pi}(\int_{\gamma_{e}(u)}\lambda+\int_{\gamma_{m}(u)}\lambda) = n_{e}a(u)+ n_{m}a_{D}(u)\end{equation}Note that under the the action of the monodromy matrix $\mathbf{M}$, $$(n_{e},n_{m})^{t}\rightarrow(\mathbf{M}^{-1})^{t}(n_{e},n_{m})^{t}$$

By the representation theory of superalgebras, it is known that the mass of a particle of charge $\gamma$ is not smaller than the norm of its central charge. It is "BPS" if this lower bound of mass is saturated \begin{equation}m=|Z_{\gamma}|\end{equation} and as a consequence this configuration preserves a fraction of the underlying supersymmetry of the theory.

The fundamental question of BPS spectra is: what are BPS particles in the theory (more precisely speaking, for which charges do there exist BPS particles and how many are there)? The answer turns out to be independent of $u$ generically but exhibits discontinuous jumps when crossing some real codimensional one hypersurfaces of $B$ defined by the condition that phases of central charges of independent charges align.  This is called wall crossing and $u$ is considered as a moduli parameter. For example, in the pure $SU(2)$ case, the hypersurface (called (marginal) stability wall, not to be confused with walls defined in section 3) is the locus$$\{u\mid \arg Z_{\gamma_{1}}= \arg Z_{\gamma_{2}}\}$$ Suppose we have a BPS particle on one side of the wall with charge $\gamma= (n_{e}, n_{m})$. The  phases of $Z_{\gamma_{1}}$ and $Z_{\gamma_{2}}$ do not align.  As a result this particle cannot be considered as two BPS particles with charge $(n_{e}, 0)$ and $(0, n_{m})$ (and masses $m_{1}$ and $m_{2}$) respectively because by the conservation of mass and BPS condition we have $$|Z_{\gamma}|= m=m_{1}+m_{2}=\mid Z_{1}\mid+\mid Z_{2}\mid$$where $Z_{1}$ and $Z_{2}$ are central charges of the two hypothetical BPS particles.  But by the additivity of central charge we have $$Z_{\gamma}=Z_{1}+Z_{2}$$These two equations cannot be simultaneously true because the phases of two central charges do not align.   Changing the moduli $u$ without touching any stability walls the above argument continues to work and we do not expect any changes of the BPS spectra. However, the contradiction argument clearly breaks down on the marginal stability wall and therefore we do expect a change of spectra when we cross a stability wall.\\

Of course, this just tells us that the spectra can change instead of how. Later we will describe a systematic way to determine the wall crossing of the BPS spectra. Nevertheless, without knowing this general method physicists determined the spectra correctly in mid 90's using some monodromy and symmetry arguments
 (see \citep{BF,Bi,SW1}). Let us describe the result for pure $SU(2)$ theory.

A stability wall in the moduli space parameterized by $u$  is a curve. There are two independent charges $\gamma_{1}, \gamma_{2}$ forming a basis of the charge lattice.  There are two stability walls. The alignment of $\gamma_{1}, \gamma_{2}$ gives one of them while the alignment of $\gamma_{1}, -\gamma_{2}$ gives the other. Let us consider the union of the two stability walls. Since it is determined by phases of periods, it has a description in term of uniformization theory \citep{AFS,Fa,Ma}.  In fact  periods of a meromorphic differential $a=Z_{A}(u)$ and $a_{D}=Z_{B}(u)$ are solutions of a Picard-Fuchs equation and they are hypergeometric functions. Their ratio satisfies a Schwartzian equation, see \citep{Y}. So the union which is the locus where the ratio is real is the pull back of an  interval by the Schwartz map.  It can also be described by the more familiar uniformization theoretical data: the modular parameter $\tau$ of elliptic curves and modular functions. In fact, let $\omega=dx/y$ be the canonical $holomorphic$ one form of the hyperelliptic curve $y^{2}=P(x)$ then \begin{equation}\partial_{u}\lambda=\omega+\mathrm{exact}\ \mathrm{form}\end{equation}so that $\tau$ defined in section 2 is the usual $\tau$ parameter of elliptic curves.\begin{equation}\tau={\int_{B}\omega\over \int_{A}\omega}\end{equation} The universal cover of the complement of the singular locus of the $u$-plane is the upper half plane. Since the three monodromies around the two singularities and the infinity of the $u$-plane generate the modular group $\Gamma(2)$ the associated modular function maps a fundamental domain of $\Gamma(2)$ to the $u$-plane and the union is the image of two arcs.  The union  is simple closed. It passes through the two singularities on the $u$ plane coordinated as $u=\pm\Lambda^{2}$, $\Lambda\in \mathbf{R}$. Since we have two nontrivial monodromies at finite places, we take two branch cuts along the real axis from the two singularities to the minus  infinity. The region outside/inside the union is also known as weak/strong coupling region (these names follow as a consequence of "asymptotic freedom"). The branch cuts divide the strong coupling region into two halves where the spectra are related by monodromies and hence it is enough to work with one half, say the one below the real axis. The result is \begin{itemize}\item In the lower half strong coupling region, the charges of possible BPS particles are $$\pm(2,-1), \pm(0,1)$$\item Across the stability wall from the above strong coupling region to the  weak coupling regions, the spectra change to $$\pm(2,0), \pm(2n,1), n\in Z$$\end{itemize} Here plus sign means they are particles while minus sign means they are "antiparticles".

The determination of BPS spectra is a problem of the four dimensional theory, but we will see that they provide the complete set of instantons that contribute to the exact form of the  hyperkahler metric of the quantum moduli space (Hitchin's moduli space) of the three dimensional compactified theory.\\

Seiberg Witten theory has been vastly generalized to allow other gauge groups and also matters (fermion fields in representations of the gauge groups). Correspondingly, we will consider moduli spaces of possibly singular solutions of Hitchin's equations with prescribed asymptotic behaviors near singularities. In fact, residues of Higgs fields are given by masses of the matters. More precisely, we add $N_{f}$ copies of a representation of the gauge group (for $SU(2)$ we usually take the spin one-half representation of $SU(2)$) and assign masses $m_{i},1\leq i\leq N_{f}$ to the fields in these representations. $N_{f}$ is called the flavor number. The mass formula is modified to \begin{equation}m=|n_{e}a+n_{m}a_{D}+\sum_{i=1}^{N_{f}}s_{i}m_{i}|\end{equation}where $s_{i}$ are integral constants called flavor charges. The transformations of charges now become integral affine symplectic transformations. In other words, under the action of $(H, \mathbf{M})$ belonging to the semi direct product of $(\mathbf{Z}^{2})^{N_{f}}$ and $SL(2,\mathbf{Z})$
$$\left(\begin{array}{c}a\\a_{D}\end{array}\right)\rightarrow \mathbf{M}\left(\begin{array}{c}a\\a_{D}\end{array}\right)+H\vec{m}=\mathbf{M}\left(\begin{array}{c}a\\a_{D}\end{array}\right)+
\sum_{i=1}^{N_{f}}m_{i}\left(\begin{array}{c}n^{i}\\n^{i}_{D}\end{array}\right)$$ \begin{equation}(n_{e},n_{m})^{t}\rightarrow (\mathbf{M}^{-1})^{t}(n_{e},n_{m})^{t}, \vec{s}\rightarrow\vec{s}-H(n_{e},n_{m})^{t}\end{equation}where $\vec{m}:=(m_{1},\cdots m_{N_{f}})^{t},\vec{s}:=(s_{1},\cdots s_{N_{f}})^{t}, H=((n^{1},n^{1}_{D}),\cdots (n^{N_{f}},n^{N_{f}}_{D}))$. Here the formulas are formulated for complex one dimensional $B$. If the dimension is higher, we just need to replace $n_{e}a+n_{m}a_{D}$ by the sum over all gauge charges and replace $SL(2,\mathbf{Z})$ by $\mathrm{Sp}(2g,\mathbf{Z})$. Geometrically it means that $m_{i}$ are residues of the Seiberg-Witten differential $\lambda$ at some singularities over the Riemann surface $C$ and $(n_{e},n_{m})$ and $\vec{s}$ are gauge charges and flavor charges  respectively so that $m=|Z_{\gamma}|$ where $\gamma\in\hat{\Gamma}$ is the sum of the gauge charge and the flavor charge.\\

It is  interesting to see that now we have two a priori different problems of instanton corrections  on a Hitchin's moduli space. One is suggested by the study of mirror symmetry of Hitchin's moduli spaces, the other is from three dimensional Seiberg Witten theory which does not involve mirror symmetry in the formulation given above. Could these two problems have equivalent answers? One of the purposes of this paper is to show  that  the answer to this question is positive.

Not every $SU(2)$ Hitchin's moduli space arises in this way. In fact, physically consistent $SU(2)$ theory requires $N_{f}\leq 4$. Another constraint is that all quantum moduli spaces $B$ of $SU(2)$ has complex dimension one. Nevertheless the notion of stability walls continue to make sense in general. Although the notion of  low energy effective theory and hence the Seiberg-Witten theory is not yet mathematically well defined, stability walls and BPS spectra can be rigorously defined and these are all we need to develop our results. A mathematical operational definition of BPS spectra will be given later. Here we write down the definition of a stability wall.
\begin{dfn}Let $\mathcal{M}$ be a Hitchin's moduli space defined in section 2. Let $\gamma_{1},\gamma_{2}\in \hat{\Gamma}$ be two charges and $B$ be the base of the Hitchin's fibration. The central charges $Z_{\gamma_{i}}(u)$ are defined in section 2. The stability wall $SW_{\gamma_{1},\gamma_{2}}(u)$ of a pair $(\gamma_{1},\gamma_{2})$ is the following real codimensional one locus\begin{equation}SW_{\gamma_{1},\gamma_{2}}(u):=\{u\mid \arg Z_{\gamma_{1}}= \arg Z_{\gamma_{2}}\}\end{equation}where $u$ is a set of  holomorphic coordinates over $B$ (usually the value of the Hitchin's map contains explicitly moduli parameters parameterizing the space of meromorphic quadratic differentials with prescribed asymptotics at singularities  and these parameters are taken to be the natural holomorphic coordinates on $B$). Note that we have suppressed the subscripts of $u$. Clearly $SW_{\gamma_{1},\gamma_{2}}(u)$ is unchanged if we exchange $\gamma_{1}$ and $\gamma_{2}$ or if we change the sign of both $\gamma_{1}$ and $\gamma_{2}$.\end{dfn}A stability wall is  codimensional one in $B$.  There could be countably many stability walls.

\section{Twistor Spaces And Instanton Corrections of Hyperkahler Metrics}

We want to write down the hyperkahler metric of a Hitchin's moduli space (interpreted as a quantum moduli space according to section 4) in a form such that it looks like a sum over BPS instantons (labeled by charges) and since we suspect that there will be nontrivial wall crossing the formula should also be designed to exhibit the phenomenon of wall crossing. Gaiotto, Moore and Neitzke realized how to do this by passing to the associated twistor space.

Any hyperkahler manifold $\mathcal{M}$ has a twistor space which is $\mathcal{M}\times CP^{1}$ with a natural complex structure together with some other data and conversely a hyperkahler manifold can be constructed from a twistor space. Let us state the theorem for the converse construction \citep{HK}.
\begin{thm}Suppose $(Z,J_{Z})$ is a $2n+1$ dimensional complex manifold together with\begin{itemize}
\item a holomorphic projection $p:Z\rightarrow CP^{1}$,
\item a holomorphic section $\Omega$ of $p^{*}(O(2))\otimes \Lambda^{2}(TF)^{*}$ which is symplectic on the fibers of $TF$ where $TF$ is the kernel of the map $dp: TZ\rightarrow TCP^{1}$,
\item a free antiholomorphic involution $\tau:Z\rightarrow Z$ which preserves $\Omega$ and $p\circ\tau=\tau^{'}\circ p$ where $\tau^{'}$ is the antipodal map of $CP^{1}$
\end{itemize}Let $\mathcal{M}$ be the set of holomorphic curves $C$ in $Z$ isomorphic to $CP^{1}$ with (same) normal bundle $2nO(1)$ and preserved by the involution. Then $\mathcal{M}$ is a hyperkahler $4n$-manifold.\end{thm}

The other direction of the twistor method is straightforward. Let ($\mathcal{M}$, $J_{1}, J_{2}, J_{3}, \omega_{1}, \omega_{2}, \omega_{3},g$) be a hyperkahler $4n$-manifold. Let $\omega_{\pm}:=\omega_{1}\pm i\omega_{2}$, then\begin{equation}\Omega(\xi):=-{i\over 2\xi}\omega_{+}+\omega_{3}-{i\over 2}\xi\omega_{-}\end{equation} is the holomorphic symplectic form in the compatible (to the metric $g$) complex structure $J_{\xi}$ defined by
$$J_{\xi}:={i(-\xi+\bar{\xi})J_{1}-(\xi+\bar{\xi})J_{2}+(1-|\xi|^{2})J_{3}\over 1+|\xi|^{2}}$$ where $\xi\in CP^{1}$ parameterizes all compatible complex structures and is called the twistor parameter. Then the tautological almost complex structure on $Z:= \mathcal{M}\times CP^{1}$ can be shown to be integrable with a tautological antiholomorphic  involution and $\Omega(\xi)$ patch together to form the holomorphic section required in the second condition in the above theorem. The fibers of the projection $\pi: Z\rightarrow \mathcal{M}$ give us holomorphic curves isomorphic to $CP^{1}$ and clearly $\mathcal{M}$ is the deformation space of such rational curves.

Hitchin's moduli spaces are not only hyperkahler but also algebraically integrable systems, i.e. there are Hitchin's fibrations. Gaiotto, Moore and Neitzke took advantage of this fibration structure and postulated the existence of a set of locally defined $\mathbf{C}^{\times}$-valued functions $\mathcal{X}_{\gamma}(u, \theta; \xi)$ where $u$ as before is the holomorphic coordinate on $B$\footnote{There is a canonical choice of $u$. We just take the value of the Hitchin's map $u=det\ \varphi$.}, "$\theta$" are angular coordinates of torus fibers and $\gamma$ is a charge (an element of the charge lattice). They should satisfy the following conditions\begin{itemize}
\item \begin{equation}\mathcal{X}_{\gamma}\mathcal{X}_{\gamma^{'}}=\mathcal{X}_{\gamma+\gamma^{'}}\end{equation}
\item There is a holomorphic Poisson bracket such that \begin{equation}\{\mathcal{X}_{\gamma}, \mathcal{X}_{\gamma^{'}}\}= \langle\gamma, \gamma^{'}\rangle\mathcal{X}_{\gamma}\mathcal{X}_{\gamma^{'}}\end{equation}
The meaning of this condition will be clear in the next section.
\item A reality condition \begin{equation}\mathcal{X}_{\gamma}(\xi)=\overline{\mathcal{X}_{-\gamma}(-1/\bar{\xi})}\end{equation}
\item They are holomorphic on the Hitchin's moduli space $\mathcal{M}$ with respect to the complex structure $J_{\xi}$ for any $\xi$
\item For any fixed $x:=(u,\theta)$, $\mathcal{X}_{\gamma}$ is holomorphic on a dense subset of $\mathbf{C}^{\times}$ (in fact the complement of a union of countably many rays). Although they are not holomorphic everywhere, the denseness is enough to guarantee the vanishing of Nijenhuis tensor and hence the almost integrable complex structure on $\mathcal{M}\times CP^{1}$ is integrable. $\mathcal{X}_{\gamma}$ is only piecewise holomorphic and we consider the discontinuous jumps across those discontinuous rays as part of our input data. There are constraints for them. This is the most subtle part of the construction and will be explained in the next section.
\item Define \begin{equation}\Omega(\xi):={1\over 8\pi^{2}R}\epsilon_{ij}{d\mathcal{X}_{\gamma_{i}}\over \mathcal{X}_{\gamma_{i}}}\wedge{d\mathcal{X}_{\gamma_{j}}\over \mathcal{X}_{\gamma_{j}}}\end{equation}where $R$ is the parameter $R$ in equation (3), $\epsilon_{ij}:= \langle\gamma_{i}, \gamma_{j}\rangle$ is the integral symplectic (intersection) pairing on the gauge charge lattice $\Gamma_{gau}\ni \gamma$ and $d$ is the exterior differential on $\mathcal{M}$. Although we define $\mathcal{X}_{\gamma}$ for all charges we only sum over a basis of gauge charges when we define $\Omega(\xi)$. It is required  that $\Omega(\xi)$ is nondegenerate and has simple poles when $\xi$ goes to zero or infinity and the discontinuous jumps of $\mathcal{X}_{\gamma}$ are such that they cancel each other and $\Omega(\xi)$ is holomorphic in $\mathbf{C}^{\times}$.
\item $$\lim_{\xi\rightarrow 0} \mathcal{X}_{\gamma}(u, \theta; \xi)\exp(-\xi^{-1}\pi R Z_{\gamma}(u))$$ exists.

\end{itemize}
By the twistor construction such data gives rise to a hyperkahler structure on $\mathcal{M}$ with $\Omega(\xi)$ as the holomorphic symplectic form and the Kahler form and hyperkahler metric can be read from it. In fact the holomorphic section in the second condition of theorem 5.1 is given by (57). The holomorphic projection is the canonical one $p:\mathcal{M}\times\ CP^{1}\rightarrow \ CP^{1}$ sending $(x,\xi)$ to $\xi$. For each $x\in \mathcal{M}$ there is a section $s_{x}: \ CP^{1}\rightarrow \mathcal{M}\times\ CP^{1}$ defined by $s_{x}(\xi)=(x,\xi)$. The normal bundle of it is shown in \citep{G1} to be $2n O(1)$. Therefore $\mathcal{M}$ is the set of rational curves in theorem  5.1. Finally the involution is $\tau(x,\xi):=(x,-1/\bar{\xi})$.

Clearly the key issue here is to describe the discontinuous rays and jumps. This will be the task of the next section where a unified treatment of wall crossing in different places are given. In this section, we will work out a simple example of this kind of ansatz and this example is also the "initial" uncorrected metric to be corrected by instantons.

Like in the Seiberg Witten theory we choose locally a splitting of the charge lattice into electric and magnetic parts $(\gamma_{e}^{i}, \gamma_{m}^{i})$ and hence $\theta=(\theta_{e}^{i},\theta_{m,i})$. Define
\begin{equation}\mathcal{X}_{\gamma}^{sf}(\xi):= \exp(\pi R\xi^{-1}Z_{\gamma}+ i\theta_{\gamma} +\pi R\xi\bar{Z}_{\gamma})\end{equation} Here $sf$ represents "semiflat". The metric defined by $\mathcal{X}_{\gamma}^{sf}(\xi)$ is called semiflat metric and is well known to be the one without any instanton corrections from the perspective of mirror symmetry\footnote{Semiflat metrics were first studied by \citep{GSVY}. For their role in mirror symmetry, see \citep{Leu}.}.
In this case, $\mathcal{X}_{\gamma}^{sf}(\xi)$ has no discontinuity and by writing down $\Omega(\xi)$ explicitly and comparing it to equation (53), one gets \begin{equation}\omega_{3}^{sf}:={i\over 2}(R(\mathrm{Im}(\tau)_{ij}da^{i}\wedge d\bar{a}^{j})+{1\over 4\pi^{2}R}((\mathrm{Im}(\tau))^{-1}_{ij}db^{i}\wedge d\bar{b}^{i})\end{equation} where $a_{i}$ and $\tau$ are defined as in section 2. $$db_{i}:=d\theta_{m,i}-\tau_{ij}d\theta_{e}^{j}$$

The first term of $\omega_{3}^{sf}$ is a constant times the Kahler form of the special Kahler metric on the affine base $B$. If we let $R\rightarrow\infty$, then after rescaling, the second term drops away and the semiflat metric of the total space of the Hitchin's fibration collapses to the special Kahler metric on $B$. This kind of degenerations of metrics is exactly what has been conjectured to be the case for degenerations of Calabi-Yau metrics in the family version of Strominger, Yau and Zaslow's formulation of mirror symmetry. In general when the torus fibration has singular fibers one expects that when the large complex limit is approached the hyperkahler\footnote{The conjecture is for general Calabi-Yau metrics. But here we are only concerned with hyperkahler spaces which means the limit on the base should be a singular special Kahler structure.} metric converges after appropriate rescaling to a singular special Kahler metric  on the base. Moreover the deviation from the  semiflat metric should be small for large $R$. Therefore if we assume that the effect of instanton corrections are exponentially suppressed by $R$ compared to the semiflat one (which is based on general physical principle), then the exact hyperkahler metric (after adding all instanton corrections) would have the properties anticipated by mirror symmetry even if the context here has a priori nothing to do with mirror symmetry.

 Let us be more specific. We impose another condition for $\mathcal{X}_{\gamma}(u, \theta; \xi)$
 \begin{itemize}
 \item We have the large $R$ asymptotic\begin{equation}\mathcal{X}_{\gamma}=\mathcal{X}_{\gamma}^{sf}(1+ \exp(-const\cdot R))\end{equation}
 \end{itemize}

It follows \citep{G1} that given a set of coordinates $\mathcal{X}_{\gamma}(u, \theta; \xi)$ satisfying all mentioned conditions, there is a hyperkahler metric $g$ on $\mathcal{M}$ with the following properties:\begin{itemize}\item It is continuous (in the next section this follows from the wall crossing formula) and smooth except for locus above singularities on the base.\item The deviation from the semiflat metric is exponentially suppressed by $R$. For a fixed $R$, it approaches the semiflat metric when $|Z_{\gamma}|\rightarrow\infty$ for all $\gamma$\footnote{This actually requires a stronger condition $\mathcal{X}_{\gamma}=\mathcal{X}_{\gamma}^{sf}(1+ \exp(-const\cdot R|Z_{\gamma}|))$.}. \item The volume of the fiber of the Hitchin's fibration is independent of $u$ and is $vol(X_{u})=(1/R)^{r}$ where $r$ is half of the rank of the charge lattice. This is a direct consequence of the fact that $\Omega(\xi)$ is the pull-back of the canonical holomorphic symplectic form on a complex torus (see the next section).\item In the complex structure $J_{3}$ the Hitchin's fibration is holomorphic and the holomorphic symplectic form is always $\omega_{+}$ independent of $R$.\end{itemize}

The discussion above also answers an important conceptual question. If there is any chance that the instanton correction of metrics of Hitchin's moduli spaces can be related to Gross and Siebert's work a Hitchin's moduli space must be put into a degenerating family. This degeneration should be either a large complex degeneration or a toric degeneration (in fact the two are conjectured to be essentially equivalent). It is not clear how to adapt the algebraic geometric definition of large complex degenerations of projective Calabi-Yau varieties to Hitchin's moduli spaces.

On the other hand Gross and Wilson considered in \citep{GW} a large complex degeneration of elliptic K3 surfaces (which are compact hyperkahler manifolds) and obtained  the same type of metric degenerations. Roughly speaking since there is a nice theory of mirror symmetry of K3 surface they use the mirror of the so-called large volume degeneration (also called the large Kahler degeneration which, according to the philosophy of mirror symmetry exchanging complex geometry and symplectic geometry, should be mirror to the large complex degeneration) to define the large complex degeneration of elliptic K3 surfaces.  The author has seen the claim in the literature that this recovers the usual definition as a family approaching a limit point with maximally unipotent monodromy although he fails to find a reference. Anyway we take the definition of Gross and Wilson as a natural one as it fits into the conjectural metric degeneration picture.

In Gross and Wilson's work the large complex degeneration can be described in the following way. The underlying complex structure in which the elliptic fibration of the K3 surfaces $X_{l}$ is fixed. Let $\omega_{l}$ be a Ricci-flat Kahler metric on $X_{l}$ with volume independent of $l$. Let $\epsilon_{l}$ be the volume of a fiber (which is assumed to be independent of the choice of the fiber) and suppose $\epsilon_{l}\rightarrow\infty$ when $l\rightarrow\infty$. Then they verified the conjectural metric collapsing picture described in the limit form of the family version of SYZ conjecture in section 3 for the sequence $(X_{l}, \epsilon_{l}\omega_{l})$.

If we compare it with the previously described behavior of Gaiotto-Moore-Neitzke's metrics it is clear that GMN's ansatz gives precisely such a family of Calabi-Yau spaces. Here $\epsilon_{l}=(1/R)^{r}$ and  $R\rightarrow\infty$ means $l\rightarrow\infty$ (for $SU(2)$ $r=1$). The volume is fixed because the top degree wedge of (fixed) $\omega_{+}$ is the holomorphic volume form. This suggests that we could view the family of Hitchin's moduli spaces parameterized by $R$ with $R\rightarrow\infty$ as the right substitute of the large complex limit.
\begin{dfn}Let $\mathcal{M}$ be the Hitchin's moduli space defined in section 2. We modify the Hitchin's equation by using (3) instead of (1) or (2). This gives us a family of Hitchin's moduli spaces denoted by $\mathcal{M}(R)$. The large complex degeneration of Hitchin's moduli spaces is the family $\mathcal{M}(R), R\rightarrow\infty$.\end{dfn}
\begin{note}Note that the complex structure of $\mathcal{M}(R)$ viewed as the moduli space of flat connections changes when $R$ changes and so does the metric. It may seem that as long as we only care about complex structures the parameter $R$ can be absorbed by scaling $\varphi$ which is part of the $\mathbf{C}^{\times}$ action described in section 2. However  recall that the $\mathbf{C}^{\times}$ action simply moves the complex structure in the space of infinitely many compatible complex structures in which the moduli interpretation is the moduli space of flat connections. So the scaling perspective is not really natural as it identifies one complex structure with a particular twistor parameter on a hyperkahler space with one with a different twistor parameter on a different hyperkahler space.  In particular changing $R$ does not preserve the canonical one which is independent on the complex structure of the Riemann surface.\end{note}
\begin{note}Ultimately, the justification of this definition comes from our ability to show the  equivalence of two instanton correction problems (see section 9). It is quite satisfying to notice that the deformation we are using is $modular$ in the sense that it is a deformation of a moduli space to a family of moduli spaces instead of just a deformation of the underlying space of a moduli space. The deformation we are studying is natural in this sense.\end{note}
\begin{note}Unlike the situation in complex geometry here we have a differentiable family instead of a complex analytic family of complex manifolds. \end{note}

\section{Kontsevich And Soibelman's Wall Crossing Formula}

Now let us answer the fundamental question left  unanswered in section 5: what are discontinuous rays and associated discontinuous jumps of Gaiotto-Moore-Neitzke's coordinates?

For  the gauge charge lattice $\Gamma_{gau}\ni \gamma_{i}$, we define $$T:=\Gamma_{gau}^{*}\otimes \mathbf{C}^{\times}$$This is a local system over $B$ and its fiber over $u$ is a complex torus $T_{u}$. Write the coordinate associated (by dual pairing) to $\gamma_{i}$ as $X_{\gamma_{i}}(u)$. There is a canonically defined holomorphic symplectic form \begin{equation}\Omega^{T}:={1\over 2}\epsilon_{ij}{dX_{\gamma_{i}}\over X_{\gamma_{i}} }\wedge {dX_{\gamma_{j}}\over X_{\gamma_{j}} }\end{equation}
$\mathcal{X}_{\gamma}$ defined in section 5 labeled by gauge charges can be viewed as pull-backs of coordinates $X_{\gamma}$ by maps $$\mathcal{X}_{u}: \mathcal{M}_{u}\rightarrow T_{u}$$ which patch together to a global map $\mathcal{X}: \mathcal{M}^{*}\rightarrow T$. Here $\mathcal{M}_{u}$ is the fiber of Hitchin's fibration over $u$ and $\mathcal{M}^{*}$ is the complement of singular fibers in $\mathcal{M}$. In other words, $$\mathcal{X}_{\gamma}(\theta)=X_{\gamma}(\mathcal{X}(\theta))$$Here we are not suggesting that we can only define $\mathcal{X}_{\gamma}$ for $u$ outside the discriminant locus of the fibration. We simply want to say something about volumes of nonsingular fibers. The restriction to $\mathcal{M}^{*}$ of the holomorphic symplectic form $\Omega(\xi)$ defined in section 5 is the pullback of the canonical one \begin{equation}\Omega(\xi)={1\over 4\pi^{2}R}\mathcal{X}^{*}(\Omega^{T})\end{equation} This condition guarantees that the volume of a nonsingular fiber is $\sim ({1\over R})^{r}$ as described in the previous section. The holomorphic Poisson bracket on $T$ is also pulled back to a bracket on $\mathcal{M}$.

Locally one can introduce a $quadratic \ refinement$, i.e. $\sigma: \hat{\Gamma}\rightarrow \mathbf{Z}_{2}$ satisfying \begin{equation}\sigma(\gamma_{1})\sigma(\gamma_{2})=(-1)^{\langle\gamma_{1}, \gamma_{2}\rangle}\sigma(\gamma_{1}+\gamma_{2})\end{equation}
 For example (this is just an demonstration that a local quadratic refinement exists. It does not mean this is the one we are going to use) one can locally split the charge lattice into a sum of the gauge lattice and the flavor lattice and split the gauge lattice into a sum of electric and magnetic charge lattices such that the electric and magnetic charges form a symplectic basis of the skew-symmetric intersection paring. Then define $$\sigma(\gamma):= (-1)^{\gamma_{e}\gamma_{m}}$$ for a charge whose associated electric and magnetic parts are $\gamma_{e}$ and $\gamma_{m}$ respectively.
 \begin{note}If there is no global assignment of a quadratic refinement then we need to introduce some global twist as explained in \citep{G1} to make up for the change of the quadratic refinement. This will force us to rethink the moduli interpretation of the total space. Fortunately, as explained in section 8, there is an a prior assignment which is global.\end{note}
 \begin{dfn}Denote the infinitesimal symplectomorphism generated by the Hamiltonian $\sigma(\gamma)X_{\gamma}$ by $e_{\gamma}$ and define a set of symplectomorphisms known as Kontsevich-Soibelman transformations:$$K_{\gamma}:=\exp(\sum_{n=1}^{\infty}{1\over n^{2}}e_{n\gamma})$$ More explicitly, $K_{\gamma}$ is given by \begin{equation}K_{\gamma}: \mathcal{X}_{\gamma^{'}}\rightarrow \mathcal{X}_{\gamma^{'}}(1-\sigma(\gamma)\mathcal{X}_{\gamma})^{\langle\gamma^{'}, \gamma\rangle}\end{equation}\end{dfn}

\begin{dfn}The discontinuous rays of $\mathcal{X}_{\gamma}$ in the twistor $CP^{1}$ minus $\{0,\infty\}$ (also called twistor $\mathbf{C}^{\times}$ from now on) are called pre-BPS rays and are labeled by charges. They are defined by\begin{equation}l_{\gamma}:=\{\xi\mid Z_{\gamma}(u)/\xi\in \mathbf{R}_{-}\}\end{equation}The definition depends on $u$, so we talk about a pre-BPS ray at $u$.\end{dfn} Now we define the jumps associated with pre-BPS rays. Let \begin{equation}S_{l}:=\prod_{\gamma\in (\Gamma_{u})_{l}}K_{\gamma}^{\Omega(\gamma; u)}\end{equation} where $l$ is a pre-BPS ray and \begin{equation}(\Gamma_{u})_{l}:=\{\gamma\mid Z_{\gamma}(u)/\xi\in \mathbf{R}_{-},\forall \xi\in l\}\end{equation} $\Omega(\gamma; u)$ (which must be an integer) is called the BPS number or BPS invariant and is meant to be a virtual counting of BPS particles with charge $\gamma$ when the moduli parameter is $u$. Physicists define  BPS numbers as  indexes obtained by taking weighted traces of certain operators over Hilbert spaces obtained by geometric quantization, see \citep{CNV,CV1,CV2,DeM,DG1}. The classical spaces to be quantized should be the moduli spaces of the corresponding BPS particles (solutions of the corresponding elliptic equations) and as far as I know $\Omega(\gamma; u)$ has not been defined rigorously. This is an important problem but the investigation of this work does not depend on the solution of this welldefinedness issue. Rather we take our results as  constraints to be satisfied by any reasonable definition of BPS numbers. In section 8 we will give a mathematical operational definition of BPS numbers which is enough for this paper's purposes.

\begin{dfn}A BPS ray at $u$ is a pre-BPS ray labeled by a charge with $\Omega(\gamma; u)\neq 0$. A charge with $\Omega(\gamma; u)\neq 0$ is called a BPS charge at $u$. The collection of all BPS charges at $u$ is called the BPS spectra at $u$. \end{dfn}
\begin{note}In section 4 we used the term BPS spectra to mean charges of BPS particles in a four dimensional gauge theory whose quantum moduli space is the base of the Hitchin's fibration of the Hitchin's moduli space. It is expected that  the BPS spectra in that section coincide exactly with the BPS spectra defined here and this is true in known examples. \end{note}
\begin{note}This definition does not tell us the value of a BPS number except that it is zero or nonzero. The values of a nonzero spectra will be assigned in section 8. The hard part of the problem is to determine whether a BPS number is nonzero. Although we define BPS charges using BPS numbers in reality this is not the approach we use to compute them. We do not obtain BPS numbers before BPS charges. Often they are determined together by the wall crossing formula. \end{note}

$S_{l_{\gamma}}$ is the discontinuous jump associated to the pre-BPS ray $l_{\gamma}$, i.e.
\begin{itemize}\item  Let $\mathcal{X}_{\gamma^{'}}^{+}$ and $\mathcal{X}_{\gamma^{'}}^{-}$ be the limits of $\mathcal{X}_{\gamma^{'}}$ as $\xi$ approaches $l_{\gamma}$ clockwise and counterclockwise respectively. Discontinuities only appear across pre-BPS rays and the discontinuous jumps are\begin{equation}\mathcal{X}_{\gamma^{'}}^{+}=S_{l_{\gamma}}\mathcal{X}_{\gamma^{'}}^{-}\end{equation}In particular $\mathcal{X}_{\gamma}$ is continuous across $l_{\gamma}$ and therefore it is well defined on $l_{\gamma}$. In $S_{l_{\gamma}}$ of  (68) we use the values of $\mathcal{X}_{\gamma}$ on $S_{l_{\gamma}}$.\end{itemize}

Since these transformations preserve the canonical symplectic form and therefore its pullback $\Omega(\xi)$ we know $\Omega(\xi)$ is holomorphic across pre-BPS rays even though $\mathcal{X}_{\gamma}$ is only piecewise analytic. The order of the product of transformations defining $S_{l}$ is unimportant as long as $u$ is not on a stability wall because then $(\Gamma_{u})_{l}$ is at most one dimensional and hence all factors commute. We do not need to define $S_{l}$ for the exceptional cases when $u$ does lie on a stability wall. We only need  the limit of ordered products  as $u$ tends to the stability wall. The wall crossing formula below guarantees that the limit is well defined.

When there are more than one rays, we can take the ordered product of $S_{l_{i}}$ according to the counterclockwise order of phases of $l_{i}$. Of course this product does depend on the order of rays. We have to describe the dependence of the above construction on the moduli parameter $u$. For any $u$ the collection of pre-BPS rays of all charge can be ordered either counterclockwise or clockwise on the twistor $\mathbf{C}^{\times}$. Without loss of generality we take the counterclockwise order. The phase of a pre-BPS ray of a given charge depends continuously on $u$.  By the definition of stability walls the following is obvious.
\begin{pro}Let $u$ vary continuously along  a path. The order of the collection of all pre-BPS rays at $u$ changes if and only if some stability wall is crossed.\end{pro}
Now let us vary $u$ without crossing any stability walls. Since physically we do not expect any change of the BPS numbers $\Omega(\gamma;u)$  we demand that the ordered product is unchanged
\footnote{There are subtleties   in this statement, see remarks below the wall crossing formula}. In fact we impose a seemingly stronger\footnote{In fact, it is not stronger because the factorization of a given composition into Kontsevich-Soibelman factors is unique if the order of pre-BPS rays is fixed.} constraint by insisting that $\Omega(\gamma; u)$ remains constant. However if we vary $u$ by letting it cross a stability wall, then the order of pre-BPS rays changes. But   we do not want  the ordered
product to change (see section 8 for the reason). This is a powerful statement as it not only means that $\Omega(\gamma; u)$ will change but also determines recursively  values of all $\Omega(\gamma; u)$ on one side of the stability wall if all $\Omega(\gamma; u)$ are given on the other side. This is Kontsevich and Soibelman's wall crossing formula proposed in \citep{KS}\footnote{But note that the context in \citep{KS} is somewhat different from ours. The construction in this paper does not depend on any motivic constructions or stability conditions of Bridgeland type.}. It is the last and the most subtle property that we require our coordinates to satisfy.
 \begin{itemize}\item The discontinuous jumps of $\mathcal{X}_{\gamma}$ must satisfy the following formula. \\
\textbf{Wall Crossing Formula}\ \ Fix two generic phases $(\theta_{-}< \theta_{+})$ on the twistor plane and define \begin{equation}S(\theta_{-},\theta_{+}; u)=\prod^{\leftarrow}_{\theta_{-}< \arg(Z_{\gamma}(u)) <\theta_{+}}K_{\gamma}^{\Omega(\gamma; u)}\end{equation} where $\prod^{\leftarrow}_{\theta_{-}< \arg(Z_{\gamma}(u)) <\theta_{+}}$ means it is a counterclockwise ordered product. We assume that the difference of the two phases is not greater than $\pi$ because any half plane already    contains the full information of BPS rays by symmetry. For any two $u, u^{'}$ not on any stability walls, suppose  that they can be connected by a path not
crossing branch cuts (but it could cross stability walls) such that along the path from   $u$ to $u^{'}$  no BPS ray crosses $\theta_{-}$ or $\theta_{+}$, then
\begin{equation}S(\theta_{-},\theta_{+}; u)=S(\theta_{-},\theta_{+}; u^{'})\end{equation}\end{itemize}
\begin{note}Although the order of BPS rays does not change if no stability wall is crossed, the set of BPS rays over which the product in (69) runs could change if $u$ is varying even if it does not cross any stability walls. For example, suppose there are two BPS charges $\gamma_{1},\gamma_{2}$ inside a stability wall
and the phase of $l_{\gamma_{2}}$ is smaller than the phase of $l_{\gamma_{1}}$. Then even if we stay inside the stability wall, it is possible that for $u$, there
are two BPS rays $l_{\gamma_{2}}, l_{\gamma_{1}}$ between $\theta_{-}$ and $\theta_{+}$ while for $u^{'}$ the two rays are instead $l_{\gamma_{1}}, l_{-\gamma_{2}}$.
The wall crossing formula certainly does not claim that $K_{\gamma_{1}}K_{\gamma_{2}}=K_{-\gamma_{2}}K_{\gamma_{1}}$. That is why we require that  no
BPS ray crosses $\theta_{-}$ or $\theta_{+}$. So if some BPS rays do cross $\theta_{-}$ or $\theta_{+}$ we just use another wall crossing formula.\end{note}
\begin{note}Remember there are nontrivial monodromies on $B$. If we connect $u$ and $u^{'}$ by a path then the path must avoid the branch cuts. So even if $u$ and $u^{'}$  are in the same chamber divided by stability walls we may not be able to connect them inside the chamber. To connect them we travel to other chambers and travel back by crossing stability walls. Along the path we may have to switch to other wall crossing formulas because some BPS rays might cross  $\theta_{-}$ or $\theta_{+}$. Of course  in the end the wall crossing formulas must respect (they do) the monodromy in the sense that the BPS spectra computed from wall crossing formulas must reproduce those determined by the monodromy action. \end{note}

\begin{note}An ordered product in the wall crossing formula is usually infinite and even infinite from both directions towards some elements in the middle. To understand its meaning we truncate the product according to the degree of charges successively. Any charge $\gamma$ can be written as a nonnegative\footnote{It is not very obvious that this is always guaranteed. So we might have to make it as an additional assumption. In the geometric construction of Gaiotto-Moore-Neitzke's coordinates in the next two sections this condition holds.  This positivity condition is a part of an interesting topic, see \citep{G3} and section 8.} (here we need $\theta_{+}- \theta_{-}\leq\pi$) linear combination of  a basis $(\gamma_{i}),1\leq i\leq r$ of $\hat{\Gamma}$ and  the degree of $\gamma$ is the sum of coefficients in the linear combination. At each stage of the truncation, the product is a finite one and the total ordered product is defined in the sense of projective limit. More precisely, we consider the algebra $F:=\mathbf{C}[\mathcal{X}_{\gamma_{1}},\cdots\mathcal{X}_{\gamma_{r}}]$ generated by $\mathcal{X}_{\gamma_{i}}$ viewed as formal variables.  We take the filtration of $F$ by ideals $I_{N}$ generated by monomials whose degrees are higher than $N$. Kontsevich-Soibelman transformations with $\theta_{-}< \arg(Z_{\gamma}(u)) <\theta_{+}$ generate a group $G_{N}(\theta_{-}, \theta_{+}; u)$ of Poisson automorphisms of $F_{N}:=F/I_{N}$ together with projections $G_{N}\rightarrow G_{N-1}$. Define $G(\theta_{-}, \theta_{+}; u)$ to be the associated projective limit and the product is in it. So we have a projective system of groups and on each projection the product is a finite one.\end{note}
\begin{note}Although the factorization of an ordered product is understood  in the projective sense this does not mean that the ordered product itself as a symplectomorphism has to be viewed in that way. On the contrary usually it is an innocent looking simple transformation. It is just that factorizing it modulo higher and higher degrees with a given order of pre-BPS rays forces us to add more and more Kontsevich-Soibelman factors.\end{note}
\begin{note}It is an algebraic fact that once the order of BPS rays is given the decomposition into Kontsevich-Soibleman factors of a given element $S$ is unique, see \citep{G1,G2,GS1}. In fact one can  make the factorization by truncations by the degree of charges. In each degree there are finitely many pre-BPS rays with an well defined order. Expand and compare both sides up to that degree one can determine the exponents of associated Kontsevich-Soibelman transformations. Moving to the next degree would introduce new factors. The algebraic procedure is actually identical to the algorithm of Kontsevich-Soibelman's theorem in section 3 because later we will identify it as a wall crossing formula. \end{note}
\begin{note}Some authors also call equation (68) describing discontinuous jumps  the wall crossing formula which may cause some confusions.  In fact  equation (68) does not involve any stability walls  and neither does it involve any walls in the sense of section 3 because walls and  BPS rays are in different spaces (although they are related, see section 9). In this paper we simply call (68) the discontinuous jumps  and the term "wall crossing" is used only when there is a stability wall being crossed.\end{note}
\begin{dfn}A set of coordinates $\mathcal{X}_{\gamma}$ satisfying all the required properties listed in this and the previous section is called Gaiotto-Moore-Neitzke coordinates. \end{dfn}

This completes our description of properties of Gaiotto-Moore-Neitzke coordinates and now the task is to build them. Finding these coordinates can be interpreted as solving an infinitely dimensional Riemann-Hilbert problem. In other words, we are to find the map $\mathcal{X}$ with prescribed discontinuous jumps and asymptotic behaviors for $\xi\rightarrow 0,\infty$. In \citep{G1}, Gaiotto, Moore and Neitzke  advocated an integral equation approach inspired by the classical  treatment of finite dimensional Riemann-Hilbert problems. It is straightforward to check that any solution $\mathcal{X}_{\gamma}(\xi)$ of the following equation has the required discontinuities across BPS rays $$\mathcal{X}_{\gamma}(\xi)=\mathcal{X}_{\gamma}(\xi)^{sf}\exp({1\over 4\pi i}\sum_{l}\int_{l}{d\xi^{'}\over \xi^{'}}{\xi^{'}+\xi\over \xi^{'}-\xi}\log{\mathcal{X}_{\gamma}(\xi^{'})\over (S_{l}\mathcal{X})_{\gamma}\xi^{'}})$$where the sum runs over all BPS rays $l$. The integral equation formulation enables them to set up an iterative approximation scheme from which the exponentially smallness of instanton corrections follows with some additional assumptions. Moreover, the integral equation above has an interpretation as the basic ansatz of an integrable field theory. This is known as the thermodynamic Bethe ansatz and has been investigated by many people in recent years, see e.g. \citep{A1,A2,Hat,Z}. Thermodynamic Bethe ansatz is also directly related to the cluster algebra structure to be formulated later. All these facts tell us that there is an integrable model hidden in the structure of the hyperkahler metric (in the instanton corrected form) of a Hitchin's moduli space (do not confuse it with the algebraically integrable system given by the Hitchin's fibration). This is clearly a promising direction but it will not be pursued further in this paper. Instead we will build these coordinates directly in a geometric way following Gaiotto, Moore and Neizrke and verify the required properties.\\

In the rest of this section let us take a closer look at the wall crossing formula. As an example, we point out the following formula \begin{equation}K_{2,-1}K_{0,1}=K_{0,1}K_{2,1}K_{4,1}\cdots K_{2,0}^{-2}\cdots K_{6,-1}K_{4,-1}K_{2,-1}\end{equation}
The charge lattice here has rank two and is identified with $\mathbf{Z}^{2}$ with the integral pairing $\langle(p, q), (p^{'},q{'})\rangle= pq^{'}-qp^{'}$. The proof of it will be discussed in section 9.

A remarkable thing happens \citep{DG2, G1}. If we go back and check the spectra of the pure $SU(2)$ theory given in section 4, we would recognize that the charges  on the left hand side of equation (71) are precisely the strong coupling spectra of BPS particles while the right hand side provides exactly the weak coupling spectra. Reversing signs gives the spectra  of antiparticles. This formula suggests that the BPS numbers of all particles with nonzero magnetic charges are all one while for the so-called "W-boson" with charge $(2,0)$ it is $-2$. The production of the BPS spectra is not a coincidence as similar wall crossing formulas give exact BPS spectra of some other gauge theories.

Now we have come to the starting point of the author's  whole investigation: another surprise! It turns out the log morphisms of scattering diagrams in section 3 are Kontsevich-Soibelman transformations after we identify  integral slopes of rays with charges and the consistency condition of scattering diagrams in Kontsevich Soibelman theorem in that section $$\theta_{s}\circ\theta_{s-1}\circ\cdots\circ\theta_{1}=\mathrm{Id}$$ is nothing but a wall crossing formula. This fact has been recognized in \citep{GPS,GP} and the authors defined the notion of a tropical vertex group as automorphisms of the torus $\mathbf{C}^{\times}\times \mathbf{C}^{\times}$ preserving the holomorphic symplectic form $\omega={dx\over x}\wedge{dy\over y}$. Log morphisms in a universal scattering diagram defined in section 3 are elements of a tropical vertex group. For example the equation (71) arises (after a change of variables and setting $t=1$) as the result of running Kontsevich and Soibelman's algorithm for the following scattering diagram$$D:=\{(\mathbf{R}(1,0), (1+tx^{-1})^{2}),(\mathbf{R}(0,1), (1+tz^{-1})^{2})\}$$See the second example in section 9.

\begin{note}If we want to identify the consistency condition as a wall crossing formula, there should be an explanation of the role of  stability walls in this construction. In Gross and Siebert's work it is not clear what a stability wall  means. The identification of the wall crossing formula and   the consistency condition is only   verified at the algebraic level. This is an important issue that we will clarify later.\end{note}

So it seems that we have three  problems involving instantons or BPS particles\begin{itemize}\item Determination of BPS spectra in four dimensional gauge theories.\item Instanton corrections of hyperkahler metrics of Hitchin's moduli spaces in three dimensional gauge theories.\item Instanton corrections of complex structures (in the form of explicit algebraic deformations of defining equations) of Calabi-Yau varieties in mirror symmetry. \end{itemize} In each of these problems, the spectra of instantons (or BPS particles) exhibit the same kind of wall crossing formulas (although for the third problem the meaning of $instantons$ and stability walls is not yet clear at the moment).  The relation between the first and the second problem is perhaps not surprising (although I believe it is not completely understood yet). After all the three dimensional theory is obtained by the four dimensional theory wrapping a circle of radius  $R$ (which is also the deformation parameter of the second problem). And it seems reasonable to expect that  instantons in the three dimension that contribute to the metric can be obtained as dimensional reductions of monopoles and dyons in the four dimension. The main result in this paper, however, is the equivalence of the second and the third problem. Note that the formulations and even the languages used in these two problems are quite different. One is a differential geometric problem while the other is completely algebraic. So this relation looks more mysterious. To make it less so we will continue our journey and introduce some geometric objects which connect both sides.

\section{Fock-Goncharov Coordinates And Quadratic Differential Foliations}
In this section a geometric construction of a set of Gaiotto-Moore-Neitzke coordinates on the Hitchin's moduli space is given. The exposition follows \citep{G2} which is based on \citep{FG1,St}.

Let $\mathcal{M}$ be our Hitchin's moduli space. We view it as a moduli space of flat connections. First, let assume there are only regular singularities $P_{i}, 1\leq i \leq l$ with regular semisimple residues ($T_{1}$). We assume $l\geq 1 $ in general and $l> 3$ if $g=0$.
\begin{dfn}Choose a triangulation  of the Riemann surface $C$ with all vertices at singularities. Let $M_{i}$ be the clockwise monodromy of flat sections around $P_{i}$. Define a decoration at $P_{i}$ to be a choice of one of the two flat eigenlines of $M_{i}$. Denote such a decorated triangulation by $T$. For an edge $E$ of $T$, we consider the two triangles bounding $E$ making up a quadrilateral $Q_{E}$ with four vertices $P_{i},1 \leq i \leq 4$ in the counterclockwise order and $E$ connecting 1 and 3. Define the Fock-Goncharov  coordinate $\mathcal{X}_{E}^{T}$ by$$\mathcal{X}_{E}^{T}:=-{(s_{1}\wedge s_{2})(s_{3}\wedge s_{4})\over(s_{2}\wedge s_{3})(s_{4}\wedge s_{1})}$$where $s_{i}$ is an element of the one dimensional  decoration at $P_{i}$ (so it is defined up to a scaling and our definition is invariant under this scaling). Since $Q_{E}$ is simply connected, $s_{i}$ can be chosen to be single-valued in it and the four eigensections are evaluated at a common point $P_{*}$ inside the quadrilateral. The value is independent of the choice of the evaluation point because it is the $SL(2, \mathrm{C})$ invariant cross ratio.\end{dfn}

\begin{note}$\mathcal{X}_{E}^{T}$ is well defined on the Zariski open set which is the complement of the locus defined by the vanishing of the denominator. It is not hard to show that it is a holomorphic coordinate on this open subset. The set of all such functions where $E$ runs over all edges of a fixed decorated triangulation is a complete set of coordinates. Moreover outside the codimension one locus where either the numerator or the denominator is zero $\mathcal{X}_{E}^{T}$ is nonzero. So we have a set of locally defined $\mathbf{C}^{\times}$ valued functions.\end{note}

If we change the decorated triangulation, then $\mathcal{X}_{E}^{T}$ changes. Any two decorated triangulations can be connected by a composition of two elementary transformations\begin{itemize}\item Flip at an edge. This means that we replace $E=E_{13}$ which connects vertices 1 and 3 in the quadrilateral containing $E$ as a diagonal edge by $E^{'}=E_{24}$ connecting vertices 2 and 4 and obtain a new triangulation $T^{'}$. \item Pop at a vertex. This simply means that we use the other possible choice of decoration at that vertex.\end{itemize}

For flips the transformations of Fock-Goncharov coordinates are $$\mathcal{X}_{E}^{T}=(\mathcal{X}_{E^{'}}^{T^{'}})^{-1}$$ and
$$\mathcal{X}_{E_{12}}^{T^{'}}=\mathcal{X}_{E_{12}}^{T}(1+\mathcal{X}_{E}^{T})$$$$ \mathcal{X}_{E_{23}}^{T^{'}}=\mathcal{X}_{E_{23}}^{T}(1+(\mathcal{X}_{E}^{T})^{-1})^{-1}
$$$$ \mathcal{X}_{E_{34}}^{T^{'}}=\mathcal{X}_{E_{34}}^{T}(1+\mathcal{X}_{E}^{T})$$
$$\mathcal{X}_{E_{41}}^{T^{'}}=\mathcal{X}_{E_{41}}^{T}(1+(\mathcal{X}_{E}^{T})^{-1})^{-1}$$We define $\langle E,E^{'}\rangle$ to be the number of faces $E$ and $E^{'}$ have in common counted with signs. The sign is positive (negative) if $E$ comes immediately before $E^{'}$ in counterclockwise (clockwise) order going around the common face. Clearly $|\langle E,E^{'}\rangle|\leq 2$. The above four equations are actually for $\langle E_{ij},E\rangle=\pm1$. In general the transformations are $$\mathcal{X}_{E}^{T}=(\mathcal{X}_{E^{'}}^{T^{'}})^{-1}$$ and \begin{equation}\mathcal{X}_{E_{ij}}^{T^{'}}=\mathcal{X}_{E_{ij}}^{T}(1+(\mathcal{X}_{E}^{T})^{-sgn(\langle E_{ij},E\rangle)})^{-\langle E_{ij},E\rangle}\end{equation}

A pop at a point $P$ can be decomposed as a composition of the following transformations: first flip all incident edges at $P$ except one. This produces a degenerate triangulation. Then pop at $P$ for this degenerate triangulation and finally flip all the flipped edges back.

Here a degenerate triangulation means that  two edges in a triangle are identified. So we have a double vertex  and the edge connecting the double vertex to itself is a loop while the double edge is an edge connecting a point on the loop (the double vertex) to another vertex $P$. To define Fock-Goncharov coordinates in this situation, we take a cover ramified at $P$ such that after taking the pull-back the triangulation is non-degenerate. Such a cover always exist and is non-unique but our definition does not depend on the choice of the cover. We pull back everything and define the Fock-Goncharov  coordinate for the degenerate edge $E$ to be the ordinary Fock-Goncharov  coordinate $\mathcal{X}_{\tilde{E}}^{T}$ on the cover where $\tilde{E}$ is any choice of the pre-images of $E$ and the definition does not depend on this choice. We need degenerate triangulations not just for pops. In fact we will construct decorated triangulations from a foliation later and that might give us a degenerate decorated triangulation. The dimensional count of Fock-Goncharov coordinates is still valid for degenerate decorated triangulations.

 Back to the problem of pops. We only need to know the transformations associated to the pop at $P$ for  degenerate triangulations and they are $$\mathcal{X}_{E}^{T^{'}}=(\mathcal{X}_{E}^{T})^{-1}$$ \begin{equation}\mathcal{X}_{E^{'}}^{T^{'}}=\mathcal{X}_{E}^{T}\mathcal{X}_{E^{'}}^{T}\end{equation} where $E$ is the degenerated edge, $E^{'}$ is the other edge in the triangle and $T^{'}$ is the new decorated degenerate triangulation obtained after the pop.

We may need to consider a process with infinitely many flips and take its limit. Consider an annular domain on the Riemann surface whose outer circle contains one singularity $P$ and inner circle contains another one $P^{'}$. We will define a sequence of triangulations. First choose two paths $E_{0\pm}$ connecting $P$ and $P^{'}$ and pointing from $P$ to $P^{'}$ such that their difference has counterclockwise winding number one around the inner circle. The two paths  form a part of a  possibly degenerate decorated  triangulation $T_{0}$. We define $T_{n+1}$ inductively by flipping $E_{n-}$ and define $E_{(n+1)+}$ to be the flip of $E_{n-}$ and $E_{(n+1)-}:=E_{n+}$. Let $\mathcal{X}_{E_{n\pm}}^{T_{n}}$ be the associated Fock-Goncharov coordinates, we then define limit Fock-Goncharov coordinates for an ideal "limit" triangulation by letting $n\rightarrow\infty$ $$\mathcal{X}_{A}^{T_{+\infty}}:= \lim_{n\rightarrow\infty}\mathcal{X}_{E_{n+}}^{T_{n}}\mathcal{X}_{E_{n-}}^{T_{n}}$$
 \begin{equation}\mathcal{X}_{B}^{T_{+\infty}}:= \lim_{n\rightarrow\infty}(\mathcal{X}_{E_{n+}}^{T_{n}})^{-n}(\mathcal{X}_{E_{n-}}^{T_{n}})^{1-n}\end{equation} One can also define $T_{-n}$ by flipping $E_{n+}$
and define $\mathcal{X}_{E_{A}}^{T_{-\infty}}$ and $\mathcal{X}_{E_{B}}^{T_{-\infty}}$ by taking limits$$\mathcal{X}_{A}^{T_{-\infty}}:= \lim_{n\rightarrow -\infty}\mathcal{X}_{E_{n+}}^{T_{n}}\mathcal{X}_{E_{n-}}^{T_{n}}$$\begin{equation}\mathcal{X}_{B}^{T_{-\infty}}:= \lim_{n\rightarrow -\infty}(\mathcal{X}_{E_{n+}}^{T_{n}})^{-n}(\mathcal{X}_{E_{n-}}^{T_{n}})^{1-n}\end{equation}
These limits exist and have been written down explicitly in \citep{G2}.

To relate the above construction to the instanton problem of Hitchin's moduli spaces (now considered as moduli spaces of flat connections), we must go back to the Hitchin's fibration.   The determinant of $\varphi$ is a quadratic differential $-\lambda^{2}$ well defined on $C$ and therefore we take the holomorphic coordinate $u$ of the base $B$ of the Hitchin' fibration (recall that $B$ can be identified as  the space of quadratic differentials) to be\footnote{In fact, the choice of $u$ in reality can be slightly different from this. Usually $\lambda^{2}$ has a fixed part and a moving part. The moving part is parameterized by $u$. Of course this does not change the discussion given below.} $u= -\lambda^{2}$. $\lambda^{2}$ has order two poles and generically has only simple zeroes. $\lambda$ is a one form  on the Riemann surface $C$ defined up to a sign but is a single valued one form on the spectral curve $S$ which is a double cover of $C$. It is  the Seiberg  Witten differential defined in section 2.

Fix an angular parameter $\vartheta$ and consider the foliation given by trajectories of $\lambda^{2}$ with phase $\vartheta$.
\begin{dfn}A trajectory of $\lambda^{2}$ with phase $\vartheta$ is a curve whose tangent vector $\partial_{t}$ satisfies $$\langle\lambda, \partial_{t}\rangle\in e^{i\vartheta}\mathrm{R}^{\times}$$ everywhere on the curve.\end{dfn}
\begin{note}A trajectory is called a WKB curve in \citep{G2}.\end{note}
There is an extensive theory of foliations given by meromorphic quadratic differentials  and the local behaviors near singularities and zeroes as well as global behaviors are known. The standard reference is Strebel' book \citep{St}. Let us summarize the results we need.

Near a point which is neither a zero nor a pole of the quadratic differential, we can straighten the foliation by choosing local coordinate $w:=\int\lambda$. Locally near an order $n$ zero, we can choose a local parameter $\zeta$ such that $\lambda^{2}$ has the representation \begin{equation}\lambda^{2}=({n+2\over 2})^{2}\zeta^{n}d\zeta^{2}\end{equation}The full angle $0\leq\arg\zeta\leq2\pi$ is divided into $n+2$ equal sectors. In particular for a simple zero (i.e. order one zero) the foliation develops three asymptotic directions surrounding and going away from the zero.

Since we only have regular singularities for the Hitchin's equations the order of poles of $\lambda^{2}$ is two. It is shown that in this case it has a local representation of the form $$\lambda^{2}={a\over\zeta^{2}}d\zeta^{2}$$ and the trajectories near the pole is either logarithmic spirals approaching the pole or radii approaching the pole or closed circles around the pole.

 Globally a trajectory belongs to one of the following cases\footnote{This classification is valid for quadratic differentials with higher order poles.}\begin{itemize}\item A generic trajectory. It is asymptotic in both directions to singular points. Generic trajectories arise in one dimensional families.\item A separating trajectory.  It is asymptotic in one direction to a simple zero and in the other direction to a singular point. Separating trajectories separate families of generic trajectories. \item A finite trajectory. It is asymptotic in both directions to a simple zero (both directions could go to the same zero) or  is closed. A finite trajectory is also called a critical trajectory.\item A divergent trajectory. It is neither closed nor approaches to a limit in one or both directions.\end{itemize}

For a generic $\vartheta$,  finite trajectories are absent and in that case Gaiotto, Moore and Neitzke showed the absence of divergent ones in out setting. We will assume the absence of finite trajectories for now to get decorated triangulations and later we will see they are instantons in the instanton correction problem.

Following Gaiotto, Moore and Neitzke, we define a decorated triangulation called WKB triangulation in the following way. \begin{dfn}We take a generic $u=-\lambda^{2}$ such that it has only simple zeroes (note that this is the generic case, the nongeneric ones are codimensional two in $B$). We choose one element from every family of generic trajectories separated by separating trajectories. They make  an ordinary triangulation of the Riemann surface. The choices of the representatives are unimportant because a triangulation is only meant to be defined up to isotopy.

In general  we want to consider quadratic differentials with higher order poles. There might be  a generic trajectory approaching the same singularity along both directions (see theorem 7.1) in which case we get a degenerate triangulation.

Near each singularity, there are two independent eigen flat sections. It is shown that one of them is exponentially small along trajectories going to the singularity while the other is exponentially large. We pick the small flat section as the decoration at the singularity. These decorations together with the triangulation define a decorated triangulation $T_{WKB}(\vartheta, \lambda^{2})$ called a WKB triangulation and therefore a set of Fock-Goncharov coordinates $\mathcal{X}_{E}^{T_{WKB}(\vartheta, \lambda^{2})}$.\end{dfn}
\begin{note}The variable $u=-\lambda^{2}$ plays two different roles in the theory. As the holomorphic coordinate of the base of the Hitchin's fibration it is a part of a set of coordinates of the total space, i.e. the Hitchin's moduli space $\mathcal{M}$ itself. But this is only implicit in the definition of Fock-Goncharov coordinates over $\mathcal{M}$ because when we define these coordinates we view the moduli space  as the moduli space of flat connections in which the Hitchin's fibration is not holomorphic. Since we are not splitting the coordinates of $\mathcal{M}$ into $(u, \theta)$\footnote{Our notations in section 5  do have such a splitting but the formulation of Gaiotto-Moore-Neitzke ansatz does not require that. We used the splitting before because we wanted to discuss the semiflat metric and the large $R$ asymptotic.} according to its fibration structure it does not seem to be easy to analyze what happens with the metric constructed in this way near the singular fibers of the Hitchin's fibration. See section 9.4 for some further discussions about the metric. On the other hand, the "moduli" parameter $u$  is an additional parameter (instead of a coordinate) of the descriptions of the moduli space of flat connections because  the WKB triangulations  depend on it. The dependence is locally constant unless some wall crossing\footnote{The precisely meaning of a wall crossing here should be a BPS wall crossing instead of a stability wall crossing, see section 9.} occurs. Therefore for the moduli parameter $u$ it is enough to pick one representative  in each chamber bounded by  BPS walls\footnote{See section 9 for the definition of BPS walls.}. That is why we   need not to consider non-generic quadratic differentials with nonsimple zeroes which form a codimensional two locus in $B$. \end{note}
\begin{note}The values of $\mathcal{X}_{E}^{T_{WKB}(\vartheta, \lambda^{2})}$ of course depend on parameters $\vartheta$, but it is the algebraic relations between some limit values of these coordinates  that will concern us when we formulate the main theorems later because we will write the defining equations of the moduli space of flat connections in terms of them.\end{note}

We can vary the phase $\vartheta$. Although for a generic $\vartheta$ there are no finite trajectories they do appear for exceptional values of $\vartheta$. If we label $\vartheta$ by rays on the complex plane then when $\vartheta$ crosses countably many exceptional rays, the decorated triangulation will change and in this way we can reproduce the transformations (jumps) of triangulation described before. \begin{itemize}\item A flip. An edge is flipped when an exceptional ray of $\vartheta$ is crossed. As $\vartheta$ goes to the exceptional ray the flipped edge degenerated to a finite trajectory connecting two simple zeroes. There is only one such finite trajectory for that exceptional value of $\vartheta$.\item A pop for a degenerate triangle. In this case as $\vartheta$ goes to the exceptional ray trajectories approaching the pole as logarithmic spirals degenerate to closed trajectories around the pole  bounded by  a finite trajectory with both directions going to the same simple zero.  The two ways of going to the pole through logarithmic spirals account for the two possible choices of decorations.\item Infinitely many flips leading to a limit configuration. In this case, the phase first passes infinitely many rays (corresponding to  flips) to reach a special configuration \footnote{The procedure of labeling by charges  described below still works in this case. See \citep{G2} for details.}with  closed trajectories  around the pole inside an annular region bounded by two finite boundary trajectories. Each boundary trajectory's  both directions go to a same simple zero (but the simple zeroes are different for the two trajectories). We can take the limit Fock-Goncharov coordinates $\mathcal{X}_{A}^{T_{+\infty}},\mathcal{X}_{B}^{T_{+\infty}}$. On the other hand, if we start from the other side of the special ray and approach it from the other direction we would also pass infinitely many rays and get limit Fock-Goncharov coordinates $\mathcal{X}_{A}^{T_{-\infty}}, \mathcal{X}_{B}^{T_{-\infty}}$. The transformations from $\mathcal{X}_{A}^{T_{+\infty}},\mathcal{X}_{B}^{T_{+\infty}}$ to $\mathcal{X}_{A}^{T_{-\infty}}, \mathcal{X}_{B}^{T_{-\infty}}$ which can be easily written down are then the "discontinuous" jumps between two ways of approaching the special ray and this operation is called a juggle. Juggles will be discussed  further later. \end{itemize} For a generic quadratic differential, the above list has exhausted all possible jumps of WKB triangulations.

A quadratic differential with only simple zeroes is also generic. But the subset of such generic quadratic differentials do not coincide with the subset of generic quadratic differentials with only three types of jumps for exceptional $\vartheta$ listed above. In fact  there are quadratic differentials with only simple zeroes but having other types of jumps.

To understand jumps of WKB triangulations for such quadratic differentials it is instructive to see how the above three types are obtained. For a jumps to occur,  finite trajectories have to appear. If we have a finite trajectory connecting two simple zeroes then we get a flip. If we have a closed trajectory the it comes with a family of closed trajectories. According to the local classification any member of this family cannot encounter a pole. But the family can surround and contract to a pole. Clearly globally the family is bounded from the other direction by other types of trajectories. They cannot be bound by a divergent one. Otherwise the closure of the divergent trajectory (which must be recurrent) is bounded by the outmost closed trajectory in the family. But this possibility is excluded in \citep{St}. The family cannot be bounded by generic or separating trajectories according to the local classification. So the only possibility is that it is a closed loop connecting several simple zeroes. Generically the boundary only meets one leading to the second type of jumps. There is also another possibility. The family of closed trajectories can be bounded from both directions. Generically both boundaries contain only one simple zeroes leading to the third type of jumps. Non-generically, the boundary or boundaries of the family could meet several simple zeroes which means we also have several flips. Later we will associate Kontsevich-Soibelman transformations to jumps of WKB triangulations and this scenario corresponds to a product of Kontsevich-Soibelman transformations associated to flips, pops and juggles since it can be decomposed into a composition of flips, pops and juggles. So even for non-generic cases as far as the wall crossing formula is concerned it is enough to consider the above three types.

\begin{note}We do not have to worry about that when finite trajectories appear the global behavior of the foliation could be wild due to the possible existence of divergent trajectories. Because we do not define Fock-Goncharov coordinates for these exceptional values of $\vartheta$. We only care about the discontinuity of limits of Fock-Goncharov coordinates defined for non-exceptional $\vartheta$.\end{note}

We will explain that the Fock-Goncharov coordinates constructed from  WKB triangulations give us a geometric realization of Gaiotto-Moore-Neitzke coordinates defined in the previous sections.

Gaiotto-Moore-Neitzke's coordinates are labeled by charges, so we need to label Fock-Goncharov coordinates by charges instead of edges. It is easy to see that every triangle in a WKB triangulation contains exactly one simple zero. Let $E$ be the edge labeling the Fock-Goncharov coordinate $\mathcal{X}_{E}^{T}$. We choose an oriented simple loop inside $Q_{E}$ surrounding the two zeroes in the two adjacent triangles and define the associated charge $\gamma_{E}$ to be the lift of the loop to the spectral curve $S$ which is a double cover of the underlying Riemann surface. Ambiguity of the sign of the cycle induced by ambiguities  of choosing orientations and one of the two sheets can be canonically fixed in the following way. Note  that $\lambda$ is a single-valued one form over $S$. We require that the positively oriented tangent vector $\partial_{t}$ of the lift of $E$ to $S$ denoted as $\hat{E}$ satisfies $$e^{-i\vartheta}\langle\lambda,\partial_{t}\rangle>0$$ The sign of the cycle $\gamma_{E}$ is fixed by $\langle\gamma_{E},\hat{E}\rangle=1$. So we can replace the labeling by $E$ by labeling by $\gamma_{E}$ (also denoted simply as $\gamma$ later in this paper). This operation respects the integral skewsymmetric pairings, i.e. $\langle\gamma_{E},\gamma_{E^{'}}\rangle=\langle E,E^{'}\rangle$.

It is easy to generalize to degenerate edges. Recall that a degenerate triangle appears when we have a generic trajectory connecting a pole to itself. In this case that pole is the double vertex and the generic trajectory is the loop edge. The double edge connects the double vertex to another pole denoted by $P$. There is a simple zero inside the degenerate triangle and two of the three separating trajectories starting from the simple zero end at the double vertex while the third one ends at the other pole $P$ (see figure 26 of \citep{G2}). Just like what we did when we defined Fock-Goncharov coordinates for degenerated triangulation we use covers to separate degenerate edges and define labeling charges. In the end the charge associated to the double edge $E$ is a loop around $P$ (so it is a flavor charge) while the charge associated to the loop edge $E^{'}$ is induced by a loop circling the simple zero inside the degenerate triangle and a simple zero outside.

It is shown that one can generate the charge lattice $\hat{\Gamma}$ of $S$ by  cycles associated to edges. We extend the definition of Fock-Goncharov coordinates to the whole lattice by the multiplicative relation $$\mathcal{X}_{\gamma_{E}}\mathcal{X}_{\gamma_{E^{'}}}=\mathcal{X}_{\gamma_{E}+\gamma_{E^{'}}}$$It is easy to show that our Fock-Goncharov coordinates satisfy the reality condition and with respect to the holomorphic Poisson bracket on the moduli space $\mathcal{M}$ $$\{\mathcal{X}_{\gamma}, \mathcal{X}_{\gamma^{'}}\}= \langle\gamma, \gamma^{'}\rangle\mathcal{X}_{\gamma}\mathcal{X}_{\gamma^{'}}$$ By WKB analysis, Gaiotto, Moore and Neitzke showed the required asymptotics for $\xi\rightarrow\infty$ \begin{equation}\mathcal{X}^{\vartheta}_{\gamma} \sim c_{\gamma} \exp({\pi R\over \xi}Z_{\gamma})\end{equation} where $c_{\gamma}$ is a constant. The proof of this asymptotic is by finding the WKB approximation of small flat sections along each edge and plugging into the definition of Fock-Goncharov coordinates. The proof also guarantees that for large $R$ (large enough such that the deviation to the WKB approximation is small enough) small flat sections on both end of an edge do not coincide (because this is the case of the WKB approximation) and hence the Fock-Goncharov coordinates are pole free in a neighborhood of any given ray in the twistor $\mathbf{C}^{\times}$. So we get \begin{itemize}\item For large enough $R$ Fock-Goncharov coordinates are piecewise holomorphic over twistor $\mathbf{C}^{\times}$.\end{itemize}Since we are only interested in large complex (large $R$) limit, this is good enough for us. Gaiotto, Moore and Neitzke also showed that $\mathcal{X}_{\gamma}^{\vartheta}$ has the required large $R$ asymptotic.

One also wants to know how to interpret $\vartheta$ in the twistor language. In the twistor $\mathbf{C}^{\times}$ we define \begin{equation}H_{\vartheta}:=\{\xi\mid \vartheta-\pi/2 < \arg\xi <\vartheta+\pi/2 \}\end{equation} The asymptotic of $\mathcal{X}^{\vartheta}_{\gamma}$ is actually for $\xi\rightarrow\infty$ within the half plane $H_{\vartheta}$. Instead of discussing the discontinuous jumps by varying $\xi$, one can equivalently discuss  the discontinuous jumps by varying $\vartheta$. Suppose we have a set of Gaiotto-Moore-Neitzke coordinates, we just define $\mathcal{X}_{\gamma}^{\vartheta}$ to be the analytic continuation of the Gaiotto-Moore-Neitzke coordinate $\mathcal{X}_{\gamma}$ starting from the central ray of the half plane $H_{\vartheta}$. $\mathcal{X}_{\gamma}^{\vartheta}$ then is holomorphic with respect to $\xi$ and jumps are associated with $\vartheta$. Conversely from $\mathcal{X}_{\gamma}^{\vartheta}$ we can divide the twistor $\mathbf{C}^{\times}$ into sectors by pre-BPS rays and define the corresponding Gaiotto-Moore-Neitzke coordinates to agree with  $\mathcal{X}_{\gamma}^{\vartheta}$ in the  sector containing the ray with phase $\vartheta$.

The pre-BPS ray for a charge $\gamma$ is now the ray defined to be the one with the direction  \begin{equation}\vartheta_{\gamma}:= \arg (-Z_{\gamma}(u))\end{equation} and conversely for a pre-BPS ray we can associate a charge. It turns out that an exceptional ray of $\vartheta$ is a pre-BPS ray and a  charge is associated to an exceptional ray of $\vartheta$. In fact, note that when $\vartheta$ is on an exceptional ray on the complex plane  finite trajectories appear. There are three cases. \begin{itemize}\item The case of a flip.  A finite trajectory connecting two zeroes appears. It is lifted to the spectral curve to be a cycle  homotopic  to the charge associated to the loop surrounding these two zeroes. Then by the definition of the central charge as a period it is clear that the phase (or anti-phase) of the central charge of that charge is the phase of the exceptional ray. So the exceptional value of $\vartheta$ picks a   charge. In this way an exceptional ray is identified with a pre-BPS ray.\item The case of a pop for a degenerate triangle. By the local behavior of critical trajectories described above for this case, it is clear the associated charge is a pure flavor charge surrounding the pole. \item The case of a limit configuration. Let $\vartheta_{c}$ be the exceptional value for the appearance of the limit configuration. Without loss of generality let us consider $T_{m},m\rightarrow\infty$ and its limit $T^{+\infty}$. Motivated by (74)(75), we define $$\gamma_{A}^{+}=\gamma_{E_{m-}}+\gamma_{E_{m+}}$$\begin{equation}\gamma_{B}^{+}=(1-m)\gamma_{E_{m-}}-m\gamma_{E_{m+}}\end{equation}Let $\gamma$ be a charge associated to an edge  away from the annular region the infinite flip process does not affect them and we have $$\lim_{\vartheta\rightarrow\vartheta_{c}}\mathcal{X}_{\gamma}^{\vartheta}=\mathcal{X}_{\gamma}^{+}$$  where $\mathcal{X_{\gamma}^{+}}$ is defined in the usual way. The  charge canonically associated to the exceptional value of $\vartheta$ (or equivalently the pre-BPS ray) in this case is defined to be $-\gamma_{A}^{+}$.  Geometrically it is induced by a loop surrounding the inner pole (the one on the inner circle of the annular region) and inner simple zero (the one that is the starting and ending of the inner finite trajectory boundary of the family of closed trajectories) and can be interpreted as the charge for the family of closed trajectories, see \citep{G2}. One can similarly deal with  $T^{-\infty}$.\end{itemize}

\begin{dfn}A geometric BPS charge at $u$ is  either a charge associated to an exceptional ray of $\vartheta$ in the first and the third cases or an anti-charge of such a charge. The set of geometric BPS charges is called the geometric BPS spectra at $u$. A geometric BPS ray is a pre-BPS ray associated to a geometric BPS charge. We may just use the names BPS charge,  BPS spectra and BPS ray because later we will show the equivalence with the definition 6.3.\end{dfn}
\begin{note}(Geometric) BPS numbers will be defined later. \end{note}
 \begin{note}Since geometric BPS charges are associated to finite trajectories these trajectories can be considered as  BPS "instantons" for our instanton correction problem. In fact they should be considered as boundaries of "M2-branes".\end{note}
 \begin{note}We do not define the charge associated to a pop in a degenerate triangle to be a geometric BPS charge for a reason explained in the next section. More discussions of these BPS charges are given in the next section and section 9.\end{note}

Clearly, the most important properties are about discontinuous transformations and the wall crossing formula. We will discuss them along with the relation to cluster algebras in the next section.\\

Till now, we have only described the construction of Fock-Goncharov coordinates for regular singularities. There is a natural extension to irregular cases. When there is an irregular singularity, we have the Stokes phenomenon. So we delete a small open disk containing this singularity and triangulate the complement. The boundary circle of the disk is decomposed into pieces by Stokes sectors and the local behaviors of foliations of quadratic differentials near  a singularity with higher than order two poles are also known. Let us quote the following proposition in Strebel's book.
\begin{thm}Let $p$ be a pole of order $n>2$. Then there are $n-2$ directions at $p$ forming equal angles and trajectories enter  these distinguished directions to go into $p$. There is a neighborhood $U$ of $p$ such that every trajectory ray which enters $U$ goes to $p$. The two rays (i.e. the two directions ) of any trajectory which stays in $U$ go to $p$ in two consecutive distinguished directions.\end{thm}
We pick a point in each of boundary pieces and consider them as vertices associated to the irregular singularity. We make an ordinary triangulation like before (the only essential difference is that we take boundary pieces as edges but the associated Fock-Goncharov coordinates will be defined to be zero). We then take the small flat section in each Stokes sector containing exactly one vertex in each of them and this defines decorations. Finally we define the Fock-Goncharov coordinates for this decorated triangulation as usual. We can label them by charges as before. When defining the labels we only use cycles associated to edges with nonzero Fock-Goncharov coordinates.

\begin{note}It seems that due to its definition, the orders of poles of $\lambda^{2}$ should always be even. But in fact the classification of trajectories applies equally well to odd order poles. Constructions of this section work regardless of the parities of orders of singularities. Allowing odd order quadratic differentials introduces two problems. First of all, from the perspective of solutions of Hitchin's equations this  means that we allow fractional exponents in the asymptotic part of the Higgs field.  It is possible to handle this kind of singularities, see \citep{W}. But it is not clear  to the author in what generality can one find a good foundational theory of the Hitchin's moduli spaces with such singularities. The reference \citep{BB} does not seem to cover all odd order cases. The second problem is that to define Fock-Goncharov coordinates one needs flat sections and it is not clear in what generality the  theory of Stokes matrices applies (it is possible that a general enough theory exists but is unknown to the author). So in section 9 we shall assume that the order is even and the leading term is regular semisimple. It is very likely that these assumption are not necessary as some examples indicate that the theory works fine in odd order cases. In fact we will meet such examples in 9.5. \end{note}

This is also a good place to explain why the author deals with only gauge group $G_{\mathrm{C}}=SL(2,\mathrm{C})$. Fock and Goncharov only defined explicit coordinates as above for $SL(n,\mathrm{C})$ and $PGL(n,\mathrm{C})$ and the latter is the Langlands dual of the former which means the associated moduli spaces are mirror to each other. By the philosophy explained in section 3, we only work with one side. We restrict attentions to $SL(2,\mathrm{C})$ because  in this case we have a thorough understanding of the foliations of differentials in the base of the Hitchin's fibration due to the extensive study of the theory of quadratic differentials in the literature.
\section{Cluster Algebras And Wall Crossing Revisited}

The transformations of charges and Fock-Goncharov coordinates can be formulated in terms of cluster algebras and cluster algebras are also closely related to Kontsevich-Soibleman's wall  crossing formula. Cluster algebras were introduced by Fomin and Zelevinsky \citep{FZ1,FZ2,FZ3}. Later Fock and Goncharov defined the notion of cluster ensembles which contains cluster algebra structures \citep{FG2}. We will give a part of the full definition which is enough for our purposes.
\begin{dfn} A seed is a datum $(\Lambda, (*, *), \{e_{i}\}, \{d_{i}\})$, where \begin{itemize}\item $\Lambda$ is an integral lattice, \item $(*, *)$ is a skewsymmetric $\mathrm{Q}$-valued bilinear form on $\Lambda$,\item $\{e_{i}\}$ is a basis of $\Lambda$, \item $\{d_{i}\}$ are positive integers asigned to $\{e_{i}\}$ and $$\varepsilon_{ij}:= (e_{i}, e_{j})d_{j}\in \mathrm{Z}$$unless $i,j \in I_{0}\times I_{0}$ where $I_{0}$ is a subset of the index set of $i$ and the basis vectors indexed by $I_{0}$ are called frozen vectors.\end{itemize} Define a torus called seed $\mathcal{X}$-torus by $$\mathcal{X}_{\Lambda}:=\mathrm{Hom} (\Lambda, \mathrm{C}^{\times})$$ An element $v\in \Lambda$ is a character of the torus and is denoted as $X_{v}$. For $e_{i}$ we use the notation $X_{i}$. The set of $X_{i}$ is called the seed $\mathcal{X}$-coordinates and there is a natural holomorphic Poisson bracket $$\{X_{v}, X_{w}\}:=(v, w)X_{v}X_{w}$$We can also define another torus known as the cluster $\mathcal{A}$-torus $$\mathcal{A}_{\Lambda}:=\mathrm{Hom}(\Lambda^{\circ}, \mathrm{C}^{\times})$$where $\Lambda^{\circ}$ is the sublattice in $\Lambda^{*}\otimes \mathrm{Q}$ spanned by $f_{i}:=d_{i}^{-1}e_{i}^{*}$. The basis $\{f_{i}\}$ is called cluster $\mathcal{A}$-coordinates and renamed as $\{A_{i}\}$. There is a symplectic form on $\mathcal{A}$ defined by $$\Omega:=(e_{i},e_{j})d\log A^{i}\wedge d\log A^{j}$$ Set $[a]_{+}:= \max (0, a)$. The seed obtained by mutation in the direction of a non-frozen vector $e_{k}$ (denoted by the symbol $\mu_{k}$) is defined to be a seed obtained by (only) replacing $e_{i}$ by $e_{i}^{'}$ $$e_{i}^{'}:= e_{i}+[\varepsilon_{ik}]_{+}e_{k}, i\neq k$$\begin{equation}e_{i}^{'}:= -e_{k}, i=k\end{equation} A mutation $\mu_{k}$ induces the following transformations (also denoted by $\mu_{k}$)
\begin{equation}
\begin{array}{ll}\mu_{k}^{*}(X_{i}^{'})&=X_{k}^{-1}, i=k\\
\mu_{k}^{*}(X_{i}^{'}) &= X_{i}(1+X_{k}^{-sgn(\varepsilon_{ik})})^{-\varepsilon_{ik}}, i\neq k\\
\end{array}
\end{equation}\begin{equation}
\begin{array}{ll}A_{k}\mu_{k}^{*}(A_{k}^{'})&=\prod_{j, \varepsilon_{kj}>0}A_{j}^{\varepsilon_{kj}}+\prod_{j, \varepsilon_{kj}<0}A_{j}^{-\varepsilon_{kj}}\\
\mu_{k}^{*}(A_{i}^{'}) &= A_{i}, i\neq k\\
\end{array}
\end{equation}These transformations are called cluster transformations. Variables $\{A_{i}\}$ together with their all mutations generate the so-called cluster algebra (with cluster transformations as relations). Later we will also call the algebra generated by dual variable $X_{i}$ together with their all mutations the cluster algebra.\end{dfn}

Now let us return to the setting of the last section. For fixed $\vartheta$ and $u$, we have a WKB triangulation and a set of charges associated to edges.
\begin{dfn}Fix $\vartheta$ and $u$. Following Gaiotto, Moore and Neitzke \citep{G3} we call those geometric BPS charges whose phases of central charges are between $\vartheta$ and $\vartheta+\pi$ positive roots. A positive root which is not a sum of other positive roots is called a simple root.\end{dfn}\begin{note}It seems that this definition depends on $\vartheta$ and $u$ but the dependence is actually weaker. We will discuss this issue in section 9.\end{note} \begin{note}It is easy to see that for a fixed $u$ and any $\vartheta$ the set of positive roots and their anti-charges is precisely the geometric BPS spectra at $u$.\end{note}
\begin{thm} \citep{G2} Fix $\vartheta$ and $u$. A complete set of simple roots is contained in  the set of $\gamma_{E}^{\vartheta}$ where $E$ runs over the set of all edges of the WKB triangulation $T_{WKB}(\vartheta,u)$.  A positive root is a sum of simple roots with nonnegative coefficients.\end{thm}Since a finite trajectory corresponding to a positive root has phase between $\vartheta$ and $\vartheta+\pi$ the intersection of this trajectory and a trajectory with phase $\vartheta$ (in particular an edge in $T_{WKB}(\vartheta,u)$) is positive. This theorem follows  from this fact, the counting of number of edges and the definition of labeling by charges.

\begin{dfn}Define a matrix by taking intersections among only simple roots $\{\gamma_{i}\}$ \begin{equation}\varepsilon_{ij}:=\langle\gamma_{i}, \gamma_{j}\rangle \Omega(\gamma_{j}; u)\end{equation}where $\Omega(\gamma_{j}; u)$ is the BPS number and $\langle\cdot\rangle$ is the intersection pairing. Since $\gamma_{E}^{\vartheta}$ generate the charge lattice we can extend this pairing to all charges. \end{dfn}

Gaiotto, Moore and Neitzke gave an a priori assignment of (geometric) BPS numbers to geometric BPS charges.
\begin{dfn}The  geometric BPS charge $\gamma$ associated to a flip is called a hypermultiplet and also denoted by $\gamma_{hyp}$.  We define the corresponding geometric BPS number by  \begin{equation}\Omega(\gamma_{hyp};u)=1\end{equation}The BPS charge associated to a juggle, i.e. $-\gamma_{A}^{+}$ is called a vectormultiplet and also denoted as $\gamma_{vec}$.  We define the geometric BPS number by \begin{equation}\Omega(\gamma_{vet};u)=-2\end{equation} Geometric BPS numbers of other charges are defined to be zero. We identify geometric BPS numbers as BPS numbers used in the definition 6.3. From now on we can drop the word "geometric".\end{dfn}
 \begin{note}One may ask what happens to charges associated to pops in degenerate triangles. We do not need  to define them or we can define them to be zero\footnote{This may not coincide with the ultimate  definition of BPS numbers as indices. The point here is that this possible discrepancy is not detectable by the wall crossing formula.} because they are irrelevant for wall crossing formulas. That is because the geometric BPS charges in this case are pure flavor charges in the radical of the intersection paring which means the associated Kontsevich-Soibelman transformations are trivial.\end{note}
 \begin{note}The philosophy here is that the notion of geometric BPS charges is the primary one and geometric BPS number are actually assigned instead of defined geometrically. Because of our definition a geometric BPS charge is identified with a BPS charge defined in  section 6. \end{note}
 There is also an a priori assignment of quadratic refinements.\begin{equation}\sigma(\gamma_{hyp})=-1,\sigma(\gamma_{vec})=1\end{equation}
 \begin{note}An important point here is that unlike the definition in section 5 this assignment is global over the whole space of $u$. See \citep{G2} for the geometric justification of this definition.\end{note}

 It is important to know the transformation formulas of charges and Fock-Goncharov coordinates labeled by charges  when we vary $\xi$ (or equivalently $\vartheta$) across a geometric BPS ray.\begin{itemize} \item For a flip from the side where $\mathrm{Im}(Z_{\gamma_{k}}/ \xi)>0$ to the side where $\mathrm{Im}(Z_{\gamma_{k}}/ \xi)<0$,the transformations of charges are $$\mu_{k}: \gamma_{k}\rightarrow \gamma_{k}^{'}:=-\gamma_{k}$$ \begin{equation}\mu_{k}: \gamma_{i}\rightarrow \gamma_{i}^{'}:=\gamma_{i}+\gamma_{k}[\varepsilon_{ik}]_{+}, i\neq k\end{equation} Transformations of Fock-Goncharov coordinates labeled by charges are obtained by taking the composition of the  transformations of Fock-Goncharov coordinates labeled by edges and the transformations of charges \begin{equation}\mathcal{X}_{\gamma_{i}}\rightarrow  \mathcal{X}_{\gamma_{i}}(1+\mathcal{X}_{\gamma_{k}})^{-<\gamma_{i}, \gamma_{k}>}\end{equation}here $\gamma_{hyp}$ is $\gamma_{k}$\footnote{Strictly speaking according to  (68) we should use the limits $\mathcal{X}_{\gamma}^{\pm}$ instead of $\mathcal{X}_{\gamma}$, but we shall use the  notations $\mathcal{X}_{\gamma}$   in the rest of the paper to make some formulas more readable. We hope it will not introduce confusions.} . \item For a juggle from $T^{-\infty}$ to $T^{+\infty}$ the transformations of charges are $$\gamma_{A}^{+}=-\gamma_{A}^{-}$$\begin{equation}\gamma_{B}^{+}=-\gamma_{B}^{-}+2\gamma_{A}^{-}\end{equation}Transformations of Fock-Goncharov coordinates labeled by charges are \begin{equation}\mathcal{X}_{\gamma}^{-}\rightarrow\mathcal{X}_{\gamma}^{-}(1-\mathcal{X}_{\gamma_{vec}}^{-})^{-2\langle\gamma,\gamma_{vet}\rangle}\end{equation} where $\gamma_{vec}:=-\gamma_{A}^{+}$. Note that $\langle\gamma_{B}^{+},\gamma_{A}^{+}\rangle=-2$.\end{itemize}

Come back to cluster algebras. Compare equation (88) (72) with equation (81) (82), we see that if we identify $e_{i}$ and $\gamma_{i}$ as well as $X$ and $\mathcal{X}$ (labeled by edges), then the transformations  of Fock-Goncharov coordinates labeled by edges under a flip are precisely  cluster transformations.
\begin{thm} Transformations of Fock-Goncharov coordinates labeled by edges obtained by crossing  BPS rays corresponding to flips are cluster transformations. The Fock-Goncharov coordinates labeled by edges generate a (dual) cluster algebra.\end{thm}

There is another crucial observation. The  equations (89) is a Kontsevich-Soibelman transformation! In fact it is \begin{equation}K_{\gamma_{k}}^{-\Omega(\gamma_{k}; u)}: \mathcal{X}_{\gamma_{i}}\rightarrow  \mathcal{X}_{\gamma_{i}}(1-\sigma(\gamma_{k})\mathcal{X}_{\gamma_{k}})^{-<\gamma_{i}, \gamma_{k}>\Omega(\gamma_{k};u)}\end{equation}Similarly (91) is also a Kontsevich-Soibelman transformation
\begin{equation}K_{\gamma_{vec}}^{\Omega(\gamma_{vec};u)}:\mathcal{X}_{\gamma}^{-}\rightarrow\mathcal{X}_{\gamma}^{-}(1-\sigma(\gamma_{vec})
\mathcal{X}_{\gamma_{vec}}^{-})^{\langle\gamma,\gamma_{vet}\rangle\Omega(\gamma_{vec};u)}\end{equation}
We also know that the Kontsevich-Soibelman transformation associated to a pop in a degenerate triangle is trivial. Therefore we get
\begin{thm}Suppose $u$ represents a quadratic differential with only simple zeroes (which is the generic case), then all transformations of Fock-Goncharov coordinates labeled by charges obtained by varying $\vartheta$ across  BPS rays (fixing $u$) are Kontsevich-Soibelman transformations or products of Kontsevich-Soibelman transformations.\end{thm}
\begin{note}The convention in section 7 of not calling the associated charge to a pop of a degenerate triangle a BPS charge is justified because the associated Kontsevich-Soibelman transformation is trivial.\end{note}

Now it is clear what the wall crossing formula means in our context. We vary the moduli $u$ and  the order of pre-BPS rays changes when $u$ crosses a stability wall. However, it is clear that the two ordered products of Kontsevich-Soibelman transformations over all rays between two fixed phases\footnote{We always assume that neither of them is an exceptional one i.e. the phase of a pre-BPS ray.} defined by \begin{equation}S(\vartheta_{-},\vartheta_{+}; u)=\prod^{\leftarrow}_{\vartheta_{-}< arg(-Z_{\gamma}(u)) <\vartheta_{+}}K_{\gamma}^{\Omega(\gamma; u)}\end{equation} is unchanged. In other words we have the following
\begin{thm}
\begin{equation}S(\vartheta_{-},\vartheta_{+}; u)=S(\vartheta_{-},\vartheta_{+}; u^{'})\end{equation} where $u$ and $u^{'}$ are two different moduli parameters in $B$ not on stability walls and  both of them are not in the codimensional two locus corresponding to quadratic differentials with nonsimple zeroes.  Suppose  that they can be connected by a path not
crossing branch cuts  such that along the path from   $u$ to $u^{'}$  no  BPS ray crosses $\vartheta_{-}$ or $\vartheta_{+}$\footnote{This assumption  is not restrictive at all. See remarks after the statement of wall crossing formula in section 6.}.\end{thm}

The proof is simple. We fix $\vartheta_{-}$($\vartheta_{+}$) and let $u$ change along the path.  Clearly this path can be assumed to avoid the locus of quadratic differentials with nonsimple zeroes. Since no phases of BPS rays could pass $\vartheta_{-}$($\vartheta_{+}$), there is  no changes of the decorated triangulation up to isotopy. So$$\mathcal{X}_{\gamma}^{\vartheta_{-}}(u)=\mathcal{X}_{\gamma}^{\vartheta_{-}}(u^{'})$$  $$\mathcal{X}_{\gamma}^{\vartheta_{+}}(u)=\mathcal{X}_{\gamma}^{\vartheta_{+}}(u^{'})$$ Since $S(\vartheta_{-},\vartheta_{+}; u)(S(\vartheta_{-},\vartheta_{+}; u^{'})) $ is the transformation mapping $\mathcal{X}_{\gamma}^{\vartheta_{-}}(u)$ ($\mathcal{X}_{\gamma}^{\vartheta_{-}}(u^{'})$)  to $\mathcal{X}_{\gamma}^{\vartheta_{+}}(u)$ ($\mathcal{X}_{\gamma}^{\vartheta_{-}}(u^{'})$), the wall crossing formula follows.
\begin{note}By the discussion in section 7 we do not have to assume that $u,u^{'}$ are generic in the sense that varying $\vartheta$ only produces the three types of jumps of WKB triangulations listed in the previous section. \end{note}
\begin{note}The BPS numbers are assigned to BPS charges. The assignment is global on $B$. On the other hand if we are given an initial assignment at a point then we can use the wall crossing formula to propagate the BPS spectra and BPS numbers across stability walls and in this way we can also calculate BPS charges and BPS numbers at other points. The two approaches are consistent because the wall crossing formula is derived by using the Kontsevich-Soibelman transformations associated to jumps and values of BPS numbers are so defined to identify jumps as Kontsevich-Soibelman transformations.  The second approach is the practical one for computations. \end{note}
\begin{note}Now we can describe  the full dependence of Fock-Goncharov coordinates. A Fock-Goncharov coordinate can be written as $\mathcal{X}^{(u_{0},\vartheta)}_{\gamma}(u,\theta;\xi)$. In the bracket $(u,\theta;\xi)$ $u$ and $\theta$ are coordinates of $\mathcal{M}$ on the base and the fiber of the Hitchin's fibration respectively. $\xi$ is the coordinate in the twistor $\mathbf{C}^{\times}$. Usually we suppress the dependence on $(u,\theta)$. Later in section 9 we will fix $\xi$. In the bracket $(u_{0},\vartheta)$ $u_{0}$ is a moduli parameter. So we do not view it as a part of coordinate system of $\mathcal{M}$. We discussed this issue before in remark 7.3. We have been denoting it by $u$ and we will continue to do so. Hopefully this would not cause any confusions. Similarly $\vartheta$ is a moduli parameter although it is also the phase of a ray in the twistor $\mathbf{C}^{\times}$. It is not the phase of $\xi$. So in particular fixing $\xi$ does not fix $\vartheta$. The word moduli here is used in the sense that  the solution of the instanton correction problem of complex structures to be discussed in section 9 depends on the choice of $(u,\vartheta)(i.e.(u_{0},\vartheta))$. From the perspective of foliations  of quadratic differentials the choice of $(u,\vartheta)$ determines a WKB triangulation. But the dependence is actually on chambers instead of points. By chambers we mean BPS chambers to be defined in section 9.\end{note}

 The wall crossing formula  has the following important consequence. Recall that the system of positive roots for $(\vartheta,u)$ and their anti-charges is the BPS spectra at $u$. The system of positive roots depends on $\vartheta$ while the BPS spectra do not. The system of positive roots is obtained by varying $\vartheta$ to $\vartheta+\pi$ with a fixed $u$  and collecting charges associated to those finite trajectories that appear during this evolution. Now we change $u$ along a path without letting  any BPS rays pass $\vartheta, \vartheta+\pi$ and without crossing any stability walls. Since the order of pre-BPS rays does not change by the unique factorization of $S(\vartheta_{-},\vartheta_{+}; u)$  the  system of positive roots does not change and neither do the BPS spectra.  If when $u$ hits some point on the path a BPS ray hits (say) $\vartheta$ then we just rotate $\vartheta$ a little bit to avoid that hitting and reduce to the no-passing situation. This is fine for determining the BPS spectra simply because the BPS spectra  do not depend on $\vartheta$. Therefore as long as no stability wall is crossed the BPS spectra do not change (of course the system of positive roots could change). There are countably many codimensional one stability walls in $B$. By a stability chamber we mean a domain (a simply connected open subset with nonempty interior) bounded by stability walls and branch cuts and having empty intersections with any stability walls or cuts. What we have just showed means that the BPS spectra are constant inside a stability chamber and therefore can be considered as data associated to the chamber. This was the point of view we used in section 4 and will be important later. However there is a tricky problem. It is possible that the stability walls are dense in (or in at least some region of) $B$ making the notion of chambers useless. Nevertheless usually we are given an initial assignment of a system of positive roots at a point and as we vary $u$ continuously most stability walls are irrelevant as the corresponding charges do not appear in the factorizations. In this case it could make sense to discuss stability chambers.

It is also very interesting to know that  the wall crossing formula can be used even when both of the factorizations on the two sides are not given. This seems unlikely but is in fact based on a simple fact. Gaiotto, Moore and Neitzke defined the following product called the spectrum generator
\begin{equation}S(\vartheta; u)=\prod^{\leftarrow}_{\vartheta< arg(-Z_{\gamma}(u)) <\vartheta+\pi}K_{\gamma}^{\Omega(\gamma; u)}\end{equation} Note that for every BPS particle either  its BPS ray lies in this half plane or its antiparticle's BPS ray does so. This means the spectrum generator captures exactly half of the BPS spectra with the other half just antiparticles. It is easy to see that  the WKB triangulation $T_{WKB}^{\vartheta+\pi}$ is obtained from $T_{WKB}^{\vartheta}$ by popping at all vertices (this operation is called an ominipop) and we know how to write down the transformations\footnote{We have shown before that the transformations of Fock-Goncharov coordinates labeled by charges are trivial for pops in degenerate triangles. But here we are considering pops in not necessarily degenerate triangles. So they are nontrivial in general.}. This transformation can be obtained without following any continuous evolution of $\vartheta$ but by its definition it is nothing but the spectrum generator. In other words, we have obtained the product without calculating any of its factors. Once we know the product we can factorize it by successive truncations according to the degree of charges. The factorization is unique provided the order of the product, i.e. the order of BPS rays is given but that is determined by the location of the moduli parameter $u$. So on each side of the stability wall, we have a factorization and in this way both sides of the wall crossing formula are constructed. Although the inductive procedure here uses essentially the  same truncation as the one used in the  inductive  proof of Theorem 3.5, there is a difference. In the latter case, one needs to know some Kontsevich-Soibelman factors (the initial data of a scattering diagram) from the beginning. \\

The wall crossing formula can also be considered as  a solution of some enumerative problems. In fact, the uniqueness of the factorization of an ordered product (given the order of (pre-)BPS rays) means that we can calculate the exponents of all Kontsevich-Soibelman factors on one side of a stability wall if all exponents are known on the other side. But an exponent is just a product of an intersection number of charges and the corresponding BPS numbers. So the wall crossing formula is a computational tool for BPS numbers. Note that BPS numbers associated to the BPS spectra are assigned   which is consistent with the wall crossing formula. The point is that once we know that there is a BPS charge then its BPS number is determined immediately but it is not easy in general to know whether a charge is BPS by following a continuous evolution. The wall crossing formula can tell us about the BPS spectra inductively.  The assignment of BPS numbers is natural. $\Omega(\gamma; u) =1$ for charges labeling flipping edges. The instantons in this case are interpreted as those finite trajectories connecting two zeroes of the quadratic differential when $\vartheta$ is on the BPS ray $l_{\gamma}$. These finite trajectories are labeled by charges just by using the charges of flipping edges and there is indeed only one finite trajectory with charge $\gamma$. So $\Omega(\gamma_{hyp}; u)=1$ gives an honest counting of critical trajectories of the quadratic differential. For juggles $\Omega(\gamma_{vec}; u)=-2$ which tells us that the result of virtual counting in this case should be $-2$.\\

Let us summarize. We have explained  the proof in \citep{G1, G2} of the following
\begin{thm}Let $\mathcal{M}(R)$ be a Hitchin's moduli space defined in section 2 with a large enough $R$. Then the Fock-Goncharov coordinates defined above for WKB triangulations $T_{WKB}(\vartheta, \lambda^{2})$ are a set of Gaiotto-Moore-Neitzke coordinates  and therefore a hyperkahler structure $\mathcal{M}(R)$ is  obtained by the twistor construction.\end{thm}

The description of the hyperkahler structure in this way is quite complicated but the effect of instanton (the BPS spectra) corrections is manifest. However it is not very explicit in the sense that   although the leading order behavior can be obtained by the WKB analysis it is not easy to extract more refined information\footnote{For example, it is not obvious to the author that the metrics defined in this way are complete. As another example it is not clear that the solution of the infinitely dimensional Riemann-Hilbert problem formulated in \citep{G2} is unique. Nevertheless if the answer was yes then it seems the iteration scheme outlined in \citep{G2} would provide an expansion of the metric beyond the leading order.}.  The point here is not that we want to investigate the differential geometry of this hyperkahler metric. What we want to do is to compare it to the instanton correction problem of complex structures in the framework of mirror symmetry.\\

 An important question to ask  is whether the instanton-corrected hyperkahler metric constructed by Fock-Goncharov coordinates and Gaiotto-Moore-Neitzke ansatz is actually the hyperkahler metric given
by the infinitely dimensional hyperkahler quotient construction of Hitchin \citep{H1}. The author would guess that the answer is yes.\\

\begin{con} For $\mathcal{M}(R)$ the Gaiotto-Moore-Neitzke metric coincides with Hitchin's hyperkahler quotient metric.\end{con}

In fact, we already know that the instanton-corrected metric has some crucial properties of the hyperkahler quotient metric such as the Hitchin's fibration is a holomorphic one for $J_{3}$.  There are at least three places where the hyperkahler metric is important in the study of mirror symmetry. First, recall that the identification of the Strominger-Yau-Zaslow mirror duality as the Langlands duality depends on the relation between Hitchin's fibration and the hyperkahler quotient metric. Second, the Fock-Goncharov coordinates are constructed for twistor parameter $\xi\in\mathbf{C}^{\times}$ of the hyperkahler quotient metric
 (because we view the moduli space as the moduli space of flat connections). Third,  we have the hyperkahler metric constructed from Fock-Goncharov coordinates
 which is of an instanton-corrected form required by mirror symmetry because of  the equivalence of two instanton correction problems to be explained in section 9. A positive answer would make all these metric  aspects of mirror symmetry of Hitchin's moduli spaces compatible with each other. A possible way to show that the two hyperkahler metrics are the same is to use an appropriate noncompact version of Yau' theorem as suggested by Seiberg and Witten (see \citep{SW2}) if one  understands the asymptotic behavior of the hyperkahler quotient metric near the infinity. This direction will not be pursued in this work.

\section{Two Instanton Correction Problems And Mirror Symmetry Through Wall Crossing}

\subsection{Outline}

Let us summarize what we have learned.

For an $SU(2)$ Hitchin's moduli space $\mathcal{M}$ on a Riemann surface $C$, one can construct an instanton-corrected hyperkahler metric from a set of coordinates using the twistor method provided that they satisfy the properties formulated by Gaiotto, Moore and Neitzke. A set of such coordinates can be constructed as Fock-Goncharov coordinates $\mathcal{X}_{\gamma}$ obtained from  WKB triangulations. A WKB triangulation is constructed from a foliation of a quadratic differential $\lambda^{2}$ on the underlying Riemann surface and depends on two parameters $\vartheta$ and $u=-\lambda^{2}$. Changing $\vartheta$ (or equivalently changing $\xi$) means that we are rotating over the twistor $CP^{1}$ and there are countably many discontinuous jumps along BPS rays which are Kontsevich-Soibelman transformations. If we vary $u$, then there are real codimensional one hypersurfaces called (marginal) stability walls such that crossing such a wall changes the order of BPS rays. However the ordered product of discontinuous transformations does not change. This statement is called the wall crossing formula and it enables us to determine BPS spectra and associated BPS numbers. The BPS spectra are our instantons in this instanton correction problem and they happen to  coincide with the BPS spectra of the four dimensional gauge theory which is a closely related but different problem. In this construction we view $\mathcal{M}$ as the moduli space of $SL(2,\mathbf{C})$-flat connections which is analytically isomorphic to the space of fundamental group representations (also denoted by $\mathcal{M}$).

On the other hand, a Hitchin's moduli space is also an example of Calabi-Yau spaces and as such it makes sense to discuss mirror symmetry. The mirror in the sense of Strominger, Yau and Zaslow is the Hitchin's moduli space over the same Riemann surface but with the gauge group changed to its Langlands dual. There are singular special
Kahler structures and in particular singular affine structures on the bases of Hitchin's fibrations of the two moduli spaces which are dual to each other in the sense of Legendre transform. The instanton correction problem asks to construct  a family of Calabi-Yau spaces whose complex structures are supposedly determined by holomorphic disks in the mirror family. In Gross and Siebert's purely algebraic approach, this is reformulated as constructing a formal toric degeneration of Calabi-Yau's from an integral singular affine structure together with a polyhedra decomposition, a polarization and a log smooth structure. The central fiber is constructed first and to get the family one deforms  affine (in the sense of algebraic geometry) pieces of the central fiber before gluing them together. One works inductively by the powers of the deformation parameter and at each order there could be inconsistencies during the gluing and we have to compose log morphisms associated to codimensional one polyhedral subsets which "correct" the gluing and hence the complex structures. These $structures$ are constructed on the same  singular affine manifolds and are supposed to encode the dual of tropical avatars of holomorphic disk instantons in the mirror. Regardless of the enumerative meaning of these $structures$, Gross and Siebert constructed them inductively and hence solved the instanton correction problem of complex structures  in an algebraic sense. \\

Now we will connect the two chains of thoughts over the same space. Let us call the two instanton correction problems the metric problem and the complex structure problem respectively.

As mentioned in section 2, according to Kapustin and Witten  the mirror
symmetry of Hitchin's moduli spaces is induced by the electric-magnetic
duality of 4d N = 4 gauge theories. So we have two instanton correction problems
over a Hitchin's moduli space. One is suggested
by gauge theory (GMN's work). The other is in
the context of mirror symmetry which according
to above is also related to gauge theory. It
would be surprising if these two problems are
unrelated. We want to show that these  two instanton correction problems are equivalent in an appropriate sense.

It may be possible to study holomorphic disk instantons
directly over Hitchin's moduli spaces
in the context of mirror symmetry. Here we will
try to use Gross-Siebert's strategy instead.

However there are some difficulties.
\begin{enumerate}
\item Objects of the two problems are quite different.
\item It is not clear how the critical
trajectories (instantons in gauge theory) can
be related to instantons in mirror symmetry.
\item It seems that the consistency conditions of log morphisms should be identified with  the wall crossing formula but the meanings of "stability walls" and instantons are not obvious in that framework. In fact since there is no central charges or even charges involved in the construction, it is not clear how to make sense of stability walls and BPS instantons labeled by charges.
\item To do mirror symmetry we need a family
of Calabi-Yau's. For Hitchin's moduli
space we need a
family which is a degeneration of moduli
spaces and is a large complex degeneration at
the same time.
\item To run Gross and Siebert's algorithm we need a log smooth structure as the input.\end{enumerate}

To formulate and prove the equivalence we rely on the following ideas
\begin{enumerate}
\item {\bf Large $R$ limit} \ As suggested in section 5 we view $R\rightarrow\infty$ as a differential geometric characterization of a large complex limit of Hitchin's moduli spaces. This provides a natural family of Hitchin's moduli spaces which could be compatible with mirror symmetry. This is conceptually important and solves the difficulty 4.
\item {\bf SYZ meets GS} \ In the limit form of SYZ conjecture the metric
degeneration is  more fundamental
than the notion of large complex limits.
In fact the mirror duality should exchange the
large complex limit with its mirror called large
volume limit (which is a Kahler degeneration)
and vice versa. The metric degeneration is
supposed to take care of both.

However in GS picture we lose the information
of Ricci
flat metrics and deal with only complex
structure degenerations. Therefore the family
version of mirror duality cannot be recovered
directly. This is reflected by the fact that the dualization
of the limit structure is defined only
for a triple consisting of the singular affine
structure, the polyhedral decomposition and
the polarization. We need to specify a log smooth
structure for the triple to be able to use GS's
solution of the reconstruction problem.

This log smooth structure cannot be obtained
by taking some kind of dualization of the log
smooth data of the original family as the dual of this log smooth structure
is supposed to recover tropical avatars
of holomorphic disks wrapping special Lagrangian
fibers which require more than the information
of complex structures of the original family.
This is consistent with the large complex limit vs
large volume limit philosophy.

Since we view large $R$ limit as the large complex
limit it makes sense to use the metrics to build
an a priori assignment of log smooth structure. This procedure solves the difficulty 5.

So we are using a hybrid version of SYZ and
GS procedures. This is  inevitable since
we want to compare the differential geometric
aspect and the algebraic geometric aspect.

Note that this is somewhat unusual because we
will try to use "instantons" on the same side
instead of "instantons" on the mirror side to do
the correction (reconstruction) problem. But
the comparison is still quite meaningful even
if one's sole purpose is to establish a mirror
duality relation. We will comment more on this point later.

\item {\bf Labeling by charges and BPS walls} \ A BPS wall associated to a charge is defined to be the locus where the values of central
charges over twistor parameters are real. We use the projections of BPS walls to build a polyhedral decomposition such that codimensional one cells are labeled by BPS charges just like BPS rays. This eventually solves the difficulty 3.

\item {\bf Wall Crossing Formula as System
of Consistency Conditions} We assign the log smooth structure such that in the end one can identify consistency conditions as wall crossing formulas after setting the deformation parameter to 1.  The assignment is obvious and is dictated by labeling by charges. We can then construct a compatible system of consistent $structures$.

\item {\bf Explicit Degenerations} To justify the naturality of previous constructions  from the perspective
of Hitchin's moduli spaces we want
to check whether the toric degeneration obtained
by Gross-Siebert's construction can recover the
Hitchin's moduli space as a generic fiber. This can be verified by the explicitly construction of the toric degeneration and the comparison to the natural way of writing down the
underlying complex manifold of $\mathcal{M}$ (viewed as the
moduli space of flat connections) in terms of
gauge invariant Fock-Goncharov coordinates
with relations the discontinuous jumps.

This step puts the whole strategy into a conceptually firm and natural framework. However it does not answer the question posed in difficulty 2 in a completely geometric way. In this sense the author views the results in this section as half conceptual and half computational.

\end{enumerate}

\subsection{Construction Of Gross-Siebert Data}

We introduce an additional parameter $R$ to the Hitchin's equation using (3) and all constructions for the metric problem now have $R$ dependence.  There is a singular special Kahler structure on the base $B$ of the Hitchin's fibration for any member of the family.  We want to use Gross and Siebert's construction over $B$ with the singular affine structure induced by the singular special Kahler structure.  We need to find the additional data making the input of Gross-Siebert's algorithm. Let us call   a singular integral affine manifold together with a polyhedral decomposition, a polarization and a positive log smooth structure (i.e. the input of Gross-Siebert's algorithm) a $Gross-Siebert$ $data$ (GS data)\footnote{In Gross and Siebert's work there is also an additional freedom of choosing the so-called open gluing data. We always take the trivial open gluing data and hence do not need to define them and do not consider them as a part of Gross-Siebert data.}.

The singular special Kahler structure on $B$ induces an integral affine structure with singularities by using its affine coordinates. More precisely, since the periods (central charges) $Z_{\gamma}(u)$ are special holomorphic coordinates we can take $$a_{i}:=\mathrm{Im}(Z_{\gamma_{i}}/ \xi)$$ for a fixed $\xi\neq0,\infty$ as natural flat affine coordinates\footnote{Of course in section 2 we used real parts of central charges, but taking the special holomorphic coordinates as $-i/\xi$ times central charges will do the job.} outside the singular locus. The affine structure is integral because the cycles $\gamma$ are integral. The singular locus $\Delta$ is the discriminant locus where some periods vanish and has codimension at least two. The monodromies are Picard-Lefschetz transformations associated to vanishing cycles.
\begin{note}When we define the affine structure $(\gamma_{i})$ are always chosen  to be a basis of gauge charges.  This gives the right dimension of $B$. So the underlying integral lattice for affine charts is  the gauge charge lattice. \end{note}

Next we need to construct a polyhedral decomposition.
\begin{dfn}Let $\tilde{B}$ be the universal cover of the complement of singular locus in $B$ denoted by $B^{*}$. Suppose $\gamma$ is a charge such that $\Omega(\gamma;u)\neq0$ for some $u$. We define the BPS wall $$W_{\gamma}:=\{(u,\xi)\mid u\in\tilde{B},\xi\in\mathbf{C}^{\times}, {Z_{\gamma}(u)\over \xi} \in \mathbf{R}\}$$We only consider charges in the BPS spectra when discussing BPS walls. Clearly $W_{\gamma}=W_{-\gamma}$.\end{dfn}\begin{note}So now there are three different kinds of "walls" in this paper: walls defined in section 3, stability walls defined in section 4 and BPS walls defined here. To avoid confusions from now on we call a wall defined in section 3 a Gross-Siebert wall (GS wall).\end{note}\begin{note}We use $\tilde{B}$ because the charges and the BPS spectra have monodromies over $B$. We can choose  branch cuts and view $u$ as $u\in B^{*}$. We will use this point of view.\end{note}\begin{note}Similar ideas appeared in \citep{G3} for different purposes.  \end{note}
\begin{note}Also in Kontsevich-Soibelman's wall crossing paper\citep{KS}
some structures similar to the projections of
BPS walls for Hitchin's moduli spaces were
considered and Kontsevich-Soibelman had the
vision that one should consider the wall crossing
formula over the base $B$. Metrics and instanton corrections were not discussed
there though the wall crossing formula seemed to be motivated by their previous work on the instanton correction problem of complex structures in the two dimension \citep{KS1}.\end{note}
$\tilde{B}\times\mathbf{C}^{\times}$ is divided by BPS walls into chambers. In other words a chamber is a connected component of \begin{equation}\Xi:=\tilde{B}\times\mathbf{C}^{\times}-\bigcup_{\gamma,\exists u, \Omega(\gamma;u)\neq0}W_{\gamma}\end{equation}Since there are several types of chambers in this paper we have to use different names. Recall that in section 3 there are chambers in a $structure$. We call them Gross-Siebert chambers (GS chambers). There are chambers in a scattering diagram divided by rays and cuts. They are called scattering chambers. There are chambers in the twistor $\mathbf{C}^{\times}$ divided by  BPS rays and they are called twistor chambers. Note that twister chambers depend on $u$. There are chambers (defined in section 6) in $B$ divided by stability walls and branch cuts. They are called stability chambers. Finally the chambers defined here are called BPS chambers. A chamber is always assumed to have nonempty interior. So it is possible that a point $(u,\xi_{0})$ which is not on any BPS wall is not contained in a BPS chamber if BPS walls can be dense. Later whenever we discuss a BPS chamber containing $(u,\xi_{0})$ we always assume that such a BPS chamber exists.

The projection of $W_{\gamma}$ to the twistor $\mathbf{C}^{\times}$ for fixed $u$ (not on a stability wall) gives us BPS lines $L_{\gamma}$ (which is the union of two BPS rays of opposite directions). Consider the projection of $(u,\xi_{0})$ to the $\xi$-plane, the twistor $\mathbf{C}^{\times}$ which is divided into twistor chambers. Recall that the angular parameter $\vartheta$ is the phase of the central ray of the half plane $H_{\vartheta}$. Let  $\vartheta=\arg\ \xi_{0}$ so that the positive roots for $(u,\vartheta)$ all lie in the half plane on the counterclockwise side of the  line  through the origin containing $\xi_{0}$.  The twistor chamber containing the ray with phase $\vartheta$ is bounded by two BPS rays.   The first ray is the first BPS ray (which has to be associated to a (minus)  simple root) that one meets if we continuously varies $\vartheta$ to $\vartheta +\pi$ while the other ray is the opposite of the last BPS ray one meets which is also associated to a simple root. Now we vary $u$. We can move from the BPS chamber containing $(u,\xi_{0})$ into another BPS chamber if and only if some BPS rays pass the ray with phase $\vartheta$ but the first such ray can only be one of the two BPS rays associated to simple roots described above.   The same argument holds for varying $\vartheta$. So we have shown the following theorem due to Gaiotto, Moore and Neitzke.
\begin{thm}BPS chambers are labeled by systems of simple roots. For a BPS chamber the associated system is given by the system of simple roots associated to a point in that BPS chamber and it does not depend on the choice of the point.\end{thm}
Similarly one can attach to a point in a BPS chamber the system of positive roots at that point.
\begin{dfn}Let $\Theta$ be a BPS chamber. The set of positive roots associated to a point in $\Theta$ is called the refined BPS spectra for that BPS chamber. The set of simple roots associated to a point in $\Theta$ is called the refined simple BPS spectra for that BPS chamber. The definitions  do not depend on the choice of the point. More generally if a point $(u,\xi_{0})$ is not necessarily contained in a BPS chamber then we just define the refined (simple) BPS spectra at that point. The definition depends only on $u$ and $\vartheta:=\arg \ \xi_{0}$. Clearly we can use $(u,\vartheta)$ instead of  $(u,\xi_{0})$. We denote the refine BPS spectra at $(u,\vartheta)$ by $\Sigma(u,\vartheta)$. \end{dfn}\begin{note}Note that the refined BPS spectra exhausts all BPS charges at $u$ that one encounters by varying $\vartheta$ to $\vartheta+\pi$. \end{note}

 Let $u$ be a point in a stability chamber. Fixing $\vartheta$ and changing $u$ inside that stability chamber may change the set of simple roots and hence the BPS chamber but the BPS spectra does not change by section 8. So if we consider BPS walls of the refined BPS spectra for some $(u,\vartheta)$  then the result does not depend on the choice of $u$ inside a stability chamber and $\vartheta$. In fact the refined BPS spectra of the BPS chamber containing $(u,\vartheta)$ contains exactly half of the BPS spectra at $u$ and $W_{\gamma}=W_{-\gamma}$. This means that we have BPS walls for the whole BPS spectra at $u$ which only depends on the stability chamber. So we can make  the following definition.
\begin{dfn}Fix a generic $\xi$. Let $\Upsilon$ be a stability chamber. Pick a point $u\in\Upsilon$ and define the system of BPS faces of $\Upsilon$ to be the collection of the projections of BPS walls of the refined BPS spectra at $(u,\vartheta)$ to $B^{*}$ ($\xi$ is fixed). Each projection is called a BPS face of $\Upsilon$ of charge $\gamma$  and is denoted by $F_{\gamma}^{\Upsilon}$. If $u$ is not necessarily contained in a stability chamber  then we simply define the BPS face for a BPS charge $\gamma$ at $u$ as the projection of the corresponding BPS wall. It is denoted by $F^{u}_{\gamma}$.  We may omit $u$ if the dependence is clear.\end{dfn}

A BPS face is a codimensional one locus. The intersection of two BPS faces $F_{\gamma_{1}}, F_{\gamma_{2}}$ is a locus (codimensional two in $B^{*}$) contained in the stability wall $SW_{\gamma_{1},\gamma_{2}}(u)$ since at the intersection the two central charges have the same phase (the phase of $\xi$).  According to remark 7.3 and 8.9 we should distinguish the space $B$ of moduli parameters $u$ in the wall crossing formula and the space $B$ which is the base of the Hitchin's fibration even if these two spaces are identical (they are both space of quadratic differentials). Since the singular affine structure is on the base of the Hitchin's moduli space it seems that we should not view a stability wall as a subset of the singular affine manifold $B$ we are dealing with here. However we can  use the canonical isomorphism between the two $B$ spaces to map stability walls to the base of the Hitchin's fibrations and from now on we can consider stability walls in the singular affine space $B$.  If we vary $\xi$ we can recover the whole stability wall. $F_{\gamma}$ is disconnected in $B^{*}$ because the singular locus in $B$ corresponding to the vanishing of charge $\gamma$ has been removed. So the real number $Z_{\gamma}(u)/\xi$ cannot be zero which disconnects $\mathbf{R}$. We can fill in the singular locus back and consider  $F_{\gamma}$ as a subset of $B$ by adding the zero to the image of  $Z_{\gamma}(u)/\xi$. From now on $F_{\gamma}$ is always understood in this sense.
\begin{note}We can fix $\xi$ because we want to consider $\mathcal{M}$ as the moduli space of flat connections whose complex structure does not depend on $\xi$. $\xi$ is assumed to be generic because we want to avoid the following two exceptional situations: an intersection of two BPS faces could be in the singular locus accidently; an intersection of two BPS faces can be contained in more than one stability walls.  In the first situation both the intersection and the singular locus are at least codimensional two. Since we know the position of singular locus is independent of $\xi$ we can avoid them by choosing a generic $\xi$. In the second situation  intersections of stability walls are codimensional two and are  independent of $\xi$. So again we can choose a generic $\xi$ to avoid this situation. But it is important to know that we cannot and do not want to avoid the situation where two or more intersections of BPS faces coincide. Finally we point out that fixing $\xi$ does not mean fixing $\vartheta$ (see remark 8.9).\end{note}

We make the following assumption to simplify  our problem.\\
\textbf{Assumption of Finite Type} \ We assume that there is at least one point $u$ such that the  BPS spectra at $u$ is a finite set. A Hitchin's moduli space satisfying this condition is called a Hitchin's moduli space of finite type. Otherwise it is called a Hitchin's moduli space of infinite type.

There are many examples of finite type Hitchin's moduli spaces. The examples in this paper are all of finite type.  We make this assumption because a Hitchin's moduli space of infinite type would mean that we glue infinitely many pieces of deformations of affine strata in Gross-Siebert's construction. It seems the only reasonable way to make sense of that is to take truncations used before to make sense of the consistency conditions (or the wall crossing formula). That way one can get a projective system of finite gluing. The main theorems in this section should have generalizations to that case.\\

Suppose the assumption of finite type is satisfied.  Take a point $u$ with finite BPS spectra and  build the system of BPS faces at $u$.

 Because of the choice of affine coordinates these BPS faces are hyperplanes in the affine coordinate system. If a BPS charge is a gauge charges this is clear because the BPS face is the zero level set of one of the linear combinations of affine coordinates. For a general charge since we have chosen a branch we can  split the charge lattice into the sum of the gauge charge lattice and the flavor charge lattice. Let $\gamma_{gau}$ and $\gamma_{fla}$ be the projection of a charge $\gamma$ to the gauge charge lattice and the flavor charge lattice respectively. Note that $Z_{\gamma_{fla}}=\int_{\gamma_{fla}}\lambda$ is  a sum of  constants times the residues (mass parameters) of $\lambda$ at some singularities which have been fixed when one formulates the moduli problem. So the central charge of a BPS charge is the central charge of the gauge charge part plus a constant and the claim follows.

 Therefore we have a system of hyperplanes giving us a polyhedral complex. However it may  not be an integral one because of the nontrivial residues.  To get an integral polyhedral complex  we can move hyperplanes labeled by charges $\gamma_{i}$ with a same gauge charge (so that they are parallel) $\gamma_{gau}$ without changing the order of their positions. It can be arranged such that the intersections with the line $\mathrm{Im}(Z_{\gamma_{gau}}(u)/\xi)$ are all rational points (i.e. all components in the coordinate system are rational). This can be done because the line is generated by an integral vector $\gamma_{gauge}$. Clearly this can be achieved by   scaling the residues of flavor charges. Since BPS charges for this particular line may not exhaust a basis of the flavor charge lattice we may need to  do this for the rest of the flavor charges. After scalings of residues (each residue could have a different scaling) we can guarantee that    all polyhedra are rational polyhedra and all vertices are rational points. We can scale residues further to clear all denominators of these rational points since there are only finitely many of them. In the end we get an integral polyhedral complex.   We call the operation of scaling residues to achieve integrality an $integral$ $scaling$ $operation$.

 It is easy to see that after this operation any intersection of two BPS faces is still contained in a stability wall for the moduli space with scaled residues. It is possible that there are finitely many coincide BPS faces. This happens when
 two charges have the same gauge charge and the residues of  the flavor parts are also the same. We therefore allow the possibility of coincide codimensional one cells. Since the whole thing is in $B$ which is isomorphic to $\mathbf{C}^{k}\simeq \mathbf{R}^{2k}$ where $k$ is half the complex dimension of $\mathcal{M}$ for each vertex in the complex we just use the trivial fan structure induced by the embedding of a neighborhood of the vertex in $\mathbf{R}^{2k}$. Thus we have obtained a polyhedral decomposition of the singular integral affine manifold $B$.

 We have proved the following theorem.
\begin{thm}After an integral scaling operation (defined in the previous paragraph) the system of BPS faces at $u$ induces a polyhedral decomposition of the integral affine structure with singularities on $B$. It is called a BPS polyhedral decomposition.\end{thm}

The collection of residues is a part of defining data of the moduli problem. These scalings therefore are perturbations of the original hyperkahler space. Since in the metric degeneration of Hitchin's moduli spaces we change $R$ which scales the residues we want to achieve an integral scaling operation  by picking a point (the choice of such points is a discrete set) from the large complex degeneration. To make sure that such a common scaling for all singularities produces an integral scaling operation we need to assume that   mass parameters (or equivalently residues, see (17)) $m_{i}$ where $i$ runs over the set of all singularities can all become  rational if we multiply them by a common real number. This condition is called the $pseudo-rationality$ condition. Clearly the set of pseudo-rational residues is dense in the space of regular semi-simple residues.

The above construction depends on the choices of $u$. But if $u$ is contained in a stability chamber then clearly it really depends on the stability chamber.

 BPS faces are labeled by some charges only up to a sign. By picking the half of the BPS spectra in the refined BPS spectra  $\Sigma(u,\vartheta)$ we can think of BPS faces as hyperplanes labeled by some charges. Moreover we can label all codimensional one cells of the polyhedral decomposition by the charges labeling BPS faces containing them. This labeling clearly depends only on the BPS chamber containing $(u,\vartheta)$.
 \begin{dfn}A BPS polyhedral decomposition together with a labeling of codimensional one cells by charges defined above is called a marked BPS polyhedral decomposition. It depends on $(u,\vartheta)$ or more precisely the BPS chamber containing $(u,\vartheta)$.\end{dfn}
 \begin{note}It is important to know that the labeling of codimensional one cells is not one to one. In fact all codimensional one cells contained in a single BPS faces are labeled by the same charge.\end{note}
 \begin{note} Later when we say "choose (pick) a pair $(u,\vartheta)$" we mean choose a pair $(u,\vartheta)$ such that $(u,\vartheta)$ is not on a BPS wall and the BPS spectra at $u$ is finite. If $(u,\vartheta)$ is contained in a BPS chamber then we mean choose that BPS chamber.\end{note}

We need to specify  a polarization $\varphi$. We want a nontrivial polarization.
\begin{dfn}Given a BPS polyhedral decomposition. A polarization $\varphi$ is nontrivial if for any vertex $v$ of the BPS polyhedral decomposition there is  a local representative $\varphi_{v}$ with $\varphi_{v}=0$ such that  the following condition is satisfied. Let $\rho_{\gamma}$ be a codimensional one cell containing $v$ labeled by a BPS charge $\gamma$.  Then the maximum of $\mathrm{ord}_{\sigma}(m_{\gamma_{gau}})$ is not equal to the minimum  $\mathrm{ord}_{\sigma}(m_{\gamma_{gau}})$ for $\sigma$ runs over all maximal dimensional cells containing $\rho_{\gamma}$. $\gamma_{gau}$ is the gauge part of the charge $\gamma$ and $m_{\gamma_{gau}}$ is defined in (99).\end{dfn}

A nontrivial polarization always exists for a given BPS polyhedral decomposition. In fact the definition just means that there is a non-trivial change passing from a maximal dimensional cell to another across a codimensional one cell which can clearly  be  arranged as $\varphi_{v}$ is only piecewise linear. There are of course many of them and we can use any of them.  The choice of $\varphi$ cannot be canonical because of the obvious fact that a deformation of the underlying variety of the Hitchin's moduli space simply as a deformation of varieties cannot be canonically determined by the Hitchin's moduli space itself. Two non-isomorphic total spaces could have two isomorphic generic fibers.

\begin{note}Recall that the deformation by changing $R$ is modular.  However we will not expect the same for the toric degeneration we are going to construct. We will be able to at least recover the Hitchin's moduli space (as the moduli space of flat connections)  as a fiber of the toric degeneration. On the other hand we do not expect to be able to embed the real family parameterized by $R$ into the complex family (the toric degeneration) due to the freedom of choosing the polarization.\end{note}

We  need to specify a log smooth structure. A log smooth structure is actually the initial data for constructing  consistent scattering diagrams and the rest are constructed inductively by imposing  consistency conditions which should be identified as essentially  wall crossing formulas. This initial data are associated to codimensional one cells and polyhedral subsets of the polyhedra decomposition. On the other hand, in the formulation of wall crossing formulas in the metric problem, the initial data (part of all Kontsevich-Soibelman factors) is associated to BPS rays.  A BPS wall projects in two directions to a BPS line and a BPS face containing codimensional one cells respectively and therefore it builds a natural correspondence between these two objects to which the initial instanton data of the metric problem and the complex structure problem are associated respectively. So we can define a log smooth structure by just using the correspondence and taking the log morphisms as the discontinuous jumps of Fock-Goncharov coordinates, i.e. \begin{equation}\theta_{\rho_{\gamma}}=\exp(\log(-f_{\rho_{\gamma}})\partial_{n_{\rho_{\gamma}}}):=K_{\gamma}\end{equation}

Let us be precise. First we build a correspondence between monoid variables and charges.  Suppose the complex dimension of the moduli space is $2n$ then the gauge charge lattice denoted by $\Gamma_{gau}$ has dimension $2n$ and an integral skew-symmetric paring which is the intersection pairing. We would like to identify it with $\mathbf{Z}^{2n}$ with the standard integral skew-symmetric paring and with a symplectic basis. But note that while the lattice of integral one cycles of the spectral curve has a symplectic basis  $\Gamma_{gau}$ may not. So in general we  only identify $\Gamma_{gau}$   with a sublattice of rank $2n$ in the lattice of all integral one cycles of the spectral curve with a symplectic basis   such that the intersection paring in $\Gamma_{gau}$ is the restriction of the standard one.

\begin{dfn}Choose a pair $(u,\vartheta)$. Since the BPS spectra is finite there are finitely many stability walls given by pairs of BPS charges at $u$. They are called primary stability walls. Chambers obtained by dividing $B$ by these primary stability walls are called primary stability chambers.\end{dfn}

Let $j$ be a joint in  the marked BPS polyhedral decomposition at $(u,\vartheta)$. It is a codimensional two cell.  We distinguish two cases.\begin{itemize}\item Non-degenerate case. In this case the joint is the intersection of exactly two BPS faces. The joint is contained in a unique primary stability wall denoted by $SW_{j}$. By the genericity of the choice of $\xi$  we can exclude the exceptional case where one of the two BPS face defining $SW_{j}$ is tangent to the stability wall. At least locally in a small neighborhood containing $j$ the stability wall $SW_{j}$ divides the neighborhood into two sides. The side containing $u$ is denoted by $Side_{u}$ while the other side is denoted by $Side_{u}^{\bot}$. Since we allow different  BPS faces with the same support (coincide BPS faces) we also allow different joints with the same support. This is harmless as the Kontsevich-Soibelman transformations of charges labeling coincide BPS faces commute.\item Degenerate case. The joint is the intersection of more than two BPS faces. By the genericity of $\xi$ the joint is still contained in just one primary stability wall. "Degenerate" does not suggest that this is non-generic in some sense.  We have this situation simply because we can have BPS charges generated by other BPS charges. Since there is a stability wall we still have  $Side_{u}$ and $Side_{u}^{\bot}$.  \end{itemize}

Let $v$ be a reference vertex in $j$, $v\in \omega\subseteq\sigma\in\mathcal{F}_{max}$. Since $\mathrm{Spec}\mathbf{C}[P_{\omega,\sigma}]\subseteq \mathrm{Spec}\mathbf{C}[P_{v}]$  to build the correspondence we only need to consider the monoid $\mathbf{C}[P_{v}]$. Recall that $P_{v}$ is generated by certain exponents $m$ and each $m$ has a projection $\bar{m}$ in $\Lambda_{v}\simeq \Gamma_{gau}$. So we get an correspondence between gauge charges $\gamma$ and  exponents with bars denoted by $\bar{m}^{'}_{\gamma}$. We put a $'$ over $\bar{m}^{'}_{\gamma}$ because the notation $\bar{m}_{\gamma}$ is reserved for a modification of $\bar{m}^{'}_{\gamma}$.

Recall that the labeling of codimensional one cells by charges labeling the BPS faces containing them is many to one. In $U_{j}=\cup_{\tau,\tau\supseteq j}Int\ \tau$ there are exactly two codimensional one cells containing $j$ and contained in a BPS face.  In $Side_{u}$ there is only one codimensional one cell for each BPS face and these codimensional one cells are ordered according to the ordering of their BPS charges in $\Sigma(u,\vartheta)$. Therefore in $Side_{u}$ each codimensional one cell can be labeled by its corresponding BPS charge and different cells are labeled differently. A codimensional one cell labeled by $\gamma$ is denoted by $\rho_{\gamma}$. It is parallel to  the zero level set of the affine coordinate labeled by the gauge charge $\gamma_{gau}$ and we  associate to it the charge $\gamma_{gau}$ and hence $\bar{m}^{'}_{\gamma}$. This association is defined up to sign because $W_{\gamma}=W_{-\gamma}$.  We fix the sign in the following way. It makes sense to talk about the positive  normal vectors of a codimensional one cell in $Side_{u}$ as a positive direction can be defined to be the direction moving along which these cells are encountered by their order in $\Sigma(u,\vartheta)$ (the choice of this direction is the same as a choice of orientation in the normal space of the joint). The sign of $\gamma_{gau}$ and hence $\bar{m}^{'}_{\gamma}$ is chosen such that $\bar{m}^{'}_{\gamma}$ is a positive normal vector. In $Side_{u}^{\bot}$ we define the positive direction to be compatible with the positive direction defined in $Side_{u}$. Each codimensional one cell in $Side_{u}^{\bot}$ is also labeled  by the BPS charge labeling the corresponding BPS face and we can define the positive  normal vectors to it.  $-\bar{m}^{'}_{\gamma}$ is a  positive normal vector to $\rho_{\gamma}$ in $Side_{u}^{\bot}$.

 In $Side_{u}$ we rotate $\bar{m}^{'}_{\gamma}$ around $j$ along the negative direction (more precisely speaking  we rotate its projection to the two dimensional normal space of $j$ along the induced negative direction  while keeping other components invariant) until it lies in $\rho_{\gamma}$ and call the vector obtained after this rotation $\bar{m}_{\gamma}$. It is uniquely defined and pointing away from $j$. In $Side_{u}^{\bot}$ we do the same thing to $-\bar{m}^{'}_{\gamma}$ and define $-\bar{m}_{\gamma}$. Clearly $-\bar{m}_{\gamma}=\bar{m}_{-\gamma}$. Note that although the same symbol $\rho_{\gamma}$ represents two codimensional one cells lying in $Side_{u}$ and $Side_{u}^{\bot}$ respectively the exponents with bars (or the corresponding  gauge charges) attached to them  are different for them.

 Since by section 8 the refined BPS spectra contains a complete set of generators of the  charge lattice  we can extend the correspondence between $\gamma$ and $\bar{m}_{\gamma}$ to the whole gauge charge lattice. This induces a multiplicative correspondence between  gauge charges and monoid variables $$\tilde{x}_{\gamma}:=z^{\bar{m}_{\gamma}}$$ where $z$ is the formal  variable converting additive relations to multiplicative relations in the monoid. For a gauge charge $\gamma$ there is a unique  monoid in $P_{v}$ as a generator. The monoid is denoted by $x_{\gamma}$ and is associated to the exponent which lies on the boundary of the cone $C_{v}$ and with its projected exponent (the exponent with bar)  $\bar{m}_{\gamma}$. Explicitly we can choose a local representation $\varphi_{v}$ of the polarization with $\varphi_{v}(v)=0$. Then \begin{equation}m_{\gamma}:=(\bar{m}_{\gamma},h),\ x_{\gamma}:=z^{(\bar{m}_{\gamma},h)}\end{equation}where $h$ is given by\begin{equation}\langle\bar{m}_{\gamma},-\lambda_{\sigma}\rangle +h=0\end{equation} for some  maximal dimensional cell $\sigma$ such that $h$ determined by (100) is maximal among all $\sigma\ni v$. $P_{v}$ is generated by finitely many monoid variables of  this form and the deformation parameter $t:=z^{(0,1)}$.  We also introduce formal multiplicative notations  $x_{\gamma_{fla}}$ for  pure flavor charges. Then we define\begin{equation}x_{\gamma}:=x_{\gamma_{gau}}x_{\gamma_{fla}}\end{equation}for a general charge $\gamma$ and we have a correspondence between $x_{\gamma}$ and $\gamma$ (and therefore $\mathcal{X}_{\gamma}$) for all charges. We  think of  $x_{\gamma_{fla}}$ as  trivial monoid variables. For computations we can think of them  as $x_{\gamma_{fla}}:=z^{(0,0)}$. So the only effect they have is that they change names of monoid variables $x_{\gamma_{gau}}$ i.e. the same monoid variable $x_{\gamma_{gau}}$ can be represented by several symbols $x_{\gamma}$ in a formula. Using labeling by general charges we can distinguish coincide codimensional one cells.

\begin{dfn} The labeling of a codimensional one cell  by the BPS charge of the BPS face containing it is called the face labeling. Let $j$ be a joint of  the marked BPS polyhedral decomposition for $(u,\vartheta)$. The labeling of a codimensional one cell containing $j$ by the  exponent with bar (or the corresponding gauge charge) is called the slab labeling.\end{dfn}

 Each cell $\rho_{\gamma}$ is considered to be the support of a slab with a slab function to be determined below. In section 3 we have explained that the log morphism attached to a slab $\rho_{\gamma}$ is \begin{equation}z^{(\bar{m}_{\gamma_{gau}^{'}},h)}\rightarrow (f_{\mathbf{\rho_{\gamma}},v})^{-\pi(\bar{m}_{\gamma_{gau}^{'}})}z^{(\bar{m}_{\gamma_{gau}^{'}},h)}\end{equation}where $\pi:\Lambda_{\sigma}\rightarrow\mathbf{Z}$ is the epimorphism with kernel $\Lambda_{\rho_{\gamma}}$ which is positive on vectors pointing from one of the two maximal dimensional cells containing $\rho_{\gamma}$ to the other  with the positive direction according to the order of elements of the refined BPS spectra. In fact $\pi(\bar{m}_{\gamma_{gau}^{'}})$ is the pairing  (standard inner product of the lattice after picking a basis) between the  positive primitive normal vector $n_{\rho_{\gamma}}$ of $\rho_{\gamma}$ and $\bar{m}_{\gamma_{gau}^{'}}$. Note that the pairing of $\gamma_{gau}$  and $\bar{m}_{\gamma_{gau}^{'}}$ is a positive integral multiple of $\pi(\bar{m}_{\gamma_{gau}^{'}})$. Let $l$ be that integer. $l$ does not depend  on $\gamma_{gau}^{'}$.

 By elementary geometry we get that in terms of the skew-symmetric integral paring (intersection pairing) $\langle\cdot\rangle$ on the lattice $\Lambda_{v}\simeq \Gamma_{gau}$, \begin{equation}-l\pi(\bar{m}_{\gamma_{gau}^{'}})=\langle\gamma_{gau}^{'},\gamma_{gau}\rangle=\langle\gamma^{'}_{gau},\gamma\rangle\end{equation}
 In fact the problem is local and two dimensional as the joint is codimensional two. So we can write $\gamma_{gau}$ as $(p,q)$ and $\gamma_{gau}^{'}$ as $(p^{'}, q^{'})$. The intersection pairing is the standard one $\langle(p_{1},q_{1}),(p_{2},q_{2})\rangle=p_{1}q_{2}-p_{2}q_{1}$ and the inner product is given by $((p_{1},q_{1}),(p_{2},q_{2}))=p_{1}p_{2}+q_{1}q_{2}$. $\bar{m}_{\gamma_{gau}^{'}}$ then is $(q^{'},-p^{'})$ and (103) follows immediately.

 So (102) is \begin{equation} x_{\gamma_{gau}^{'}}\rightarrow x_{\gamma_{gau}^{'}}(f_{\mathbf{\rho_{\gamma},},v})^{l^{-1}\langle\gamma_{gau}^{'},\gamma\rangle}\end{equation}This can be extended in the trivial way to \begin{equation}x_{\gamma^{'}}\rightarrow x_{\gamma^{'}}(f_{\mathbf{\rho_{\gamma}},v})^{l^{-1}\langle\gamma^{'},\gamma\rangle}\end{equation}for any charges $\gamma$ and $\gamma^{'}$.

 On the other hand $\gamma$ also labels Fock-Goncharov coordinates and Kontsevich-Soibleman transformations  are  \begin{equation}K_{\gamma}^{\Omega(\gamma;u)}:\mathcal{X}_{\gamma^{'}}\rightarrow \mathcal{X}_{\gamma^{'}}(1-\sigma(\gamma)\mathcal{X}_{\gamma})^{\Omega(\gamma;u)\langle\gamma^{'},\gamma\rangle}\end{equation} Here we assume that we have chosen $(u,\vartheta)$ and we vary $\vartheta$ to $\vartheta+\pi$. So $\gamma$ runs over all refined BPS charges labeling $\rho_{\gamma}$ in $Side_{u}$. If we make the identification $\mathcal{X}_{\gamma}\rightarrow x_{\gamma}$ then (106) becomes (105) if we identify $f_{\rho_{\gamma},v}$ by the following equation. \begin{equation}f_{\rho_{\gamma},v}=(1-\sigma(\gamma)x_{\gamma})^{l\Omega(\gamma;u)}\end{equation}
We still need to define $f_{\rho_{\gamma},v}$ attached to $\rho_{\gamma}$ in $Side_{u}^{\bot}$. Note that with respect to the positive direction in $Side_{u}^{\bot}$ the cells are encountered in the order from $\vartheta+\pi$ to $\vartheta$. Since we want to get $(K_{\gamma}^{\Omega(\gamma;u)})^{-1}$ in this case $f_{\rho_{\gamma},v}$ attached to $\rho_{\gamma}$ in $Side_{u}^{\bot}$ must be still $$f_{\rho_{\gamma},v}=(1-\sigma(\gamma)x_{\gamma})^{l\Omega(\gamma;u)}$$
It is important to know that even though the two cells labeled by $\gamma$ in $U_{j}$ have different attached monoid variables $x_{\gamma}$ and $x_{-\gamma}$ (note that $x_{\gamma}\neq -x_{-\gamma}$ although $\bar{m}_{\gamma_{gau}}$=$-\bar{m}_{-\gamma_{gau}}$) the slab functions attached to them are identical. Slab functions are labeled in the face labeling. That is why we sometimes use the ambiguous face labeling $\rho_{\gamma}$.

If there are several coincide codimensional one cells then we have several slab functions attached to the same support of several slabs. It is easy to see that the product of Kontsevich-Soibelman transformations associated to them does not depend on the order of the product and is induced by  the product of these slab function.  These facts  follow from the following obvious properties of Kontsevich-Soibelman transformations.\begin{pro}If the intersection pairing between  two charges is zero then the two associated Kontsevich-Soibelman transformations commute. The Kontsevich-Soibelman transformation associated to a pure flavor charge is the identity. Any Kontsevich-Soibelman transformation acts on a Fock-Goncharov coordinate labeled by a pure flavor charge trivially.\end{pro}

 We have defined slab functions for a joint with a reference vertex.  If we change the reference vertex $v$ we get parallel transports of exponents in  slab functions and slab functions have to satisfy the formula (35). So  for each support of slabs we simply pick an initial  reference  vertex and use (107) to define the slab function with respect to that vertex and  then use (35) to define the slab function on the same support with respect to other reference vertices. By the definition of $m^{\rho}_{v^{'}v}$ if we swap two vertices the formula still holds and the definition is also consistent  along a loop avoiding $\Delta$. As explained in section 3 this definition guarantees that the log morphism attached to a slab does not depend on the choice of reference vertices.

 After that we use (36) to define $f_{e}$. This is in fact unnecessary because we might take slab functions inducing the log smooth structure instead of the log smooth structure itself as the input of our construction (which is actually more convenient as we shall see). Nevertheless for completeness let us verify that this produces a log smooth structure. (23) in the definition of log smooth structures is just (35). Now consider all codimensional one cells containing a joint $j$.   For each codimensional one cell with  positive primitive normal vectors $n_{\rho_{\gamma}}=\check{d}_{\rho_{\gamma}}$ there is another one with the opposite positive primitive normal vectors. They are labeled by $x_{\gamma}$ and $x_{-\gamma}$ respectively in the slab labeling. $f_{e}$ associated to these two cells may be different even if $f_{\rho_{\gamma}}$ are the same. But the difference arises only from the change of linear functions of the polarization across codimensional one cells and are contributed by exponents with nonzero components in the normal directions to $j$. When we restrict them to $V(j)$   those components are killed. So the restrictions of $f_{e}$ associated to the two cells containing $j$ are the same. This means (22) holds.  Therefore the collection of $f_{e}$ defines a log smooth structure. The only way that a negative exponent can appear in (107) is that the corresponding BPS charge is a vectormultiplet. Because of the assumption of finite type and the fact that one needs to cross infinitely many flips to reach a vectormultiplet  the exponent in (107) is always positive. Then the log smooth structure is clearly positive as all slab functions are pole free.
 \begin{dfn}The positive log smooth structure defined above is called the BPS log smooth structure at $(u,\vartheta)$. It really depends on the BPS chamber containing $(u,\vartheta)$. The collection of a BPS polyhedral decomposition (denoted by $\mathcal{F}_{BPS}$), a nontrivial polarization for this decomposition and a BPS log smooth structure (denoted by $Log_{BPS}$) is called a BPS Gross-Siebert data and is denoted by $GS_{BPS}(u,\vartheta)$.\end{dfn}

 Let us summarize. Let $\mathcal{M}$ be our moduli space. More precisely it is the moduli space of $SL(2,\mathbf{C})$ flat connections over a Riemann surface $C$ with $l$ singularities with prescribed singular parts (15) at singularities modulo the  gauge equivalence. We assume that for every singularity  $T_{n}$ is regular semi-simple.  $l>0$ and if the underlying Riemann surface is $CP^{1}$ then the number $l$ should  make the dimension of $\mathcal{M}$ positive (e.g. if all singularities are regular then $l>3$). $\mathcal{M}$ is also the moduli space of solutions of Hitchin's equations (1) or (2) with singularities with prescribed singular parts  (local models of abelian singularities and no fractional exponents, see section 2 and remark 7.10)  modulo gauge equivalence and as such it is endowed with a hyperkahler structure. $\mathcal{M}$ as the moduli space of flat connections is realized in complex structures parameterized by $\xi\in\mathbf{C}^{\times}$.  We also assume that the set of residues is pseudo-rational which can  be achieved by an arbitrarily small perturbation of residues of the original moduli problem. Let $R$ be a positive real number large enough such that Fock-Goncharov coordinates are piecewise holomorphic (pole free) in the twistor $\mathbf{C}^{\times}$ and let $\mathcal{M}(R)$ be the corresponding moduli space using the modified equations  (3). The moduli interpretation as a moduli space of flat connections does not depend on $R$. BPS spectra are also independent of $R$. We assume that $\mathcal{M}$ (and hence $\mathcal{M}(R)$) is  of finite type.  We pick $(u,\vartheta)$ such that at $u$ the BPS spectra is finite. We have proved the following theorem (it is actually a definition-theorem).
 \begin{thm}Fix a generic $\xi$. After  choosing an $\mathcal{M}(R)$ which provides an integral scaling operation of $\mathcal{M}$, over the singular integral affine manifold $B$ which is the base of the Hitchin's fibration  there is a BPS Gross-Siebert data.   The choice of the nontrivial polarization is not unique. The marked BPS polyhedral decomposition and the BPS log smooth structure depend on $(u,\vartheta)$. If $(u,\vartheta)$ is contained in a BPS chamber then $\mathcal{F}_{BPS}$ and $Log_{BPS}$ depend only on the choice of the BPS chamber. \end{thm}

 \noindent {\bf Example 1}\ \ In this example we consider a moduli space of Hitchin's equation with one irregular singularity at $\infty$ over the Riemann sphere $CP^{1}$. Strictly speaking this is not an example for theorem 9.4 because the order of quadratic differentials is odd.   Nevertheless this example is the simplest case with a nontrivial wall crossing formula. As explained in section 7  the construction of Fock-Goncharov coordinates  is still valid and so is the construction of BPS Gross-Siebert data. Moreover we expect that the results in  section 9 generalize to  odd order cases once the corresponding Hitchin's moduli spaces have expected properties.

 We need to specify the asymptotic behaviors of solutions of Hitchin's equations near $z=\infty$. Let $A_{0}$ and $\varphi_{0}$ be singular parts of the connection $A$ and the Higgs field $\varphi$ near $\infty$ respectively. So as $z\rightarrow \infty$ $$A\sim A_{0}, \varphi\sim \varphi_{0}$$

We prescribe the asymptotic behaviors by first choosing a quadratic differential $\lambda^{2}=P_{N}(z)dz^{2}$ where $P_{N}(z)$ is a monic polynomial of degree $N$. Define $\Delta(z)$ to be the singular part  of the expansion of $\sqrt{P_{N}}$ (up to a sign) near $\infty$ $$\sqrt{P_{N}(z)}=\Delta(z)+ o(z^{-1})$$ So $\lambda^{2}\sim\Delta(z)^{2}$. When $N$ is even $\Delta(z)$ does not contain non-integral powers and we can  diagonalize the Higgs field $\varphi_{0}$ (note the the trace must be zero since it is in $sl(2)$)$$\varphi_{0}=\left(\begin{array}{cc}\Delta(z) & 0\\0 & -\Delta(z) \end{array}\right)$$ so that $-\lambda^{2}\sim det \ \varphi_{0}, z\rightarrow \infty$ as required by the definition of Hitchin's fibration. This prescribes the asymptotic behavior of the Higgs field. In our case, however,  $N=3$. Since it is odd we cannot diagonalize $\varphi_{0}$. Instead for odd $N$ one can find a complex gauge transformation to put it into the form (see \citep{G2})$$\varphi_{0}=\Delta(z)\left(\begin{array}{cc}0 & (\bar{z}/z)^{1/4}\\(z/\bar{z})^{1/4} & 0 \end{array}\right)$$We prescribe $A_{0}$ compatibly $$A_{0}={1\over 8}\left(\begin{array}{cc}1 & 0\\0 & -1 \end{array}\right)({dz\over z}-{d\bar{z}\over \bar{z}})$$

The example we now begin to  study is called Argyres-Douglas theory. We choose a  quadratic differential by choosing $$P_{3}(z)=z^{3}-3\Lambda^{2}z$$ Here $\Lambda$ is a positive real constant. As just explained we use it to prescribe  asymptotic behaviors of solutions of the Hitchin's equations and this define the  Hitchin's moduli space. Now considered the base $B$ of Hitchin's fibrations which contains $z^{3}-3\Lambda^{2}z$ and its deformations. Deformations can only be of the following form $$\lambda^{2}=(z^{3}-3\Lambda^{2}z+u)dz^{2}$$where $u$ is a complex number. Otherwise the corresponding behaviors determined by $\lambda^{2}$ would not be the prescribed ones. This means $B$ is complex one dimensional and $u$ is the moduli parameter.

The spectral curve is clearly an elliptic curve whose Jacobian is itself. Since the Jacobian of $CP^{1}$ is trivial, the Prym is the trivial quotient of the Jacobian of the elliptic curve. Therefore the fiber are spectral curves themselves. Each of them has two cycles $\gamma_{1}, \gamma_{2}$ and they generate the gauge charge lattice with $\langle\gamma_{1},\gamma_{2}\rangle=1$. The flavor charge lattice is trivial in this example. The singularities of the affine base $B$ are the zero locus of  the discriminant of $P_{3}(z)$. There are two of them given by $u_{1}=-2\Lambda^{3}$ and $u_{2}=2\Lambda^{3}$. According to the Picard-Lefschetz transformation, we can choose the two nontrivial counterclockwise monodromies around $u_{1}, u_{2}$ to be $$\left(\begin{array}{cc}1 & 1\\0 & 1 \end{array}\right) \ \mathrm{and}\ \left(\begin{array}{cc}1 & 0\\-1 & 1 \end{array}\right)$$The union of stability walls is given by $$W:= \{u\mid \arg Z_{\gamma_{1}}(u)=\pm \arg Z_{\gamma_{2}}(u)\}$$It is a simple closed curve passing though $u_{1}, u_{2}$. It is symmetric with respect to the reflection by the real axis (the line connecting $u_{1}$ and  $u_{2}$).  We take two branch cuts along the real axis from $u_{1}$ to $-\infty$ and from $u_{2}$ to $+\infty$ respectively. \\

For solutions of Hitchin's equations there is a single singularity at $z=\infty$ with order seven and hence it is an irregular one. It has five Stokes sectors. By the prescription given in section 7, we obtain a WKB triangulation for a generic $\vartheta$. This triangulation (of the complement of a small disk containing $\infty$) is a triangulation of a pentagon. There are three triangles with each containing a simple zero of $\lambda^{2}$. According to the identification of charges (cycles) and edges, there are two charges associated to the two inner edges. Fock-Goncharov coordinates for outer edges are defined to be zero according to section 7, so we only need to study inner edges.

 The stability chamber inside the union of stability walls is also called the strong coupling region. Let us start from the point $u=0$ where it is clear that the three zeroes of the quadratic differential is colinear. The two associated charges (up to a sign) are $\pm\gamma_{1}$ and  $\pm\gamma_{2}$. Since a refined BPS spectra which covers exactly half of the BPS spectra is generated by simple roots which are associated to these two charge we know the BPS spectra at $u=0$ is $(\pm\gamma_{1},\pm\gamma_{2})$. Let us assume the phase of $l_{\gamma_{1}}$ is smaller than the phase of $l_{\gamma_{2}}$. The BPS spectra remain to be the same as along as we stay inside this stability chamber. So we can assume $u$ is any value inside $W$ in the following discussion. Now we vary $\vartheta$ from $\vartheta$ to $\vartheta+\pi$. Varying $\vartheta$ generically means we change the triangulation by isotopy and  it is clear that we meet two flips of edges labeled by $(\pm)\gamma_{1}$ (only one of them) and  $(\pm)\gamma_{2}$ (only one of them). Without loss of generality we assume that the two flips are associated to $\gamma_{1}$ and $\gamma_{2}$ respectively. So the (ordered) refined BPS spectra at $(u,\vartheta)$ is $(\gamma_{1},\gamma_{2})$. This is also the refined simple BPS spectra.

 We take the projections of the BPS walls to the $u$ plane for a fixed generic $\xi$ and obtain two BPS faces $$\{u\mid\mathrm{ Im}(Z_{\gamma_{1}}(u)/\xi)=0)\}, \{u\mid \mathrm{Im}(Z_{\gamma_{2}}(u)/\xi)=0)\}$$In the corresponding affine coordinates, these are two lines intersecting at a point on the stability wall. The two lines pass through two singularities respectively. Without loss of generality we assume  that intersection point is on the upper half of the stability wall. We  identify $\bar{m}_{\gamma_{1}},\bar{m}_{\gamma_{2}}$ with $(-1,0),(0,-1)$ in $\mathbf{Z}^{2}$ (so $\gamma_{1}$ and $\gamma_{2}$ are $(0,-1)$ and $(1,0)$ respectively). The two rays generated by $\bar{m}_{\gamma_{1}}$ and $\bar{m}_{\gamma_{2}}$ are positioned as the negative "$\mathrm{x}$-axis" and negative "$\mathrm{y}$-axis" and the singularity with counterclockwise monodromy $\left(\begin{array}{cc}1 & 0\\-1 & 1 \end{array}\right)$ is on the negative part of the "$\mathrm{y}$-axis" and singularity with the counterclockwise monodromy $\left(\begin{array}{cc}1 & 1\\0 & 1 \end{array}\right)$ is on the negative part of the "$\mathrm{x}$-axis"\footnote{This choice is made such that our picture matches the one in \citep{GS3} where the consistency condition of this example is  discussed.}.

Now we have a polyhedral decomposition where the intersection is the only vertex $v$. The two lines give us  four rays which are considered as codimensional one cells. They are denoted as $\rho_{i}, 1\leq i\leq 4$ in the counterclockwise order such that $\rho_{1}$ is the positive part of the $\mathrm{x}$-axis. They decompose the plane into four quadrants which are maximal dimensional cells. These cells are denoted by $\sigma_{i}, 1\leq i\leq 4$ in the counterclockwise order such that $\sigma_{1}$ is the first quadrant. The fan structure at the vertex is the obvious one induced by the standard one of the plane.

Following Gross and Siebert we pick the polarization $\varphi$ to be the piecewise linear function   $$\varphi\mid_{\sigma_{1}}=\mathrm{x+y},\varphi\mid_{\sigma_{2}}=\mathrm{y},\varphi\mid_{\sigma_{3}}=\mathrm{0},\varphi\mid_{\sigma_{4}}=\mathrm{x}$$ So in particular $\varphi(1,0)=1, \varphi(0,1)=1, \varphi(-1,0)=0, \varphi(0,-1)=0$.

We specify the log smooth structure (the initial data for a $structure$) by using Kontsevich-Soibelman transformations to slabs. In this case, these are the four rays starting from the vertex and they are labeled by charges $\gamma_{1}, \gamma_{2}$. The orientation of $B$  is the one under which the composition order of Kontsevich-Soibelman transformations for $u$ inside the union is the counterclockwise order over $B$ while for  $u$ outside the union it is the clockwise order.

First we introduce monoid variables. We define the deformation parameter $t=z^{(0,1)}$ where  $0=(0,0)$ is the zero in the underlying rank two charge lattice $\mathbf{Z}^{2}$  of the problem. We also define\footnote{There might be a clash of notations here. Note the difference between $x,y$ and $\mathrm{x}, \mathrm{y}$. Also in the third identity, $z$ has different meanings.} \begin{equation}x:=z^{((1,0),1)}, y:=z^{((-1,0),0)},z:=z^{((0,1),1)},w:=z^{((0,-1),0)}\end{equation} Here $(1,0), (-1,0), (0,1), (0,-1)$  are the primitive generator of the positive (negative) $\mathrm{x}$-axis and the primitive generator of the positive (negative)  $\mathrm{y}$-axis respectively and as such are the (projections of) exponents or gauge charges. Since we identify $\bar{m}_{\gamma_{1}},\bar{m}_{\gamma_{2}}$ with $(-1,0),(0,-1)$,  $\{x,y, z,w\}$ are $\{x_{-\gamma_{1}},x_{\gamma_{1}},x_{-\gamma_{2}},x_{\gamma_{2}}\}$ respectively. $\{x,y,z,w,t\}$ generate $P_{v}$.

Clearly $xy=t, zw=t$ and these are relations in the ring $P_{\rho_{i}\rightarrow\rho_{i}, \sigma_{i/i-1}}$. Here the last components of the exponents ($1, 0,1,0$ for $x,y,z,w$ respectively) are chosen according to the four values of $\varphi$ given above. $\varphi$ is picked so that the power of $t$ in the relation $xy=t, zw=t$ is one. If we replace it by some powers of $t$ by changing $\varphi$, the defining equations obtained in the end will also be changed by replacing $t$ by its powers.

With these notations, we can write down the slab functions (which induce the log smooth structure) induced by Kontsevich-Soibelman transformations associated to BPS charges $\gamma_{1},\gamma_{2}$. They are\footnote{We have used the fact that the quadratic refinements are one because we have flips.} \begin{equation}f_{\rho_{1},v}=f_{\rho_{3},v}=1+x^{-1}t=1+y, f_{\rho_{2},v}=f_{\rho_{4},v}=1+z^{-1}t=1+w\end{equation}

\subsection{Equivalence Of Two Instanton Correction Problems}

Suppose we are given a BPS Gross Siebert data $GS_{BPS}(u,\vartheta)$. As explained in  section 3 the log smooth structure already gives us a $structure$ consistent to order zero. We want to show that it is the order zero part of a compatible system of consistent $structures$. There are two equivalent ways to do that. We can forget about the wall crossing formula and just follow the procedure described in section 3. Or we can take advantages of the wall crossing formula by building a system of $structures$ from it and then verify that the system is a compatible system of consistent $structures$. Since such a system is determined by the initial data i.e. the log smooth structure the two approaches yield the same result. Here we will follow the second route which is easier. \\

\textbf{Step I}: Universal $Structures$

 Consider the initial $structure$ which consists of only slabs. In other words every codimensional cell $\rho_{\gamma}$ in the face labeling in the marked  BPS polyhedral decomposition is the support of a slab and the slab function is $f_{\rho_{\gamma}}$. It is clear that a joint $j$ can only be a codimensional two joint. Pick a reference vertex $v\in j$. The consistency discussed here does not depend on this choice. The elements of the  refined BPS spectra $\Sigma(u,\vartheta)$ are ordered according to their BPS rays' counterclockwise order in the twistor $CP^{1}$. We denote this ordered finite set by $(\gamma_{1},\gamma_{2},\cdots, \gamma_{r})$. We also have a finite set of primary stability walls \begin{equation}(SW_{\gamma_{i},\gamma_{l}}), i\neq l, 1\leq i, l\leq r\end{equation} dividing $B$ into finitely many primary stability chambers. Since the set $\Sigma(u.\vartheta)$ has exhausted all charges appearing in the refined BPS spectra other stability walls are not relevant for wall crossings involving only these charges. Of course after we apply the wall crossing at these primary stability walls more BPS charges will appear which introduces more stability walls. But as long as we truncate the wall crossing formulas involved  then at any given order (or degree) there are only finitely many BPS charges and we  continue having stability chambers obtained by dividing $B$.  Note that  a stability chamber at certain order could be destroyed in the next order. Recall that in a small neighborhood of $j$ there are two sides $Side_{u}$ and $Side_{u}^{\bot}$. The intersections of these two sides with  the normal space of $j$ are denoted by the same notations. In the normal space they can not be connected by a continuous path without intersecting some curve which is the intersection of the normal space  and a primary stability wall. Otherwise following the path the order of pre-BPS rays can not change contradicting the definition of $Side_{u}$ and $Side_{u}^{\bot}$. We orient the two dimensional normal space such that in $Side_{u}$ the order of the codimensional one cells labeled by refined BPS charges  is the counterclockwise order. Note that in $Side_{u}^{\bot}$ the order is clockwise with respect to this orientation.

Let us consider the non-degenerate case first.  $j$ is contained in a unique primary stability wall. Assume it is $SW_{\gamma_{i},\gamma_{l}}$. Clearly in the normal space denoted by $Q_{j}$ we have obtained a scattering diagram at $v$ for $j$ by taking the projections of codimensional one cells containing $j$
(they are labeled by $\gamma_{i},\gamma_{l}$) together with slab functions as cuts. $\gamma_{i}$ and $\gamma_{l}$ each in the face labeling labels precisely two codimensional one cells (slabs). The two slabs with the same charge must lie in different sides. Recall that the associated log morphisms (or Kontsevich-Soibelman transformations) on the two slabs with the same charge are inverse to each other.  Now we can apply the wall crossing  formula at this stability wall and obtain an identity \begin{equation}K_{\gamma_{l}}K_{\gamma_{i}}=S(\vartheta^{j}_{-},\vartheta^{j}_{+};u^{'})=K_{\gamma_{i}}S^{'}(\vartheta^{j}_{-},\vartheta^{j}_{+};u^{'}) K_{\gamma_{l}}\end{equation}where $S(\vartheta^{j}_{-},\vartheta^{j}_{+};u^{'})$ is the ordered product  of all Kontsevich-Soibelman transformations over refined BPS charges at $(u^{'},\vartheta)$ between $\vartheta^{j}_{-}$ and $\vartheta^{j}_{+}$ on the other side of the stability wall. $S^{'}(\vartheta^{j}_{-},\vartheta^{j}_{+};u^{'})$ is defined by the second equality. Here we have assumed that between $\vartheta^{j}_{-}$ and $\vartheta^{j}_{+}$ and at $u$ one only encounters $\gamma_{i},\gamma_{l}$ which can always be arranged due to the finiteness of the BPS spectra. $u^{'}$ is on the other side of the stability wall and is close enough to $u$ such that no BPS ray passes $\vartheta^{j}_{-}$ or $\vartheta^{j}_{+}$ when we fix them and move from $u$ to $u^{'}$.\begin{note}The positions and exponents of $K_{\gamma_{i}},K_{\gamma_{l}}$ reflect the inverse relation of Kontsevich-Soibelman transformations of slabs with same charges. It is consistent with the uniqueness of factorization. In fact since $\gamma_{i},\gamma_{l}$ are clearly simple roots by truncating to the lowest nontrivial degree we see that we must include these two factors.  \end{note}

For each charge $\gamma$ that appears in $S^{'}(\vartheta^{j}_{-},\vartheta^{j}_{+};u^{'})$ we add a ray $\mathbf{t}_{\gamma}$ starting from $v$ (which is identified as the origin of $Q_{j}$) and with underlying line $\bar{\bar{\rho_{\gamma}}}$. We require that the direction of the ray must be pointing to the side of the stability wall containing $u^{'}$.

We can do the same thing in a degenerate case. The only difference is that we use the wall crossing formula which reverses completely the order of elements of $\Sigma(u,\vartheta)$ that appear in $Side_{u}$ at the joint and add rays labeled by BPS charges appearing in the wall crossing formula.

If a BPS charge appears on the side containing $u$ but disappears on the other side then we modify the slab function by setting it to 1. This is a possible scenario. For example in   the pentagon relation (134)   if we start from the side with $K_{\gamma_{1}}K_{\gamma_{1}+\gamma_{2}}K_{\gamma_{2}}$ then we need to modify the the slab function attached to the ray contained in the BPS face $\rho_{\gamma_{1}+\gamma_{2}}$ and lying on the other side of the stability wall to be 1. If the slabs function of a cut (corresponding to a slab) is modified to 1 in $Side_{u}^{\bot}$ we  delete that cut (slab) in $Side_{u}^{\bot}$. It is possible that there are  rays whose supports coincide with the support of a cut. Clearly when we compose the ordered product in $Side_{u}^{\bot}$ the order of log morphisms attached to this cut (slab) and these rays is ambiguous. This is fine because the product is invariant with different orders of them by proposition 9.3. This unordered product of Kontsevich-Soibelmann transformations can  be replaced by a single one. Instead of adding these rays we modify the slab function to be the product of the original slab function and the ray  functions. The Kontsevich-Soibelman transformation associated to the modification is precisely the  unordered product of Kontsevich-Soibelmann transformations. We always do such a modification whenever there are rays coincide with a cut. This is consistent with the requirement  that a wall and a slab are not allowed to share the same support.

The following propositions is  trivial by our construction.
\begin{pro}The rays added to the scattering diagram according to BPS charges appearing in the wall crossing formula are projections of BPS faces at $u^{'}$\end{pro}

We let $\mathbf{p}_{\gamma}$ be the wall $(j-\mathbf{R}_{\geq0}\bar{m}_{\gamma})\cap\sigma$ where $\sigma$ is the maximal dimensional cell whose projection to $Q_{j}$ is the scattering chamber containing the ray $\mathbf{t}_{\gamma}$. Note that $\mathbf{p}_{\gamma}$ is a half plane of the BPS face of $\gamma$. Here we assume that the choice of $\bar{m}_{\gamma}$ is such that the projection of the wall is the ray $\mathbf{t}_{\gamma}$ (other wise we use $\bar{m}_{-\gamma}$). We attach the function $f_{\rho_{\gamma}}$ to $\mathbf{p}_{\gamma}$. In this way $\mathbf{p}_{\gamma}$ becomes a Gross-Siebert wall.

\begin{dfn}We add Gross-Siebert walls labeled by BPS charges in the wall crossing formula to the initial $structure$  for all joints of the marked BPS polyhedral decomposition and modify slab function if necessary  to obtain a new $structure$. After doing that  these Gross-Siebert walls and slabs may have new intersections which are contained in other stability walls. We then do the above algorithm again. In other words we build a scattering diagram for each joint  of the new structure. Then we do wall crossing calculations again for each joint.  After that we build Gross-Siebert walls labeled by BPS charges in the wall crossing calculations and modify slab functions accordingly if necessary.  We repeat this algorithm and go on like this. The collection of all Gross-Siebert walls and slabs whose attached log morphisms (or equivalently functions inducing these log morphisms) are determined by corresponding Kontsevich-Soibelman transformations of BPS charges is called the universal $structure$ of $GS_{BPS}$ at $(u,\vartheta)$. We use $\aleph_{uni}(u,\vartheta)$ to denote it.\end{dfn}

\textbf{Step II}: Truncations

Refined simple BPS spectra and refine BPS spectra are determined by finite trajectories one encounters between $\vartheta$ and $\vartheta+\pi$. We define refined (simple) BPS spectra for $(\vartheta^{j}_{-},\vartheta^{j}_{+})$ using the same definition except that we replace $\vartheta$ and $\vartheta+\pi$ by $\vartheta^{j}_{-}$ and $\vartheta^{j}_{+}$.
\begin{pro}In the non-degenerate case  $(\gamma_{i},\gamma_{l})$ is the refined simple BPS spectra for $(\vartheta^{j}_{-},\vartheta^{j}_{+})$ on both sides of the stability wall $SW_{\gamma_{i},\gamma_{l}}$.\end{pro}
\textbf{Proof}\ \ The statement is trivial on the side with only two charges $(\gamma_{i},\gamma_{l})$. For the other side we consider the rank two sublattice generated by $(\gamma_{i},\gamma_{l})$ with the induced intersection paring. Everything in the formalism of wall crossing formula explained in section 6 works for this sublattice and we do the wall crossing computation algebraically within this sublattice. In other words on the other side we truncate the product given by $K_{\gamma_{l}}K_{\gamma_{i}}$ according to the degree of positive linear combinations of $\gamma_{i},\gamma_{l}$. The factorization obtained in this way must coincide with the factorization of  $K_{\gamma_{l}}K_{\gamma_{i}}$ in the full lattice because of the uniqueness of factorization. Since this factorization is obtained in the sublattice by the above truncation procedure it is clear that charges it contains are all positive linear combinations of $\gamma_{i},\gamma_{l}$. So $(\gamma_{i},\gamma_{l})$ is also the refined simple BPS spectra for $(\vartheta^{j}_{-},\vartheta^{j}_{+})$ on the other side.\\

For degenerate case  we can use the (finite) refined simple BPS spectra at $(u,\vartheta)$ which positively generate the whole refined BPS spectra which includes the refined BPS spectra for $(\vartheta^{j}_{-},\vartheta^{j}_{+})$. So a subset of the refined simple BPS spectra at $(u,\vartheta)$ is the refined simple BPS spectra for $(\vartheta^{j}_{-},\vartheta^{j}_{+})$. Moreover the refined simple BPS spectra at $(u^{'},\vartheta)$  (which is the set of charges associated to all edges in the WKB triangulation $T_{WKB}^{\vartheta,u^{'}}$) is actually the  refined simple BPS spectra at $(u,\vartheta)$ because no BPS ray passes $\vartheta$ or $\vartheta+\pi$. The point is that like the non-degenerate case the refine simple BPS spectra for $(\vartheta^{j}_{-},\vartheta^{j}_{+})$ is determined on the side containing $u$.  As explained in section 6 we use the degree induced by the refined simple BPS spectra to truncate the wall crossing formula.\\

 Recall that the log morphism induced by the Kontsevich-Soibelman transformation $K_{\gamma}$ is given by $$x_{\gamma^{'}}\rightarrow x_{\gamma^{'}}(1-\sigma(\gamma)x_{\gamma})^{\langle\gamma^{'},\gamma\rangle}$$where $x_{\gamma}=z^{(\bar{m}_{\gamma_{gau}},h)}$. The deformation parameter $t$ is $t=z^{(0,1)}$.  $x_{\gamma}$ does not contain powers of $t$ explicitly (the order of $t$ is zero). However using  $x_{\gamma}$ to write down the slab function is not a canonical choice. We can change it to an expression in terms of other monoid variables using relations in the canonical thickening.  So a slab function or a wall function can contain  $t$ explicitly and in term of different sets of monoid variables we may have different $t$-orders for each term. We define the $t$-order of a slab function or a wall function to be the   $t$-order of its lowest order nonconstant terms in the formal Taylor expansion in the expression of monoid variables  which makes this lowest $t$-order largest. Then it makes sense to compute products of log morphisms modulo powers of $t$.   A single log morphism is viewed as a trivial product of log morphisms. The truncation of a product (composition) of log morphisms to the order $k$ is defined to be the truncation  modulo $t^{k+1}$ in the expression defining  the $t$-order. We can  truncate each log morphism to order $k$   before we compose them (and we may need to truncate again after that).

 The $t$-orders depend on the polarization and therefore the truncation by degrees and the truncation by powers of $t$ are not correspondent naturally. However the following proposition is enough for us.
 \begin{pro}Given a joint $j$. For any $k$ there exists a large enough $N$ such that the degree $N$ truncation contains all such BPS charges  that the  log morphisms associated to the rays and cuts labeled by them are nontrivial modulo $t^{k+1}$. Moreover   there are only finitely many of rays and cuts whose associated log morphisms are nontrivial modulo $t^{k+1}$.\end{pro}
 \textbf{Proof}\ \  First we show that the $t$-order of a nontrivial $f_{\rho_{\gamma}}$ is strictly positive. Without loss of generality suppose that in the slab labeling the cell $\rho_{\gamma}$ is labeled by a bar exponent whose exponent is $m_{\gamma_{gau}}=(\bar{m}_{\gamma_{gau}},h)$ then the order of $t$ of $f_{\rho_{\gamma}}=(1-\sigma(\gamma)x_{\gamma})^{l\Omega(\gamma;u)}$ in terms of $x_{\gamma}$ is clearly zero. Since the polarization is nontrivial by proposition 3.2 the order of $m_{\gamma_{gau}}$ strictly increases if we change the maximal dimensional cells form one that contains $\bar{m}_{\gamma_{gau}}$ to  one that contains $-\bar{m}_{\gamma_{gau}}$ . That means the $h$-component  used to define $m_{\gamma_{gau}}$ and $x_{\gamma}$ makes the left hand side of (100) positive for the maximal dimensional cell containing $-\bar{m}_{\gamma_{gau}}$. Therefore the $t$-order of $x_{\gamma}$ is positive in terms of $x_{-\gamma}$.

 On the other hand note that the refined simple BPS spectra at either side  of the joint (equivalently either side of the wall crossing formula) is finite. A charge appeared in the wall crossing formula is a positive linear combination of elements of the refined simple BPS spectra on the relevant side. Since  a slab function or a wall function associated to each such element has a positive $t$-order, if the degree (i.e. the sum of coefficients of the positive linear combination) is large enough then any ray labeled by a charge with degree larger  than that has associated wall function whose  $t$-order is higher than $k$. So the truncation by a large degree includes  all such BPS charges  that the  log morphisms associated to the rays and cuts labeled by them are nontrivial modulo $t^{k+1}$.

 Finally the finiteness statement in the proposition follows for  the finiteness of the refined simple BPS spectra at either side  of the joint (equivalently either side of the wall crossing formula).
 \begin{pro}The scattering diagram obtained by adding rays and modifying slab functions according to the wall crossing formula is a consistent scattering diagram to any order. \end{pro}
 \textbf{Proof}\ \ Because all log morphisms are Kontsevich-Soibelman transformations under the correspondence $x_{\gamma}\rightarrow\mathcal{X}_{\gamma}$ the wall crossing formula formulated as the fact the the ordered product along a loop is the identity in the truncated (by the degree) and projective sense becomes that condition that the ordered product of log morphisms along a loop containing the origin in the normal space of a joint is the identity in the truncated (by the degree) and projective sense. Since for any order $k$ rays and cuts whose attached log morphisms are nontrivial modulo $t^{k+1}$ are all included in a large enough degree truncation clearly the ordered product is the identity modulo $t^{k+1}$ to any order $k$.

\begin{note}The definition of a scattering diagram in section 3 the set of cuts and rays is required to be a finite set. We can generalize the definition by allowing a possibly infinite set of rays and cuts. We simply need to interpret  a possibly infinite product of log morphisms in the projective sense. Then it is clear that the wall crossing formula is a system of consistency conditions at $j$.\end{note}
 \begin{pro}Suppose $\mathbf{p}_{\gamma}$ is a Gross-Siebert wall obtained in the wall crossing calculation at certain joint of the universal structure. Assume the log morphism attached to it is trivial modulo $t^{k+1}$. Let $j$ be a joint obtained by intersecting $\mathbf{p}_{\gamma}$ and another Gross-Siebert wall or slab. If the joint $j$ is non-degenerate then all log morphisms associated to added rays associated to BPS charges in the wall crossing formula at $j$ are trivial modulo $t^{k+1}$. If the joint $j$ is degenerate we can consider added rays in the wall crossing without $\mathbf{p}_{\gamma}$ and denote that set by $\Pi(j-\mathbf{p}_{\gamma})$. Denote the set of added rays at $j$ with $\mathbf{p}_{\gamma}$ by $\Pi(j)$. Then log morphisms attached to elements in $\Pi(j)-\Pi(j-\mathbf{p}_{\gamma})$ are trivial modulo $t^{k+1}$.\end{pro}
 \textbf{Proof}\ \ For the non-degenerate $j$ we use the fact that $\gamma$ and the other charge labeling the other codimensional one cell of the $structure$ containing $j$ form the system of simple roots at both sides of the stability wall containing $j$. Since that means the BPS charge $\gamma^{'}$ of an added ray has $\gamma$ degree greater than or equal to that of $\gamma$, $x_{\gamma^{'}}$ has greater or same $t$-order and the proposition follows. The proof for the degenerate cases is analogous.

 \begin{dfn}Fix $k$. In the universal $structure$ $\aleph_{uni}(u,\vartheta)$, at each joint we truncate the wall crossing formula there by a large enough degree such that the order $k$ truncation  is included. We then check all rays and cuts in the associated scattering diagram in the degree truncation and delete those whose associated log morphisms are trivial modulo $t^{k+1}$.  For each joint we construct the Gross-Siebert walls from the rays that remain. We construct slabs from cuts that remain. The collection of all such Gross-Siebert walls and slabs is a $structure$. We call it the order $k$ truncated $structure$ and denote it by $\aleph_{k}(u,\vartheta)$.\end{dfn}
Proposition 9.9 guarantees that the  definition is a consistent one. If a ray is deleted then all rays from wall crossing calculations  induced by the Gross-Siebert wall associated to this ray should not exist any more. However in the definition this is the case if and only if all such rays' associated log morphisms are trivial modulo $t^{k+1}$. This is true by proposition 9.9.

\begin{thm}The collection of order $k$ truncated $structures$ is a compatible system of consistent $structures$. The wall crossing formula at a joint is equivalent to a collection of compatible consistency conditions of scattering diagrams.\end{thm}
\textbf{Proof}\ \ First we show that the order $k$ truncated $structure$ $\aleph_{k}(u,\vartheta)$ is consistent for each $k$. We have the wall crossing formula at each joint $j$ of $\aleph_{k}(u,\vartheta)$. We mod out both sides of the formula by $t^{k+1}$ and get an identity. Clearly by definition each side of the identity is also the product (modulo $t^{k+1}$) of log morphisms which are nontrivial modulo $t^{k+1}$ and therefore is the product of log morphisms attached to Gross-Siebert walls and slabs of $\aleph_{k}(u,\vartheta)$. This means that $\aleph_{k}(u,\vartheta)$ is consistent to order $k$.
Now we consider the compatibility of $\aleph_{k}(u,\vartheta)$ and $\aleph_{k+1}(u,\vartheta)$. The first condition follows directly from  the definition of $\aleph_{k}(u,\vartheta)$. As for the second condition note that in wall crossings we only add Gross-Siebert walls but not slabs. So the only  case to check is that we modify the slab functions. If a slab function is modified  to 1 for a slab after a wall crossing then the slab is deleted so that there is no intersection between the old slab in $\aleph_{k}(u,\vartheta)$ and the (nonexistent) new one in $\aleph_{k+1}(u,\vartheta)$ with different slab functions.  If the slab function is modified by multiplying  it by some ray (wall) functions then $t$-orders of these functions are  $k+1$ because they are in $\aleph_{k+1}(u,\vartheta)$ but not in $\aleph_{k}(u,\vartheta)$. Any attached function has  a constant term which is always 1.  So  the new slab function which is defined to be the product of the old slab function and ray (wall) functions will be the same modulo $t^{k+1}$.  The second condition follows. Hence we have a compatible system of consistent $structures$. The statement about the wall crossing formula is self-evident.
\begin{note}In section 3 we mentioned that to achieve the consistency of a $structure$ at a codimensional two joints we may have to normalize the log smooth structures which means the constant term of $f_{\rho_{\gamma}}$ is 1. Here our log smooth structure is automatically normalized.\end{note}

As a corollary of theorem 9.10, theorem 9.4 and proposition 3.3 we immediately get
\begin{thm}Let $\mathcal{M}(R)$  be a Hitchin's moduli space described in section 9.1. We assume the settings in the statement of theorem 9.4. Then there is a compatible system of consistent $structures$ associated to the BPS Gross-Siebert data $GS_{BPS}$ at $(u,\vartheta)$ and a  formal toric degeneration of Calabi-Yau varieties inducing the  BPS Gross-Siebert data can be constructed.\end{thm}

We have shown the equivalence of instanton data of the two instanton correction problems in a good sense (see section 9.4 for further discussions). We also want to see what this tells us about the relation of  the two instanton problems themselves. For that purpose we need to find the explicit degeneration.

Since we use unbounded cells each fiber in the toric degeneration is an affine variety. We expect that the Hitchin's moduli space as an affine variety is contained in this family.

First of all let us find a realization  the Hitchin's moduli space viewed as the moduli space of flat connections.  Recall that  the set of Fock-Goncharov coordinates labeled by charges associated to edges for a WKB triangulation provides a set of gauge invariant functions which is also a complete set of coordinates of $\mathcal{M}$. We want to use them as variables in the ideal defining $\mathcal{M}$ as an affine variety (but see remark 9.14). We need to enlarge the set of variables and then find relations between them.

Pick a pair $(u,\vartheta)$ and consider the associated marked BPS polyhedral decomposition. Elements of the refined BPS spectra $\Sigma(u,\vartheta)$ are ordered. The order determines a positive direction around each joint.

   Let $j$ be a joint of the marked BPS polyhedral decomposition for $(u,\vartheta)$. We use the slab labeling for slabs containing $j$.  Suppose  a BPS charge $\gamma_{i+1}$ labels a slab $\rho_{\gamma_{i+1}}$. If in $Side_{u}$ $\rho_{\gamma_{i+1}}$ is reached by which we mean that its BPS ray is reached by $\vartheta$ then we have the Kontsevich-Soibelman  transformations
$$K_{\gamma_{i+1}}^{\Omega(\gamma_{i+1};u)}: \mathcal{X}_{\gamma}\rightarrow
\mathcal{X}_{\gamma}(1-\sigma(\gamma_{i+1})
\mathcal{X}_{\gamma_{i+1}})^{\Omega(\gamma_{i+1};u)\langle\gamma,\gamma_{i+1}\rangle}$$ In particular let $\gamma_{i}$ be the  BPS charge of the last BPS ray encountered in the positive direction before  $l_{\gamma_{i+1}}$ and $\gamma_{i+2}$ be the first after $l_{\gamma_{i+1}}$ in $Side_{u}$. According to our convention there is no joint with multiple slabs. The slab function is \begin{equation}f_{\rho_{\gamma_{i+1}}}=(1-\sigma(\gamma_{i+1})
\mathcal{X}_{\gamma_{i+1}})^{l_{i+1}\Omega(\gamma_{i+1};u)}\end{equation}Define
$$a_{i}=\langle\gamma_{i},\gamma_{i+1}\rangle=-l_{i+1}\pi(\bar{m}_{\gamma_{i}^{gau}})$$
\begin{equation}b_{i+2}=-\langle\gamma_{i+2},\gamma_{i+1}\rangle=l_{i+1}\pi(\bar{m}_{\gamma_{i+2}^{gau}})\end{equation}
Note that $a_{i},b_{i+2}>0$ and $l_{i+1}$ is the integer defined by the above equation and is not $l_{\gamma_{i+1}}$. It is easy to see that $b_{i+2}\gamma_{i}+a_{i}\gamma_{i+2}$ is a charge such that its component in the direction orthogonal (in the inner product) to the BPS face labeled by $\gamma_{i+1}$ is zero. It is possible that $\gamma_{i+1}$ is the first one in $Side_{u}$ in which case we replace $\gamma_{i}$ by the charge which labels (in the slab labeling) the last one before $\gamma_{i+1}$ in $U_{j}$. Similarly if $\gamma_{i+1}$ is the last one in $Side_{u}$ then we replace $\gamma_{i+2}$ by the first one after $\gamma_{i+1}$ in $U_{j}$.

We have $$ \mathcal{X}_{\gamma_{i}}\rightarrow f_{\rho_{\gamma_{i+1}}}^{-\pi(\bar{m}_{\gamma_{i}^{gau}})}\mathcal{X}_{\gamma_{i}}, \ \mathcal{X}_{\gamma_{i+2}}\rightarrow f_{\rho_{\gamma_{i+1}}}^{-\pi(\bar{m}_{\gamma_{i+2}^{gau}})}\mathcal{X}_{\gamma_{i+2}}$$We define a two components form of $\mathcal{X}_{\gamma_{i}}$ by \begin{equation}\bar{\mathcal{X}}_{\gamma_{i}}:= (\mathcal{X}_{\gamma_{i}},f_{\rho_{\gamma_{i+1}}}^{-\pi(\bar{m}_{\gamma_{i}^{gau}})}\mathcal{X}_{\gamma_{i}})\end{equation} where the first component is the restriction of the Fock-Goncharov coordinate $\mathcal{X}_{\gamma_{i}}$ to the sector in the twistor $\mathbf{C}^{\times}$ bounded by BPS rays $l_{\gamma_{i}}$ and $l_{\gamma_{i+1}}$. Then the second component is the restriction to the sector bounded by $l_{\gamma_{i+1}}$ and $l_{\gamma_{i+2}}$. This is a local expression in the sector bounded  $l_{\gamma_{i}}$ and $l_{\gamma_{i+2}}$ of a global one. In fact for a discontinuity labeled by $\gamma_{s}$ we just define  more generally \begin{equation}\bar{\mathcal{X}}_{\gamma_{i}}:= (\mathcal{X}_{\gamma_{i}},f_{\rho_{\gamma_{s}}}^{-\pi^{s}(\bar{m}_{\gamma_{i}^{gau}})}\mathcal{X}_{\gamma_{i}})\end{equation} in the sector bound by $l_{\gamma_{s-1}}$ and $l_{\gamma_{s+1}}$ where $\gamma_{s-1}$ and $\gamma_{s+1}$ are the last slab charge before $\gamma_{s}$ and the first slab charge after $\gamma_{s}$ respectively. Here we add a superscript $s$ to restore the dependence of $-\pi(\bar{m}_{\gamma_{i}^{gau}})$ on $\gamma_{s}$ suppressed in (113). Clearly the definitions in different sectors patch together consistently and we get a global  definition  of $\bar{\mathcal{X}}_{\gamma_{i}}$ over sectors bounded by BPS rays labeled by charges we meet in $Side_{u}$ labeling slabs containing $j$. In other words we consider the restriction of $\mathcal{X}^{(u,\vartheta)}_{\gamma_{i}}$ with fixed $u$ to  the range of $\vartheta$ such that only the refine BPS charges containing $j$ in $Side_{u}$ are contained in this range. Then each component of (115) extends to  this globally defined $\mathcal{X}^{(u,\vartheta)}_{\gamma_{i}}$.

Now in the sector bounded by $l_{\gamma_{i}}$ and $l_{\gamma_{i+2}}$ \begin{equation}\bar{\mathcal{X}}_{\gamma_{i+2}}:= (f_{\rho_{\gamma_{i+1}}}^{\pi(\bar{m}_{\gamma_{i+2}^{gau}})}\mathcal{X}_{\gamma_{i+2}},\mathcal{X}_{\gamma_{i+2}})\end{equation} We use this representation to set $\mathcal{X}_{\gamma_{i+2}}$ be in the same twistor chamber of $f^{-\pi(\bar{m}_{\gamma_{i}^{gau}})}_{\rho_{\gamma_{i+1}}}\mathcal{X}_{i}$.

By (114) and (116) we have $$\bar{\mathcal{X}}_{\gamma_{i}}^{b_{i+2}}\bar{\mathcal{X}}_{\gamma_{i+2}}^{a_{i}}=((1-\sigma(\gamma_{i+1})\mathcal{X}_{\gamma_{i+1}})^{a_{i}b_{i+2}\Omega(\gamma_{i+1};u)},(1-\sigma(\gamma_{i+1})\mathcal{X}_{\gamma_{i+1}})^{a_{i}b_{i+2}\Omega(\gamma_{i+1};u)})\mathcal{X}_{s_{\gamma_{i+1}}}$$
    $\mathcal{X}_{s_{\gamma_{i+1}}}$ is defined by $$\mathcal{X}_{s_{\gamma_{i+1}}}:=\mathcal{X}_{\gamma_{i}}^{b_{i+2}}\mathcal{X}_{\gamma_{i+2}}^{a_{i}}$$
The symbol $s_{\gamma_{i+1}}$ indicates that the charge is contained in the BPS face $\rho_{\gamma_{i+1}}$. Since the discontinuous jump of $\mathcal{X}_{\gamma_{i+1}}$ across the ray $l_{\gamma_{i+1}}$ is trivial the two component form of $\mathcal{X}_{\gamma_{i+1}}$ for the sector bounded by $l_{\gamma_{i}}$ and $l_{\gamma_{i+2}}$ is $$\bar{\mathcal{X}}_{\gamma_{i+1}}=(\mathcal{X}_{\gamma_{i+1}},\mathcal{X}_{\gamma_{i+1}})$$The same formula holds for $\mathcal{X}_{s_{\gamma_{i+1}}}$. Therefore we can write the above equation in the following form $$\bar{\mathcal{X}}_{\gamma_{i}}^{b_{i+2}}\bar{\mathcal{X}}_{\gamma_{i+2}}^{a_{i}}=(1-\sigma(\gamma_{i+1})\bar{\mathcal{X}}_{\gamma_{i+1}})^{a_{i}b_{i+2}\Omega(\gamma_{i+1};u)}\bar{\mathcal{X}}_{s_{\gamma_{i+1}}}$$
where $1$  and $\sigma(\gamma_{i+1})$ are understood in the two component sense. Note that $\bar{\mathcal{X}}_{s_{\gamma_{i+1}}}$ is a power of $\bar{\mathcal{X}}_{\gamma_{i+1}}$.

Here we are really considering some limit values (in the twistor $\mathbf{C}^{\times}$) of Fock-Goncharov coordinates (see the footnote after (89)) which is natural because we only care about the moduli space $\mathcal{M}$ instead of its twistor space and therefore must eliminate the twistor parameter part of Fock-Goncharov coordinates. Like many other constructions the limits are determined by $(u,\vartheta)$ or the BPS chamber containing it.\\

Now we want to show that $a_{i}, b_{i+2}$ in the above equation are equal.
\begin{thm}Let $\bar{\mathcal{X}}_{\gamma_{i}}$, $\bar{\mathcal{X}}_{\gamma_{i+1}}$, $\bar{\mathcal{X}}_{\gamma_{i+2}}$ and $\bar{\mathcal{X}}_{s_{\gamma_{i+1}}}$ be symbols defined above, then $a_{i}=b_{i+2}$ and (note that now  $\mathcal{X}_{s_{\gamma_{i+1}}}:=\mathcal{X}_{\gamma_{i}}\mathcal{X}_{\gamma_{i+2}}$)
\begin{equation}\bar{\mathcal{X}}_{\gamma_{i}}\bar{\mathcal{X}}_{\gamma_{i+2}}=(1-\sigma(\gamma_{i+1})\bar{\mathcal{X}}_{\gamma_{i+1}})^{a_{i}\Omega(\gamma_{i+1};u)}\bar{\mathcal{X}}_{s_{\gamma_{i+1}}}\end{equation}
\end{thm}
\noindent {\bf Proof}\ \ We simply want to show that  the contractions to the normal direction of $\rho_{\gamma_{i+1}}$ of $\gamma_{i}$ and $\gamma_{i+2}$ are opposite. If the joint is non-degenerate then this is true because $\gamma_{i}$ and $\gamma_{i+2}$ are opposite charges.

Now let us assume the joint is degenerate.  Let us denote  the first charge in $Side_{u}$ by $\gamma_{1}$ and last one by $\gamma_{2}$. Clearly $\langle\gamma_{1},\gamma_{2}\rangle\neq0$. Consider the set of  BPS charges in $Side_{u}$ of the degenerate joint and denote this set by $L$. 
Clearly any charge in $L$ is a nonnegative  rational combination  of $\gamma_{1}$ and $\gamma_{2}$ (because we are in codimensional two). The combination is in fact integral because the BPS charges $\gamma_{1},\gamma_{2}$ are geometrically  represented by loops around simple zeroes and so are all the other BPS charges. When one varies $\vartheta$ there is no way to produce fractional loops. On the Riemann surface $\gamma_{1}$ is represented by a loop circling two simple zeros which we denote by $O_{1}$ and $O_{2}$. They are contained in two triangles denoted by $T_{1}$  and $T_{2}$. $\gamma_{2}$ is represented by a loop circling two simple zeros which we denote by $O_{2}$ and $O_{3}$ (there must be a common simple zero otherwise the intersection would be zero). They are contained in two triangles denoted by $T_{2}$  and $T_{3}$. We will exhaust all possible scenarios of the geometric relation between $\gamma_{1}$ and $\gamma_{2}$.
\begin{itemize}
\item $O_{1}\neq O_{3}$. \ \  The loop circling the two simple zeroes $O_{3}$ and $O_{1}$  (with the proper orientation) represents the charge $\gamma_{1}+\gamma_{2}$ in $Side_{u}$.  Therefore we know that $L$ contains $\gamma_{1},\gamma_{2}$ and possibly also $\gamma_{1}+\gamma_{2}$.   Geometric BPS charges correspond to finite trajectories connecting simple zeroes. There is no finite trajectory connecting more than two points in $(O_{1},O_{2},O_{3})$ at any given critical value of $\vartheta$. In fact such a hypothetical finite trajectory would have different phases (which violates the constant phase condition) in different segments separated by the simple zeroes other than the starting and the ending point. Therefore there is no way to geometrically realize $m\gamma_{1}+n\gamma_{2},\ m>1 or\ n>1$\footnote{One may wonder the possibility of the realization of some multiples of just one finite trajectory such as $m\gamma_{1}$. But this is still impossible as a finite trajectory connecting two simple zeroes only crosses the corresponding edge once.}.  This tells us that for a degenerate joint $L$ only contains $\gamma_{1},\gamma_{2}$ and $\gamma_{1}+\gamma_{2}$.
\item $O_{1}=O_{3}$ and one of  $T_{1}$ and $T_{2}$ is non-degenerate.\ \ In this case the two triangles must have two commons edges corresponding to $\gamma_{1}, \gamma_{2}$. Since one triangle is non-degenerate  these two common edges have three vertices which means they can not degenerate to finite trajectories at the same time. So  we do not have nontrivial linear combinations of $\gamma_{1}, \gamma_{2}$.
\item $O_{1}=O_{3}$, both $T_{1}$ and $T_{2}$ are degenerate.\ \ If we want the two edges corresponding to the two charges to be able to degenerate to finite trajectories at the same time they must have two common vertices. The two edges then are the two double edges in the two degenerate triangles. This brings us to the scenario of infinite flips described in  section 7 with the two loop edges as the two circles bounding the annular region. There are two possibilities. If the winding of the sum of the two edges around the inner circle is zero we get two flips corresponding to $\gamma_{1},\gamma_{2}$ which gives us a non-degenerate joint. Otherwise we have the process of infinite flips. By the assumption of finite type this scenario has been excluded.
\end{itemize}
By the above classification the possibilities of three consecutive charges in $Side_{u}$ of a joint can only be $(\gamma_{i},\gamma_{i+1},\gamma_{i+2}=-\gamma_{i})$, $(\gamma_{i},\gamma_{i+1}=\gamma_{i}+\gamma_{i+2},\gamma_{i+2})$, $(\gamma_{i},\gamma_{i+1},\gamma_{i+2}=\gamma_{i+1}+(-\gamma_{i}))$ and $(\gamma_{i}=-\gamma_{i+2}+\gamma_{i+1},\gamma_{i+1},\gamma_{i+2})$. The theorem follows immediately.\\

We can do the same thing in $Side_{u}^{\bot}$ with the positive direction around the joint replaced by the negative direction around the joint. We use the wall crossing formula to add Gross-Siebert walls (equivalently BPS rays) and modify slab functions (equivalently associated Kontsevich-Soibelman transformations). Of course in $Side_{u}^{\bot}$ we could have infinitely many Gross-Siebert walls corresponding to infinitely many sectors in the twistor $\mathbf{C}^{\times}$. But we only consider $\bar{\mathcal{X}}_{\gamma}$ labeled by the refined BPS spectra at $(u,\vartheta)$. So there are only finitely many relations. We need to define them by incorporating the effects of all Kontsevich-Soibelman transformations associated to added BPS rays (or equivalently Gross-Siebert walls). If the slab function of a slab is modified to 1 in $Side_{u}^{\bot}$ we delete that slab and the corresponding BPS ray. There are Gross-Siebert walls between two slabs in $Side_{u}^{\bot}$ labeled (in the face labeling) by two adjacent charges in the refined BPS spectra at $(u,\vartheta)$. \begin{itemize} \item If in $Side_{u}^{\bot}$ $\rho_{\gamma_{i+1}}$ is a slab such that the associated Kontsevich-Soibelman transformation of the corresponding BPS rays $l_{\gamma_{i+1}}$ is a transformation induced by a flip, then in (114) the second component is meant to be the restriction of the global $\mathcal{X}_{\gamma_{i}}$ to the  sector bounded by $l_{\gamma_{i+1}}$ and the first BPS ray one encounters when rotating  $l_{\gamma_{i+1}}$ in the positive direction of  the twistor $\mathbf{C}^{\times}$ induced by the order of the BPS spectra in this new stability chamber (note that this corresponds to the negative direction around $j$). This means that the first component is the restriction to the sector bounded by $l_{\gamma_{i+1}}$ and the first BPS ray one encounters when rotating $l_{\gamma_{i+1}}$ to $l_{\gamma_{i+2}}$ in the negative direction of the twistor $\mathbf{C}^{\times}$. Similarly in (116) we match the two components with the two in (114).\item Suppose in $Side_{u}^{\bot}$ $\rho_{\gamma_{i+1}}$ is a slab such that the associated Kontsevich-Soibelman transformation of the corresponding BPS rays $l_{\gamma_{i+1}}$ is a transformation induced by a juggle then it is surrounded by infinitely many Gross-Siebert walls. Equivalently the corresponding BPS ray $l_{\gamma_{i+1}}$ is surrounded by infinitely many BPS rays.  The second component of (114) is the limit Fock-Goncharov coordinate labeled by $\gamma_{i}$ from the negative direction of the twistor $\mathbf{C}^{\times}$. Then the first component is  the limit Fock-Goncharov coordinate labeled by $\gamma_{i}$ from the positive direction of the twistor $\mathbf{C}^{\times}$. \item Suppose in $Side_{u}^{\bot}$ $\rho_{\gamma_{i+1}}$ is a slab coincide with several Gross-Siebert walls such that the associated Kontsevich-Soibelman transformation of the corresponding BPS rays $l_{\gamma_{i+1}}$ is a transformation induced by a composition of a juggle and several flips\footnote{This composition does not depend on the order of composing the transformations.}, we change $f_{\rho_{\gamma_{i+1}}}$ in (112) to the product of $f_{\rho_{\gamma_{i+1}}}$ with the wall functions  attached to the coincide Gross-Siebert walls.  After that we follow the same procedure in the previous cases.\end{itemize} Then by the definition and the wall crossing formula the same equation (117) holds for slabs in $Side_{u}^{\bot}$. In this way we have a relation (117) for each slab containing $j$. Note that although two slabs with opposite normal vectors are labeled by the same charge in the face labeling they are labeled in the notation of (117) differently. Of course the labeling used in (117) can be identified with the slab labeling by charges if we indicate which slab in $Side_{u}^{\bot}$ is the opposite one to a given one in $Side_{u}$.

 We call the relation defined by  (117) a Fock-Goncharov relation. We have defined them for a joint with a reference vertex.  If we change the reference vertex  charges are related by parallel transport of charges and $f_{\rho_{\gamma},v}$ transforms by (35). The log morphism attached to a slab and hence the corresponding Kontsevich-Soibelman transformation associated to the labeling charge in the face labeling do not depend on the choice of the reference vertex.

 We have defined Fock-Goncharov relations according to the structure of the marked BPS polyhedral decomposition but in fact we can equivalently define them according to data in the metric problem. To realize the Hitchin's moduli space  via Fock-Goncharov relations, it is natural from the perspective of the metric problem to collect   relations associated to all discontinuous jumps across all BPS rays  that appear.  Since we want to compare it with a generic fiber of the degeneration constructed by Gross and Siebert's method it would be convenient to equivalently formulate this collection in terms of data in the complex structure problem. From the perspective of complex structure problem this means we use the following definition.
\begin{dfn}Let $i+1$ runs over all slabs in the marked BPS polyhedral decomposition at $(u,\vartheta)$. We call the ideal generated by all such Fock-Goncharov relations the Fock-Goncharov ideal $I_{FG}(u,\vartheta)$ and the defining equations are called the Fock-Goncharov realization at $(u,\vartheta)$ of the Hitchin's moduli space. \end{dfn}

This construction gives us the correction dimension. It follows form the fact that  the number of independent Fock-Goncharov coordinates for a fixed $(u,\vartheta)$ is equal to the complex dimension of $\mathcal{M}$.

\begin{note}In GMN's ansatz Fock-Goncharov coordinates are valued in $\mathbf{C}^{\times}$. So a Fock-Goncharov realization really gives a variety in a product of copies of $\mathbf{C}^{\times}$. However one can also extend Fock-Goncharov coordinates to be valued in $\mathbf{C}$ (see remark 7.1) and therefore it is also natural to view the Fock-Goncharov realization as an affine variety in a product of copies of $\mathbf{C}$. What is the geometric meaning of this affine variety? Since the Hitchin's moduli space as the moduli space of flat connections is analytically isomorphic to the moduli space of fundamental group representations which is an  affine variety and Fock-Goncharov coordinates are gauge invariant holomorphic coordinates on the moduli space of flat connections the following conjecture is likely to be true.
\begin{con} The Hitchin's moduli space as an affine variety in a Fock-Goncharov realization is algebraically isomorphic to the underlying affine variety of the moduli space of fundamental group representations. \end{con}
But remember that the moduli space of flat connections is not algebraically isomorphic to the moduli space of fundamental group representations. In particular we should not consider the moduli space of flat connections as an affine variety. 
It is shown that traces of monodromy matrices are generated by Fock-Goncharov coordinates (see the appendix A of \citep{G2}). Since we can use traces and their relations to write down the underlying variety of the moduli space of fundamental group representations it would be interesting to  see what the relation between this realization and a Fock-Goncharov realization is.\end{note}

 Now let us consider the degeneration given by Gross and Siebert's construction from a compatible system of consistent $structures$. We will use the following proposition  proved in \citep{GS1} by Gross and Siebert.

  Suppose $\rho=\sigma_{+}\cap\sigma_{-}$ is a codimensional one cell with two adjacent maximal dimensional cells $\sigma_{+}$ and $\sigma_{-}$. $v\in\omega\subseteq\rho$. Denote $P_{x}$ by $P$ for $x\in Int(\omega)-\Delta$. Let $R_{\sigma_{-}}:=R_{\omega\rightarrow\sigma_{-},\sigma_{-}}^{k},R_{\sigma_{+}}:=R_{\omega\rightarrow\sigma_{+},\sigma_{+}}^{k},R_{\rho}:=R_{\omega\rightarrow\rho,\sigma_{+}}^{k}$ denote  the coordinate rings  of $k$-th order canonical thickenings associated to  $\sigma_{-}$, $\sigma_{+}$ and $\rho$ respectively and $f_{\rho}$ is the slab function attached to $\rho$. Here we have assumed that the reference maximal dimensional cell for $\omega\rightarrow\rho$ is $\sigma_{+}$ so that the gluing of $R_{\sigma_{-}}$ and $R_{\sigma_{+}}$ along $R_{\rho}$ is given by the fiber product with respect to the canonical quotient homomorphism $R_{\sigma_{+}}\rightarrow R_{\rho}$ and the homomorphism $R_{\sigma_{-}}\rightarrow R_{\rho}$ which is the composition of the canonical homomorphism and the log morphism  induced by the change of chamber from $\sigma_{-}$ to $\sigma_{+}$. Let $F_{\rho}$ be the set of  monoid variable whose projections are contained in $\Lambda_{\rho}$. By choosing appropriate coordinates and local representative of the polarization the elements of the set of monoid variable  $P+F_{\rho}$ can be put into a standard form $\Lambda_{\rho}+S_{e},e\geq1$ where $S_{e}$ is rank two (because $\Lambda_{\rho}$ is rank $n-1$) and is generated by $(-a,e),(a,0)$ and $(0,1)$.  Here the second components are the $h$-components in $(\bar{m},h)$. The first components are $primitive$ generators in $\Lambda_{\rho}^{\bot}$ identified with $\mathbf{Z}$ (one can set $a=1$ if one likes) and we have suppressed the components in $\Lambda_{\rho}$. $e$ is the increase of the piecewise integral linear function of the local representative of the polarization along $\Lambda_{\rho}^{\bot}$ by  $a$ units. Let $x:=z^{(-a,e)}, y=z^{(a,0)}, t=z^{(0,1)}$. Let $R_{-},R_{+},R_{\cap}$ be the localizations of $R_{\sigma_{-}},R_{\sigma_{+}},R_{\rho}$ at $\{z^{m}\}_{m\in F_{\rho}}$ respectively. Then explicitly  $$R_{-}= \mathbf{C}[\Lambda_{\rho}][x,y,t]/(xy-t^{e},y^{\beta}t^{\gamma}\mid \beta e+\gamma\geq k+1)$$  $$R_{+}= \mathbf{C}[\Lambda_{\rho}][x,y,t]/(xy-t^{e},x^{\alpha}t^{\gamma}\mid \alpha e+\gamma\geq k+1)$$
\begin{equation}R_{\cap}= \mathbf{C}[\Lambda_{\rho}][x,y,t]/(xy-t^{e},x^{\alpha}y^{\beta}t^{\gamma}\mid max\{\alpha,\beta\}e+\gamma\geq k+1)_{f_{\rho}}\end{equation}The log morphism associated to $\rho$ is\begin{equation}x\rightarrow f_{\rho}^{a}x,\ y\rightarrow f_{\rho}^{-a}y\end{equation}

Let $R_{+}\rightarrow R_{\cap}$ be the canonical quotient homomorphism followed by the localization at $f_{\rho}$ and $R_{-}\rightarrow R_{\cap}$ be the quotient morphism (followed by the localization) composed with the log morphism. The gluing is the fiber product of these two homomorphisms. \begin{thm}(Lemma 2.34 of \citep{GS1}) The fiber product $R_{-}\times_{R_{\cap}}R_{+}$ denoted by $R_{\cup}$  is \begin{equation}R_{\cup}:=\mathbf{C}[\Lambda_{\rho}][X,Y,t]/(XY-f_{\rho}^{a}t^{e},t^{k+1})\end{equation} The map $R_{\cup}\rightarrow R_{-}\times_{R_{\cap}}R_{+}$ which is an isomorphism is given by \begin{equation}X\rightarrow(x,f_{\rho}^{a}x),\ Y\rightarrow(f_{\rho}^{a}y,y)\end{equation}\end{thm}
\textbf{Proof}\ \ As $\mathbf{C}[\Lambda_{\rho}][t]/(t^{k+1})$ modules $R_{\pm}$ and $R_{\cap}$ are generated by $1, x^{i},y^{l}$. Elements $g_{\pm}\in R_{\pm}$ can be written as $g_{-}=\sum_{i\geq0}a_{i}x^{i}+h_{-}(y,t)$ and $g_{+}=\sum_{l\geq0}b_{l}y^{l}+h_{+}(x,t)$ with $h_{\pm}(0,t)=0$. $(g_{-},g_{+})\in R_{-}\times_{R_{\cap}}R_{+}$ if and only if $$a_{0}=b_{0},\  h_{-}(y,t)=\sum_{l>0}b_{l}f_{\rho}^{al}y^{l}, \ h_{+}(x,t)=\sum_{i>0}a_{i}f_{\rho}^{ai}x^{i}$$ in $R_{\cap}$. Then $(g_{-},g_{+})$ is clearly the image of $\sum_{i\geq0}a_{i}X^{i}+\sum_{l\geq0}b_{l}Y^{l}\in R_{\cup}$. This shows the surjectivity. The injectivity follows from the fact that $R_{\cup}$ is a free $\mathbf{C}[\Lambda_{\rho}][t]/(t^{k+1})$ module generated by $1,X^{i},Y^{l},i,l>0$.\\

Let us come back to the marked BPS polyhedral decomposition. Let $j$ be a joint.  Denote the slab labeled by $\gamma_{i}$ in the slab labeling by $slab_{\gamma_{i}}$. Let $\mathcal{C}_{i}$ and $\mathcal{C}_{i+1}$ denote  the Gross-Siebert chambers whose boundaries containing ($slab_{\gamma_{i}}$,$slab_{\gamma_{i+1}}$) and ($slab_{\gamma_{i+1}}$,$slab_{\gamma_{i+2}}$) respectively.  $\mathcal{C}_{i}$ and $\mathcal{C}_{i+1}$ are maximal dimensional cells in the BPS polyhedral decomposition playing the role of $\sigma_{-}$ and $\sigma_{+}$ with $\mathcal{C}_{i}\cap\mathcal{C}_{i+1}=\rho:=\rho_{\gamma_{i+1}}$.

We calculate $\pi(\bar{m}_{\gamma_{i}}),\pi(\bar{m}_{\gamma_{i+2}})$ for the epimorphism $\pi$ from the gauge lattice to $\mathbf{Z}$ with kernel $\Lambda_{\rho_{\gamma_{i+1}}}$. Let the $h$-components according to the polarization be $(h_{i},h_{i+2})$.  Finally set $x:=x_{\gamma_{i}}^{\bot}=z^{(\pi(\bar{m}_{\gamma_{i}}),h_{i})},y:=x_{\gamma_{i+2}}^{\bot}=z^{(\pi(\bar{m}_{\gamma_{i+2}}),h_{i+2})}$. In this notation we have suppressed the components along $\Lambda_{\rho_{\gamma_{i+1}}}$. Then according to theorem 9.12 $$x_{\gamma_{i}}^{\bot}x_{\gamma_{i+2}}^{\bot}=t^{e_{i+1}}$$$$x_{\gamma_{i}}x_{\gamma_{i+2}}=x_{s_{\rho_{\gamma_{i+1}}}}t^{e_{i+1}}$$ where $x_{s_{\rho_{\gamma_{i+1}}}}\in\Lambda_{\rho_{\gamma_{i+1}}}$ and $e_{i+1}$ is a positive integer determined by $h_{i},h_{i+2}$. In theorem 9.13 $x,y$ are supposed  to be associated to   primitive normal generators of $\rho$ (which generate lattice points in $\Lambda_{\rho}^{\bot}$) so that they are generators of the coordinated rings. Here our generators are also primitive because $\gamma_{i},\gamma_{i+2}$ are primitive charges meaning that $\gamma_{i}$ ($\gamma_{i+2}$) is  not a multiple of another  charge. This follows from their geometric construction as geometric BPS charges because we do not allow multiple loops around a given finite trajectory connecting two simple zeroes.

Since the slab function attached to $\rho_{\gamma_{i+1}}$ is given by $f_{\rho_{\gamma_{i+1}}}=(1-\sigma(\gamma_{i+1})x_{\gamma_{i+1}})^{l_{i+1}\Omega(\gamma_{i+1};u)}$, using theorem 9.13 and then putting back the components along $\Lambda_{\rho_{\gamma_{i+1}}}$ (which are not affected by the log morphisms) we know that the gluing of canonical thickenings of the affine strata labeled by maximal dimensional cells $\mathcal{C}_{i}$ and $\mathcal{C}_{i+1}$ along  the canonical thickening of the strata labeled by $\rho_{\gamma_{i+1}}$ gives us the equation \begin{equation}X_{\gamma_{i}}X_{\gamma_{i+2}}=(1-\sigma(\gamma_{i+1})x_{\gamma_{i+1}})^{a_{i}\Omega(\gamma_{i+1};u)}x_{s_{\rho_{\gamma_{i+1}}}}t^{e_{i+1}}\end{equation} where \begin{equation}X_{\gamma_{i}}=(x_{\gamma_{i}},f_{\rho_{\gamma_{i+1}}}^{a_{i}}x_{\gamma_{i}}),\ X_{\gamma_{i+2}}=(f^{a_{i}}_{\rho_{\gamma_{i+1}}}x_{\gamma_{i+2}},x_{\gamma_{i+2}})\end{equation}

 $X_{\gamma_{i}}$ is the lift of $x_{\gamma_{i}}$ in the fiber product obtained by the  gluing of $\mathcal{C}_{i}$ and $\mathcal{C}_{i+1}$. Since we know that the gluing of all Gross-Siebert chambers in a consistent $structure$ is consistent, there is a well defined lift not only on this glued piece but also on the whole space and we still denote it by $X_{\gamma_{i}}$. In fact  a gluing is either a change of strata where no log morphism is introduced or a change of chambers where a log morphism is composed between canonical thickenings. A log morphism in our case is attached to  a slab $\rho_{\gamma_{l}}$ and is given by $$x_{\gamma_{i}}\rightarrow x_{\gamma_{i}}f_{\rho_{\gamma_{l}}}^{-\pi(\bar{m}^{gau}_{\gamma_{i}})}$$ Gluing of chambers produces (122) in the relevant strata. So (123) is really the local expression of $X_{\gamma_{i}}$.

Also note that $(1-\sigma(\gamma_{i+1})x_{\gamma_{i+1}})^{a_{i}\Omega(\gamma_{i+1};u)}t^{e_{i+1}}$ in (122) is really $$((1-\sigma(\gamma_{i+1})x_{\gamma_{i+1}})^{a_{i}\Omega(\gamma_{i+1};u)},(1-\sigma(\gamma_{i+1})x_{\gamma_{i+1}})^{a_{i}\Omega(\gamma_{i+1};u)})t^{e_{i+1}}$$ On the other hand $$((1-\sigma(\gamma_{i+1})x_{\gamma_{i+1}})^{a_{i}\Omega(\gamma_{i+1};u)},(1-\sigma(\gamma_{i+1})x_{\gamma_{i+1}})^{\Omega(a_{i}\gamma_{i+1};u)})$$ is the lift of $(1-\sigma(\gamma_{i+1})x_{\gamma_{i+1}})^{a_{i}\Omega(\gamma_{i+1};u)}$. Hence $(1-\sigma(\gamma_{i+1})x_{\gamma_{i+1}})^{a_{i}\Omega(\gamma_{i+1};u)}t^{e_{i+1}}$ in (122) is really \begin{equation}(1-\sigma(\gamma_{i+1})X_{\gamma_{i+1}})^{a_{i}\Omega(\gamma_{i+1};u)}t^{e_{i+1}}\end{equation}This is analogous to the meaning of (117).

We have the following proposition.\begin{pro}The gluing of all $k$-th order canonical thickenings associated to  all strata labeled by cells in the marked BPS polyhedral decomposition for $(u,\vartheta)$  is given by  the ideal generated by \begin{equation} (X_{\gamma_{i}}X_{\gamma_{i+2}}-(1-\sigma(\gamma_{i+1})
X_{\gamma_{i+1}})^{a_{i}\Omega(\gamma_{i+1};u)}x_{s_{\rho_{\gamma_{i+1}}}}t^{e_{i+1}},t^{k+1})\end{equation} where $X_{\gamma_{i}}$ is the lift of monoid variable $x_{\gamma_{i}}$ to the total space whose existence is guaranteed by theorem 9.11. $i+1$ runs over all slabs (suppose the cardinality of the set of slabs is $s$) in the marked BPS polyhedral decomposition.\end{pro}
\textbf{Proof}\ \ According to section 3 the gluing consists of two types: changes of strata and changes of chambers. Changes of strata do not produce new relations or change relations of monoid variables in  canonical thickenings   because they are just the embedding relations between canonical thickenings of affine strata i.e. all morphisms are canonical quotient homomorphisms. Changes of (Gross-Siebert) chambers produce fiber products with log morphisms composed.  Then the proposition follows from (122) and (123).
\begin{note}We use $x_{s_{\rho_{\gamma_{i+1}}}}$ instead of $X_{s_{\rho_{\gamma_{i+1}}}}$ in (125). Unlike the situation in (117) $x_{s_{\rho_{\gamma_{i+1}}}}$ cannot be immediately considered as the lift (two component form) of a monoid variable  in the gluing. That is simply because the $h$-component of it is zero. Therefore one only need to absorb some powers of $t$ to get  $h$ defined in (100). In other words we can replace $x_{s_{\rho_{\gamma_{i+1}}}}$ by $X_{s_{\rho_{\gamma_{i+1}}}}$  if we are willing to lose some power of $t$. Just like the situation of (117) $X_{s_{\rho_{\gamma_{i+1}}}}$ is a power of $X_{\gamma_{i+1}}$ since the charge $\gamma_{i+1}$ is primitive. So it is not a  new generator.\end{note}
\begin{note}Not all relations in the ideal $I$ are necessary. We will see the redundancy of some of them demanded by the wall crossing formula. \end{note}
\begin{thm}Let $\mathcal{M}(\mathcal{M}(R))$ be a Hitchin's moduli space described in section 9.1. Choose a pair $(u,\vartheta)$. Let $GS_{BPS}$ be a BPS Gross-Siebert data at $(u,\vartheta)$. Let $k$ be large enough. Then the Hitchin's moduli space in the Fock-Goncharov realization at $(u,\vartheta)$ is a generic fiber of the degeneration over $\mathrm{Spec}\mathbf{C}[t]/(t^{k+1})$ from Gross-Siebert's construction. \end{thm}
\textbf{Proof}\ \ Let $k$ be large enough one can assume that $k+1$ is larger than all exponents $e_{i+1}$ in (125).
Set $t=1$ we get  the defining equation of a generic fiber of the degeneration. On the other hand by the monoid-charge correspondence we make the following change of variable\begin{equation}\mathcal{X}_{\gamma}\rightarrow x_{\gamma}, \ \bar{\mathcal{X}_{\gamma}}\rightarrow X_{\gamma}\end{equation}Then the Fock-Goncharov ideal of the Hitchin's moduli space is mapped to the ideal of that generic fiber and vice versa.\\

In the above theorem we glue canonical thickenings of all Gross-Siebert chambers. The consistency around a joint is implicit. Now we will make it explicit for non-degenerate joints.

Let $j$ be a non-degenerate joint contained in the unique stability wall $SW_{j}$. Suppose $SW_{j}=SW_{\gamma_{1},\gamma_{2}}$. So  there are three Gross-Siebert chambers having nonempty intersections with $Side_{u}$. They are bounded by $(-\gamma_{2},\gamma_{1})$, $(\gamma_{1},\gamma_{2})$ and $(\gamma_{2},-\gamma_{1})$ respectively. Let us denote them by $\mathcal{C}_{1}$, $\mathcal{C}_{2}$,$\mathcal{C}_{3}$ respectively. There is one more Gross-Siebert chamber denoted by $\mathcal{C}_{4}$.  It is contained in $Side_{u}^{\bot}$. In $Side_{u}^{\bot}$ we glue $\mathcal{C}_{1}$ and $\mathcal{C}_{4}$ as well as $\mathcal{C}_{4}$ and $\mathcal{C}_{3}$. In $Side_{u}$ we glue $\mathcal{C}_{1}$ and $\mathcal{C}_{2}$ as well as $\mathcal{C}_{2}$ and $\mathcal{C}_{3}$. The gluing of $\mathcal{C}_{1}$ and $\mathcal{C}_{2}$ gives us the ideal \begin{equation}(X_{-\gamma_{2}}X_{\gamma_{2}}-(1-\sigma(\gamma_{1})
X_{\gamma_{1}})^{a_{-\gamma_{2}}\Omega(\gamma_{1};u)}x_{s_{\rho_{\gamma_{1}}}}t^{e_{\gamma_{1}}},t^{k+1})\end{equation}

 We will show that gluing $\mathcal{C}_{4}$ and $\mathcal{C}_{3}$ produces the same relations.  The relevant monoids are still $x_{-\gamma_{2}}$ and $x_{\gamma_{2}}$. But they are in different Gross-Siebert chambers now and therefore must be obtained from  $x_{-\gamma_{2}}$ and $x_{\gamma_{2}}$ in the first gluing by using log morphisms.

We have changes of strata which are surjective quotient homomorphisms$$ R^{k}_{j\rightarrow\mathcal{C}_{1},\mathcal{C}_{1}}\rightarrow R^{k}_{j\rightarrow\rho_{\gamma_{1}},\mathcal{C}_{1}}\rightarrow R^{k}_{j\rightarrow j,\mathcal{C}_{1}}$$\begin{equation} R^{k}_{j\rightarrow\mathcal{C}_{4},\mathcal{C}_{4}}\rightarrow R^{k}_{j\rightarrow\rho_{-\gamma_{1}},\mathcal{C}_{4}}\rightarrow R^{k}_{j\rightarrow j,\mathcal{C}_{4}}\end{equation}
$$R^{k}_{j\rightarrow\rho_{\gamma_{1}},\mathcal{C}_{1}}\rightarrow R^{k}_{\rho_{\gamma_{1}}\rightarrow\rho_{\gamma_{1}},\mathcal{C}_{1}}$$
\begin{equation}R^{k}_{j\rightarrow\rho_{-\gamma_{1}},\mathcal{C}_{4}}\rightarrow R^{k}_{\rho_{-\gamma_{1}}\rightarrow\rho_{-\gamma_{1}},\mathcal{C}_{4}}\end{equation}
We also have $$ R^{k}_{\rho_{\gamma_{1}}\rightarrow\mathcal{C}_{1},\mathcal{C}_{1}}\rightarrow R^{k}_{\rho_{\gamma_{1}}\rightarrow\rho_{\gamma_{1}},\mathcal{C}_{1}}$$
\begin{equation}R^{k}_{\rho_{-\gamma_{1}}\rightarrow\mathcal{C}_{4},\mathcal{C}_{4}}
\rightarrow R^{k}_{\rho_{-\gamma_{1}}\rightarrow\rho_{-\gamma_{1}},\mathcal{C}_{4}}\end{equation}

 The gluing of $\mathcal{C}_{1}$ and $\mathcal{C}_{2}$ is the fiber product of $R^{k}_{\rho_{\gamma_{1}}\rightarrow\mathcal{C}_{1},\mathcal{C}_{1}}\rightarrow R^{k}_{\rho_{\gamma_{1}}\rightarrow\rho_{\gamma_{1}},\mathcal{C}_{1}}$ and $R^{k}_{\rho_{\gamma_{1}}\rightarrow\mathcal{C}_{2},\mathcal{C}_{2}}\rightarrow R^{k}_{\rho_{\gamma_{1}}\rightarrow\rho_{\gamma_{1}},\mathcal{C}_{2}}$ with the log morphism \begin{equation}\theta_{\rho_{\gamma_{1}}}: R^{k}_{\rho_{\gamma_{1}}\rightarrow\rho_{\gamma_{1}},\mathcal{C}_{1}}\rightarrow R^{k}_{\rho_{\gamma_{1}}\rightarrow\rho_{\gamma_{1}},\mathcal{C}_{2}}\end{equation}
$\theta_{\rho_{\gamma_{1}}}$ is actually induced from the log morphism \begin{equation}\theta_{\rho_{\gamma_{1}}}: R^{k}_{j\rightarrow\rho_{\gamma_{1}},\mathcal{C}_{1}}\rightarrow R^{k}_{j\rightarrow\rho_{\gamma_{1}},\mathcal{C}_{2}}\end{equation} We move $x_{\gamma_{2}}$ and $x_{-\gamma_{2}}$ clockwise from the $slab_{\gamma_{1}}$ to $slab_{-\gamma_{1}}$ crossing all Gross-Siebert walls and slabs between them and we obtain the monoid variables for the gluing of $\mathcal{C}_{4}$ and $\mathcal{C}_{3}$ along $slab_{-\gamma_{1}}$ with $f_{\rho_{\gamma_{1}}}$ attached to it. We still use $x_{\gamma_{2}}$ and $x_{-\gamma_{2}}$ to denote  these monoid variables as they are still components of the local expressions of the global elements $X_{\gamma_{2}}$ and $X_{-\gamma_{2}}$. Similarly we move $x_{\gamma_{2}}$ and $x_{-\gamma_{2}}$ counterclockwise from the $slab_{\gamma_{1}}$ to $slab_{-\gamma_{1}}$ crossing all Gross-Siebert walls and slabs between them and we obtain the monoid variables for the gluing of $\mathcal{C}_{4}$ and $\mathcal{C}_{3}$ along $slab_{-\gamma_{1}}$.
The log morphism between them across $slab_{\gamma_{1}}$ counterclockwise  is \begin{equation}x_{-\gamma_{2}}\rightarrow f_{\rho_{\gamma_{1}}}^{a}x_{-\gamma_{2}},\ x_{\gamma_{2}}\rightarrow f_{\rho_{\gamma_{1}}}^{-b}x_{\gamma_{2}}\end{equation}We have suppressed the indices of $a,b$. Because the ordered composition along a loop is an identity and all Gross Siebert walls and slabs except $slab_{\gamma_{1}}$ and $slab_{-\gamma_{1}}$ have been crossed it is clear that the log morphisms between $x_{-\gamma_{2}}$ and $x_{\gamma_{2}}$ across clockwise $slab_{-\gamma_{1}}$ must be also given by (133). Also note that $$x_{s_{\rho_{\gamma_{1}}}}=x_{s_{\rho_{\gamma_{1}}}}, \ e_{\gamma_{1}}=e_{-\gamma_{1}}$$The second relation is obvious while the first one really means that we take the composition of actions by log morphisms on $x_{s_{\rho_{\gamma_{1}}}}$ (so here for $x_{s_{\rho_{\gamma_{1}}}}$ we use the same convention of notations for $x_{\gamma_{2}}$ and $x_{-\gamma_{2}}$).  Note  that the composition of actions by log morphisms on $x_{s_{\rho_{\gamma_{1}}}}$ along a loop is the identity and the log morphisms across $slab_{\pm\gamma_{1}}$ are trivial for $x_{s_{\rho_{\gamma_{1}}}}$.

Therefore the gluing of $\mathcal{C}_{4}$ and $\mathcal{C}_{3}$ yields the same relation. Similarly the gluing of $\mathcal{C}_{2}$ and $\mathcal{C}_{3}$ yields the same relation obtained by gluing $\mathcal{C}_{1}$ and $\mathcal{C}_{4}$. This shows the global consistency as well as the redundancy of the relations in proposition 9.14 explicitly. We have proved the following proposition. \begin{pro}For a non-degenerate joint  all relations can be obtained by gluing  Gross-Siebert chambers from only  one side of the primary stability wall.\end{pro}

\subsection{Consequences And Discussions}

\noindent {\bf 1. \ Equivalence of Instanton Data}\\

We view  theorem 9.10, 9.11 and 9.15 as an equivalence between the metric instanton data  and the complex structure instanton data associated to the metric problem and the complex structure problem respectively.  The labeling by charges   and construction of  BPS log smooth structures and truncated $structures$ build correspondence between the metric and the complex structure instanton data in such a way that discontinuities of the metric instanton data are identified with log morphisms of the gluing of deformations in the complex structure problem while the wall crossing formula  is identified with  consistency conditions of gluing. Labeling by charges    gives  instanton meanings to Gross and Siebert's "corrections". In fact since slabs and Gross-Siebert walls are labeled by BPS charges one can think of the corrections (log morphisms) attached to them as being associated to BPS instantons which are finite trajectories with BPS charges.

Geometrically this identification is nontrivial. On the metric side the wall crossing formula is a computational tool for the enumerative problem of critical (i.e. finite) trajectories of quadratic differential foliations while on the complex structure side the consistency condition is an obstruction of the deformation problem. It is not easy to see a priori how these two problem can be related. Of course one may object by saying that the instanton correction problem of  complex structures  seems to be artificially set up to get the identification. However we have shown that it is a natural thing to do because eventually the toric degeneration obtained is a  degeneration of the Hitchin's moduli space (in its Fock-Goncharov realization) viewed as the moduli space of flat connections. \\

Here is another interesting observation. The  ordered factorizations in the wall crossing formula and the system of consistency conditions can  be inductively calculated by taking truncations. This is the algebraic way of deriving or proving them.   On the other hand a wall crossing as proved in section 8 can be derived by following a continuous variation of  $\vartheta$ and collecting the BPS charges and their Kontsevich-Soibelman transformations along the way. These two ways  yield the same answer in the end but they are actually  different in the process. The first method works by truncation and at each stage it is  possible that not all of the rays which appear in the end have appeared. At the next stage new rays can appear and can appear on both sides of existing rays which means  that these rays  do not pop out in their natural order given by the order of the refined BPS spectra. The algorithm is inductive and a closed formula of the factorization is not guaranteed. The second method can give the closed formula if one can follow the changes of decorated triangulations. This is certainly challenging in general but in some examples we can do that. Also in the second method rays appear in their natural order in the refined BPS spectra  following  the continuous variation of $\vartheta$. In this sense the identification of the wall crossing formulas in the metric problem with a system of consistent conditions is a nontrivial result relating two different mechanisms of incorporating instanton corrections. \\

\noindent {\bf 2.  \ Metric Instanton Corrections in Mirror Symmetry}\\

The solution of the metric instanton problem by Gaiotto, Moore and Neitzke is in the context of gauge theory and instantons are critical trajectories which are physically expected to be boundaries of some branes. Therefore it is not   mathematically a priori clear that it also gives instanton corrected form of the Calabi-Yau (hyperkahler) metric required in mirror symmetry.

Although the equivalence proved above by itself does not tell anything new about the description of the hyperkahler metric based on Gaiotto-Moore-Neitzke ansatz  it answers the question whether the GMN's hyperkahler metric is given by the instanton corrections required by mirror symmetry positively because the solution of the complex structure instanton correction problem given by Gross and Siebert is for mirror symmetry.

It would be nice if one can actually count holomorphic
disks in the mirror Hitchin's moduli
spaces and check the match from that point
of view.

It seems that the metric instanton correction problem  in mirror symmetry has rarely been considered in the literature.  On interesting paper that might be relevant to Hitchin's moduli spaces if \citep{Ch}. It  studies the Ooguri-Vafa metric \citep{OV} which can be considered as
 a simple local model around a single affine singularity without wall
 crossing\footnote{The title of \citep{Ch} contains "wall crossing", but it really means discontinuous jumps in our sense.}.\\

Since the complex structure and the compatible symplectic structure  determine the metric we can say that  there is no instanton corrections to the symplectic structure now that intanton corrections to complex
structures and instanton corrections to
metrics are equivalent. This fact is a general belief in the field but the author could not find a completely convincing argument. Here we have an example for which this matter is settled.\\

This equivalence also provides us a set of examples of fiberwise compact Calabi-Yau's (i.e. Hitchin's moduli spaces) for which instanton corrections to complex structures and metrics
can in principle be calculated. We just pick one point $u$ (a quadratic differential)
and follow the evolution of the decorated
triangulations when $\vartheta$ changes to $\vartheta + \pi$. Collect
all critical trajectories together with their BPS
charges one encounters and the wall crossing
formula will determine the rest.\\

The author feels that one of the most important parts in the whole picture is using the twistor method to transform the metric problem to a problem of holomorphic functions. To appreciate this point suppose we want to handle the Calabi-Yau metric using the coordinates on the Calabi-Yau manifold itself which seems to be the only choice in general cases. We would imagine that we first identify the semiflat part and after that manage to add the instanton corrections. Now suppose we want to compare this to the complex structure instanton correction problem. Naturally we would want to work with complex coordinates (e.g. Kodaira-Spencer theory). There does exist a very interesting heuristic proposal due to Fukaya \citep{Fu} which uses classical deformation  theory and other tools to deal with the complex structure problem (but not the metric problem) along this line. However that is not the approach of Gross-Siebert which does not construct a semiflat complex structure first and then deform it (see the introduction of \citep{GS1} for a brief history of their ideas).  The point is that if we want to compare the metric problem to the complex structure problem in Gross-Siebert's approach we had better have a way which does not use the coordinates of Calabi-Yau itself. For Hitchin's moduli spaces which are hyperkahler we use the twistor description which  has a holomorphic nature making the comparison to the complex structure problem much more straightforward. But this also tells us that for general Calabi-Yau's or even Calabi-Yau threefolds where a twistor-like description is absent the strategy in this paper will not work and we really need some new ideas.\\

\noindent {\bf 3.  \ Toward a Mirror Theorem}\\

The equivalence has an unexpected implication: no matter what the instantons in mirror symmetry of Hitchin's moduli spaces are their enumerative geometry must be equivalent to the enumerative geometry of critical trajectories of foliations and there should be a geometric way to see that. This would imply something nontrivial in mirror symmetry.

In fact the equivalence has set up a link between an enumerative problem and a deformation problem. This is very similar in spirit to the well known mirror formula connecting the enumerative problem of Gromov-Witten invariants and the deformation problem of calculating periods. The theorem should  also be interpreted as a mirror theorem. Again one could object by saying that while Gromov-Witten invariants and periods are obtained on Calabi-Yau varieties mirror to each other the equivalence here  identifies things on the same Hitchin's moduli space. However as explained in section 3 the instanton data of the complex structure problem is expected to be the dual data of some enumerative data in the mirror Hitchin's moduli space. So the theorem here is really half of the following conjectural full mirror statement connecting two enumerative problems of instantons on two (families of) Hitchin's moduli spaces mirror to each other.\\

\begin{con} Let $\mathcal{M}$ be an $SU(2)$ Hitchin's moduli space with prescribed singularities and $\hat{\mathcal{M}}$ its SYZ mirror which is an $PGL(2)$ Hitchin's moduli space with prescribed singularities\footnote{Strictly speaking the SYZ mirror symmetry of Hitchin's moduli spaces with Langlands dual gauge groups has not be extended in complete mathematical rigors to include prescribed singularities. But it is very likely to be true and is probably known to some experts. See \citep{GW} and \citep{W} for physicists' treatments.}. Then the enumerative problem of holomorphic disks wrapping special Lagrangian fibers in $\hat{\mathcal{M}}$\footnote{Of course this problem has to be properly formulated first. Perhaps we should only count those wrapping singular fibers?} is equivalent to the enumerative problem of critical trajectories of quadratic differential foliations on the Riemann surface.\end{con}

The mystery of this mirror conjecture is that it identifies an enumerative problem on a Hitchin's moduli space with an enumerative problem  on the Riemann surface. It is not clear geometrically how this could be true. \\

\noindent {\bf 4. \ Wall Crossings}\\

The meanings of charges, central charges and wall crossing  with respect
to stability walls are clarified in Gross-Siebert's
construction. This is evident from the theorem 9.10 and the construction of BPS Gross-Siebert data.\\

However there is a more subtle kind of wall crossing phenomenon. We call it $wall$ $crossing$ $of$ $degenerations$.

The construction of the Gross-Siebert data depends on the choice of $(u,\vartheta)$ or more generally the BPS chamber containing the pair. The affine realization of Hitchin's moduli space in terms of Fock-Goncharov relations also has the same dependency and is  not intrinsic.  Different realizations could yield an isomorphic variety. It seems that we should expect that they always yield isomorphic variety. So it is important to understand what happens if we change $(u,\vartheta)$.

Holding $u$ fixed while changing $\vartheta$ is easy to understand. This operation changes the refined BPS spectra but does not change the BPS spectra. So the BPS polyhedral decomposition is fixed and we just change the labeling of codimensional one cells. The BPS log smooth structure is changed accordingly. All joints are preserved and at each joint we use a new wall crossing formula according to the prescription in remark 6.4. Then we still get a universal $structure$ and a compatible system of consistent $structures$. Finally the defining ideals of the degeneration and the Hitchin's moduli space are changed by relabeling variables.

Holding $\vartheta$ fixed while changing $u$ is more complicated. This should be considered as the wall crossing of mirror degenerations. If a primary stability wall for the initial $u$ is crossed the BPS spectra would change which changes the BPS polyhedral decomposition. So the new Fock-Goncharov realization of $\mathcal{M}$ will be different and can not be obtained by simply renaming variables. In fact we now have different numbers of variables and relations. Nevertheless because of the wall crossing formula  at least in some examples one can show that the new degeneration is obtained by the old one by adding new relations for new variables without changing old relations between old variables so that the generic fibers for different $u$ are naturally isomorphic. So in this sense the wall crossing is realized manifestly as the change of numbers of variables and relations without changing the underlying variety. Note that in general the BPS spectra is infinite which means we have infinitely many cells to glue. As in the wall crossing formula we have to take truncations to a given order to get finiteness. So in general after a wall crossing we really have a projective system of degenerations and the wall crossing of degenerations is understood in the truncated and projective sense. It seems we need to use the machinery of formal schemes. The author believes that fully clarifying the meaning of wall crossing of degenerations in general is an important problem. In this paper we will just describe it in examples.\\

Yet  another interesting direction
is that the Kontsevich-Soibelman wall crossing formula
is actually designed to describe the wall crossings of some enumerative
problems of stable objects of certain
Bridgeland type stability conditions in certain
triangulated categories \citep{KS}. The possibility has
also been speculated by Gaiotto-Moore-Neitzke in \citep{G1}. More spectacularly, the whole picture appears to have deep relations with the study of entropy and microstates of some black holes from which the Kontsevich-Soibelman wall crossing formula can be derived \citep{A3}.\\

\noindent {\bf 5. \ SYZ vs GS}\\

The compatibility of the metric side and
algebraic side of the equivalence and the production
of toric degeneration from large $R$ degeneration
should be considered as a check of
the compatibility of the differential geometric
limit form of SYZ mirror conjecture and its
algebraic geometric version (Gross-Siebert's
version).\\

\noindent {\bf 6. \ Degenerations of Hitchin's Moduli Spaces}\\

The instanton correction problem of complex structures in this paper is understood in an algebraic sense. That is, it is in the form of explicit deformations of algebraic defining equations.   Usually in a problem of mirror symmetry a degeneration is given a priori and the task is to construct its mirror. However for Hitchin's moduli space it is not clear how to do that  and even if we had one it may not be appropriate for mirror symmetry.  Introducing $R$-deformation in section 5 is a promising step for metric aspects of mirror symmetry but that does not give us an algebraic degeneration. The results proved in this section provides a way to construct such a degeneration. It is built from some input data (Fock-Goncharov coordinates and BPS spectra) which have natural geometric meanings in the moduli interpretation of the hyperkahler space.  And because this degeneration is obtained by running Gross-Siebert algorithm it is automatically a degeneration needed by mirror symmetry. It would be interesting to see if these degenerations have any significance in other contexts. \\

What is the relation between the  deformation parameter for Hitchin's equations (namely $R$) and the  deformation parameter $t$ for the same moduli space in Gross and Siebert's construction? The following conjecture is  quite plausible given our construction.\\

\begin{con} The central fiber of Gross-Siebert's toric degeneration of Hitchin's moduli spaces coincides with the large complex limit point for the large complex degeneration of Hitchin's moduli space.\end{con}

Note that in the definition 5.1 we defined the large complex degeneration as the large $R$ family but we did not define the limit point. So to make sense of the conjecture we should first define it.  This is also a nontrivial check of the compatibility of the SYZ picture and GS picture mentioned above. Finally it is  natural to ask whether the central fiber $t=0$ (or $R=\infty$) can be given also as a moduli space of some kind of degenerate objects.\\

So an intuitive geometric picture of two degenerations in the moduli space of Hitchin's moduli spaces is the following:

There is a moduli space $\Re$ of hyperkahler structures on the underlying  manifold  of an $SU(2)$ Hitchin's moduli space $\mathcal{M}$\footnote{We may want to perturb the moduli problem slightly if necessary. This corresponds to achieve the pseudo-rationality condition in section 9.2.}. The global structure of $\Re$ is not clear. But we have a real family of Hitchin's moduli spaces with changing hyperkahler structures containing $\mathcal{M}$ and approaching a large complex point in $\Re$. Each element of this family that is close enough to the large complex point is endowed with a hyperkahler metric which is given exactly and has incorporated all instanton corrections required by mirror symmetry. There is a discrete set of  elements in this family with the following property. Fix one of them (this corresponds to choosing an integral scaling operation in section 9.2) which is close enough to the large complex point (to get large enough $R$, see section 7) and we can construct a complex family of Calabi-Yau varieties (the toric degeneration). This complex family solves the algebraic geometric version of the instanton correction problem in mirror symmetry. The complex family is not uniquely determined by the original large complex family and the choice of a fixed hyperkahler structure (fixed $R$). It depends on an additional choice (the choice of a polarization). When we say the complex family we mean we have chosen a polarization. However the instanton data associated to it is always equivalent to the instanton data for the metrics. Moreover it always contains an element which is $canonically$ isomorphic to the original Hitchin's moduli space (in its Fock-Goncharov realization) whose complex structure corresponds to the  complex structure of the  moduli space of $SL(2,\mathrm{C})$ flat connections. The complex family also goes to the large complex point (assume the previous conjecture is true). The (real) large complex family is not expected to be embedded in the complex family. \\

Unlike the (real) large complex family it is not clear under which conditions each element of the complex family is isomorphic to a Hitchin's moduli space. Since the toric degeneration is simply a complex family of complex manifolds it is not so easy to exclude the somewhat perverse possibility that an element as a complex manifold is isomorphic to a Hitchin's moduli space with one of its compatible complex structures without fixing the twistor parameter\footnote{However this point of view is not natural form the perspective of the full hyperkahler geometry of Hitchin's moduli spaces, see remark 5.1.}. The question is also complicated by the  fact that the complex family is  given as explicit deformations of ideals while the complex structures of Hitchin's moduli spaces that we are talking about can not be easily extracted from the defining ideal. There are two issues here. First of all this is clearly related to the  issue of the analytic but nonalgebraic isomorphism between the moduli space of flat connections and the moduli space of fundamental group representations, see remark 9.14. Unfortunately the author's knowledge on this issue is not enough for him to determine what this would imply. Second we are facing the highly nontrivial issue of converting the moduli information in terms of deformations of ideals to the moduli information in terms of other means and vice versa.\\

This second point deserves further remarks. To appreciate the non-triviality of the issue let us take a look at the Legendre family of elliptic curves.
$$y^{2}=x(x-1)(x-\lambda)$$ This is what we meant by deformations of the defining ideal. On the other hand we can view an elliptic curve as a complex torus and as such the moduli can be labeled by a point $\tau$ (ratio of periods) in the upper-half plane. This is what we meant by  other means. The relation between $\lambda$ and $\tau$ is given by the elliptic modular lambda function $\lambda(\tau)=16q^{1/2}-128q+704q^{3/2}+\cdots$ where $q=exp(2\pi i\tau)$.

Back to our situation. The Gross-Siebert approach gives us deformations of ideals while the hyperkahler structure (Ricci-flat metric) and  complex structures of the Hitchin's moduli space are given by more intrinsic means. For the example of elliptic curves in the previous paragraph the more intrinsic means is the periods approach as one can write $dz=dx+\tau dy$ and the Kahler form is $dz\wedge d\bar{z}$. Of course for a compact Calabi-Yau manifold or a Hitchin's moduli space we do not know how to write down the Ricci-flat metric in this way but the example of elliptic curves suggests that perhaps we need to convert the deformation of ideals to the deformation of periods. In the situation of Hitchin's moduli spaces, however, we do not know any analogues of  modular functions.
So once again (after the discussion in 9.5.2) we see that we are having a clash of two perspectives associated to the two types of degenerations. This is really one of the deep problems in mirror symmetry.  This paper offers some insights using the twistor space as a bridge but the situation is still largely unclear.\\

Finally we notice that $t$ is complex and $R$ is real. Is there a complexification of a large $R$ degeneration? It seems the $\mathbf{C}^{\times}$ action in section 2 provides such a complexification. But what is its role in the family version of SYZ mirror symmetry?\\

\subsection{Examples}
\noindent {\bf Example 1 \ (Continued)}
Let us consider the stability chamber inside the union of stability walls first. Let $u$ be a point in this stability chamber (strong coupling region). The BPS spectra at $u$ is $(\pm\gamma_{1},\pm\gamma_{2})$. So  the  ordered product of Kontsevich-Soibelman  transformations associated to the refined BPS spectra at $(u,\vartheta)$ is either $K_{\gamma_{2}}K_{\gamma_{1}}$ or $K_{-\gamma_{1}}K_{\gamma_{2}}$ or $K_{-\gamma_{2}}K_{-\gamma_{1}}$ or $K_{\gamma_{1}}K_{-\gamma_{2}}$. Without loss of generality let us take $K_{\gamma_{2}}K_{\gamma_{1}}$.

Now we let $u$ cross a stability wall from the strong coupling region to the weak coupling region. If we follow a continuous evolution of $\vartheta$ by drawing pictures carefully (for these pictures and many more, see \citep{G2}) or using computer we would see that while in the strong coupling region there are only two flips in the nearby weak coupling region there are three flips of edges labeled successively by $\gamma_{2}, \gamma_{1}+\gamma_{2}$ and $\gamma_{1}$. This order is actually mandatary to us without following the evolution because we know  that  the phases of BPS rays of $\gamma_{1}$ and $\gamma_{2}$ have switched. So the ordered product should be $$K_{\gamma_{1}}K_{\gamma_{1}+\gamma_{2}}K_{\gamma_{2}}$$ The wall crossing formula in this case says \begin{equation}K_{\gamma_{2}}K_{\gamma_{1}}=K_{\gamma_{1}}K_{\gamma_{1}+\gamma_{2}}K_{\gamma_{2}}\end{equation} Note that we have used the fact that $\Omega(\gamma; u)$ is always one for flips. Once we know the closed formula, it is straightforward to verify it if one finds the proof given above is not rigorous enough (it is rigorous). Since the wall crossing formula is local we can locally identify the charge lattice with $\mathbf{Z}^{2}$ such that $\bar{m}_{\gamma_{1}}=(-1,0), \bar{m}_{\gamma_{2}}=(0,-1)$ (so $\gamma_{1}= (0,-1), \gamma_{2}=(1,0)$) and also identify $\mathcal{X}_{\gamma_{1}}$ and $\mathcal{X}_{\gamma_{2}}$ as a complete set of  independent variables  satisfying the multiplicative relation and the Poisson bracket relation.  Then we compose automorphisms on both sides applied to both $\mathcal{X}_{\gamma_{1}}$ and $\mathcal{X}_{\gamma_{2}}$ and see that they coincide which is sufficient to deduce the formula.

The derivation is symmetric with respect to the change of the role of $\gamma_{1}$ and $\gamma_{2}$. In fact, varying from $\vartheta+\pi$ to $\vartheta$ gives us
\begin{equation}K_{\gamma_{1}}K_{\gamma_{2}}=K_{\gamma_{2}}K_{\gamma_{1}+\gamma_{2}}K_{\gamma_{1}}\end{equation}Similarly we have more identities such as \begin{equation}K_{-\gamma_{2}}K_{\gamma_{1}}=K_{\gamma_{1}}K_{\gamma_{1}-\gamma_{2}}K_{-\gamma_{2}}\end{equation}

We have obtained the formula by following the continuous evolution of $\vartheta$ or at least knowing the BPS spectra of charges. But in fact we can derive it without doing so. We just use Kontsevich-Soibelman's theorem in section 3. We produce a scattering diagram on the plane in the following way. We use the set up of example 1 in section 9.1. Define $$D:=\{(\mathbf{R}(1,0),(1+t\mathrm{x})), (\mathbf{R}(0,1),(1+t\mathrm{y}))\}$$Then by the definition given in section 3, the associated log automorphisms are $$K_{1}: \mathrm{x}\rightarrow \mathrm{x}, \mathrm{y}\rightarrow \mathrm{y}(1+t\mathrm{x}),\ \mathrm{for \ ray} \ \mathbf{R}_{\leq 0}(1,0)$$
$$K_{2}: \mathrm{x}\rightarrow \mathrm{x}/(1+t\mathrm{y}), \mathrm{y}\rightarrow \mathrm{y},\ \mathrm{for \ ray} \ \mathbf{R}_{\leq 0}(0,1)$$
$$K_{1}^{-1}:\mathrm{x}\rightarrow \mathrm{x}, \mathrm{y}\rightarrow \mathrm{y}/(1+t\mathrm{x}),\ \mathrm{for \ ray} \ \mathbf{R}_{\geq 0}(1,0)$$
\begin{equation}K_{2}^{-1}:\mathrm{x}\rightarrow \mathrm{x}(1+t\mathrm{y}), \mathrm{y}\rightarrow \mathrm{y},\ \mathrm{for \ ray} \ \mathbf{R}_{\geq 0}(0,1)\end{equation}
where the automorphisms are taken when one crosses the ray counterclockwise. After setting $t=1$ the above transformations are identified as the Kontsevich-Soibelmann transformations. In fact, $$\mathrm{x}\rightarrow\mathcal{X}_{-\gamma_{1}},\ \mathrm{y}\rightarrow\mathcal{X}_{-\gamma_{2}}$$\begin{equation}t\rightarrow 1, K_{i}\rightarrow K_{\gamma_{i}}\end{equation} Now if we follow a loop starting from the first quadrant counterclockwise then the counterclockwise ordered product of automorphisms is $$K_{1}^{-1}K_{2}K_{1}K_{2}^{-1}\neq 1$$ So we have to follow the procedure in the proof of the theorem to add rays such that the new diagram is consistent. In this simple case, the consistency in the first order persists to higher orders and the result is that we need to add only one ray.  It is $(\mathbf{R}(1,1),(1+t^{2}\mathrm{xy}))$ whose associated counterclockwise automorphism is $$K_{1+2}^{-1}: \mathrm{x}\rightarrow \mathrm{x}/(1+t^{2}\mathrm{xy}), \mathrm{y}\rightarrow \mathrm{y}(1+t^{2}\mathrm{xy})$$ so that  following the loop counterclockwise we have $$K_{1+2}^{-1}K_{1}^{-1}K_{2}K_{1}K_{2}^{-1}= 1$$ which is nothing but (134) by using (138).

Although we have got the same answer in the end, the second approach is  based on a different mechanism. In the first approach we follow the order of BPS rays as we meet them along a continuous evolution to derive  the wall crossing formula which selects the BPS spectra. In the second approach one does not need to use stability walls and even if one puts the stability wall into the picture one does not meet all BPS rays because some of them (ray (1,1)) are not known until the wall crossing formula is obtained and they can appear without respecting the order of BPS rays.\\

As promised in section 8 there is a third way to derive the formula. Recall the definition of a spectrum generator given there.  The cumulative result of the variation from $\vartheta$ to $\vartheta+\pi$ is an omnipop whose associated transformation $S$ in our example is (see \citep{G2} for the derivation) $$\mathcal{X}_{\gamma_{1}}\rightarrow \mathcal{X}_{\gamma_{1}}(1+\mathcal{X}_{\gamma_{2}})$$$$\mathcal{X}_{\gamma_{2}}\rightarrow\mathcal{X}_{\gamma_{1}}
(1+\mathcal{X}_{\gamma_{1}}+\mathcal{X}_{\gamma_{1}}\mathcal{X}_{\gamma_{2}})^{-1}$$We assume the phase of $l_{\gamma_{1}}$ is large than the phase of $l_{\gamma_{2}}$ following the convention used there. The author wants to emphasize that $S$ is obtained without knowing any of the Kontsevich-Soibleman factors  or even any charges appearing in either side of the wall crossing formula! However due to the strong constraint that the wall crossing formula imposes this is already enough. We are seeking a decomposition of the form $$S=\prod_{m,n\geq 0}K_{m\gamma_{1}+n\gamma_{2}}^{\Omega(m\gamma_{1}+n\gamma_{2};u)}$$Inside the stability walls, since the phase of $l_{\gamma_{1}}$ is large than the phase of $l_{\gamma_{2}}$, the product order must be the increasing order of $m/n$ from $0$ to $\infty$. We then  truncate the product successively by the degree and then take the projective limit. As expected, truncation up to order two (hence involving only $K_{\gamma_{1}}$ and $K_{\gamma_{2}}$) gives us the expected decomposition $K_{\gamma_{1}}K_{\gamma_{2}}$ which by induction can be shown to persist to all higher orders. Therefore we have derived the left hand side of (135). Working outside the stability wall by decomposing in the reverse order we get the right hand side.\\

Identity (134) (or its variations like (135) (136)) is called the pentagon identity and is the simplest nontrivial wall crossing formula. It has been encountered by many authors. The example presented here was described in \citep{G2} (for the Fock-Goncharov coordinates part). The scattering diagram part follows \citep{GPS}. \\

Next we are going to consider the corresponding complex structure problem. This example has been studied in \citep{GS3}.

We have an initial $structure$ with only one joint $v$ which is the origin. A $structure$ in the two dimension with only one joint reduces to a scattering diagrams.  It is  the scattering diagram we just described  in this section and the wall crossing gives the additional Gross-Siebert walls  that have to be inserted. It is nothing but the ray $(1,1)$ together with its log automorphism.

 We have canonical thickenings\footnote{We define $i$ modulo 4.} such as  $$R_{v\rightarrow\sigma_{i},\sigma_{i}}^{k}, R_{v\rightarrow\rho_{i},\sigma_{i}}^{k},R_{v\rightarrow\rho_{i},\sigma_{i-1}}^{k}, R_{v\rightarrow v,\sigma_{i-1}}^{k},$$ $$ R_{\rho_{i}\rightarrow\sigma_{i},\sigma_{i}}^{k}, R_{\rho_{i}\rightarrow\sigma_{i-1},\sigma_{i-1}}^{k}, R_{\rho_{i}\rightarrow\rho_{i},\sigma_{i}}^{k},R_{\rho_{i}\rightarrow\sigma_{i-1},\sigma_{i-1}}^{k}$$They are glued in two ways. The first type is change of strata. In other words, for $\tau^{'}\subseteq\tau$, we have $R_{v\rightarrow \tau,\sigma_{i}}^{k}\rightarrow R_{v\rightarrow \tau^{'},\sigma_{i}}^{k}$. There is no need to compose automorphisms. The second type is change of chambers. For example we glue $R_{\rho_{i}\rightarrow\sigma_{i},\sigma_{i}}^{k}$ and $R_{\rho_{i}\rightarrow\sigma_{i-1},\sigma_{i-1}}^{k}$ by identifying $R_{\rho_{i}\rightarrow\rho_{i},\sigma_{i}}^{k}$ and  $R_{\rho_{i}\rightarrow\rho_{i},\sigma_{i-1}}^{k}$. To identify $R_{\rho_{i}\rightarrow\rho_{i},\sigma_{i}}^{k}$ and  $R_{\rho_{i}\rightarrow\rho_{i},\sigma_{i-1}}^{k}$ we need a parallel transport from $\sigma_{i-1}$ to $\sigma_{i}$. Suppose we only have the one axis and only one singularity which is the one on that axis, then the naive gluing is already inconsistent because of the nontrivial monodromy which means that the two different ways of parallel transports bypassing the singularity from different sides give different identifications. The gluing is therefore not well defined. The automorphism attached to the cell $\rho_{i}$  induced by $f_{\rho_{i},v}$ makes this gluing consistent (in the absence of the other singularity), see \citep{GS3}. Similar things happen for the other singularity assuming the absence of this one. However, when both of the two singularities appear, due to their interaction ("scattering"), the slab functions given here that once guarantee consistency separately do not make the gluing consistent any more. The new consistency condition of the gluing (which is not induced by monodromies) is the one given in the definition 3.19. So moving along a loop around the only codimensional two cell $v$, the composition of morphisms induced by the slab functions between these affine pieces $$R_{v\rightarrow \sigma_{1}, \sigma_{1}}^{k}\rightarrow R_{v\rightarrow \sigma_{2}, \sigma_{2}}^{k}\rightarrow R_{v\rightarrow \sigma_{3}, \sigma_{3}}^{k}\rightarrow R_{v\rightarrow \sigma_{4}, \sigma_{4}}^{k}\rightarrow R_{v\rightarrow \sigma_{1}, \sigma_{1}}^{k}$$must be  the identity. The  composition of the log automorphisms in terms of generators are  \footnote{The direction of the loop is counterclockwise starting from the first quadrant.}$$x\rightarrow (1+wy)^{-1}x, y\rightarrow (1+wy)y$$ $$z\rightarrow (1+wy)z, w\rightarrow (1+wy)^{-1}w$$ which is not the identity. To compensate that, one only need to add a ray (a Gross-Siebert wall) $p:=\mathbf{R}_{\geq 0}(1,1)$ with $f_{p, v}:=(1+wy)=(1+t^{2}x^{-1}z^{-1})$. This is precisely the result presented before by a change of variables.

Note that the union of stability walls separates the five rays (one Gross-Siebert wall and four slabs). Inside it we have two rays associated to two charges $\gamma_{1},\gamma_{2}$ and outside it we have three rays  with the order of rays labeled by $\gamma_{1}$ and $\gamma_{2}$ reversed and a new one added according to the wall crossing formula. By its slope the new ray is clearly the projection of the BPS wall associated to the charge $\gamma_{1}+\gamma_{2}$.   This is the wall crossing interpretation of the consistency condition in Gross and Siebert's construction.\\

Let us  continue and try to find the defining equations of the degeneration. $p$ divides $\sigma_{1}$ into two chambers denoted by $\mathrm{u}_{1}, \mathrm{u}_{2}$ clockwise. Note that in section 3's notations $\sigma_{\mathrm{u}_{1}}=\sigma_{\mathrm{u}_{2}}=\sigma_{1}$. The degeneration is obtained by gluing $R_{\rho_{2}\rightarrow\sigma_{1},\sigma_{\mathrm{u}_{1}}}^{k}$ and $R_{\rho_{2}\rightarrow\sigma_{2},\sigma_{2}}^{k}$, $R_{\rho_{3}\rightarrow\sigma_{2},\sigma_{2}}^{k}$ and $R_{\rho_{3}\rightarrow\sigma_{3},\sigma_{3}}^{k}$, $R_{\rho_{4}\rightarrow\sigma_{3},\sigma_{3}}^{k}$ and $R_{\rho_{4}\rightarrow\sigma_{4},\sigma_{4}}^{k}$, $R_{\rho_{1}\rightarrow\sigma_{4},\sigma_{4}}^{k}$ and $R_{\rho_{1}\rightarrow\sigma_{1},\sigma_{\mathrm{u}_{2}}}^{k}$, and finally $R_{\rho_{1}\rightarrow\sigma_{1},\sigma_{\mathrm{u}_{2}}}^{k}$ and $R_{\rho_{2}\rightarrow\sigma_{1},\sigma_{\mathrm{u}_{1}}}^{k}$. The treatments of the second and the third  gluing are almost the same. The only difference is just a renaming of variables. So let us consider the second gluing $R_{\rho_{3}\rightarrow\sigma_{2},\sigma_{2}}^{k}$ and $R_{\rho_{3}\rightarrow\sigma_{3},\sigma_{3}}^{k}$. Using the definitions in section 3, one gets $$R_{\rho_{3}\rightarrow\sigma_{2},\sigma_{2}}^{k}=\mathbf{C}[t]/(t^{k+1})[\Lambda_{\rho_{3}}][z,w]/(zw-t, w^{k+1})$$  Similarly we have $$R_{\rho_{3}\rightarrow\sigma_{3},\sigma_{3}}^{k}=\mathbf{C}[t]/(t^{k+1})[\Lambda_{\rho_{3}}][z,w]/(zw-t, z^{k+1})$$ There are canonical quotient homomorphisms from $R_{\rho_{3}\rightarrow\sigma_{2},\sigma_{2}}^{k}$ and $R_{\rho_{3}\rightarrow\sigma_{3},\sigma_{3}}^{k}$ to $R_{\rho_{3}\rightarrow\rho_{3},\sigma_{2}}^{k}$ and $R_{\rho_{3}\rightarrow\rho_{3},\sigma_{3}}^{k}$ respectively (followed by localizations) and the gluing is obtained by identifying $R_{\rho_{3}\rightarrow\rho_{3},\sigma_{2}}^{k}$ and $R_{\rho_{3}\rightarrow\rho_{3},\sigma_{3}}^{k}$ by the log automorphism. $R_{\rho_{3}\rightarrow\rho_{3},\sigma_{2}}^{k}$ and $R_{\rho_{3}\rightarrow\rho_{3},\sigma_{3}}^{k}$ are the same ring $\mathbf{C}[\Lambda_{\rho_{3}}][t]/(t^{k+1})[z,w][zw-t, z^{k+1},w^{k+1}]$ and the log automorphism is $$z\rightarrow f_{\rho_{3},v}z=(1+y)z, w\rightarrow f_{\rho_{3},v}^{-1}w=(1+y)^{-1}w$$
The fiber product is isomorphic to $$\mathbf{C}[t]/(t^{k+1})[\Lambda_{\rho_{3}}][Z,W][ZW-(1+y)t]$$by the homomorphism $$Z\rightarrow(z,(1+y)z), W\rightarrow((1+y)w,w)$$By the same argument the third gluing gives $$\mathbf{C}[t]/(t^{k+1})[\Lambda_{\rho_{4}}][X,Y][XY-(1+w)t]$$ $$X\rightarrow((1+w)x,x), Y\rightarrow(y,(1+w)y)$$  The first gluing and the fourth gluing reproduce these two gluing results by proposition 9.16 because the joint is non-degenerate. As for the fifth gluing, it is not really a gluing of different affine pieces along a substrata. All such gluing has been done. In fact it is induced by  two changes of of strata  from $\rho\rightarrow \sigma$ to $v\rightarrow\sigma$ and a change of chambers in the middle. The nontrivial part is the change of chambers  induced by the division of $\sigma_{1}$ into two chambers by the wall $p$. $$R_{v\rightarrow\sigma_{1},\sigma_{\mathrm{u}_{1}}}^{k}\rightarrow R_{v\rightarrow\sigma_{1},\sigma_{\mathrm{u}_{2}}}^{k}$$ Since $\sigma_{1}=\sigma_{\mathrm{u}_{1}}=\sigma_{\mathrm{u}_{2}}$, the prescription in section 3 says that the transformation is just the log automorphism induced by $f_{p,v}$ which clearly preserves the relations $xy=t, zw=t$. 

The above calculation tells us that  the toric degeneration produced by our construction is  $$\mathrm{Spec}\ \mathbf{C}[t]/(t^{k+1})[X,Y,Z,W,t][XY-(1+W)t, ZW-(1+Y)t]$$From the explicit presentation of the ring it is clear that the  $k+1$-th order ring is naturally compatible to the $k$-th order ring simply by taking the quotient homomorphisms $\mathbf{C}[t]/(t^{k+2})\rightarrow \mathbf{C}[t]/(t^{k+1})$.  In particular we do not need to change any other  relations and  log morphisms. This is the compatibility of consistent $k$-th and $k+1$-th $structures$ in the sense of Gross and Siebert. So we can send $k$ to infinity and know that the defining equations of the total degeneration over $\mathrm{Spec}\ \mathbf{C}[[t]]$ are $XY=(1+W)t, ZW=(1+Y)t$ and the  fiber over $t=1$ is   \begin{equation}XY=(1+W), ZW=(1+Y)\end{equation}

This is an intersection of two degree two hypersurfaces in $\mathbf{C}^{4}$ and therefore is an affine Calabi-Yau variety.\\

On the other hand, the same equations are also the equations of the Hitchin's moduli space in terms of Fock-Goncharov coordinates associated to two independent charges/edges. Although this has been proved before for general cases the proof ignores the role of cluster transformations as we only use Kontsevich-Soibelman transformations before. It is interesting  to derive it in a way making cluster transformations transparent.  Denote the two Fock-Goncharov coordinates associated to two edges for the initial WKB triangulation at $(u,\vartheta)$ as $x_{1}, y_{1}$. First let us consider the case when the moduli parameter $u$ is inside the stability walls (the strong coupling region) and vary $\vartheta$ to $\vartheta+\pi$. Then we encounter two flips.  According to the cluster transformations of two successive flips, we define $$y_{2}=x_{1}^{-1}, x_{2}=y_{1}(1+x_{1}), y_{3}=x_{2}^{-1}, x_{3}=y_{2}(1+x_{2})$$ We immediately get $$x_{1}x_{3}=1+x_{2}=1+{1\over y_{3}}$$ \begin{equation}{1\over y_{1}}{1\over y_{3}}=1+{1\over y_{2}}=1+x_{1}\end{equation} It is important to track the charge/edge labels. $x_{1}, x_{3}$ are labeled by the same charge (up to a sign)  with different signs (or equivalently labeled by edges) and the same goes for $y_{1}^{-1}, y_{3}^{-1}$ for the other charge\footnote{Note that $y_{i}^{-1}$ is the Fock-Goncharov coordinate with the negative charge of $y_{i}$ and therefore it is the coordinate labeled by charge with the effect of the mutation of charges incorporated (see  section 8). Since the Kontsevich-Soibelman transformations differ from cluster transformation by the mutations of charges, using $y_{i}^{-1}$ instead of $y_{i}$ is consistent.}.  These two equations (with redundant but natural variables) define the Hitchin's moduli space as an affine variety in a way which is democratic to all BPS charges in the strong coupling region.

The equations are identical to the equation (139) by the identification $$X\rightarrow x_{3}, Y\rightarrow x_{1}, Z\rightarrow {1\over y_{1}}, W\rightarrow {1\over y_{3}}$$and therefore the Hitchin's moduli space is indeed embedded into the toric degeneration constructed by Gross and Siebert's algorithm.  The identification  keeps  track of labeling by charges.\\

We can do the same thing for the moduli region outside the union of stability walls (the weak coupling region). There are three flips corresponding to three charge $\gamma_{1}, \gamma_{1}+\gamma_{2}, \gamma_{2}$ of the BPS spectra in this region. The corresponding defining equations are  given by $$y_{n+1}=x_{n}^{-1}, x_{n+1}=y_{n}(1+x_{n}), \ n=1,2,3.$$ Although we now have two more equations the underlying varieties (as varieties in a product of copies of $\mathbf{C}^{\times}$, see remark 9.14) are isomorphic via the canonical map $$(x_{1},x_{2},x_{3},x_{4},y_{1},y_{2},y_{3},y_{4})\rightarrow(x_{1},x_{2},x_{3},y_{1},y_{2},y_{3})$$ with inverse map $$(x_{1},x_{2},x_{3},y_{1},y_{2},y_{3})\rightarrow(x_{1},x_{2},x_{3},y_{3}(1+x_{3}),y_{1},y_{2},y_{3}, x_{3}^{-1})$$Eliminating some variables, we get that the ideal of the variety is generated by \begin{equation}x_{1}x_{3}-(1+x_{2}), x_{2}x_{4}-(1+x_{3}), x_{4}x_{1}-(1+x_{5})\end{equation}where $x_{5}:=y_{4}(1+x_{4})$. These variables are  labeled by BPS charges in that stability chamber.

From the perspective of Gross and Siebert's approach to the complex structure problem, this just means that we can start from other BPS polyhedral decompositions and everything works consistently. We define a polyhedral decomposition using projections of BPS walls labeled by $\gamma_{1}, \gamma_{1}+\gamma_{2}, \gamma_{2}$ together with log morphisms according to Kontsevich-Soibelman transformations  associated to the three independent flips labeled by these three charges in this region. Then we would get an inconsistent scattering diagram with six cuts (there is a cut with slope $(-1,-1)$). After the wall crossing calculation the slab function attached to the cut with slope $(-1,-1)$ is modified to 1 and is therefore deleted.  The obtained consistent $structure$ is  the one determined above with five slabs.

Now let us construct the degeneration from this new BPS polyhedral decomposition arising from the weak coupling region. It is instructive to see how the construction is consistent to the construction arising  from the strong coupling region.  Note that unlike the situation inside the stability wall this time we have five maximal dimensional cells instead of four. Let  same symbols $x,y,z,w,t$ to denote the same monoid variable as before. We define $p_{0}:=z^{(1,1,0)}$. The monoid variable associated to the fifth ray $\rho=\mathbf{R}_{\geq 0}(1,1)$ (with slab function $f_{\rho}=(1+wy)=(1+t^{2}x^{-1}z^{-1})$) is \begin{equation}p=z^{(1,1,2)}=p_{0}t^{2}\end{equation}Let $x_{s_{\rho_{1}}}:=z^{(1,0,0)}$. Clearly \begin{equation}p=zx_{s_{\rho_{1}}}t\end{equation}

By the exactly same calculations for the strong coupling region the gluing of the second quadrant and the third quadrant gives us the equation $ZW=(1+Y)t$. The gluing of the part of the first quadrant below $\rho$ and the fourth quadrant  gives us  \begin{equation}PW=(1+Y)x_{s_{\rho_{1}}}t^{2}\end{equation}Using  (143) the equation (144) becomes $ZW=(1+Y)t$.

 Similarly we get $XY=(1+W)t$ by gluing either  the third quadrant and the fourth quadrant or the part of the first quadrant above $\rho$ and the second quadrant. We are left with the gluing of  the part of the first quadrant above $\rho$ and the part below it.  This gluing is given by the ideal generated by \begin{equation}XZ-(1+t^{2}x^{-1}z^{-1})p_{0}t^{2}\end{equation} Because of the relations $x^{-1}z^{-1}p_{0}t^{2}=1$ and $p=p_{0}t^{2}$, (145) is $$XZ-(t^{2}+P)$$So the gluing of all maximal cells gives us \begin{equation}XY=(1+W)t, ZW=(1+Y)t, XZ=(t^{2}+P)\end{equation}and setting $t=1$ one gets $$XY=(1+W), ZW=(1+Y), XZ=(1+P)$$which recovers (141) obtained from Fock-Goncharov coordinates in the weak coupling region. This variety is an intersection of three degree two hypersurfaces in $\mathbf{C}^{6}$  and therefore is an affine Calabi-Yau variety\footnote{It is in $\mathbf{C}^{6}$ instead of $\mathbf{C}^{5}$. We want our variables to be labeled by all BPS charges at $u$ in the slab labeling and therefore we must keep the variable labeled by $-\gamma_{1}-\gamma_{2}$ even if its corresponding slab has been deleted in the construction of the universal $structure$.}.\\

One can  unify the above two complementary descriptions by comparing the WKB triangulation $T_{WKB}(\theta,u)$ and $T_{WKB}(\theta+\pi,u^{'})$ where $u$ and $u^{'}$ are close enough\footnote{They are close enough to avoid the crossing of $\vartheta$, $\vartheta+\pi$ by BPS rays' phases.} and are inside and outside of the union of stability walls respectively. The corresponding wall crossing formula is $K_{\gamma_{2}}^{-1}K_{\gamma_{1}}^{-1}K_{\gamma_{2}}K_{\gamma_{1}+\gamma_{2}}K_{\gamma_{1}}=1$ which tells us that we can also follow a loop $(\vartheta,u^{'})\rightarrow(\vartheta+\pi,u^{'})\rightarrow(\vartheta+\pi,u)\rightarrow(\vartheta,u)$ where $5=3+2$ flips are encountered in the order of the factors in the wall crossing formula. Therefore we define recursively $$y_{n+1}=x_{n}^{-1}, x_{n+1}=y_{n}(1+x_{n})$$ Take a look at the picture of flips or the wall crossing formula, you would agree that these seemingly infinitely many variables are actually periodic with period five  which can be easily verified by algebra. So  the democratic (to all BPS charges for the whole moduli region) way of writing  the ideal is $$I:=(x_{n-1}x_{n+1}-(1+x_{n}))$$ with $x_{n}=x_{n+5}$. This is of course consistent to the above two descriptions in the weak and strong coupling region.

Now we can give a more explicit interpretation of  "wall crossing"  in the complex structure problem. The natural defining ideals obtained via the Fock-Goncharov relations are in terms of variables labeled by BPS charges and as such they exhibit wall crossing phenomenon. The wall crossing formula guarantees the consistency of different descriptions.

The relations $y_{n+1}=x_{n}^{-1}, x_{n+1}=y_{n}(1+x_{n})$ are known as Zamolodchikov's Y-system. It is related to the thermodynamic Bethe Ansatz mentioned in section 6. Periodicity in Y-systems is a beautiful story and has been studied by many people, see for example \citep{FZ4}. However, most wall crossing formulas give rise to non-periodic relations.\\\\

\noindent{\bf Example 2}\ \ Here we describe a true example of the cases (quadratic differentials with even order poles) studied in the main results of this section. Since it is very similar to the previous example and everything is a simple analogue to its counterpart the presentation will be very brief. We only describe the wall crossing formula. This example is also from \citep{G2}.

The Riemann surface $C$ is still $CP^{1}$. The space $B$ is complex one dimensional and consists of the following quadratic differentials $$\lambda^{2}=(z^{4}+4\Lambda^{2}z^{2}+2mz+u)dz^{2}$$ where $u$ parameterizes $B$ and both $\Lambda$ and $m$ are constants. There is an order eight pole at the infinity and clearly $m$ is the mass parameter. For simplicity let us set $m=0$.

There are four simple zeroes on $C$. There are two singular points on $B$ given by $u=0$ (multiplicity one) and $u=4\Lambda^{4}$ (multiplicity two).

$S_{u}$ is an elliptic curve with two punctures lying over the infinity. The charge lattice $\hat{\Gamma}\simeq\mathbf{Z}^{3}$ after choosing branch cuts and it contains a one dimensional flavor charge lattice. $\hat{\Gamma}$ is generated by three charges $\gamma_{1}, \gamma_{2}, \gamma_{3}$ such that $\gamma_{2}+\gamma_{3}$ is a pure flavor charge and $$\langle\gamma_{2}, \gamma_{3}\rangle=0, \langle\gamma_{3}, \gamma_{1}\rangle=\langle\gamma_{1}, \gamma_{2}\rangle=1$$ Since the residue $m$ is zero, for the pure flavor charge $\gamma_{2}+\gamma_{3}$ $$Z_{\gamma_{2}}+Z_{\gamma_{3}}=0$$
This tells us that there are only two stability walls. One is given by the alignment of $Z_{\gamma_{1}},Z_{\gamma_{2}}$ while the other is given by the alignment of $Z_{\gamma_{1}},Z_{\gamma_{3}}$. The union is a closed curve passing through the two singularities.

It is easy to check that one possible wall crossing formula is $$K_{\gamma_{1}}K_{\gamma_{2}}K_{-\gamma_{3}}=
K_{\gamma_{2}}K_{-\gamma_{3}}K_{\gamma_{1}+\gamma_{2}-\gamma_{3}}K_{\gamma_{1}+\gamma_{2}}K_{\gamma_{1}-\gamma_{3}}K_{\gamma_{1}}$$ where the left hand side is for the stability chamber inside the union. Note that here we have multiple joints with the same support.

It is interesting to see how the stability walls and wall crossing formulas split when we allow a nonzero $m$. Details can be found in \citep{G2}.\\\\

\noindent{\bf Example 3}\ \ We  want to consider an example with infinitely many jumps in the wall crossing formula.

The metric problem part of this example is considered in \citep{G2} and the BPS spectra is identified with the BPS spectra  of the pure $SU(2)$ gauge theory. The corresponding scattering diagram is determined in \citep{GPS}.

The underlying Riemann surface is still $CP^{1}$. The quadratic differentials are$$\lambda^{2}=({\Lambda^{2}\over z^{3}}+{2u\over z^{2}}+{\Lambda^{2}\over z})dz^{2}$$
As before $\Lambda$ is a real positive constant and $u$ is the moduli parameter on the complex one dimensional affine base of the Hitchin's fibration. There are two order three irregular poles of $\lambda^{2}$ at $z=0$ and $z=\infty$. This form of $\lambda^{2}$ is chosen to maintain the asymptotic behaviors of solutions of Hitchin's equation prescribed by an element of this one parameter family of quadratic differentials. There are two simple zeroes which collide at $u_{1}=-\Lambda^{2}$ and $u_{2}=\Lambda^{2}$. $u_{1}$ and $u_{2}$ are the singularities of the affine structure. The charge lattice has rank two with basis denoted by $\gamma_{1}$ and $\gamma_{2}$ such that $\gamma_{1}$($\gamma_{2}$) is the vanishing cycle at $u_{1}$($u_{2}$). There are two stability walls and their union is a simple closed curve passing through $u_{1}$ and $u_{2}$.

Let us consider the WKB triangulations for $u$ staying inside the union. Like in the first example, the BPS spectra (charges) here are $\gamma_{1},\gamma_{2}$ (up to a sign). According to the description  of local behaviors of trajectories near singularities in section 7,  since the order of (say) the singularity at zero is three there is a single Stokes ray and hence a trajectory connecting the singularity\footnote{More precisely the endpoint is on the boundary of a small open disk containing zero. But since there is only one, we can identify it with the singularity itself.} to itself. This means that we have degenerate triangles. There are also generic trajectories connecting the two singularities and these trajectories arise in two one parameter families separated by the two simple zeroes. Pictures can be found on the page 123 in \citep{G2}. From the picture it is easy to see that $$\langle\gamma_{1},\gamma_{2}\rangle=2$$

Varying $\vartheta$ to $\vartheta+\pi$ inside that stability chamber we encounter two flips. The relevant critical finite trajectories are two trajectories connecting the two zeroes and they are on different sides of the singularity at zero\footnote{Of course, one can also use the other singularity.}. So the ordered product of Kontsevich-Soibelman transformations is\footnote{Assuming without loss of generality that the BPS ray of $\gamma_{1}$ has larger phase. } $$K_{\gamma_{1}}K_{\gamma_{2}}$$ with
$$K_{\gamma_{1}}: \mathcal{X}_{\gamma_{1}}\rightarrow\mathcal{X}_{\gamma_{1}}, \ \mathcal{X}_{\gamma_{2}}\rightarrow\mathcal{X}_{\gamma_{2}}(1+\mathcal{X}_{\gamma_{1}})^{-2}$$ $$K_{\gamma_{2}}: \mathcal{X}_{\gamma_{2}}\rightarrow\mathcal{X}_{\gamma_{2}}, \ \mathcal{X}_{\gamma_{1}}\rightarrow\mathcal{X}_{\gamma_{1}}(1+\mathcal{X}_{\gamma_{2}})^{2}$$

Now let us move to the other stability chamber. We meet the flip labeled by $\gamma_{1}$ first. Then we meet infinitely many flips in the scenario of reaching a limit configuration as described in section 7. In fact, the relevant BPS rays are rays associated to $((n+1)\gamma_{1}+n\gamma_{2})$ where $n$ is a nonnegative integer\footnote{Now that the phase of the BPS ray of $\gamma_{1}$ is smaller than the phase of the BPS ray of $\gamma_{2}$ outside the union of stability walls, this is an counterclockwise order when $n\rightarrow \infty$}. To see this fact just notice that the region between the two trajectories connecting the two singularities to themselves is an annular region and together with the other two trajectories this is the initial configuration of the infinite flip scenario in section 7. The loop formed by joining the two generic trajectories has winding number one around the inner circle. Flipping once increases the winding number by one  and therefore we have the above BPS spectra. The similar thing does not happen inside the union of stability walls because in that case initially the winding number of the loop formed by joining  the  two generic trajectories is zero.

This infinite sequence of BPS rays converges to the "limit" ray with charge $\gamma_{1}+\gamma_{2}$. If we vary the angular phase $\vartheta$ starting from a phase larger than the phase of the BPS ray of $\gamma_{2}$ clockwise instead of counterclockwise as we have been doing, then we have another infinite sequence of flips with BPS charges $(n+1)\gamma_{2}+n\gamma_{1}$ also converging to $\gamma_{1}+\gamma_{2}$. Therefore we are in the situation of having to take a juggle. In other words we compose infinitely many flips labeled by $((n+1)\gamma_{1}+n\gamma_{2})$ (which is  in a $finite$ range of the variation of $\vartheta$) succeeded with a jump from the limit Fock-Goncharov coordinates $\mathcal{X}^{+}_{\gamma}$ to the limit Fock-Goncharov coordinates $\mathcal{X}^{-}_{\gamma}$. Then we continue increasing $\vartheta$ and pass through infinitely many flips labeled by $(n+1)\gamma_{2}+n\gamma_{1}$ and compose the associated Kontsevich-Soibelman tranformations counterclockwise. While maybe one can make sense of the composition of the first infinite sequence in the ordinary limit sense one certainly cannot do that for the second infinite sequence as they are to be composed backwards from the limit. The total infinite composition therefore has to be understood in the truncated and projective limit sense. To write down the wall crossing formula, we need to know the transformations associated to juggles. According to section 8 in our example it is $K_{\gamma_{1}+\gamma_{2}}^{-2}$. So the wall crossing formula is \begin{equation}K_{\gamma_{1}}K_{\gamma_{2}}=K_{\gamma_{2}}K_{\gamma_{1}+2\gamma_{2}}K_{2\gamma_{1}+3\gamma_{3}}\cdots K_{\gamma_{1}+\gamma_{2}}^{-2}\cdots K_{3\gamma_{1}+2\gamma_{2}}K_{2\gamma_{1}+\gamma_{2}}K_{\gamma_{1}}\end{equation}\begin{note}If we identify the charge lattice with
 $\mathbf{Z}^{2}$ with the integral pairing $\langle(p, q), (p^{'},q{'})\rangle= pq^{'}-qp^{'}$ and pick the basis as $\gamma_{1}=(2,-1)$ and $\gamma_{2}=(0,1)$, then we recover the equation $$K_{2,-1}K_{0,1}=K_{0,1}K_{2,1}K_{4,1}\cdots K_{2,0}^{-2}\cdots K_{6,-1}K_{4,-1}K_{2,-1}$$in section 6 which determines the BPS spectrum of the pure $SU(2)$ gauge theory.\end{note}

For the corresponding complex structure problem, the steps of determining the singular integral affine structure with a BPS polyhedral decomposition, a log smooth structure and the slab functions inducing it are almost identical to the first example since we also have only two charges in this case. We use the same polarization. We use the same notations used in example 1. The slab functions are given by  $$f_{\rho_{1},v}=f_{\rho_{3},v}=(1+y)^{2}=(1+tx^{-1})^{2}, \ f_{\rho_{2},v}=f_{\rho_{4},v}=(1+w)^{2}=(1+tz^{-1})^{2}$$ where every symbol has the same meaning as in the first example. In particular $$xy=t, zw=t$$ Here we have used a trick to avoid setting up new notations. Instead of picking two charges with intersection number 2 we stick to the old notations $\gamma_{1}$ and $\gamma_{2}$ (with intersection number 1) but we make the compensation by raising the power of the slab function from 1 to 2 so that the log morphisms stay the same.

By the wall crossing formula  the new rays (Gross-Siebert walls) one has to add to make a compatible system of consistent $structures$ are $$\{(\mathbf{R}_{\geq0}(n+1,n),(1+y^{n+1}w^{n})^{2}),(\mathbf{R}_{\geq0}(n,n+1),(1+y^{n}w^{n+1})^{2}), (\mathbf{R}_{\geq0}(1,1),(1-wy)^{-4}) \}$$

It is instructive to run the algorithm of truncating by powers of $t$ and see how one can obtain the above result. We start from $k=1$. So we work over $\mathrm{Spec}\ C[[t]]/(t^{2})$. The composition of four log automorphisms associated to the four initial rays is the identity modulo $t^{2}$. In fact starting from the cell $\sigma_{1}$ it is given by\footnote{Note that we mod out by $t^{2}$ after each log automorphism is composed instead of doing that after composing all automorphisms. This reduces considerably the amount of calculations.} $$x\rightarrow x(1+z^{-1}t)^{2}=x(1+2z^{-1}t)\ \mathrm{mod} \ t^{2}\rightarrow$$ $$\rightarrow x(1+2z^{-1}(1+x^{-1}t)^{-2}t)=x(1+2z^{-1}t)\ \mathrm{mod} \ t^{2}\rightarrow $$ $$\rightarrow x(1+z^{-1}t)^{-2} (1+2z^{-1}t)=x \ \mathrm{mod} \ t^{2}\rightarrow  x \ \mathrm{mod} \ t^{2}$$and similar results for $y,z,w$. So it is consistent over $\mathrm{Spec}\ C[t]/(t^{2})$. Now let $k=2$, we get $$x\rightarrow x(1+2z^{-1}t+z^{-2}t^{2})\ \mathrm{mod}\ t^{3}\rightarrow x(1+2z^{-1}t -4z^{-1}x^{-1}t^{2}+z^{-2}t^{2})\ \mathrm{mod}\ t^{3}$$ $$\rightarrow x(1-4z^{-1}x^{-1}t^{2})\ \mathrm{mod}\ t^{3}\rightarrow x(1-4z^{-1}x^{-1}t^{2})\ \mathrm{mod}\ t^{3}$$ and similar results for $y,z,w$. The nontrivial automorphisms are canceled modulo $t^{3}$ by adding a ray $\rho$ with slope $(1,1)$  with attached function\footnote{One can check that after crossing this ray we get $$x(1+4z^{-1}x^{-1}t^{2})(1-4z^{-1}x^{-1}t^{2})\ \mathrm{mod} \ t^{3}=x \ \mathrm{mod} \ t^{3}$$} $f_{\rho}=(1-t^{2}x^{-1}z^{-1})^{-4}=(1-wy)^{-4}$.

We can continue and according to the power of $t$ we get $(\mathbf{R}(n+1,n),(1+y^{n+1}w^{n})^{2})$ and $(\mathbf{R}(n,n+1),(1+y^{n}w^{n+1})^{2})$ at the $2n+1$-th order. Let $k\rightarrow\infty$ and we are done. Note that the power of $t$ is the same of the sum of negative degrees of $x,z$ indicating that the truncation here is the same as the degree truncation.  Also note that for each $k$, only finitely many rays are added while in the derivation by listing critical trajectories it is obtained by one strike by moving $\vartheta$ to $\vartheta+\pi$ without taking any truncations.

One can repeat the gluing algorithm to  get the defining equations of the Hitchin's moduli space. Now we have four maximal dimensional cells in the polyhedral decomposition and infinitely many chambers and we know the gluing in the first quadrant (the "infinite" region) are just changes of chambers such that  the log automorphisms attached to Gross-Siebert walls guarantee the gluing consistency. The gluing of thickennings  $R_{\rho_{2}\rightarrow\sigma_{1},\sigma_{\mathrm{u}_{1}}}^{k}$ and $R_{\rho_{2}\rightarrow\sigma_{2},\sigma_{2}}^{k}$, etc proceed  analogously as calculating fiber products with log automorphisms composed.

 By proposition 9.14 we know that the ideal of the toric degeneration  is generated by \begin{equation}XY-(1+W)^{2}t, ZW-(1+Y)^{2}t\end{equation} because $f_{\rho_{i}}$ are $(1+w)^{2}$ and $(1+y)^{2}$.

Set $t=1$ and we get \begin{equation}XY=(1+W)^{2}, ZW=(1+Y)^{2}\end{equation} As the intersection of two degree two hypersurfaces in $\mathbf{C}^{4}$ it is an affine Calabi-Yau variety.

 Choose a WKB triangulation corresponding a moduli parameter inside the union of stability walls  and choose $x_{1}, y_{1}$ as Fock-Goncharov coordinates labeled by the two nondegenerate edges. Then according to the cluster transformations under the two flips encountered by changing $\vartheta$ to $\vartheta+\pi$ we define $$y_{2}=x_{1}^{-1}, x_{2}=y_{1}(1+x_{1})^{2}$$ $$y_{3}=x_{2}^{-1}, x_{3}=y_{2}(1+x_{2})^{2}$$ and  we get the defining equations $$x_{1}x_{3}=(1+x_{2})^{2}=(1+{1\over y_{3}})^{2}$$ $${1\over y_{1}}{1\over y_{3}}=(1+{1\over y_{2}})^{2}=(1+x_{1})^{2}$$which are equivalent to (149) via the identification of the variables labeled by the same charges.

 On the other hand, if we are in the other stability chamber then we have a projective system of coordinate rings. For each $k$ the consistent scattering diagram at the order $2k+1$ has finitely many rays $$\{(\mathbf{R}(n+1,n),(1+y^{n+1}w^{n})^{2}),(\mathbf{R}(n,n+1),(1+y^{n}w^{n+1})^{2}),(\mathbf{R}(1,1),(1-wy)^{-4}), n\leq k\}$$  We choose  BPS faces associated to these charges and we have finitely many maximal dimensional cells in the corresponding polyhedral decomposition. So in this case we have a  system of degenerations over $\mathrm{Spec}\mathbf{C}[t]/(t^{k+1})$ for varying $k$.\\

 What if we want to do the continuous evolution along a loop in the space of pairs of $(\vartheta,u)$ crossing the stability wall? Infinitely many rays/charges will be encountered between the ray $(1,1)$ and the ray $(k+1,k)$ as well as the $(1,1)$ and $(k,k+1)$. If we move only along one direction we cannot crossing the ray (1,1) without doing truncations which spoils the order of rays. So we move in both directions. In other words, we can start from a phase between the phase of BPS rays of $\gamma_{1}$ and $\gamma_{2}$ at a point inside the stability walls and then follow the loop in both directions to cross the stability walls and enter the "upper" and "lower" infinite regions of rays. So in either direction we can take the limit of ordered product without taking truncations and in the end we figure out the "juggle" transformation between the two limits. For example, the wall crossing formula (147) can be understood  in the following way. We let \begin{equation}K_{\gamma_{2}}^{-1}K_{\gamma_{1}}^{-1}K_{\gamma_{2}}K_{\gamma_{1}+2\gamma_{2}}K_{2\gamma_{1}+3\gamma_{3}}\cdots \end{equation}act on $\mathcal{X}_{\gamma_{1}}, \mathcal{X}_{\gamma_{2}}$ by putting  $\mathcal{X}_{\gamma_{i}}$ to the left of $K_{\gamma_{2}}^{-1}$ and try to find the limit.\footnote{This is just an awkward way to  say that we are really calculating $$\cdots K_{\gamma_{1}+2\gamma_{2}}^{-1}K_{\gamma_{2}}^{-1}K_{\gamma_{1}}K_{\gamma_{2}}\mathcal{X}_{\gamma_{i}}$$ in the usual sense. The arrangement of (150) has the virtue of making the derivation of the wall crossing formula easier to visualize.}  The we try to find the limit $$\cdots K_{3\gamma_{1}+2\gamma_{2}}K_{2\gamma_{1}+\gamma_{2}}K_{\gamma_{1}}\mathcal{X}_{\gamma_{i}}$$Finally we verify that the transformation $K_{\gamma_{1}+\gamma_{2}}^{-2}$ brings the second limit to the first one. In this way we can derive the wall crossing formula without doing any truncations. See the appendix of \citep{G1} for detailed calculations of this example\footnote{But be aware that the signs are different there due to the different assignments of quadratic refinements.}.

 So we have three ways of deriving an explicit wall crossing formula. By following a continuous evolution together with this algebraic trick, by inductive truncations given an initial  factorization in a primary stability chamber and by using the spectrum generator. They yield the same result.\\

\end{document}